\documentclass[a4paper,12pt,leqno]{article}%{amsart}
\usepackage{amssymb}
\usepackage{amsmath}
\usepackage{array}
\usepackage{multirow}
\usepackage{theorem}

\theoremstyle{change}
\theorembodyfont{\itshape}

\newcommand{\Lg}{\mbox{$\mathfrak g$}}
\newcommand{\Ll}{\mbox{$\mathfrak l$}}
\newcommand{\Lh}{\mbox{$\mathfrak h$}}
\newcommand{\Lk}{\mbox{$\mathfrak k$}}
\newcommand{\Lp}{\mbox{$\mathfrak p$}}
\newcommand{\La}{\mbox{$\mathfrak a$}}

\newcommand{\Lc}{\mbox{$\mathfrak c$}}
\newcommand{\Lt}{\mbox{$\mathfrak t$}}

\newcommand{\Pf}{{\em Proof}. }
\newcommand{\EPf}{\hfill$\square$}

\newcommand{\inn}[2]{\mbox{$<\!#1\,,#2\!>$}}
\newcommand{\her}[2]{\mbox{$<\!\!<\!#1\,,#2\!>\!\!>$}}
\newcommand{\ad}[1]{\mbox{$\mbox{ad}_{#1}$}}

\newcommand{\SU}[1]{\mbox{$\mathbf{SU}(#1)$}}
\newcommand{\U}[1]{\mbox{$\mathbf{U}(#1)$}}
\newcommand{\SP}[1]{\mbox{$\mathbf{Sp}(#1)$}}
\newcommand{\SO}[1]{\mbox{$\mathbf{SO}(#1)$}}
\newcommand{\OG}[1]{\mbox{$\mathbf{O}(#1)$}}
\newcommand{\Spin}[1]{\mbox{$\mathbf{Spin}(#1)$}}
\newcommand{\A}[1]{\mbox{$\mathbf{A}_{#1}$}}
\newcommand{\B}[1]{\mbox{$\mathbf{B}_{#1}$}}
\newcommand{\C}[1]{\mbox{$\mathbf{C}_{#1}$}}
\newcommand{\D}[1]{\mbox{$\mathbf{D}_{#1}$}}
\newcommand{\G}{\mbox{$\mathbf{G}_2$}}
\newcommand{\E}[1]{\mbox{$\mathbf{E}_{#1}$}}
\newcommand{\F}{\mbox{$\mathbf{F}_4$}}

\newtheorem{thm}{Theorem}

\newtheorem{rmk}[thm]{Remark}
\newtheorem{cor}[thm]{Corollary}
\newtheorem{prop}[thm]{Proposition}
\newtheorem{lem}[thm]{Lemma}
\newtheorem{exs}[thm]{Examples}

\theoremstyle{plain}
\newtheorem{claim}{\sc Claim}

\begingroup\makeatletter\ifx\SetFigFont\undefined
\def\x#1#2#3#4#5#6#7\relax{\def\x{#1#2#3#4#5#6}}%
\expandafter\x\fmtname xxxxxx\relax \def\y{splain}%
\ifx\x\y   % LaTeX or SliTeX?
\gdef\SetFigFont#1#2#3{%
  \ifnum #1<17\tiny\else \ifnum #1<20\small\else
  \ifnum #1<24\normalsize\else \ifnum #1<29\large\else
  \ifnum #1<34\Large\else \ifnum #1<41\LARGE\else
     \huge\fi\fi\fi\fi\fi\fi
  \csname #3\endcsname}%
\else
\gdef\SetFigFont#1#2#3{\begingroup
  \count@#1\relax \ifnum 25<\count@\count@25\fi
  \def\x{\endgroup\@setsize\SetFigFont{#2pt}}%
  \expandafter\x
    \csname \romannumeral\the\count@ pt\expandafter\endcsname
    \csname @\romannumeral\the\count@ pt\endcsname
  \csname #3\endcsname}%
\fi
\fi\endgroup

\newcommand{\Avii}[7]{
\begin{picture}(7552,895)(1025,-4136)
\thinlines
\put(1201,-3961){\circle{336}}
\put(2401,-3961){\circle{336}}
\put(3601,-3961){\circle{336}}
\put(4801,-3961){\circle{336}}
\put(6001,-3961){\circle{336}}
\put(7201,-3961){\circle{336}}
\put(8401,-3961){\circle{336}}
\put(1369,-3961){\line( 1, 0){864}}
\put(2569,-3961){\line( 1, 0){864}}
\put(3769,-3961){\line( 1, 0){864}}
\put(4969,-3961){\line( 1, 0){864}}
\put(6169,-3961){\line( 1, 0){864}}
\put(7369,-3961){\line( 1, 0){864}}
\put(1141,-3436){\makebox(0,0)[lb]{\smash{\SetFigFont{6}{7.2}{rm}#1}}}
\put(2341,-3436){\makebox(0,0)[lb]{\smash{\SetFigFont{6}{7.2}{rm}#2}}}
\put(3541,-3436){\makebox(0,0)[lb]{\smash{\SetFigFont{6}{7.2}{rm}#3}}}
\put(4741,-3436){\makebox(0,0)[lb]{\smash{\SetFigFont{6}{7.2}{rm}#4}}}
\put(5941,-3436){\makebox(0,0)[lb]{\smash{\SetFigFont{6}{7.2}{rm}#5}}}
\put(7141,-3436){\makebox(0,0)[lb]{\smash{\SetFigFont{6}{7.2}{rm}#6}}}
\put(8341,-3436){\makebox(0,0)[lb]{\smash{\SetFigFont{6}{7.2}{rm}#6}}}
\end{picture}
}
\begingroup\makeatletter\ifx\SetFigFont\undefined
\def\x#1#2#3#4#5#6#7\relax{\def\x{#1#2#3#4#5#6}}%
\expandafter\x\fmtname xxxxxx\relax \def\y{splain}%
\ifx\x\y   % LaTeX or SliTeX?
\gdef\SetFigFont#1#2#3{%
  \ifnum #1<17\tiny\else \ifnum #1<20\small\else
  \ifnum #1<24\normalsize\else \ifnum #1<29\large\else
  \ifnum #1<34\Large\else \ifnum #1<41\LARGE\else
     \huge\fi\fi\fi\fi\fi\fi
  \csname #3\endcsname}%
\else
\gdef\SetFigFont#1#2#3{\begingroup
  \count@#1\relax \ifnum 25<\count@\count@25\fi
  \def\x{\endgroup\@setsize\SetFigFont{#2pt}}%
  \expandafter\x
    \csname \romannumeral\the\count@ pt\expandafter\endcsname
    \csname @\romannumeral\the\count@ pt\endcsname
  \csname #3\endcsname}%
\fi
\fi\endgroup

\newcommand{\Avi}[6]{
\begin{picture}(6352,895)(1025,-4136)
\thinlines
\put(1201,-3961){\circle{336}}
\put(2401,-3961){\circle{336}}
\put(3601,-3961){\circle{336}}
\put(4801,-3961){\circle{336}}
\put(6001,-3961){\circle{336}}
\put(7201,-3961){\circle{336}}
\put(1369,-3961){\line( 1, 0){864}}
\put(2569,-3961){\line( 1, 0){864}}
\put(3769,-3961){\line( 1, 0){864}}
\put(4969,-3961){\line( 1, 0){864}}
\put(6169,-3961){\line( 1, 0){864}}
\put(1141,-3436){\makebox(0,0)[lb]{\smash{\SetFigFont{6}{7.2}{rm}#1}}}
\put(2341,-3436){\makebox(0,0)[lb]{\smash{\SetFigFont{6}{7.2}{rm}#2}}}
\put(3541,-3436){\makebox(0,0)[lb]{\smash{\SetFigFont{6}{7.2}{rm}#3}}}
\put(4741,-3436){\makebox(0,0)[lb]{\smash{\SetFigFont{6}{7.2}{rm}#4}}}
\put(5941,-3436){\makebox(0,0)[lb]{\smash{\SetFigFont{6}{7.2}{rm}#5}}}
\put(7141,-3436){\makebox(0,0)[lb]{\smash{\SetFigFont{6}{7.2}{rm}#6}}}
\end{picture}
}
\begingroup\makeatletter\ifx\SetFigFont\undefined
\def\x#1#2#3#4#5#6#7\relax{\def\x{#1#2#3#4#5#6}}%
\expandafter\x\fmtname xxxxxx\relax \def\y{splain}%
\ifx\x\y   % LaTeX or SliTeX?
\gdef\SetFigFont#1#2#3{%
  \ifnum #1<17\tiny\else \ifnum #1<20\small\else
  \ifnum #1<24\normalsize\else \ifnum #1<29\large\else
  \ifnum #1<34\Large\else \ifnum #1<41\LARGE\else
     \huge\fi\fi\fi\fi\fi\fi
  \csname #3\endcsname}%
\else
\gdef\SetFigFont#1#2#3{\begingroup
  \count@#1\relax \ifnum 25<\count@\count@25\fi
  \def\x{\endgroup\@setsize\SetFigFont{#2pt}}%
  \expandafter\x
    \csname \romannumeral\the\count@ pt\expandafter\endcsname
    \csname @\romannumeral\the\count@ pt\endcsname
  \csname #3\endcsname}%
\fi
\fi\endgroup

\newcommand{\Av}[5]{
\begin{picture}(5152,895)(1025,-4136)
\thinlines
\put(1201,-3961){\circle{336}}
\put(2401,-3961){\circle{336}}
\put(3601,-3961){\circle{336}}
\put(4801,-3961){\circle{336}}
\put(6001,-3961){\circle{336}}
\put(1369,-3961){\line( 1, 0){864}}
\put(2569,-3961){\line( 1, 0){864}}
\put(3769,-3961){\line( 1, 0){864}}
\put(4969,-3961){\line( 1, 0){864}}
\put(1141,-3436){\makebox(0,0)[lb]{\smash{\SetFigFont{6}{7.2}{rm}#1}}}
\put(2341,-3436){\makebox(0,0)[lb]{\smash{\SetFigFont{6}{7.2}{rm}#2}}}
\put(3541,-3436){\makebox(0,0)[lb]{\smash{\SetFigFont{6}{7.2}{rm}#3}}}
\put(4741,-3436){\makebox(0,0)[lb]{\smash{\SetFigFont{6}{7.2}{rm}#4}}}
\put(5941,-3436){\makebox(0,0)[lb]{\smash{\SetFigFont{6}{7.2}{rm}#5}}}
\end{picture}
}
\begingroup\makeatletter\ifx\SetFigFont\undefined
\def\x#1#2#3#4#5#6#7\relax{\def\x{#1#2#3#4#5#6}}%
\expandafter\x\fmtname xxxxxx\relax \def\y{splain}%
\ifx\x\y   % LaTeX or SliTeX?
\gdef\SetFigFont#1#2#3{%
  \ifnum #1<17\tiny\else \ifnum #1<20\small\else
  \ifnum #1<24\normalsize\else \ifnum #1<29\large\else
  \ifnum #1<34\Large\else \ifnum #1<41\LARGE\else
     \huge\fi\fi\fi\fi\fi\fi
  \csname #3\endcsname}%
\else
\gdef\SetFigFont#1#2#3{\begingroup
  \count@#1\relax \ifnum 25<\count@\count@25\fi
  \def\x{\endgroup\@setsize\SetFigFont{#2pt}}%
  \expandafter\x
    \csname \romannumeral\the\count@ pt\expandafter\endcsname
    \csname @\romannumeral\the\count@ pt\endcsname
  \csname #3\endcsname}%
\fi
\fi\endgroup

\newcommand{\Aiv}[4]{
\begin{picture}(3952,895)(1025,-4136)
\thinlines
\put(1201,-3961){\circle{336}}
\put(2401,-3961){\circle{336}}
\put(3601,-3961){\circle{336}}
\put(4801,-3961){\circle{336}}
\put(1369,-3961){\line( 1, 0){864}}
\put(2569,-3961){\line( 1, 0){864}}
\put(3769,-3961){\line( 1, 0){864}}
\put(1141,-3436){\makebox(0,0)[lb]{\smash{\SetFigFont{6}{7.2}{rm}#1}}}
\put(2341,-3436){\makebox(0,0)[lb]{\smash{\SetFigFont{6}{7.2}{rm}#2}}}
\put(3541,-3436){\makebox(0,0)[lb]{\smash{\SetFigFont{6}{7.2}{rm}#3}}}
\put(4741,-3436){\makebox(0,0)[lb]{\smash{\SetFigFont{6}{7.2}{rm}#4}}}
\end{picture}
}

\begingroup\makeatletter\ifx\SetFigFont\undefined
\def\x#1#2#3#4#5#6#7\relax{\def\x{#1#2#3#4#5#6}}%
\expandafter\x\fmtname xxxxxx\relax \def\y{splain}%
\ifx\x\y   % LaTeX or SliTeX?
\gdef\SetFigFont#1#2#3{%
  \ifnum #1<17\tiny\else \ifnum #1<20\small\else
  \ifnum #1<24\normalsize\else \ifnum #1<29\large\else
  \ifnum #1<34\Large\else \ifnum #1<41\LARGE\else
     \huge\fi\fi\fi\fi\fi\fi
  \csname #3\endcsname}%
\else
\gdef\SetFigFont#1#2#3{\begingroup
  \count@#1\relax \ifnum 25<\count@\count@25\fi
  \def\x{\endgroup\@setsize\SetFigFont{#2pt}}%
  \expandafter\x
    \csname \romannumeral\the\count@ pt\expandafter\endcsname
    \csname @\romannumeral\the\count@ pt\endcsname
  \csname #3\endcsname}%
\fi
\fi\endgroup

\newcommand{\Aiii}[3]{
\begin{picture}(2752,895)(1025,-4136)
\thinlines
\put(1201,-3961){\circle{336}}
\put(2401,-3961){\circle{336}}
\put(3601,-3961){\circle{336}}
\put(1369,-3961){\line( 1, 0){864}}
\put(2569,-3961){\line( 1, 0){864}}
\put(1141,-3436){\makebox(0,0)[lb]{\smash{\SetFigFont{6}{7.2}{rm}#1}}}
\put(2341,-3436){\makebox(0,0)[lb]{\smash{\SetFigFont{6}{7.2}{rm}#2}}}
\put(3541,-3436){\makebox(0,0)[lb]{\smash{\SetFigFont{6}{7.2}{rm}#3}}}
\end{picture}
}

\begingroup\makeatletter\ifx\SetFigFont\undefined
\def\x#1#2#3#4#5#6#7\relax{\def\x{#1#2#3#4#5#6}}%
\expandafter\x\fmtname xxxxxx\relax \def\y{splain}%
\ifx\x\y   % LaTeX or SliTeX?
\gdef\SetFigFont#1#2#3{%
  \ifnum #1<17\tiny\else \ifnum #1<20\small\else
  \ifnum #1<24\normalsize\else \ifnum #1<29\large\else
  \ifnum #1<34\Large\else \ifnum #1<41\LARGE\else
     \huge\fi\fi\fi\fi\fi\fi
  \csname #3\endcsname}%
\else
\gdef\SetFigFont#1#2#3{\begingroup
  \count@#1\relax \ifnum 25<\count@\count@25\fi
  \def\x{\endgroup\@setsize\SetFigFont{#2pt}}%
  \expandafter\x
    \csname \romannumeral\the\count@ pt\expandafter\endcsname
    \csname @\romannumeral\the\count@ pt\endcsname
  \csname #3\endcsname}%
\fi
\fi\endgroup

\newcommand{\Aii}[2]{
\begin{picture}(1552,895)(1025,-4136)
\thinlines
\put(1201,-3961){\circle{336}}
\put(2401,-3961){\circle{336}}
\put(1369,-3961){\line( 1, 0){864}}
\put(1141,-3436){\makebox(0,0)[lb]{\smash{\SetFigFont{6}{7.2}{rm}#1}}}
\put(2341,-3436){\makebox(0,0)[lb]{\smash{\SetFigFont{6}{7.2}{rm}#2}}}
\end{picture}
}

\begingroup\makeatletter\ifx\SetFigFont\undefined
\def\x#1#2#3#4#5#6#7\relax{\def\x{#1#2#3#4#5#6}}%
\expandafter\x\fmtname xxxxxx\relax \def\y{splain}%
\ifx\x\y   % LaTeX or SliTeX?
\gdef\SetFigFont#1#2#3{%
  \ifnum #1<17\tiny\else \ifnum #1<20\small\else
  \ifnum #1<24\normalsize\else \ifnum #1<29\large\else
  \ifnum #1<34\Large\else \ifnum #1<41\LARGE\else
     \huge\fi\fi\fi\fi\fi\fi
  \csname #3\endcsname}%
\else
\gdef\SetFigFont#1#2#3{\begingroup
  \count@#1\relax \ifnum 25<\count@\count@25\fi
  \def\x{\endgroup\@setsize\SetFigFont{#2pt}}%
  \expandafter\x
    \csname \romannumeral\the\count@ pt\expandafter\endcsname
    \csname @\romannumeral\the\count@ pt\endcsname
  \csname #3\endcsname}%
\fi
\fi\endgroup

\newcommand{\Ai}[1]{
\begin{picture}(352,895)(1025,-4136)
\thinlines
\put(1201,-3961){\circle{336}}
\put(1141,-3436){\makebox(0,0)[lb]{\smash{\SetFigFont{6}{7.2}{rm}\mbox{$#1$}}}}
\end{picture}
}

\begingroup\makeatletter\ifx\SetFigFont\undefined
\def\x#1#2#3#4#5#6#7\relax{\def\x{#1#2#3#4#5#6}}%
\expandafter\x\fmtname xxxxxx\relax \def\y{splain}%
\ifx\x\y   % LaTeX or SliTeX?
\gdef\SetFigFont#1#2#3{%
  \ifnum #1<17\tiny\else \ifnum #1<20\small\else
  \ifnum #1<24\normalsize\else \ifnum #1<29\large\else
  \ifnum #1<34\Large\else \ifnum #1<41\LARGE\else
     \huge\fi\fi\fi\fi\fi\fi
  \csname #3\endcsname}%
\else
\gdef\SetFigFont#1#2#3{\begingroup
  \count@#1\relax \ifnum 25<\count@\count@25\fi
  \def\x{\endgroup\@setsize\SetFigFont{#2pt}}%
  \expandafter\x
    \csname \romannumeral\the\count@ pt\expandafter\endcsname
    \csname @\romannumeral\the\count@ pt\endcsname
  \csname #3\endcsname}%
\fi
\fi\endgroup

\newcommand{\Bviii}[8]{
\begin{picture}(8752,895)(1025,-4136)
\thinlines
\put(1201,-3961){\circle{336}}
\put(2401,-3961){\circle{336}}
\put(3601,-3961){\circle{336}}
\put(4801,-3961){\circle{336}}
\put(6001,-3961){\circle{336}}
\put(7201,-3961){\circle{336}}
\put(8401,-3961){\circle{336}}
\put(9601,-3961){\circle*{336}}
\put(1369,-3961){\line( 1, 0){864}}
\put(2569,-3961){\line( 1, 0){864}}
\put(3769,-3961){\line( 1, 0){864}}
\put(4969,-3961){\line( 1, 0){864}}
\put(6169,-3961){\line( 1, 0){864}}
\put(7369,-3961){\line( 1, 0){864}}
\put(8569,-3911){\line( 1, 0){864}}
\put(8569,-4011){\line( 1, 0){864}}
\put(1141,-3436){\makebox(0,0)[lb]{\smash{\SetFigFont{6}{7.2}{rm}#1}}}
\put(2341,-3436){\makebox(0,0)[lb]{\smash{\SetFigFont{6}{7.2}{rm}#2}}}
\put(3541,-3436){\makebox(0,0)[lb]{\smash{\SetFigFont{6}{7.2}{rm}#3}}}
\put(4741,-3436){\makebox(0,0)[lb]{\smash{\SetFigFont{6}{7.2}{rm}#4}}}
\put(5941,-3436){\makebox(0,0)[lb]{\smash{\SetFigFont{6}{7.2}{rm}#5}}}
\put(7141,-3436){\makebox(0,0)[lb]{\smash{\SetFigFont{6}{7.2}{rm}#6}}}
\put(8341,-3436){\makebox(0,0)[lb]{\smash{\SetFigFont{6}{7.2}{rm}#7}}}
\put(9541,-3436){\makebox(0,0)[lb]{\smash{\SetFigFont{6}{7.2}{rm}#8}}}
\end{picture}
}
\begingroup\makeatletter\ifx\SetFigFont\undefined
\def\x#1#2#3#4#5#6#7\relax{\def\x{#1#2#3#4#5#6}}%
\expandafter\x\fmtname xxxxxx\relax \def\y{splain}%
\ifx\x\y   % LaTeX or SliTeX?
\gdef\SetFigFont#1#2#3{%
  \ifnum #1<17\tiny\else \ifnum #1<20\small\else
  \ifnum #1<24\normalsize\else \ifnum #1<29\large\else
  \ifnum #1<34\Large\else \ifnum #1<41\LARGE\else
     \huge\fi\fi\fi\fi\fi\fi
  \csname #3\endcsname}%
\else
\gdef\SetFigFont#1#2#3{\begingroup
  \count@#1\relax \ifnum 25<\count@\count@25\fi
  \def\x{\endgroup\@setsize\SetFigFont{#2pt}}%
  \expandafter\x
    \csname \romannumeral\the\count@ pt\expandafter\endcsname
    \csname @\romannumeral\the\count@ pt\endcsname
  \csname #3\endcsname}%
\fi
\fi\endgroup

\newcommand{\Bvii}[7]{
\begin{picture}(7552,895)(2225,-4136)
\thinlines
\put(2401,-3961){\circle{336}}
\put(3601,-3961){\circle{336}}
\put(4801,-3961){\circle{336}}
\put(6001,-3961){\circle{336}}
\put(7201,-3961){\circle{336}}
\put(8401,-3961){\circle{336}}
\put(9601,-3961){\circle*{336}}
\put(2569,-3961){\line( 1, 0){864}}
\put(3769,-3961){\line( 1, 0){864}}
\put(4969,-3961){\line( 1, 0){864}}
\put(6169,-3961){\line( 1, 0){864}}
\put(7369,-3961){\line( 1, 0){864}}
\put(8569,-3911){\line( 1, 0){864}}
\put(8569,-4011){\line( 1, 0){864}}
\put(2341,-3436){\makebox(0,0)[lb]{\smash{\SetFigFont{6}{7.2}{rm}#1}}}
\put(3541,-3436){\makebox(0,0)[lb]{\smash{\SetFigFont{6}{7.2}{rm}#2}}}
\put(4741,-3436){\makebox(0,0)[lb]{\smash{\SetFigFont{6}{7.2}{rm}#3}}}
\put(5941,-3436){\makebox(0,0)[lb]{\smash{\SetFigFont{6}{7.2}{rm}#4}}}
\put(7141,-3436){\makebox(0,0)[lb]{\smash{\SetFigFont{6}{7.2}{rm}#5}}}
\put(8341,-3436){\makebox(0,0)[lb]{\smash{\SetFigFont{6}{7.2}{rm}#6}}}
\put(9541,-3436){\makebox(0,0)[lb]{\smash{\SetFigFont{6}{7.2}{rm}#7}}}
\end{picture}
}
\begingroup\makeatletter\ifx\SetFigFont\undefined
\def\x#1#2#3#4#5#6#7\relax{\def\x{#1#2#3#4#5#6}}%
\expandafter\x\fmtname xxxxxx\relax \def\y{splain}%
\ifx\x\y   % LaTeX or SliTeX?
\gdef\SetFigFont#1#2#3{%
  \ifnum #1<17\tiny\else \ifnum #1<20\small\else
  \ifnum #1<24\normalsize\else \ifnum #1<29\large\else
  \ifnum #1<34\Large\else \ifnum #1<41\LARGE\else
     \huge\fi\fi\fi\fi\fi\fi
  \csname #3\endcsname}%
\else
\gdef\SetFigFont#1#2#3{\begingroup
  \count@#1\relax \ifnum 25<\count@\count@25\fi
  \def\x{\endgroup\@setsize\SetFigFont{#2pt}}%
  \expandafter\x
    \csname \romannumeral\the\count@ pt\expandafter\endcsname
    \csname @\romannumeral\the\count@ pt\endcsname
  \csname #3\endcsname}%
\fi
\fi\endgroup

\newcommand{\Bvi}[6]{
\begin{picture}(6352,895)(3425,-4136)
\thinlines
\put(3601,-3961){\circle{336}}
\put(4801,-3961){\circle{336}}
\put(6001,-3961){\circle{336}}
\put(7201,-3961){\circle{336}}
\put(8401,-3961){\circle{336}}
\put(9601,-3961){\circle*{336}}
\put(3769,-3961){\line( 1, 0){864}}
\put(4969,-3961){\line( 1, 0){864}}
\put(6169,-3961){\line( 1, 0){864}}
\put(7369,-3961){\line( 1, 0){864}}
\put(8569,-3911){\line( 1, 0){864}}
\put(8569,-4011){\line( 1, 0){864}}
\put(3541,-3436){\makebox(0,0)[lb]{\smash{\SetFigFont{6}{7.2}{rm}#1}}}
\put(4741,-3436){\makebox(0,0)[lb]{\smash{\SetFigFont{6}{7.2}{rm}#2}}}
\put(5941,-3436){\makebox(0,0)[lb]{\smash{\SetFigFont{6}{7.2}{rm}#3}}}
\put(7141,-3436){\makebox(0,0)[lb]{\smash{\SetFigFont{6}{7.2}{rm}#4}}}
\put(8341,-3436){\makebox(0,0)[lb]{\smash{\SetFigFont{6}{7.2}{rm}#5}}}
\put(9541,-3436){\makebox(0,0)[lb]{\smash{\SetFigFont{6}{7.2}{rm}#6}}}
\end{picture}
}
\begingroup\makeatletter\ifx\SetFigFont\undefined
\def\x#1#2#3#4#5#6#7\relax{\def\x{#1#2#3#4#5#6}}%
\expandafter\x\fmtname xxxxxx\relax \def\y{splain}%
\ifx\x\y   % LaTeX or SliTeX?
\gdef\SetFigFont#1#2#3{%
  \ifnum #1<17\tiny\else \ifnum #1<20\small\else
  \ifnum #1<24\normalsize\else \ifnum #1<29\large\else
  \ifnum #1<34\Large\else \ifnum #1<41\LARGE\else
     \huge\fi\fi\fi\fi\fi\fi
  \csname #3\endcsname}%
\else
\gdef\SetFigFont#1#2#3{\begingroup
  \count@#1\relax \ifnum 25<\count@\count@25\fi
  \def\x{\endgroup\@setsize\SetFigFont{#2pt}}%
  \expandafter\x
    \csname \romannumeral\the\count@ pt\expandafter\endcsname
    \csname @\romannumeral\the\count@ pt\endcsname
  \csname #3\endcsname}%
\fi
\fi\endgroup

\newcommand{\Bv}[5]{
\begin{picture}(5152,895)(4625,-4136)
\thinlines
\put(4801,-3961){\circle{336}}
\put(6001,-3961){\circle{336}}
\put(7201,-3961){\circle{336}}
\put(8401,-3961){\circle{336}}
\put(9601,-3961){\circle*{336}}
\put(4969,-3961){\line( 1, 0){864}}
\put(6169,-3961){\line( 1, 0){864}}
\put(7369,-3961){\line( 1, 0){864}}
\put(8569,-3911){\line( 1, 0){864}}
\put(8569,-4011){\line( 1, 0){864}}
\put(4741,-3436){\makebox(0,0)[lb]{\smash{\SetFigFont{6}{7.2}{rm}#1}}}
\put(5941,-3436){\makebox(0,0)[lb]{\smash{\SetFigFont{6}{7.2}{rm}#2}}}
\put(7141,-3436){\makebox(0,0)[lb]{\smash{\SetFigFont{6}{7.2}{rm}#3}}}
\put(8341,-3436){\makebox(0,0)[lb]{\smash{\SetFigFont{6}{7.2}{rm}#4}}}
\put(9541,-3436){\makebox(0,0)[lb]{\smash{\SetFigFont{6}{7.2}{rm}#5}}}
\end{picture}
}
\begingroup\makeatletter\ifx\SetFigFont\undefined
\def\x#1#2#3#4#5#6#7\relax{\def\x{#1#2#3#4#5#6}}%
\expandafter\x\fmtname xxxxxx\relax \def\y{splain}%
\ifx\x\y   % LaTeX or SliTeX?
\gdef\SetFigFont#1#2#3{%
  \ifnum #1<17\tiny\else \ifnum #1<20\small\else
  \ifnum #1<24\normalsize\else \ifnum #1<29\large\else
  \ifnum #1<34\Large\else \ifnum #1<41\LARGE\else
     \huge\fi\fi\fi\fi\fi\fi
  \csname #3\endcsname}%
\else
\gdef\SetFigFont#1#2#3{\begingroup
  \count@#1\relax \ifnum 25<\count@\count@25\fi
  \def\x{\endgroup\@setsize\SetFigFont{#2pt}}%
  \expandafter\x
    \csname \romannumeral\the\count@ pt\expandafter\endcsname
    \csname @\romannumeral\the\count@ pt\endcsname
  \csname #3\endcsname}%
\fi
\fi\endgroup

\newcommand{\Biv}[4]{
\begin{picture}(3952,895)(5825,-4136)
\thinlines
\put(6001,-3961){\circle{336}}
\put(7201,-3961){\circle{336}}
\put(8401,-3961){\circle{336}}
\put(9601,-3961){\circle*{336}}
\put(6169,-3961){\line( 1, 0){864}}
\put(7369,-3961){\line( 1, 0){864}}
\put(8569,-3911){\line( 1, 0){864}}
\put(8569,-4011){\line( 1, 0){864}}
\put(5941,-3436){\makebox(0,0)[lb]{\smash{\SetFigFont{6}{7.2}{rm}#1}}}
\put(7141,-3436){\makebox(0,0)[lb]{\smash{\SetFigFont{6}{7.2}{rm}#2}}}
\put(8341,-3436){\makebox(0,0)[lb]{\smash{\SetFigFont{6}{7.2}{rm}#3}}}
\put(9541,-3436){\makebox(0,0)[lb]{\smash{\SetFigFont{6}{7.2}{rm}#4}}}
\end{picture}
}
\begingroup\makeatletter\ifx\SetFigFont\undefined
\def\x#1#2#3#4#5#6#7\relax{\def\x{#1#2#3#4#5#6}}%
\expandafter\x\fmtname xxxxxx\relax \def\y{splain}%
\ifx\x\y   % LaTeX or SliTeX?
\gdef\SetFigFont#1#2#3{%
  \ifnum #1<17\tiny\else \ifnum #1<20\small\else
  \ifnum #1<24\normalsize\else \ifnum #1<29\large\else
  \ifnum #1<34\Large\else \ifnum #1<41\LARGE\else
     \huge\fi\fi\fi\fi\fi\fi
  \csname #3\endcsname}%
\else
\gdef\SetFigFont#1#2#3{\begingroup
  \count@#1\relax \ifnum 25<\count@\count@25\fi
  \def\x{\endgroup\@setsize\SetFigFont{#2pt}}%
  \expandafter\x
    \csname \romannumeral\the\count@ pt\expandafter\endcsname
    \csname @\romannumeral\the\count@ pt\endcsname
  \csname #3\endcsname}%
\fi
\fi\endgroup

\newcommand{\Biii}[3]{
\begin{picture}(2752,895)(7025,-4136)
\thinlines
\put(7201,-3961){\circle{336}}
\put(8401,-3961){\circle{336}}
\put(9601,-3961){\circle*{336}}
\put(7369,-3961){\line( 1, 0){864}}
\put(8569,-3911){\line( 1, 0){864}}
\put(8569,-4011){\line( 1, 0){864}}
\put(7141,-3436){\makebox(0,0)[lb]{\smash{\SetFigFont{6}{7.2}{rm}#1}}}
\put(8341,-3436){\makebox(0,0)[lb]{\smash{\SetFigFont{6}{7.2}{rm}#2}}}
\put(9541,-3436){\makebox(0,0)[lb]{\smash{\SetFigFont{6}{7.2}{rm}#3}}}
\end{picture}
}
\begingroup\makeatletter\ifx\SetFigFont\undefined
\def\x#1#2#3#4#5#6#7\relax{\def\x{#1#2#3#4#5#6}}%
\expandafter\x\fmtname xxxxxx\relax \def\y{splain}%
\ifx\x\y   % LaTeX or SliTeX?
\gdef\SetFigFont#1#2#3{%
  \ifnum #1<17\tiny\else \ifnum #1<20\small\else
  \ifnum #1<24\normalsize\else \ifnum #1<29\large\else
  \ifnum #1<34\Large\else \ifnum #1<41\LARGE\else
     \huge\fi\fi\fi\fi\fi\fi
  \csname #3\endcsname}%
\else
\gdef\SetFigFont#1#2#3{\begingroup
  \count@#1\relax \ifnum 25<\count@\count@25\fi
  \def\x{\endgroup\@setsize\SetFigFont{#2pt}}%
  \expandafter\x
    \csname \romannumeral\the\count@ pt\expandafter\endcsname
    \csname @\romannumeral\the\count@ pt\endcsname
  \csname #3\endcsname}%
\fi
\fi\endgroup

\begingroup\makeatletter\ifx\SetFigFont\undefined
\def\x#1#2#3#4#5#6#7\relax{\def\x{#1#2#3#4#5#6}}%
\expandafter\x\fmtname xxxxxx\relax \def\y{splain}%
\ifx\x\y   % LaTeX or SliTeX?
\gdef\SetFigFont#1#2#3{%
  \ifnum #1<17\tiny\else \ifnum #1<20\small\else
  \ifnum #1<24\normalsize\else \ifnum #1<29\large\else
  \ifnum #1<34\Large\else \ifnum #1<41\LARGE\else
     \huge\fi\fi\fi\fi\fi\fi
  \csname #3\endcsname}%
\else
\gdef\SetFigFont#1#2#3{\begingroup
  \count@#1\relax \ifnum 25<\count@\count@25\fi
  \def\x{\endgroup\@setsize\SetFigFont{#2pt}}%
  \expandafter\x
    \csname \romannumeral\the\count@ pt\expandafter\endcsname
    \csname @\romannumeral\the\count@ pt\endcsname
  \csname #3\endcsname}%
\fi
\fi\endgroup

\begingroup\makeatletter\ifx\SetFigFont\undefined
\def\x#1#2#3#4#5#6#7\relax{\def\x{#1#2#3#4#5#6}}%
\expandafter\x\fmtname xxxxxx\relax \def\y{splain}%
\ifx\x\y   % LaTeX or SliTeX?
\gdef\SetFigFont#1#2#3{%
  \ifnum #1<17\tiny\else \ifnum #1<20\small\else
  \ifnum #1<24\normalsize\else \ifnum #1<29\large\else
  \ifnum #1<34\Large\else \ifnum #1<41\LARGE\else
     \huge\fi\fi\fi\fi\fi\fi
  \csname #3\endcsname}%
\else
\gdef\SetFigFont#1#2#3{\begingroup
  \count@#1\relax \ifnum 25<\count@\count@25\fi
  \def\x{\endgroup\@setsize\SetFigFont{#2pt}}%
  \expandafter\x
    \csname \romannumeral\the\count@ pt\expandafter\endcsname
    \csname @\romannumeral\the\count@ pt\endcsname
  \csname #3\endcsname}%
\fi
\fi\endgroup

\begingroup\makeatletter\ifx\SetFigFont\undefined
\def\x#1#2#3#4#5#6#7\relax{\def\x{#1#2#3#4#5#6}}%
\expandafter\x\fmtname xxxxxx\relax \def\y{splain}%
\ifx\x\y   % LaTeX or SliTeX?
\gdef\SetFigFont#1#2#3{%
  \ifnum #1<17\tiny\else \ifnum #1<20\small\else
  \ifnum #1<24\normalsize\else \ifnum #1<29\large\else
  \ifnum #1<34\Large\else \ifnum #1<41\LARGE\else
     \huge\fi\fi\fi\fi\fi\fi
  \csname #3\endcsname}%
\else
\gdef\SetFigFont#1#2#3{\begingroup
  \count@#1\relax \ifnum 25<\count@\count@25\fi
  \def\x{\endgroup\@setsize\SetFigFont{#2pt}}%
  \expandafter\x
    \csname \romannumeral\the\count@ pt\expandafter\endcsname
    \csname @\romannumeral\the\count@ pt\endcsname
  \csname #3\endcsname}%
\fi
\fi\endgroup

\begingroup\makeatletter\ifx\SetFigFont\undefined
\def\x#1#2#3#4#5#6#7\relax{\def\x{#1#2#3#4#5#6}}%
\expandafter\x\fmtname xxxxxx\relax \def\y{splain}%
\ifx\x\y   % LaTeX or SliTeX?
\gdef\SetFigFont#1#2#3{%
  \ifnum #1<17\tiny\else \ifnum #1<20\small\else
  \ifnum #1<24\normalsize\else \ifnum #1<29\large\else
  \ifnum #1<34\Large\else \ifnum #1<41\LARGE\else
     \huge\fi\fi\fi\fi\fi\fi
  \csname #3\endcsname}%
\else
\gdef\SetFigFont#1#2#3{\begingroup
  \count@#1\relax \ifnum 25<\count@\count@25\fi
  \def\x{\endgroup\@setsize\SetFigFont{#2pt}}%
  \expandafter\x
    \csname \romannumeral\the\count@ pt\expandafter\endcsname
    \csname @\romannumeral\the\count@ pt\endcsname
  \csname #3\endcsname}%
\fi
\fi\endgroup

\begingroup\makeatletter\ifx\SetFigFont\undefined
\def\x#1#2#3#4#5#6#7\relax{\def\x{#1#2#3#4#5#6}}%
\expandafter\x\fmtname xxxxxx\relax \def\y{splain}%
\ifx\x\y   % LaTeX or SliTeX?
\gdef\SetFigFont#1#2#3{%
  \ifnum #1<17\tiny\else \ifnum #1<20\small\else
  \ifnum #1<24\normalsize\else \ifnum #1<29\large\else
  \ifnum #1<34\Large\else \ifnum #1<41\LARGE\else
     \huge\fi\fi\fi\fi\fi\fi
  \csname #3\endcsname}%
\else
\gdef\SetFigFont#1#2#3{\begingroup
  \count@#1\relax \ifnum 25<\count@\count@25\fi
  \def\x{\endgroup\@setsize\SetFigFont{#2pt}}%
  \expandafter\x
    \csname \romannumeral\the\count@ pt\expandafter\endcsname
    \csname @\romannumeral\the\count@ pt\endcsname
  \csname #3\endcsname}%
\fi
\fi\endgroup

\begingroup\makeatletter\ifx\SetFigFont\undefined
\def\x#1#2#3#4#5#6#7\relax{\def\x{#1#2#3#4#5#6}}%
\expandafter\x\fmtname xxxxxx\relax \def\y{splain}%
\ifx\x\y   % LaTeX or SliTeX?
\gdef\SetFigFont#1#2#3{%
  \ifnum #1<17\tiny\else \ifnum #1<20\small\else
  \ifnum #1<24\normalsize\else \ifnum #1<29\large\else
  \ifnum #1<34\Large\else \ifnum #1<41\LARGE\else
     \huge\fi\fi\fi\fi\fi\fi
  \csname #3\endcsname}%
\else
\gdef\SetFigFont#1#2#3{\begingroup
  \count@#1\relax \ifnum 25<\count@\count@25\fi
  \def\x{\endgroup\@setsize\SetFigFont{#2pt}}%
  \expandafter\x
    \csname \romannumeral\the\count@ pt\expandafter\endcsname
    \csname @\romannumeral\the\count@ pt\endcsname
  \csname #3\endcsname}%
\fi
\fi\endgroup

\newcommand{\Ciii}[3]{
\begin{picture}(2752,895)(7025,-4136)
\thinlines
\put(7201,-3961){\circle*{336}}
\put(8401,-3961){\circle*{336}}
\put(9601,-3961){\circle{336}}
\put(7369,-3961){\line( 1, 0){864}}
\put(8569,-3911){\line( 1, 0){864}}
\put(8569,-4011){\line( 1, 0){864}}
\put(7141,-3436){\makebox(0,0)[lb]{\smash{\SetFigFont{6}{7.2}{rm}#1}}}
\put(8341,-3436){\makebox(0,0)[lb]{\smash{\SetFigFont{6}{7.2}{rm}#2}}}
\put(9541,-3436){\makebox(0,0)[lb]{\smash{\SetFigFont{6}{7.2}{rm}#3}}}
\end{picture}
}
\begingroup\makeatletter\ifx\SetFigFont\undefined
\def\x#1#2#3#4#5#6#7\relax{\def\x{#1#2#3#4#5#6}}%
\expandafter\x\fmtname xxxxxx\relax \def\y{splain}%
\ifx\x\y   % LaTeX or SliTeX?
\gdef\SetFigFont#1#2#3{%
  \ifnum #1<17\tiny\else \ifnum #1<20\small\else
  \ifnum #1<24\normalsize\else \ifnum #1<29\large\else
  \ifnum #1<34\Large\else \ifnum #1<41\LARGE\else
     \huge\fi\fi\fi\fi\fi\fi
  \csname #3\endcsname}%
\else
\gdef\SetFigFont#1#2#3{\begingroup
  \count@#1\relax \ifnum 25<\count@\count@25\fi
  \def\x{\endgroup\@setsize\SetFigFont{#2pt}}%
  \expandafter\x
    \csname \romannumeral\the\count@ pt\expandafter\endcsname
    \csname @\romannumeral\the\count@ pt\endcsname
  \csname #3\endcsname}%
\fi
\fi\endgroup

\newcommand{\Cii}[2]{
\begin{picture}(1552,895)(8225,-4136)
\thinlines
\put(8401,-3961){\circle*{336}}
\put(9601,-3961){\circle{336}}
\put(8569,-3911){\line( 1, 0){864}}
\put(8569,-4011){\line( 1, 0){864}}
\put(8341,-3436){\makebox(0,0)[lb]{\smash{\SetFigFont{6}{7.2}{rm}#1}}}
\put(9541,-3436){\makebox(0,0)[lb]{\smash{\SetFigFont{6}{7.2}{rm}#2}}}
\end{picture}
}

\begingroup\makeatletter\ifx\SetFigFont\undefined
\def\x#1#2#3#4#5#6#7\relax{\def\x{#1#2#3#4#5#6}}%
\expandafter\x\fmtname xxxxxx\relax \def\y{splain}%
\ifx\x\y   % LaTeX or SliTeX?
\gdef\SetFigFont#1#2#3{%
  \ifnum #1<17\tiny\else \ifnum #1<20\small\else
  \ifnum #1<24\normalsize\else \ifnum #1<29\large\else
  \ifnum #1<34\Large\else \ifnum #1<41\LARGE\else
     \huge\fi\fi\fi\fi\fi\fi
  \csname #3\endcsname}%
\else
\gdef\SetFigFont#1#2#3{\begingroup
  \count@#1\relax \ifnum 25<\count@\count@25\fi
  \def\x{\endgroup\@setsize\SetFigFont{#2pt}}%
  \expandafter\x
    \csname \romannumeral\the\count@ pt\expandafter\endcsname
    \csname @\romannumeral\the\count@ pt\endcsname
  \csname #3\endcsname}%
\fi
\fi\endgroup

\begingroup\makeatletter\ifx\SetFigFont\undefined
\def\x#1#2#3#4#5#6#7\relax{\def\x{#1#2#3#4#5#6}}%
\expandafter\x\fmtname xxxxxx\relax \def\y{splain}%
\ifx\x\y   % LaTeX or SliTeX?
\gdef\SetFigFont#1#2#3{%
  \ifnum #1<17\tiny\else \ifnum #1<20\small\else
  \ifnum #1<24\normalsize\else \ifnum #1<29\large\else
  \ifnum #1<34\Large\else \ifnum #1<41\LARGE\else
     \huge\fi\fi\fi\fi\fi\fi
  \csname #3\endcsname}%
\else
\gdef\SetFigFont#1#2#3{\begingroup
  \count@#1\relax \ifnum 25<\count@\count@25\fi
  \def\x{\endgroup\@setsize\SetFigFont{#2pt}}%
  \expandafter\x
    \csname \romannumeral\the\count@ pt\expandafter\endcsname
    \csname @\romannumeral\the\count@ pt\endcsname
  \csname #3\endcsname}%
\fi
\fi\endgroup

\newcommand{\Dviii}[8]{\parbox[c]{4.8cm}{$
\begin{picture}(7551,2200)(2125,-5050)
\thinlines
\put(2401,-3961){\circle{336}}
\put(3601,-3961){\circle{336}}
\put(4801,-3961){\circle{336}}
\put(6001,-3961){\circle{336}}
\put(7201,-3961){\circle{336}}
\put(8401,-3961){\circle{336}}
\put(9253,-3109){\circle{336}}
\put(9253,-4813){\circle{336}}
\put(2569,-3961){\line( 1, 0){870}}
\put(3796,-3961){\line( 1, 0){870}}
\put(4969,-3961){\line( 1, 0){870}}
\put(6169,-3961){\line( 1, 0){870}}
\put(7369,-3961){\line( 1, 0){870}}
\put(8519,-3843){\line( 1, 1){615}}
\put(8519,-4079){\line( 1, -1){615}}
\put(2341,-3436){\makebox(0,0)[lb]{\smash{\SetFigFont{6}{7.2}{rm}#1}}}
\put(3541,-3436){\makebox(0,0)[lb]{\smash{\SetFigFont{6}{7.2}{rm}#2}}}
\put(4741,-3436){\makebox(0,0)[lb]{\smash{\SetFigFont{6}{7.2}{rm}#3}}}
\put(5941,-3436){\makebox(0,0)[lb]{\smash{\SetFigFont{6}{7.2}{rm}#4}}}
\put(7141,-3451){\makebox(0,0)[lb]{\smash{\SetFigFont{6}{7.2}{rm}#5}}}
\put(8341,-3511){\makebox(0,0)[lb]{\smash{\SetFigFont{6}{7.2}{rm}#6}}}
\put(9676,-3211){\makebox(0,0)[lb]{\smash{\SetFigFont{6}{7.2}{rm}#7}}}
\put(9676,-4936){\makebox(0,0)[lb]{\smash{\SetFigFont{6}{7.2}{rm}#8}}}
\end{picture}$}
}
\begingroup\makeatletter\ifx\SetFigFont\undefined
\def\x#1#2#3#4#5#6#7\relax{\def\x{#1#2#3#4#5#6}}%
\expandafter\x\fmtname xxxxxx\relax \def\y{splain}%
\ifx\x\y   % LaTeX or SliTeX?
\gdef\SetFigFont#1#2#3{%
  \ifnum #1<17\tiny\else \ifnum #1<20\small\else
  \ifnum #1<24\normalsize\else \ifnum #1<29\large\else
  \ifnum #1<34\Large\else \ifnum #1<41\LARGE\else
     \huge\fi\fi\fi\fi\fi\fi
  \csname #3\endcsname}%
\else
\gdef\SetFigFont#1#2#3{\begingroup
  \count@#1\relax \ifnum 25<\count@\count@25\fi
  \def\x{\endgroup\@setsize\SetFigFont{#2pt}}%
  \expandafter\x
    \csname \romannumeral\the\count@ pt\expandafter\endcsname
    \csname @\romannumeral\the\count@ pt\endcsname
  \csname #3\endcsname}%
\fi
\fi\endgroup

\begingroup\makeatletter\ifx\SetFigFont\undefined
\def\x#1#2#3#4#5#6#7\relax{\def\x{#1#2#3#4#5#6}}%
\expandafter\x\fmtname xxxxxx\relax \def\y{splain}%
\ifx\x\y   % LaTeX or SliTeX?
\gdef\SetFigFont#1#2#3{%
  \ifnum #1<17\tiny\else \ifnum #1<20\small\else
  \ifnum #1<24\normalsize\else \ifnum #1<29\large\else
  \ifnum #1<34\Large\else \ifnum #1<41\LARGE\else
     \huge\fi\fi\fi\fi\fi\fi
  \csname #3\endcsname}%
\else
\gdef\SetFigFont#1#2#3{\begingroup
  \count@#1\relax \ifnum 25<\count@\count@25\fi
  \def\x{\endgroup\@setsize\SetFigFont{#2pt}}%
  \expandafter\x
    \csname \romannumeral\the\count@ pt\expandafter\endcsname
    \csname @\romannumeral\the\count@ pt\endcsname
  \csname #3\endcsname}%
\fi
\fi\endgroup

\newcommand{\Dvi}[6]{\parbox[c]{3.3cm}{$
\begin{picture}(5151,2200)(4525,-5050)
\thinlines
\put(4801,-3961){\circle{336}}
\put(6001,-3961){\circle{336}}
\put(7201,-3961){\circle{336}}
\put(8401,-3961){\circle{336}}
\put(9253,-3109){\circle{336}}
\put(9253,-4813){\circle{336}}
\put(4969,-3961){\line( 1, 0){870}}
\put(6169,-3961){\line( 1, 0){870}}
\put(7369,-3961){\line( 1, 0){870}}
\put(8519,-3843){\line( 1, 1){615}}
\put(8519,-4079){\line( 1, -1){615}}
\put(4741,-3436){\makebox(0,0)[lb]{\smash{\SetFigFont{6}{7.2}{rm}#1}}}
\put(5941,-3436){\makebox(0,0)[lb]{\smash{\SetFigFont{6}{7.2}{rm}#2}}}
\put(7141,-3436){\makebox(0,0)[lb]{\smash{\SetFigFont{6}{7.2}{rm}#3}}}
\put(8341,-3436){\makebox(0,0)[lb]{\smash{\SetFigFont{6}{7.2}{rm}#4}}}
\put(9676,-3211){\makebox(0,0)[lb]{\smash{\SetFigFont{6}{7.2}{rm}#5}}}
\put(9676,-4936){\makebox(0,0)[lb]{\smash{\SetFigFont{6}{7.2}{rm}#6}}}
\end{picture}$}
}
\begingroup\makeatletter\ifx\SetFigFont\undefined
\def\x#1#2#3#4#5#6#7\relax{\def\x{#1#2#3#4#5#6}}%
\expandafter\x\fmtname xxxxxx\relax \def\y{splain}%
\ifx\x\y   % LaTeX or SliTeX?
\gdef\SetFigFont#1#2#3{%
  \ifnum #1<17\tiny\else \ifnum #1<20\small\else
  \ifnum #1<24\normalsize\else \ifnum #1<29\large\else
  \ifnum #1<34\Large\else \ifnum #1<41\LARGE\else
     \huge\fi\fi\fi\fi\fi\fi
  \csname #3\endcsname}%
\else
\gdef\SetFigFont#1#2#3{\begingroup
  \count@#1\relax \ifnum 25<\count@\count@25\fi
  \def\x{\endgroup\@setsize\SetFigFont{#2pt}}%
  \expandafter\x
    \csname \romannumeral\the\count@ pt\expandafter\endcsname
    \csname @\romannumeral\the\count@ pt\endcsname
  \csname #3\endcsname}%
\fi
\fi\endgroup

\newcommand{\Dv}[5]{\parbox[c]{2.8cm}{$
\begin{picture}(3951,2200)(5725,-5050)
\thinlines
\put(6001,-3961){\circle{336}}
\put(7201,-3961){\circle{336}}
\put(8401,-3961){\circle{336}}
\put(9253,-3109){\circle{336}}
\put(9253,-4813){\circle{336}}
\put(6169,-3961){\line( 1, 0){870}}
\put(7369,-3961){\line( 1, 0){870}}
\put(8519,-3843){\line( 1, 1){615}}
\put(8519,-4079){\line( 1, -1){615}}
\put(5941,-3436){\makebox(0,0)[lb]{\smash{\SetFigFont{6}{7.2}{rm}#1}}}
\put(7141,-3451){\makebox(0,0)[lb]{\smash{\SetFigFont{6}{7.2}{rm}#2}}}
\put(8341,-3511){\makebox(0,0)[lb]{\smash{\SetFigFont{6}{7.2}{rm}#3}}}
\put(9676,-3211){\makebox(0,0)[lb]{\smash{\SetFigFont{6}{7.2}{rm}#4}}}
\put(9676,-4936){\makebox(0,0)[lb]{\smash{\SetFigFont{6}{7.2}{rm}#5}}}
\end{picture}$}
}
\begingroup\makeatletter\ifx\SetFigFont\undefined
\def\x#1#2#3#4#5#6#7\relax{\def\x{#1#2#3#4#5#6}}%
\expandafter\x\fmtname xxxxxx\relax \def\y{splain}%
\ifx\x\y   % LaTeX or SliTeX?
\gdef\SetFigFont#1#2#3{%
  \ifnum #1<17\tiny\else \ifnum #1<20\small\else
  \ifnum #1<24\normalsize\else \ifnum #1<29\large\else
  \ifnum #1<34\Large\else \ifnum #1<41\LARGE\else
     \huge\fi\fi\fi\fi\fi\fi
  \csname #3\endcsname}%
\else
\gdef\SetFigFont#1#2#3{\begingroup
  \count@#1\relax \ifnum 25<\count@\count@25\fi
  \def\x{\endgroup\@setsize\SetFigFont{#2pt}}%
  \expandafter\x
    \csname \romannumeral\the\count@ pt\expandafter\endcsname
    \csname @\romannumeral\the\count@ pt\endcsname
  \csname #3\endcsname}%
\fi
\fi\endgroup

\newcommand{\Div}[4]{\parbox[c]{1.8cm}{$
\begin{picture}(2751,2200)(6925,-5050)
\thinlines
\put(7201,-3961){\circle{336}}
\put(8401,-3961){\circle{336}}
\put(9253,-3109){\circle{336}}
\put(9253,-4813){\circle{336}}
\put(7369,-3961){\line( 1, 0){870}}
\put(8519,-3843){\line( 1, 1){615}}
\put(8519,-4079){\line( 1, -1){615}}
\put(7141,-3451){\makebox(0,0)[lb]{\smash{\SetFigFont{6}{7.2}{rm}#1}}}
\put(8341,-3511){\makebox(0,0)[lb]{\smash{\SetFigFont{6}{7.2}{rm}#2}}}
\put(9676,-3211){\makebox(0,0)[lb]{\smash{\SetFigFont{6}{7.2}{rm}#3}}}
\put(9676,-4936){\makebox(0,0)[lb]{\smash{\SetFigFont{6}{7.2}{rm}#4}}}
\end{picture}$}
}

\begingroup\makeatletter\ifx\SetFigFont\undefined
\def\x#1#2#3#4#5#6#7\relax{\def\x{#1#2#3#4#5#6}}%
\expandafter\x\fmtname xxxxxx\relax \def\y{splain}%
\ifx\x\y   % LaTeX or SliTeX?
\gdef\SetFigFont#1#2#3{%
  \ifnum #1<17\tiny\else \ifnum #1<20\small\else
  \ifnum #1<24\normalsize\else \ifnum #1<29\large\else
  \ifnum #1<34\Large\else \ifnum #1<41\LARGE\else
     \huge\fi\fi\fi\fi\fi\fi
  \csname #3\endcsname}%
\else
\gdef\SetFigFont#1#2#3{\begingroup
  \count@#1\relax \ifnum 25<\count@\count@25\fi
  \def\x{\endgroup\@setsize\SetFigFont{#2pt}}%
  \expandafter\x
    \csname \romannumeral\the\count@ pt\expandafter\endcsname
    \csname @\romannumeral\the\count@ pt\endcsname
  \csname #3\endcsname}%
\fi
\fi\endgroup

\newcommand{\Fiv}[4]{
\begin{picture}(3952,895)(1025,-4136)
\thinlines
\put(1201,-3961){\circle{336}}
\put(2401,-3961){\circle{336}}
\put(3601,-3961){\circle*{336}}
\put(4801,-3961){\circle*{336}}
\put(1369,-3961){\line( 1, 0){864}}
\put(2569,-3911){\line( 1, 0){864}}
\put(2569,-4011){\line( 1, 0){864}}
\put(3769,-3961){\line( 1, 0){864}}
\put(1141,-3436){\makebox(0,0)[lb]{\smash{\SetFigFont{6}{7.2}{rm}\mbox{$#1$}}}}
\put(2341,-3436){\makebox(0,0)[lb]{\smash{\SetFigFont{6}{7.2}{rm}\mbox{$#2$}}}}
\put(3541,-3436){\makebox(0,0)[lb]{\smash{\SetFigFont{6}{7.2}{rm}\mbox{$#3$}}}}
\put(4741,-3436){\makebox(0,0)[lb]{\smash{\SetFigFont{6}{7.2}{rm}\mbox{$#4$}}}}
\end{picture}
}
\begingroup\makeatletter\ifx\SetFigFont\undefined
\def\x#1#2#3#4#5#6#7\relax{\def\x{#1#2#3#4#5#6}}%
\expandafter\x\fmtname xxxxxx\relax \def\y{splain}%
\ifx\x\y   % LaTeX or SliTeX?
\gdef\SetFigFont#1#2#3{%
  \ifnum #1<17\tiny\else \ifnum #1<20\small\else
  \ifnum #1<24\normalsize\else \ifnum #1<29\large\else
  \ifnum #1<34\Large\else \ifnum #1<41\LARGE\else
     \huge\fi\fi\fi\fi\fi\fi
  \csname #3\endcsname}%
\else
\gdef\SetFigFont#1#2#3{\begingroup
  \count@#1\relax \ifnum 25<\count@\count@25\fi
  \def\x{\endgroup\@setsize\SetFigFont{#2pt}}%
  \expandafter\x
    \csname \romannumeral\the\count@ pt\expandafter\endcsname
    \csname @\romannumeral\the\count@ pt\endcsname
  \csname #3\endcsname}%
\fi
\fi\endgroup

\newcommand{\Gii}[2]{
\begin{picture}(1552,895)(1025,-4136)
\thinlines
\put(1201,-3961){\circle*{336}}
\put(2401,-3961){\circle{336}}
\put(1339,-3861){\line( 1, 0){924}}
\put(1369,-3961){\line( 1, 0){864}}
\put(1339,-4061){\line( 1, 0){924}}
\put(1141,-3436){\makebox(0,0)[lb]{\smash{\SetFigFont{6}{7.2}{rm}\mbox{$#1$}}}}
\put(2341,-3436){\makebox(0,0)[lb]{\smash{\SetFigFont{6}{7.2}{rm}\mbox{$#2$}}}}
\end{picture}
}
\begingroup\makeatletter\ifx\SetFigFont\undefined
\def\x#1#2#3#4#5#6#7\relax{\def\x{#1#2#3#4#5#6}}%
\expandafter\x\fmtname xxxxxx\relax \def\y{splain}%
\ifx\x\y   % LaTeX or SliTeX?
\gdef\SetFigFont#1#2#3{%
  \ifnum #1<17\tiny\else \ifnum #1<20\small\else
  \ifnum #1<24\normalsize\else \ifnum #1<29\large\else
  \ifnum #1<34\Large\else \ifnum #1<41\LARGE\else
     \huge\fi\fi\fi\fi\fi\fi
  \csname #3\endcsname}%
\else
\gdef\SetFigFont#1#2#3{\begingroup
  \count@#1\relax \ifnum 25<\count@\count@25\fi
  \def\x{\endgroup\@setsize\SetFigFont{#2pt}}%
  \expandafter\x
    \csname \romannumeral\the\count@ pt\expandafter\endcsname
    \csname @\romannumeral\the\count@ pt\endcsname
  \csname #3\endcsname}%
\fi
\fi\endgroup

\newcommand{\Evi}[6]{\raisebox{2.4mm}{\parbox[c]{3.4cm}{$
\begin{picture}(5300,2200)(950,-5400)
\thinlines
\put(1201,-3961){\circle{336}}
\put(2401,-3961){\circle{336}}
\put(3601,-3961){\circle{336}}
\put(4801,-3961){\circle{336}}
\put(6001,-3961){\circle{336}}
\put(3601,-5161){\circle{336}}
\put(1369,-3961){\line( 1, 0){864}}
\put(2569,-3961){\line( 1, 0){864}}
\put(3769,-3961){\line( 1, 0){864}}
\put(4969,-3961){\line( 1, 0){864}}
\put(3601,-4129){\line( 0, -1){864}}
\put(1141,-3436){\makebox(0,0)[lb]{\smash{\SetFigFont{6}{7.2}{rm}\mbox{$#1$}}}}
\put(2341,-3436){\makebox(0,0)[lb]{\smash{\SetFigFont{6}{7.2}{rm}\mbox{$#3$}}}}
\put(3541,-3436){\makebox(0,0)[lb]{\smash{\SetFigFont{6}{7.2}{rm}\mbox{$#4$}}}}
\put(4741,-3436){\makebox(0,0)[lb]{\smash{\SetFigFont{6}{7.2}{rm}\mbox{$#5$}}}}
\put(5941,-3436){\makebox(0,0)[lb]{\smash{\SetFigFont{6}{7.2}{rm}\mbox{$#6$}}}}
\put(4024,-5263){\makebox(0,0)[lb]{\smash{\SetFigFont{6}{7.2}{rm}\mbox{$#2$}}}}
\end{picture}$}}
}

\begingroup\makeatletter\ifx\SetFigFont\undefined
\def\x#1#2#3#4#5#6#7\relax{\def\x{#1#2#3#4#5#6}}%
\expandafter\x\fmtname xxxxxx\relax \def\y{splain}%
\ifx\x\y   % LaTeX or SliTeX?
\gdef\SetFigFont#1#2#3{%
  \ifnum #1<17\tiny\else \ifnum #1<20\small\else
  \ifnum #1<24\normalsize\else \ifnum #1<29\large\else
  \ifnum #1<34\Large\else \ifnum #1<41\LARGE\else
     \huge\fi\fi\fi\fi\fi\fi
  \csname #3\endcsname}%
\else
\gdef\SetFigFont#1#2#3{\begingroup
  \count@#1\relax \ifnum 25<\count@\count@25\fi
  \def\x{\endgroup\@setsize\SetFigFont{#2pt}}%
  \expandafter\x
    \csname \romannumeral\the\count@ pt\expandafter\endcsname
    \csname @\romannumeral\the\count@ pt\endcsname
  \csname #3\endcsname}%
\fi
\fi\endgroup

\newcommand{\Evii}[7]{\raisebox{2.4mm}{\parbox[c]{4.2cm}{$
\begin{picture}(6500,2200)(950,-5400)
\thinlines
\put(1201,-3961){\circle{336}}
\put(2401,-3961){\circle{336}}
\put(3601,-3961){\circle{336}}
\put(4801,-3961){\circle{336}}
\put(6001,-3961){\circle{336}}
\put(7201,-3961){\circle{336}}
\put(3601,-5161){\circle{336}}
\put(1369,-3961){\line( 1, 0){864}}
\put(2569,-3961){\line( 1, 0){864}}
\put(3769,-3961){\line( 1, 0){864}}
\put(4969,-3961){\line( 1, 0){864}}
\put(6169,-3961){\line( 1, 0){864}}
\put(3601,-4129){\line( 0, -1){864}}
\put(1141,-3436){\makebox(0,0)[lb]{\smash{\SetFigFont{6}{7.2}{rm}\mbox{$#1$}}}}
\put(2341,-3436){\makebox(0,0)[lb]{\smash{\SetFigFont{6}{7.2}{rm}\mbox{$#3$}}}}
\put(3541,-3436){\makebox(0,0)[lb]{\smash{\SetFigFont{6}{7.2}{rm}\mbox{$#4$}}}}
\put(4741,-3436){\makebox(0,0)[lb]{\smash{\SetFigFont{6}{7.2}{rm}\mbox{$#5$}}}}
\put(5941,-3436){\makebox(0,0)[lb]{\smash{\SetFigFont{6}{7.2}{rm}\mbox{$#6$}}}}
\put(7141,-3436){\makebox(0,0)[lb]{\smash{\SetFigFont{6}{7.2}{rm}\mbox{$#7$}}}}
\put(4024,-5263){\makebox(0,0)[lb]{\smash{\SetFigFont{6}{7.2}{rm}\mbox{$#2$}}}}
\end{picture}$}}
}
\begingroup\makeatletter\ifx\SetFigFont\undefined
\def\x#1#2#3#4#5#6#7\relax{\def\x{#1#2#3#4#5#6}}%
\expandafter\x\fmtname xxxxxx\relax \def\y{splain}%
\ifx\x\y   % LaTeX or SliTeX?
\gdef\SetFigFont#1#2#3{%
  \ifnum #1<17\tiny\else \ifnum #1<20\small\else
  \ifnum #1<24\normalsize\else \ifnum #1<29\large\else
  \ifnum #1<34\Large\else \ifnum #1<41\LARGE\else
     \huge\fi\fi\fi\fi\fi\fi
  \csname #3\endcsname}%
\else
\gdef\SetFigFont#1#2#3{\begingroup
  \count@#1\relax \ifnum 25<\count@\count@25\fi
  \def\x{\endgroup\@setsize\SetFigFont{#2pt}}%
  \expandafter\x
    \csname \romannumeral\the\count@ pt\expandafter\endcsname
    \csname @\romannumeral\the\count@ pt\endcsname
  \csname #3\endcsname}%
\fi
\fi\endgroup

\newcommand{\Eviii}[8]{\raisebox{2.4mm}{\parbox[c]{5cm}{$
\begin{picture}(7700,2200)(950,-5400)
\thinlines
\put(1201,-3961){\circle{336}}
\put(2401,-3961){\circle{336}}
\put(3601,-3961){\circle{336}}
\put(4801,-3961){\circle{336}}
\put(6001,-3961){\circle{336}}
\put(7201,-3961){\circle{336}}
\put(8401,-3961){\circle{336}}
\put(3601,-5161){\circle{336}}
\put(1369,-3961){\line( 1, 0){864}}
\put(2569,-3961){\line( 1, 0){864}}
\put(3769,-3961){\line( 1, 0){864}}
\put(4969,-3961){\line( 1, 0){864}}
\put(6169,-3961){\line( 1, 0){864}}
\put(7369,-3961){\line( 1, 0){864}}
\put(3601,-4129){\line( 0, -1){864}}
\put(1141,-3436){\makebox(0,0)[lb]{\smash{\SetFigFont{6}{7.2}{rm}\mbox{$#1$}}}}
\put(2341,-3436){\makebox(0,0)[lb]{\smash{\SetFigFont{6}{7.2}{rm}\mbox{$#3$}}}}
\put(3541,-3436){\makebox(0,0)[lb]{\smash{\SetFigFont{6}{7.2}{rm}\mbox{$#4$}}}}
\put(4741,-3436){\makebox(0,0)[lb]{\smash{\SetFigFont{6}{7.2}{rm}\mbox{$#5$}}}}
\put(5941,-3436){\makebox(0,0)[lb]{\smash{\SetFigFont{6}{7.2}{rm}\mbox{$#6$}}}}
\put(7141,-3436){\makebox(0,0)[lb]{\smash{\SetFigFont{6}{7.2}{rm}\mbox{$#7$}}}}
\put(8341,-3436){\makebox(0,0)[lb]{\smash{\SetFigFont{6}{7.2}{rm}\mbox{$#8$}}}}
\put(4024,-5263){\makebox(0,0)[lb]{\smash{\SetFigFont{6}{7.2}{rm}\mbox{$#2$}}}}
\end{picture}$}}
}

\begingroup\makeatletter\ifx\SetFigFont\undefined
\def\x#1#2#3#4#5#6#7\relax{\def\x{#1#2#3#4#5#6}}%
\expandafter\x\fmtname xxxxxx\relax \def\y{splain}%
\ifx\x\y   % LaTeX or SliTeX?
\gdef\SetFigFont#1#2#3{%
  \ifnum #1<17\tiny\else \ifnum #1<20\small\else
  \ifnum #1<24\normalsize\else \ifnum #1<29\large\else
  \ifnum #1<34\Large\else \ifnum #1<41\LARGE\else
     \huge\fi\fi\fi\fi\fi\fi
  \csname #3\endcsname}%
\else
\gdef\SetFigFont#1#2#3{\begingroup
  \count@#1\relax \ifnum 25<\count@\count@25\fi
  \def\x{\endgroup\@setsize\SetFigFont{#2pt}}%
  \expandafter\x
    \csname \romannumeral\the\count@ pt\expandafter\endcsname
    \csname @\romannumeral\the\count@ pt\endcsname
  \csname #3\endcsname}%
\fi
\fi\endgroup

\newcommand{\An}[3]{
\begin{picture}(3952,895)(1025,-4136)
\thinlines
\put(1201,-3961){\circle{336}}
\put(2401,-3961){\circle{336}}
\put(4771,-3961){\circle{336}}
\put(3436,-3961){\circle*{30}}
\put(3586,-3961){\circle*{30}}
\put(3736,-3961){\circle*{30}}
\put(1369,-3961){\line( 1, 0){870}}
\put(2569,-3961){\line( 1, 0){600}}
\put(4003,-3961){\line( 1, 0){600}}
\put(1141,-3436){\makebox(0,0)[lb]{\smash{\SetFigFont{6}{7.2}{rm}\mbox{$#1$}}}}
\put(2341,-3436){\makebox(0,0)[lb]{\smash{\SetFigFont{6}{7.2}{rm}\mbox{$#2$}}}}
\put(4741,-3436){\makebox(0,0)[lb]{\smash{\SetFigFont{6}{7.2}{rm}\mbox{$#3$}}}}
\end{picture}
}
\begingroup\makeatletter\ifx\SetFigFont\undefined
\def\x#1#2#3#4#5#6#7\relax{\def\x{#1#2#3#4#5#6}}%
\expandafter\x\fmtname xxxxxx\relax \def\y{splain}%
\ifx\x\y   % LaTeX or SliTeX?
\gdef\SetFigFont#1#2#3{%
  \ifnum #1<17\tiny\else \ifnum #1<20\small\else
  \ifnum #1<24\normalsize\else \ifnum #1<29\large\else
  \ifnum #1<34\Large\else \ifnum #1<41\LARGE\else
     \huge\fi\fi\fi\fi\fi\fi
  \csname #3\endcsname}%
\else
\gdef\SetFigFont#1#2#3{\begingroup
  \count@#1\relax \ifnum 25<\count@\count@25\fi
  \def\x{\endgroup\@setsize\SetFigFont{#2pt}}%
  \expandafter\x
    \csname \romannumeral\the\count@ pt\expandafter\endcsname
    \csname @\romannumeral\the\count@ pt\endcsname
  \csname #3\endcsname}%
\fi
\fi\endgroup
\newcommand{\Anl}[3]{
\begin{picture}(3952,895)(1025,-4136)
\thinlines
\put(4801,-3961){\circle{336}}
\put(3601,-3961){\circle{336}}
\put(1201,-3961){\circle{336}}
\put(2536,-3961){\circle*{30}}
\put(2386,-3961){\circle*{30}}
\put(2236,-3961){\circle*{30}}
\put(3769,-3961){\line( 1, 0){870}}
\put(2803,-3961){\line( 1, 0){600}}
\put(1369,-3961){\line( 1, 0){600}}
\put(4741,-3436){\makebox(0,0)[lb]{\smash{\SetFigFont{6}{7.2}{rm}\mbox{$#3$}}}}
\put(3541,-3436){\makebox(0,0)[lb]{\smash{\SetFigFont{6}{7.2}{rm}\mbox{$#2$}}}}
\put(1141,-3436){\makebox(0,0)[lb]{\smash{\SetFigFont{6}{7.2}{rm}\mbox{$#1$}}}}
\end{picture}
}
\begingroup\makeatletter\ifx\SetFigFont\undefined
\def\x#1#2#3#4#5#6#7\relax{\def\x{#1#2#3#4#5#6}}%
\expandafter\x\fmtname xxxxxx\relax \def\y{splain}%
\ifx\x\y   % LaTeX or SliTeX?
\gdef\SetFigFont#1#2#3{%
  \ifnum #1<17\tiny\else \ifnum #1<20\small\else
  \ifnum #1<24\normalsize\else \ifnum #1<29\large\else
  \ifnum #1<34\Large\else \ifnum #1<41\LARGE\else
     \huge\fi\fi\fi\fi\fi\fi
  \csname #3\endcsname}%
\else
\gdef\SetFigFont#1#2#3{\begingroup
  \count@#1\relax \ifnum 25<\count@\count@25\fi
  \def\x{\endgroup\@setsize\SetFigFont{#2pt}}%
  \expandafter\x
    \csname \romannumeral\the\count@ pt\expandafter\endcsname
    \csname @\romannumeral\the\count@ pt\endcsname
  \csname #3\endcsname}%
\fi
\fi\endgroup

\newcommand{\Ans}[2]{
\begin{picture}(2752,895)(1025,-4136)
\thinlines
\put(1201,-3961){\circle{336}}
\put(3601,-3961){\circle{336}}
\put(2236,-3961){\circle*{30}}
\put(2386,-3961){\circle*{30}}
\put(2536,-3961){\circle*{30}}
\put(1369,-3961){\line( 1, 0){600}}
\put(2803,-3961){\line( 1, 0){600}}
\put(1141,-3436){\makebox(0,0)[lb]{\smash{\SetFigFont{6}{7.2}{rm}\mbox{$#1$}}}}
\put(3541,-3436){\makebox(0,0)[lb]{\smash{\SetFigFont{6}{7.2}{rm}\mbox{$#2$}}}}
\end{picture}
}

\begingroup\makeatletter\ifx\SetFigFont\undefined
\def\x#1#2#3#4#5#6#7\relax{\def\x{#1#2#3#4#5#6}}%
\expandafter\x\fmtname xxxxxx\relax \def\y{splain}%
\ifx\x\y   % LaTeX or SliTeX?
\gdef\SetFigFont#1#2#3{%
  \ifnum #1<17\tiny\else \ifnum #1<20\small\else
  \ifnum #1<24\normalsize\else \ifnum #1<29\large\else
  \ifnum #1<34\Large\else \ifnum #1<41\LARGE\else
     \huge\fi\fi\fi\fi\fi\fi
  \csname #3\endcsname}%
\else
\gdef\SetFigFont#1#2#3{\begingroup
  \count@#1\relax \ifnum 25<\count@\count@25\fi
  \def\x{\endgroup\@setsize\SetFigFont{#2pt}}%
  \expandafter\x
    \csname \romannumeral\the\count@ pt\expandafter\endcsname
    \csname @\romannumeral\the\count@ pt\endcsname
  \csname #3\endcsname}%
\fi
\fi\endgroup

\newcommand{\Ani}[4]{
\begin{picture}(5152,895)(1025,-4136)
\thinlines
\put(1201,-3961){\circle{336}}
\put(2401,-3961){\circle{336}}
\put(4771,-3961){\circle{336}}
\put(5971,-3961){\circle{336}}
\put(3436,-3961){\circle*{30}}
\put(3586,-3961){\circle*{30}}
\put(3736,-3961){\circle*{30}}
\put(1369,-3961){\line( 1, 0){870}}
\put(2569,-3961){\line( 1, 0){600}}
\put(4003,-3961){\line( 1, 0){600}}
\put(4969,-3961){\line( 1, 0){870}}
\put(1141,-3436){\makebox(0,0)[lb]{\smash{\SetFigFont{6}{7.2}{rm}\mbox{$#1$}}}}
\put(2341,-3436){\makebox(0,0)[lb]{\smash{\SetFigFont{6}{7.2}{rm}\mbox{$#2$}}}}
\put(4741,-3436){\makebox(0,0)[lb]{\smash{\SetFigFont{6}{7.2}{rm}\mbox{$#3$}}}}
\put(5941,-3436){\makebox(0,0)[lb]{\smash{\SetFigFont{6}{7.2}{rm}\mbox{$#4$}}}}
\end{picture}
}
\begingroup\makeatletter\ifx\SetFigFont\undefined
\def\x#1#2#3#4#5#6#7\relax{\def\x{#1#2#3#4#5#6}}%
\expandafter\x\fmtname xxxxxx\relax \def\y{splain}%
\ifx\x\y   % LaTeX or SliTeX?
\gdef\SetFigFont#1#2#3{%
  \ifnum #1<17\tiny\else \ifnum #1<20\small\else
  \ifnum #1<24\normalsize\else \ifnum #1<29\large\else
  \ifnum #1<34\Large\else \ifnum #1<41\LARGE\else
     \huge\fi\fi\fi\fi\fi\fi
  \csname #3\endcsname}%
\else
\gdef\SetFigFont#1#2#3{\begingroup
  \count@#1\relax \ifnum 25<\count@\count@25\fi
  \def\x{\endgroup\@setsize\SetFigFont{#2pt}}%
  \expandafter\x
    \csname \romannumeral\the\count@ pt\expandafter\endcsname
    \csname @\romannumeral\the\count@ pt\endcsname
  \csname #3\endcsname}%
\fi
\fi\endgroup

\newcommand{\Anii}[7]{
\begin{picture}(9952,895)(1025,-4136)
\thinlines
\put(1201,-3961){\circle{336}}
\put(2401,-3961){\circle{336}}
\put(3436,-3961){\circle*{30}}
\put(3586,-3961){\circle*{30}}
\put(3736,-3961){\circle*{30}}
\put(4771,-3961){\circle{336}}
\put(5971,-3961){\circle{336}}
\put(7171,-3961){\circle{336}}
\put(8206,-3961){\circle*{30}}
\put(8356,-3961){\circle*{30}}
\put(8506,-3961){\circle*{30}}
\put(9571,-3961){\circle{336}}
\put(10771,-3961){\circle{336}}

\put(1369,-3961){\line( 1, 0){864}}
\put(2569,-3961){\line( 1, 0){600}}
\put(4003,-3961){\line( 1, 0){600}}
\put(4969,-3961){\line( 1, 0){864}}
\put(6169,-3961){\line( 1, 0){864}}
\put(7369,-3961){\line( 1, 0){600}}
\put(8803,-3961){\line( 1, 0){600}}
\put(9769,-3961){\line( 1, 0){864}}

\put(1141,-3436){\makebox(0,0)[lb]{\smash{\SetFigFont{6}{7.2}{rm}\mbox{$#1$}}}}
\put(2341,-3436){\makebox(0,0)[lb]{\smash{\SetFigFont{6}{7.2}{rm}\mbox{$#2$}}}}
\put(4741,-3436){\makebox(0,0)[lb]{\smash{\SetFigFont{6}{7.2}{rm}\mbox{$#3$}}}}
\put(5941,-3436){\makebox(0,0)[lb]{\smash{\SetFigFont{6}{7.2}{rm}\mbox{$#4$}}}}
\put(7141,-3436){\makebox(0,0)[lb]{\smash{\SetFigFont{6}{7.2}{rm}\mbox{$#5$}}}}
\put(9541,-3436){\makebox(0,0)[lb]{\smash{\SetFigFont{6}{7.2}{rm}\mbox{$#6$}}}}
\put(10741,-3436){\makebox(0,0)[lb]{\smash{\SetFigFont{6}{7.2}{rm}\mbox{$#7$}}}}
\end{picture}
}
\begingroup\makeatletter\ifx\SetFigFont\undefined
\def\x#1#2#3#4#5#6#7\relax{\def\x{#1#2#3#4#5#6}}%
\expandafter\x\fmtname xxxxxx\relax \def\y{splain}%
\ifx\x\y   % LaTeX or SliTeX?
\gdef\SetFigFont#1#2#3{%
  \ifnum #1<17\tiny\else \ifnum #1<20\small\else
  \ifnum #1<24\normalsize\else \ifnum #1<29\large\else
  \ifnum #1<34\Large\else \ifnum #1<41\LARGE\else
     \huge\fi\fi\fi\fi\fi\fi
  \csname #3\endcsname}%
\else
\gdef\SetFigFont#1#2#3{\begingroup
  \count@#1\relax \ifnum 25<\count@\count@25\fi
  \def\x{\endgroup\@setsize\SetFigFont{#2pt}}%
  \expandafter\x
    \csname \romannumeral\the\count@ pt\expandafter\endcsname
    \csname @\romannumeral\the\count@ pt\endcsname
  \csname #3\endcsname}%
\fi
\fi\endgroup

\newcommand{\Aniii}[8]{
\begin{picture}(11152,895)(1025,-4136)
\thinlines
\put(1201,-3961){\circle{336}}
\put(2401,-3961){\circle{336}}
\put(3436,-3961){\circle*{30}}
\put(3586,-3961){\circle*{30}}
\put(3736,-3961){\circle*{30}}
\put(4771,-3961){\circle{336}}
\put(5971,-3961){\circle{336}}
\put(7171,-3961){\circle{336}}
\put(8371,-3961){\circle{336}}
\put(9406,-3961){\circle*{30}}
\put(9556,-3961){\circle*{30}}
\put(9706,-3961){\circle*{30}}
\put(10771,-3961){\circle{336}}
\put(11971,-3961){\circle{336}}

\put(1369,-3961){\line( 1, 0){864}}
\put(2569,-3961){\line( 1, 0){600}}
\put(4003,-3961){\line( 1, 0){600}}
\put(4969,-3961){\line( 1, 0){864}}
\put(6169,-3961){\line( 1, 0){864}}
\put(7369,-3961){\line( 1, 0){864}}
\put(8569,-3961){\line( 1, 0){600}}
\put(10003,-3961){\line( 1, 0){600}}
\put(10969,-3961){\line( 1, 0){864}}

\put(1141,-3436){\makebox(0,0)[lb]{\smash{\SetFigFont{6}{7.2}{rm}\mbox{$#1$}}}}
\put(2341,-3436){\makebox(0,0)[lb]{\smash{\SetFigFont{6}{7.2}{rm}\mbox{$#2$}}}}
\put(4741,-3436){\makebox(0,0)[lb]{\smash{\SetFigFont{6}{7.2}{rm}\mbox{$#3$}}}}
\put(5941,-3436){\makebox(0,0)[lb]{\smash{\SetFigFont{6}{7.2}{rm}\mbox{$#4$}}}}
\put(7141,-3436){\makebox(0,0)[lb]{\smash{\SetFigFont{6}{7.2}{rm}\mbox{$#5$}}}}
\put(8341,-3436){\makebox(0,0)[lb]{\smash{\SetFigFont{6}{7.2}{rm}\mbox{$#6$}}}}
\put(10741,-3436){\makebox(0,0)[lb]{\smash{\SetFigFont{6}{7.2}{rm}\mbox{$#7$}}}}
\put(11941,-3436){\makebox(0,0)[lb]{\smash{\SetFigFont{6}{7.2}{rm}\mbox{$#8$}}}}
\end{picture}
}
\begingroup\makeatletter\ifx\SetFigFont\undefined
\def\x#1#2#3#4#5#6#7\relax{\def\x{#1#2#3#4#5#6}}%
\expandafter\x\fmtname xxxxxx\relax \def\y{splain}%
\ifx\x\y   % LaTeX or SliTeX?
\gdef\SetFigFont#1#2#3{%
  \ifnum #1<17\tiny\else \ifnum #1<20\small\else
  \ifnum #1<24\normalsize\else \ifnum #1<29\large\else
  \ifnum #1<34\Large\else \ifnum #1<41\LARGE\else
     \huge\fi\fi\fi\fi\fi\fi
  \csname #3\endcsname}%
\else
\gdef\SetFigFont#1#2#3{\begingroup
  \count@#1\relax \ifnum 25<\count@\count@25\fi
  \def\x{\endgroup\@setsize\SetFigFont{#2pt}}%
  \expandafter\x
    \csname \romannumeral\the\count@ pt\expandafter\endcsname
    \csname @\romannumeral\the\count@ pt\endcsname
  \csname #3\endcsname}%
\fi
\fi\endgroup

\newcommand{\Bn}[4]{
\begin{picture}(5152,895)(1025,-4136)
\thinlines
\put(1201,-3961){\circle{336}}
\put(2401,-3961){\circle{336}}
\put(4801,-3961){\circle{336}}
\put(6001,-3961){\circle*{336}}
\put(3436,-3946){\circle*{30}}
\put(3586,-3946){\circle*{30}}
\put(3736,-3946){\circle*{30}}
\put(1369,-3961){\line( 1, 0){864}}
\put(2569,-3961){\line( 1, 0){600}}
\put(4003,-3961){\line( 1, 0){600}}
\put(4969,-3911){\line( 1, 0){864}}
\put(4969,-4011){\line( 1, 0){864}}
\put(1141,-3436){\makebox(0,0)[lb]{\smash{\SetFigFont{6}{7.2}{rm}\mbox{$#1$}}}}
\put(2341,-3436){\makebox(0,0)[lb]{\smash{\SetFigFont{6}{7.2}{rm}\mbox{$#2$}}}}
\put(4741,-3436){\makebox(0,0)[lb]{\smash{\SetFigFont{6}{7.2}{rm}\mbox{$#3$}}}}
\put(5941,-3436){\makebox(0,0)[lb]{\smash{\SetFigFont{6}{7.2}{rm}\mbox{$#4$}}}}
\end{picture}
}
\begingroup\makeatletter\ifx\SetFigFont\undefined
\def\x#1#2#3#4#5#6#7\relax{\def\x{#1#2#3#4#5#6}}%
\expandafter\x\fmtname xxxxxx\relax \def\y{splain}%
\ifx\x\y   % LaTeX or SliTeX?
\gdef\SetFigFont#1#2#3{%
  \ifnum #1<17\tiny\else \ifnum #1<20\small\else
  \ifnum #1<24\normalsize\else \ifnum #1<29\large\else
  \ifnum #1<34\Large\else \ifnum #1<41\LARGE\else
     \huge\fi\fi\fi\fi\fi\fi
  \csname #3\endcsname}%
\else
\gdef\SetFigFont#1#2#3{\begingroup
  \count@#1\relax \ifnum 25<\count@\count@25\fi
  \def\x{\endgroup\@setsize\SetFigFont{#2pt}}%
  \expandafter\x
    \csname \romannumeral\the\count@ pt\expandafter\endcsname
    \csname @\romannumeral\the\count@ pt\endcsname
  \csname #3\endcsname}%
\fi
\fi\endgroup

\newcommand{\Bns}[3]{
\begin{picture}(3952,895)(1025,-4136)
\thinlines
\put(1201,-3961){\circle{336}}
\put(3601,-3961){\circle{336}}
\put(4801,-3961){\circle*{336}}
\put(2236,-3946){\circle*{30}}
\put(2386,-3946){\circle*{30}}
\put(2536,-3946){\circle*{30}}
\put(1369,-3961){\line( 1, 0){600}}
\put(2803,-3961){\line( 1, 0){600}}
\put(3769,-3911){\line( 1, 0){864}}
\put(3769,-4011){\line( 1, 0){864}}
\put(1141,-3436){\makebox(0,0)[lb]{\smash{\SetFigFont{6}{7.2}{rm}\mbox{$#1$}}}}
\put(3541,-3436){\makebox(0,0)[lb]{\smash{\SetFigFont{6}{7.2}{rm}\mbox{$#2$}}}}
\put(4741,-3436){\makebox(0,0)[lb]{\smash{\SetFigFont{6}{7.2}{rm}\mbox{$#3$}}}}
\end{picture}
}
\begingroup\makeatletter\ifx\SetFigFont\undefined
\def\x#1#2#3#4#5#6#7\relax{\def\x{#1#2#3#4#5#6}}%
\expandafter\x\fmtname xxxxxx\relax \def\y{splain}%
\ifx\x\y   % LaTeX or SliTeX?
\gdef\SetFigFont#1#2#3{%
  \ifnum #1<17\tiny\else \ifnum #1<20\small\else
  \ifnum #1<24\normalsize\else \ifnum #1<29\large\else
  \ifnum #1<34\Large\else \ifnum #1<41\LARGE\else
     \huge\fi\fi\fi\fi\fi\fi
  \csname #3\endcsname}%
\else
\gdef\SetFigFont#1#2#3{\begingroup
  \count@#1\relax \ifnum 25<\count@\count@25\fi
  \def\x{\endgroup\@setsize\SetFigFont{#2pt}}%
  \expandafter\x
    \csname \romannumeral\the\count@ pt\expandafter\endcsname
    \csname @\romannumeral\the\count@ pt\endcsname
  \csname #3\endcsname}%
\fi
\fi\endgroup

\newcommand{\Bni}[5]{
\begin{picture}(6352,895)(1025,-4136)
\thinlines
\put(1201,-3961){\circle{336}}
\put(2401,-3961){\circle{336}}
\put(3601,-3961){\circle{336}}
\put(6001,-3961){\circle{336}}
\put(7201,-3961){\circle*{336}}

\put(1369,-3961){\line( 1, 0){864}}
\put(2569,-3961){\line( 1, 0){864}}
\put(3769,-3961){\line( 1, 0){600}}
\put(5203,-3961){\line( 1, 0){600}}
\put(6169,-3911){\line( 1, 0){864}}
\put(6169,-4011){\line( 1, 0){864}}

\put(4636,-3946){\circle*{30}}
\put(4786,-3946){\circle*{30}}
\put(4936,-3946){\circle*{30}}

\put(1141,-3436){\makebox(0,0)[lb]{\smash{\SetFigFont{6}{7.2}{rm}\mbox{$#1$}}}}
\put(2341,-3436){\makebox(0,0)[lb]{\smash{\SetFigFont{6}{7.2}{rm}\mbox{$#2$}}}}
\put(3541,-3436){\makebox(0,0)[lb]{\smash{\SetFigFont{6}{7.2}{rm}\mbox{$#3$}}}}
\put(5941,-3436){\makebox(0,0)[lb]{\smash{\SetFigFont{6}{7.2}{rm}\mbox{$#4$}}}}
\put(7141,-3436){\makebox(0,0)[lb]{\smash{\SetFigFont{6}{7.2}{rm}\mbox{$#5$}}}}

\end{picture}
}
\begingroup\makeatletter\ifx\SetFigFont\undefined
\def\x#1#2#3#4#5#6#7\relax{\def\x{#1#2#3#4#5#6}}%
\expandafter\x\fmtname xxxxxx\relax \def\y{splain}%
\ifx\x\y   % LaTeX or SliTeX?
\gdef\SetFigFont#1#2#3{%
  \ifnum #1<17\tiny\else \ifnum #1<20\small\else
  \ifnum #1<24\normalsize\else \ifnum #1<29\large\else
  \ifnum #1<34\Large\else \ifnum #1<41\LARGE\else
     \huge\fi\fi\fi\fi\fi\fi
  \csname #3\endcsname}%
\else
\gdef\SetFigFont#1#2#3{\begingroup
  \count@#1\relax \ifnum 25<\count@\count@25\fi
  \def\x{\endgroup\@setsize\SetFigFont{#2pt}}%
  \expandafter\x
    \csname \romannumeral\the\count@ pt\expandafter\endcsname
    \csname @\romannumeral\the\count@ pt\endcsname
  \csname #3\endcsname}%
\fi
\fi\endgroup

\newcommand{\Bnii}[6]{
\begin{picture}(7552,895)(1025,-4136)
\thinlines
\put(1201,-3961){\circle{336}}
\put(2401,-3961){\circle{336}}
\put(3601,-3961){\circle{336}}
\put(4801,-3961){\circle{336}}
\put(7201,-3961){\circle{336}}
\put(8401,-3961){\circle*{336}}

\put(1369,-3961){\line( 1, 0){864}}
\put(2569,-3961){\line( 1, 0){864}}
\put(3769,-3961){\line( 1, 0){864}}

\put(4903,-3961){\line( 1, 0){600}}
\put(6403,-3961){\line( 1, 0){600}}
\put(7369,-3911){\line( 1, 0){864}}
\put(7369,-4011){\line( 1, 0){864}}

\put(5836,-3946){\circle*{30}}
\put(5986,-3946){\circle*{30}}
\put(6136,-3946){\circle*{30}}

\put(1141,-3436){\makebox(0,0)[lb]{\smash{\SetFigFont{6}{7.2}{rm}\mbox{$#1$}}}}
\put(2341,-3436){\makebox(0,0)[lb]{\smash{\SetFigFont{6}{7.2}{rm}\mbox{$#2$}}}}
\put(3541,-3436){\makebox(0,0)[lb]{\smash{\SetFigFont{6}{7.2}{rm}\mbox{$#3$}}}}
\put(4741,-3436){\makebox(0,0)[lb]{\smash{\SetFigFont{6}{7.2}{rm}\mbox{$#4$}}}}
\put(7141,-3436){\makebox(0,0)[lb]{\smash{\SetFigFont{6}{7.2}{rm}\mbox{$#5$}}}}
\put(8341,-3436){\makebox(0,0)[lb]{\smash{\SetFigFont{6}{7.2}{rm}\mbox{$#6$}}}}

\end{picture}
}
\begingroup\makeatletter\ifx\SetFigFont\undefined
\def\x#1#2#3#4#5#6#7\relax{\def\x{#1#2#3#4#5#6}}%
\expandafter\x\fmtname xxxxxx\relax \def\y{splain}%
\ifx\x\y   % LaTeX or SliTeX?
\gdef\SetFigFont#1#2#3{%
  \ifnum #1<17\tiny\else \ifnum #1<20\small\else
  \ifnum #1<24\normalsize\else \ifnum #1<29\large\else
  \ifnum #1<34\Large\else \ifnum #1<41\LARGE\else
     \huge\fi\fi\fi\fi\fi\fi
  \csname #3\endcsname}%
\else
\gdef\SetFigFont#1#2#3{\begingroup
  \count@#1\relax \ifnum 25<\count@\count@25\fi
  \def\x{\endgroup\@setsize\SetFigFont{#2pt}}%
  \expandafter\x
    \csname \romannumeral\the\count@ pt\expandafter\endcsname
    \csname @\romannumeral\the\count@ pt\endcsname
  \csname #3\endcsname}%
\fi
\fi\endgroup

\newcommand{\Bniii}[7]{
\begin{picture}(8752,895)(1025,-4136)
\thinlines
\put(1201,-3961){\circle{336}}
\put(2401,-3961){\circle{336}}
\put(3601,-3961){\circle{336}}
\put(4801,-3961){\circle{336}}
\put(7201,-3961){\circle{336}}
\put(8401,-3961){\circle{336}}
\put(9601,-3961){\circle*{336}}

\put(1369,-3961){\line( 1, 0){864}}
\put(2569,-3961){\line( 1, 0){864}}
\put(3769,-3961){\line( 1, 0){864}}

\put(4903,-3961){\line( 1, 0){600}}
\put(6403,-3961){\line( 1, 0){600}}

\put(7369,-3961){\line( 1, 0){864}}

\put(8569,-3911){\line( 1, 0){864}}
\put(8569,-4011){\line( 1, 0){864}}

\put(5836,-3946){\circle*{30}}
\put(5986,-3946){\circle*{30}}
\put(6136,-3946){\circle*{30}}

\put(1141,-3436){\makebox(0,0)[lb]{\smash{\SetFigFont{6}{7.2}{rm}\mbox{$#1$}}}}
\put(2341,-3436){\makebox(0,0)[lb]{\smash{\SetFigFont{6}{7.2}{rm}\mbox{$#2$}}}}
\put(3541,-3436){\makebox(0,0)[lb]{\smash{\SetFigFont{6}{7.2}{rm}\mbox{$#3$}}}}
\put(4741,-3436){\makebox(0,0)[lb]{\smash{\SetFigFont{6}{7.2}{rm}\mbox{$#4$}}}}
\put(7141,-3436){\makebox(0,0)[lb]{\smash{\SetFigFont{6}{7.2}{rm}\mbox{$#5$}}}}
\put(8341,-3436){\makebox(0,0)[lb]{\smash{\SetFigFont{6}{7.2}{rm}\mbox{$#6$}}}}
\put(9541,-3436){\makebox(0,0)[lb]{\smash{\SetFigFont{6}{7.2}{rm}\mbox{$#7$}}}}

\end{picture}
}
\begingroup\makeatletter\ifx\SetFigFont\undefined
\def\x#1#2#3#4#5#6#7\relax{\def\x{#1#2#3#4#5#6}}%
\expandafter\x\fmtname xxxxxx\relax \def\y{splain}%
\ifx\x\y   % LaTeX or SliTeX?
\gdef\SetFigFont#1#2#3{%
  \ifnum #1<17\tiny\else \ifnum #1<20\small\else
  \ifnum #1<24\normalsize\else \ifnum #1<29\large\else
  \ifnum #1<34\Large\else \ifnum #1<41\LARGE\else
     \huge\fi\fi\fi\fi\fi\fi
  \csname #3\endcsname}%
\else
\gdef\SetFigFont#1#2#3{\begingroup
  \count@#1\relax \ifnum 25<\count@\count@25\fi
  \def\x{\endgroup\@setsize\SetFigFont{#2pt}}%
  \expandafter\x
    \csname \romannumeral\the\count@ pt\expandafter\endcsname
    \csname @\romannumeral\the\count@ pt\endcsname
  \csname #3\endcsname}%
\fi
\fi\endgroup

\newcommand{\Bniv}[8]{
\begin{picture}(9952,895)(1025,-4136)
\thinlines
\put(1201,-3961){\circle{336}}
\put(2401,-3961){\circle{336}}
\put(3601,-3961){\circle{336}}
\put(4801,-3961){\circle{336}}
\put(7201,-3961){\circle{336}}
\put(8401,-3961){\circle{336}}
\put(9601,-3961){\circle{336}}
\put(10801,-3961){\circle*{336}}

\put(1369,-3961){\line( 1, 0){864}}
\put(2569,-3961){\line( 1, 0){864}}
\put(3769,-3961){\line( 1, 0){864}}

\put(4903,-3961){\line( 1, 0){600}}
\put(6403,-3961){\line( 1, 0){600}}

\put(7369,-3961){\line( 1, 0){864}}

\put(8569,-3961){\line( 1, 0){864}}

\put(9769,-3911){\line( 1, 0){864}}
\put(9769,-4011){\line( 1, 0){864}}

\put(5836,-3946){\circle*{30}}
\put(5986,-3946){\circle*{30}}
\put(6136,-3946){\circle*{30}}

\put(1141,-3436){\makebox(0,0)[lb]{\smash{\SetFigFont{6}{7.2}{rm}\mbox{$#1$}}}}
\put(2341,-3436){\makebox(0,0)[lb]{\smash{\SetFigFont{6}{7.2}{rm}\mbox{$#2$}}}}
\put(3541,-3436){\makebox(0,0)[lb]{\smash{\SetFigFont{6}{7.2}{rm}\mbox{$#3$}}}}
\put(4741,-3436){\makebox(0,0)[lb]{\smash{\SetFigFont{6}{7.2}{rm}\mbox{$#4$}}}}
\put(7141,-3436){\makebox(0,0)[lb]{\smash{\SetFigFont{6}{7.2}{rm}\mbox{$#5$}}}}
\put(8341,-3436){\makebox(0,0)[lb]{\smash{\SetFigFont{6}{7.2}{rm}\mbox{$#6$}}}}
\put(9541,-3436){\makebox(0,0)[lb]{\smash{\SetFigFont{6}{7.2}{rm}\mbox{$#7$}}}}
\put(10741,-3436){\makebox(0,0)[lb]{\smash{\SetFigFont{6}{7.2}{rm}\mbox{$#8$}}}}

\end{picture}
}
\begingroup\makeatletter\ifx\SetFigFont\undefined
\def\x#1#2#3#4#5#6#7\relax{\def\x{#1#2#3#4#5#6}}%
\expandafter\x\fmtname xxxxxx\relax \def\y{splain}%
\ifx\x\y   % LaTeX or SliTeX?
\gdef\SetFigFont#1#2#3{%
  \ifnum #1<17\tiny\else \ifnum #1<20\small\else
  \ifnum #1<24\normalsize\else \ifnum #1<29\large\else
  \ifnum #1<34\Large\else \ifnum #1<41\LARGE\else
     \huge\fi\fi\fi\fi\fi\fi
  \csname #3\endcsname}%
\else
\gdef\SetFigFont#1#2#3{\begingroup
  \count@#1\relax \ifnum 25<\count@\count@25\fi
  \def\x{\endgroup\@setsize\SetFigFont{#2pt}}%
  \expandafter\x
    \csname \romannumeral\the\count@ pt\expandafter\endcsname
    \csname @\romannumeral\the\count@ pt\endcsname
  \csname #3\endcsname}%
\fi
\fi\endgroup

\newcommand{\Cn}[4]{
\begin{picture}(5152,895)(1025,-4136)
\thinlines
\put(1201,-3961){\circle*{336}}
\put(2401,-3961){\circle*{336}}
\put(4801,-3961){\circle*{336}}
\put(6001,-3961){\circle{336}}
\put(3436,-3946){\circle*{30}}
\put(3586,-3946){\circle*{30}}
\put(3736,-3946){\circle*{30}}
\put(1369,-3961){\line( 1, 0){864}}
\put(2569,-3961){\line( 1, 0){600}}
\put(4003,-3961){\line( 1, 0){600}}
\put(4969,-3911){\line( 1, 0){864}}
\put(4969,-4011){\line( 1, 0){864}}
\put(1141,-3436){\makebox(0,0)[lb]{\smash{\SetFigFont{6}{7.2}{rm}\mbox{$#1$}}}}
\put(2341,-3436){\makebox(0,0)[lb]{\smash{\SetFigFont{6}{7.2}{rm}\mbox{$#2$}}}}
\put(4741,-3436){\makebox(0,0)[lb]{\smash{\SetFigFont{6}{7.2}{rm}\mbox{$#3$}}}}
\put(5941,-3436){\makebox(0,0)[lb]{\smash{\SetFigFont{6}{7.2}{rm}\mbox{$#4$}}}}
\end{picture}
}
\begingroup\makeatletter\ifx\SetFigFont\undefined
\def\x#1#2#3#4#5#6#7\relax{\def\x{#1#2#3#4#5#6}}%
\expandafter\x\fmtname xxxxxx\relax \def\y{splain}%
\ifx\x\y   % LaTeX or SliTeX?
\gdef\SetFigFont#1#2#3{%
  \ifnum #1<17\tiny\else \ifnum #1<20\small\else
  \ifnum #1<24\normalsize\else \ifnum #1<29\large\else
  \ifnum #1<34\Large\else \ifnum #1<41\LARGE\else
     \huge\fi\fi\fi\fi\fi\fi
  \csname #3\endcsname}%
\else
\gdef\SetFigFont#1#2#3{\begingroup
  \count@#1\relax \ifnum 25<\count@\count@25\fi
  \def\x{\endgroup\@setsize\SetFigFont{#2pt}}%
  \expandafter\x
    \csname \romannumeral\the\count@ pt\expandafter\endcsname
    \csname @\romannumeral\the\count@ pt\endcsname
  \csname #3\endcsname}%
\fi
\fi\endgroup

\newcommand{\Cns}[3]{
\begin{picture}(3952,895)(1025,-4136)
\thinlines
\put(1201,-3961){\circle*{336}}
\put(3601,-3961){\circle*{336}}
\put(4801,-3961){\circle{336}}
\put(2236,-3946){\circle*{30}}
\put(2386,-3946){\circle*{30}}
\put(2536,-3946){\circle*{30}}
\put(1369,-3961){\line( 1, 0){600}}
\put(2803,-3961){\line( 1, 0){600}}
\put(3769,-3911){\line( 1, 0){864}}
\put(3769,-4011){\line( 1, 0){864}}
\put(1141,-3436){\makebox(0,0)[lb]{\smash{\SetFigFont{6}{7.2}{rm}\mbox{$#1$}}}}
\put(3541,-3436){\makebox(0,0)[lb]{\smash{\SetFigFont{6}{7.2}{rm}\mbox{$#2$}}}}
\put(4741,-3436){\makebox(0,0)[lb]{\smash{\SetFigFont{6}{7.2}{rm}\mbox{$#3$}}}}
\end{picture}
}
\begingroup\makeatletter\ifx\SetFigFont\undefined
\def\x#1#2#3#4#5#6#7\relax{\def\x{#1#2#3#4#5#6}}%
\expandafter\x\fmtname xxxxxx\relax \def\y{splain}%
\ifx\x\y   % LaTeX or SliTeX?
\gdef\SetFigFont#1#2#3{%
  \ifnum #1<17\tiny\else \ifnum #1<20\small\else
  \ifnum #1<24\normalsize\else \ifnum #1<29\large\else
  \ifnum #1<34\Large\else \ifnum #1<41\LARGE\else
     \huge\fi\fi\fi\fi\fi\fi
  \csname #3\endcsname}%
\else
\gdef\SetFigFont#1#2#3{\begingroup
  \count@#1\relax \ifnum 25<\count@\count@25\fi
  \def\x{\endgroup\@setsize\SetFigFont{#2pt}}%
  \expandafter\x
    \csname \romannumeral\the\count@ pt\expandafter\endcsname
    \csname @\romannumeral\the\count@ pt\endcsname
  \csname #3\endcsname}%
\fi
\fi\endgroup

\newcommand{\Cni}[5]{
\begin{picture}(6352,895)(1025,-4136)
\thinlines
\put(1201,-3961){\circle*{336}}
\put(2401,-3961){\circle*{336}}
\put(3601,-3961){\circle*{336}}
\put(6001,-3961){\circle*{336}}
\put(7201,-3961){\circle{336}}

\put(1369,-3961){\line( 1, 0){864}}
\put(2569,-3961){\line( 1, 0){864}}
\put(3769,-3961){\line( 1, 0){600}}
\put(5203,-3961){\line( 1, 0){600}}
\put(6169,-3911){\line( 1, 0){864}}
\put(6169,-4011){\line( 1, 0){864}}

\put(4636,-3946){\circle*{30}}
\put(4786,-3946){\circle*{30}}
\put(4936,-3946){\circle*{30}}

\put(1141,-3436){\makebox(0,0)[lb]{\smash{\SetFigFont{6}{7.2}{rm}#1}}}
\put(2341,-3436){\makebox(0,0)[lb]{\smash{\SetFigFont{6}{7.2}{rm}#2}}}
\put(3541,-3436){\makebox(0,0)[lb]{\smash{\SetFigFont{6}{7.2}{rm}#3}}}
\put(5941,-3436){\makebox(0,0)[lb]{\smash{\SetFigFont{6}{7.2}{rm}#4}}}
\put(7141,-3436){\makebox(0,0)[lb]{\smash{\SetFigFont{6}{7.2}{rm}#5}}}

\end{picture}
}
\begingroup\makeatletter\ifx\SetFigFont\undefined
\def\x#1#2#3#4#5#6#7\relax{\def\x{#1#2#3#4#5#6}}%
\expandafter\x\fmtname xxxxxx\relax \def\y{splain}%
\ifx\x\y   % LaTeX or SliTeX?
\gdef\SetFigFont#1#2#3{%
  \ifnum #1<17\tiny\else \ifnum #1<20\small\else
  \ifnum #1<24\normalsize\else \ifnum #1<29\large\else
  \ifnum #1<34\Large\else \ifnum #1<41\LARGE\else
     \huge\fi\fi\fi\fi\fi\fi
  \csname #3\endcsname}%
\else
\gdef\SetFigFont#1#2#3{\begingroup
  \count@#1\relax \ifnum 25<\count@\count@25\fi
  \def\x{\endgroup\@setsize\SetFigFont{#2pt}}%
  \expandafter\x
    \csname \romannumeral\the\count@ pt\expandafter\endcsname
    \csname @\romannumeral\the\count@ pt\endcsname
  \csname #3\endcsname}%
\fi
\fi\endgroup

\newcommand{\Cnii}[6]{
\begin{picture}(7552,895)(1025,-4136)
\thinlines
\put(1201,-3961){\circle*{336}}
\put(2401,-3961){\circle*{336}}
\put(3601,-3961){\circle*{336}}
\put(4801,-3961){\circle*{336}}
\put(7201,-3961){\circle*{336}}
\put(8401,-3961){\circle{336}}

\put(1369,-3961){\line( 1, 0){864}}
\put(2569,-3961){\line( 1, 0){864}}
\put(3769,-3961){\line( 1, 0){864}}

\put(4903,-3961){\line( 1, 0){600}}
\put(6403,-3961){\line( 1, 0){600}}
\put(7369,-3911){\line( 1, 0){864}}
\put(7369,-4011){\line( 1, 0){864}}

\put(5836,-3946){\circle*{30}}
\put(5986,-3946){\circle*{30}}
\put(6136,-3946){\circle*{30}}

\put(1141,-3436){\makebox(0,0)[lb]{\smash{\SetFigFont{6}{7.2}{rm}#1}}}
\put(2341,-3436){\makebox(0,0)[lb]{\smash{\SetFigFont{6}{7.2}{rm}#2}}}
\put(3541,-3436){\makebox(0,0)[lb]{\smash{\SetFigFont{6}{7.2}{rm}#3}}}
\put(4741,-3436){\makebox(0,0)[lb]{\smash{\SetFigFont{6}{7.2}{rm}#4}}}
\put(7141,-3436){\makebox(0,0)[lb]{\smash{\SetFigFont{6}{7.2}{rm}#5}}}
\put(8341,-3436){\makebox(0,0)[lb]{\smash{\SetFigFont{6}{7.2}{rm}#6}}}

\end{picture}
}
\begingroup\makeatletter\ifx\SetFigFont\undefined
\def\x#1#2#3#4#5#6#7\relax{\def\x{#1#2#3#4#5#6}}%
\expandafter\x\fmtname xxxxxx\relax \def\y{splain}%
\ifx\x\y   % LaTeX or SliTeX?
\gdef\SetFigFont#1#2#3{%
  \ifnum #1<17\tiny\else \ifnum #1<20\small\else
  \ifnum #1<24\normalsize\else \ifnum #1<29\large\else
  \ifnum #1<34\Large\else \ifnum #1<41\LARGE\else
     \huge\fi\fi\fi\fi\fi\fi
  \csname #3\endcsname}%
\else
\gdef\SetFigFont#1#2#3{\begingroup
  \count@#1\relax \ifnum 25<\count@\count@25\fi
  \def\x{\endgroup\@setsize\SetFigFont{#2pt}}%
  \expandafter\x
    \csname \romannumeral\the\count@ pt\expandafter\endcsname
    \csname @\romannumeral\the\count@ pt\endcsname
  \csname #3\endcsname}%
\fi
\fi\endgroup

\newcommand{\Dn}[5]{\parbox[c]{3.4cm}{$
\begin{picture}(5151,2200)(4525,-5050)
\thinlines
\put(4801,-3961){\circle{336}}
\put(6001,-3961){\circle{336}}
\put(8401,-3961){\circle{336}}
\put(9253,-3109){\circle{336}}
\put(9253,-4813){\circle{336}}
\put(7036,-3946){\circle*{30}}
\put(7186,-3946){\circle*{30}}
\put(7336,-3946){\circle*{30}}
\put(4969,-3961){\line( 1, 0){864}}
\put(6169,-3961){\line( 1, 0){600}}
\put(7603,-3961){\line( 1, 0){600}}
\put(8519,-3843){\line( 1, 1){615}}
\put(8519,-4079){\line( 1, -1){615}}
\put(4741,-3436){\makebox(0,0)[lb]{\smash{\SetFigFont{6}{7.2}{rm}\mbox{$#1$}}}}
\put(5941,-3436){\makebox(0,0)[lb]{\smash{\SetFigFont{6}{7.2}{rm}\mbox{$#2$}}}}
\put(8341,-3436){\makebox(0,0)[lb]{\smash{\SetFigFont{6}{7.2}{rm}\mbox{$#3$}}}}
\put(9676,-3211){\makebox(0,0)[lb]{\smash{\SetFigFont{6}{7.2}{rm}\mbox{$#4$}}}}
\put(9676,-4936){\makebox(0,0)[lb]{\smash{\SetFigFont{6}{7.2}{rm}\mbox{$#5$}}}}
\end{picture}$}
}
\begingroup\makeatletter\ifx\SetFigFont\undefined
\def\x#1#2#3#4#5#6#7\relax{\def\x{#1#2#3#4#5#6}}%
\expandafter\x\fmtname xxxxxx\relax \def\y{splain}%
\ifx\x\y   % LaTeX or SliTeX?
\gdef\SetFigFont#1#2#3{%
  \ifnum #1<17\tiny\else \ifnum #1<20\small\else
  \ifnum #1<24\normalsize\else \ifnum #1<29\large\else
  \ifnum #1<34\Large\else \ifnum #1<41\LARGE\else
     \huge\fi\fi\fi\fi\fi\fi
  \csname #3\endcsname}%
\else
\gdef\SetFigFont#1#2#3{\begingroup
  \count@#1\relax \ifnum 25<\count@\count@25\fi
  \def\x{\endgroup\@setsize\SetFigFont{#2pt}}%
  \expandafter\x
    \csname \romannumeral\the\count@ pt\expandafter\endcsname
    \csname @\romannumeral\the\count@ pt\endcsname
  \csname #3\endcsname}%
\fi
\fi\endgroup

\newcommand{\Dns}[4]{\parbox[c]{2.5cm}{$
\begin{picture}(3951,2200)(4525,-5050)
\thinlines
\put(4801,-3961){\circle{336}}
\put(7201,-3961){\circle{336}}
\put(8053,-3109){\circle{336}}
\put(8053,-4813){\circle{336}}
\put(5836,-3946){\circle*{30}}
\put(5986,-3946){\circle*{30}}
\put(6136,-3946){\circle*{30}}
\put(4969,-3961){\line( 1, 0){600}}
\put(6403,-3961){\line( 1, 0){600}}
\put(7319,-3843){\line( 1, 1){615}}
\put(7319,-4079){\line( 1, -1){615}}
\put(4741,-3436){\makebox(0,0)[lb]{\smash{\SetFigFont{6}{7.2}{rm}\mbox{$#1$}}}}
\put(7141,-3436){\makebox(0,0)[lb]{\smash{\SetFigFont{6}{7.2}{rm}\mbox{$#2$}}}}
\put(8476,-3211){\makebox(0,0)[lb]{\smash{\SetFigFont{6}{7.2}{rm}\mbox{$#3$}}}}
\put(8476,-4936){\makebox(0,0)[lb]{\smash{\SetFigFont{6}{7.2}{rm}\mbox{$#4$}}}}
\end{picture}$}
}
\begingroup\makeatletter\ifx\SetFigFont\undefined
\def\x#1#2#3#4#5#6#7\relax{\def\x{#1#2#3#4#5#6}}%
\expandafter\x\fmtname xxxxxx\relax \def\y{splain}%
\ifx\x\y   % LaTeX or SliTeX?
\gdef\SetFigFont#1#2#3{%
  \ifnum #1<17\tiny\else \ifnum #1<20\small\else
  \ifnum #1<24\normalsize\else \ifnum #1<29\large\else
  \ifnum #1<34\Large\else \ifnum #1<41\LARGE\else
     \huge\fi\fi\fi\fi\fi\fi
  \csname #3\endcsname}%
\else
\gdef\SetFigFont#1#2#3{\begingroup
  \count@#1\relax \ifnum 25<\count@\count@25\fi
  \def\x{\endgroup\@setsize\SetFigFont{#2pt}}%
  \expandafter\x
    \csname \romannumeral\the\count@ pt\expandafter\endcsname
    \csname @\romannumeral\the\count@ pt\endcsname
  \csname #3\endcsname}%
\fi
\fi\endgroup

\newcommand{\Dni}[6]{\parbox[c]{4.1cm}{$
\begin{picture}(6351,2200)(4525,-5050)
\thinlines
\put(4801,-3961){\circle{336}}
\put(6001,-3961){\circle{336}}
\put(7201,-3961){\circle{336}}

\put(9601,-3961){\circle{336}}
\put(10453,-3109){\circle{336}}
\put(10453,-4813){\circle{336}}
\put(8236,-3946){\circle*{30}}
\put(8386,-3946){\circle*{30}}
\put(8536,-3946){\circle*{30}}
\put(4969,-3961){\line( 1, 0){864}}
\put(6169,-3961){\line( 1, 0){864}}
\put(7369,-3961){\line( 1, 0){600}}
\put(8803,-3961){\line( 1, 0){600}}
\put(9719,-3843){\line( 1, 1){615}}
\put(9719,-4079){\line( 1, -1){615}}
\put(4741,-3436){\makebox(0,0)[lb]{\smash{\SetFigFont{6}{7.2}{rm}#1}}}
\put(5941,-3436){\makebox(0,0)[lb]{\smash{\SetFigFont{6}{7.2}{rm}#2}}}
\put(7141,-3436){\makebox(0,0)[lb]{\smash{\SetFigFont{6}{7.2}{rm}#3}}}
\put(9541,-3436){\makebox(0,0)[lb]{\smash{\SetFigFont{6}{7.2}{rm}#4}}}
\put(10876,-3211){\makebox(0,0)[lb]{\smash{\SetFigFont{6}{7.2}{rm}#5}}}
\put(10876,-4936){\makebox(0,0)[lb]{\smash{\SetFigFont{6}{7.2}{rm}#6}}}
\end{picture}$}
}
\begingroup\makeatletter\ifx\SetFigFont\undefined
\def\x#1#2#3#4#5#6#7\relax{\def\x{#1#2#3#4#5#6}}%
\expandafter\x\fmtname xxxxxx\relax \def\y{splain}%
\ifx\x\y   % LaTeX or SliTeX?
\gdef\SetFigFont#1#2#3{%
  \ifnum #1<17\tiny\else \ifnum #1<20\small\else
  \ifnum #1<24\normalsize\else \ifnum #1<29\large\else
  \ifnum #1<34\Large\else \ifnum #1<41\LARGE\else
     \huge\fi\fi\fi\fi\fi\fi
  \csname #3\endcsname}%
\else
\gdef\SetFigFont#1#2#3{\begingroup
  \count@#1\relax \ifnum 25<\count@\count@25\fi
  \def\x{\endgroup\@setsize\SetFigFont{#2pt}}%
  \expandafter\x
    \csname \romannumeral\the\count@ pt\expandafter\endcsname
    \csname @\romannumeral\the\count@ pt\endcsname
  \csname #3\endcsname}%
\fi
\fi\endgroup

\newcommand{\Dnii}[7]{\parbox[c]{4.8cm}{$
\begin{picture}(7551,2200)(4525,-5050)
\thinlines
\put(4801,-3961){\circle{336}}
\put(6001,-3961){\circle{336}}
\put(7201,-3961){\circle{336}}
\put(8401,-3961){\circle{336}}
\put(10801,-3961){\circle{336}}
\put(11653,-3109){\circle{336}}
\put(11653,-4813){\circle{336}}
\put(9436,-3946){\circle*{30}}
\put(9586,-3946){\circle*{30}}
\put(9736,-3946){\circle*{30}}
\put(4969,-3961){\line( 1, 0){864}}
\put(6169,-3961){\line( 1, 0){864}}
\put(7369,-3961){\line( 1, 0){864}}
\put(8569,-3961){\line( 1, 0){600}}
\put(10003,-3961){\line( 1, 0){600}}
\put(10919,-3843){\line( 1, 1){615}}
\put(10919,-4079){\line( 1, -1){615}}
\put(4741,-3436){\makebox(0,0)[lb]{\smash{\SetFigFont{6}{7.2}{rm}#1}}}
\put(5941,-3436){\makebox(0,0)[lb]{\smash{\SetFigFont{6}{7.2}{rm}#2}}}
\put(7141,-3436){\makebox(0,0)[lb]{\smash{\SetFigFont{6}{7.2}{rm}#3}}}
\put(8341,-3436){\makebox(0,0)[lb]{\smash{\SetFigFont{6}{7.2}{rm}#4}}}
\put(10741,-3436){\makebox(0,0)[lb]{\smash{\SetFigFont{6}{7.2}{rm}#5}}}
\put(12076,-3211){\makebox(0,0)[lb]{\smash{\SetFigFont{6}{7.2}{rm}#6}}}
\put(12076,-4936){\makebox(0,0)[lb]{\smash{\SetFigFont{6}{7.2}{rm}#7}}}
\end{picture}$}
}
\begingroup\makeatletter\ifx\SetFigFont\undefined
\def\x#1#2#3#4#5#6#7\relax{\def\x{#1#2#3#4#5#6}}%
\expandafter\x\fmtname xxxxxx\relax \def\y{splain}%
\ifx\x\y   % LaTeX or SliTeX?
\gdef\SetFigFont#1#2#3{%
  \ifnum #1<17\tiny\else \ifnum #1<20\small\else
  \ifnum #1<24\normalsize\else \ifnum #1<29\large\else
  \ifnum #1<34\Large\else \ifnum #1<41\LARGE\else
     \huge\fi\fi\fi\fi\fi\fi
  \csname #3\endcsname}%
\else
\gdef\SetFigFont#1#2#3{\begingroup
  \count@#1\relax \ifnum 25<\count@\count@25\fi
  \def\x{\endgroup\@setsize\SetFigFont{#2pt}}%
  \expandafter\x
    \csname \romannumeral\the\count@ pt\expandafter\endcsname
    \csname @\romannumeral\the\count@ pt\endcsname
  \csname #3\endcsname}%
\fi
\fi\endgroup

\newcommand{\Dniii}[8]{\parbox[c]{5.55cm}{$
\begin{picture}(8751,2200)(1125,-5050)
\thinlines
\put(1201,-3961){\circle{336}}
\put(2401,-3961){\circle{336}}
\put(3601,-3961){\circle{336}}
\put(4801,-3961){\circle{336}}

\put(5836,-3946){\circle*{30}}
\put(5986,-3946){\circle*{30}}
\put(6136,-3946){\circle*{30}}

\put(7201,-3961){\circle{336}}
\put(8401,-3961){\circle{336}}

\put(9253,-3109){\circle{336}}
\put(9253,-4813){\circle{336}}

\put(1369,-3961){\line( 1, 0){870}}
\put(2569,-3961){\line( 1, 0){870}}
\put(3796,-3961){\line( 1, 0){870}}

\put(4969,-3961){\line( 1, 0){600}}
\put(6403,-3961){\line( 1, 0){600}}

\put(7369,-3961){\line( 1, 0){870}}
\put(8519,-3843){\line( 1, 1){615}}
\put(8519,-4079){\line( 1, -1){615}}

\put(1141,-3436){\makebox(0,0)[lb]{\smash{\SetFigFont{6}{7.2}{rm}\mbox{$#1$}}}}
\put(2341,-3436){\makebox(0,0)[lb]{\smash{\SetFigFont{6}{7.2}{rm}\mbox{$#2$}}}}
\put(3541,-3436){\makebox(0,0)[lb]{\smash{\SetFigFont{6}{7.2}{rm}\mbox{$#3$}}}}
\put(4741,-3436){\makebox(0,0)[lb]{\smash{\SetFigFont{6}{7.2}{rm}\mbox{$#4$}}}}
\put(7141,-3436){\makebox(0,0)[lb]{\smash{\SetFigFont{6}{7.2}{rm}\mbox{$#5$}}}}
\put(8341,-3436){\makebox(0,0)[lb]{\smash{\SetFigFont{6}{7.2}{rm}\mbox{$#6$}}}}
\put(9676,-3211){\makebox(0,0)[lb]{\smash{\SetFigFont{6}{7.2}{rm}\mbox{$#7$}}}}
\put(9676,-4936){\makebox(0,0)[lb]{\smash{\SetFigFont{6}{7.2}{rm}\mbox{$#8$}}}}
\end{picture}$}
}
\begingroup\makeatletter\ifx\SetFigFont\undefined
\def\x#1#2#3#4#5#6#7\relax{\def\x{#1#2#3#4#5#6}}%
\expandafter\x\fmtname xxxxxx\relax \def\y{splain}%
\ifx\x\y   % LaTeX or SliTeX?
\gdef\SetFigFont#1#2#3{%
  \ifnum #1<17\tiny\else \ifnum #1<20\small\else
  \ifnum #1<24\normalsize\else \ifnum #1<29\large\else
  \ifnum #1<34\Large\else \ifnum #1<41\LARGE\else
     \huge\fi\fi\fi\fi\fi\fi
  \csname #3\endcsname}%
\else
\gdef\SetFigFont#1#2#3{\begingroup
  \count@#1\relax \ifnum 25<\count@\count@25\fi
  \def\x{\endgroup\@setsize\SetFigFont{#2pt}}%
  \expandafter\x
    \csname \romannumeral\the\count@ pt\expandafter\endcsname
    \csname @\romannumeral\the\count@ pt\endcsname
  \csname #3\endcsname}%
\fi
\fi\endgroup

\newcommand{\Dniv}[9]{\parbox[c]{6.35cm}{$
\begin{picture}(9951,2200)(1125,-5050)
\thinlines
\put(1201,-3961){\circle{336}}
\put(2401,-3961){\circle{336}}
\put(3601,-3961){\circle{336}}
\put(4801,-3961){\circle{336}}

\put(5836,-3946){\circle*{30}}
\put(5986,-3946){\circle*{30}}
\put(6136,-3946){\circle*{30}}

\put(7201,-3961){\circle{336}}
\put(8401,-3961){\circle{336}}
\put(9601,-3961){\circle{336}}

\put(10453,-3109){\circle{336}}
\put(10453,-4813){\circle{336}}

\put(1369,-3961){\line( 1, 0){870}}
\put(2569,-3961){\line( 1, 0){870}}
\put(3796,-3961){\line( 1, 0){870}}

\put(4969,-3961){\line( 1, 0){600}}
\put(6403,-3961){\line( 1, 0){600}}

\put(7369,-3961){\line( 1, 0){870}}
\put(8569,-3961){\line( 1, 0){870}}
\put(9719,-3843){\line( 1, 1){615}}
\put(9719,-4079){\line( 1, -1){615}}

\put(1141,-3436){\makebox(0,0)[lb]{\smash{\SetFigFont{6}{7.2}{rm}\mbox{$#1$}}}}
\put(2341,-3436){\makebox(0,0)[lb]{\smash{\SetFigFont{6}{7.2}{rm}\mbox{$#2$}}}}
\put(3541,-3436){\makebox(0,0)[lb]{\smash{\SetFigFont{6}{7.2}{rm}\mbox{$#3$}}}}
\put(4741,-3436){\makebox(0,0)[lb]{\smash{\SetFigFont{6}{7.2}{rm}\mbox{$#4$}}}}
\put(7141,-3436){\makebox(0,0)[lb]{\smash{\SetFigFont{6}{7.2}{rm}\mbox{$#5$}}}}
\put(8341,-3436){\makebox(0,0)[lb]{\smash{\SetFigFont{6}{7.2}{rm}\mbox{$#6$}}}}
\put(9541,-3436){\makebox(0,0)[lb]{\smash{\SetFigFont{6}{7.2}{rm}\mbox{$#7$}}}}
\put(10876,-3211){\makebox(0,0)[lb]{\smash{\SetFigFont{6}{7.2}{rm}\mbox{$#8$}}}}
\put(10876,-4936){\makebox(0,0)[lb]{\smash{\SetFigFont{6}{7.2}{rm}\mbox{$#9$}}}}
\end{picture}$}
}

\begin{document}

\title{\bf Representations of compact Lie groups and\\
the osculating spaces of their orbits}
\author{\sc Claudio Gorodski\footnote{\emph{Alexander von Humboldt Research
      Fellow} at the University of Cologne during the completion of this
      work.}\ \
      and Gudlaugur Thorbergsson}
%\date{\it September 2000}
\date{}
\maketitle

\vfill

\parbox[t]{7cm}{\footnotesize\sc Instituto de Matem\'atica e Estat\'\i stica\\
                Universidade de S\~ao Paulo\\
                Rua do Mat\~ao, 1010\\
                S\~ao Paulo, SP 05508-900\\
                Brazil\\ \hfill\\
                E-mail: {\tt gorodski@ime.usp.br}}{}\hfill{}
\parbox[t]{7cm}{\footnotesize\sc Mathematisches Institut\\
                Universit\"at zu K\"oln\\
                Weyertal 86-90\\
                50931 K\"oln\\
                Germany\\ \hfill\\
                E-mail: {\tt gthorber@mi.uni-koeln.de}}\hfil{}
\vspace{.5cm}

\thispagestyle{empty}
\pagebreak

\thispagestyle{empty}
\tableofcontents
\pagebreak

\thispagestyle{empty}

{\ }
\vspace{4cm}\relax

\begin{center}
\bf\Large Abstract
\end{center}

Several classes of irreducible orthogonal representations of compact Lie
groups that
are of importance in Differential Geometry have
the property that the second osculating space of all of their nontrivial
orbits
coincide with the representation space. We say
that representations with this property are of class $\mathcal O^2$. Our
approach in the present paper will be to find
restrictions on the class $\mathcal O^2$ and  then apply them to classify
variationally complete and taut representations. The known classifications
of cohomogeneity one and two orthogonal representations 
and more generally of polar representations 
will also 
follow easily. 

\vfill

\footnotetext{2000 {\em Mathematics Subject Classification}: primary, 
57S15; secondary, 53C30, 53C40, 53C42.}

\pagebreak

%name of the file: prolog.tex

\section{Introduction}

Several classes of irreducible orthogonal representations of compact Lie
groups that
are of importance in Differential Geometry have
the property that the second osculating space of all of their nontrivial
orbits
coincide with the representation space. We say
that representations with this property are of class $\mathcal O^2$. Our
approach in the present paper will be to find
restrictions on the class $\mathcal O^2$ and  then apply them to classify
variationally complete and taut representations. The known classifications
of cohomogeneity one and two orthogonal representations 
and more generally of polar representations 
will also 
follow easily.  We refer to Section~\ref{sec:review} for
the
definitions of concepts used in this introduction.\medskip

Variationally complete actions were
introduced by Bott in \cite{Bott}.  Bott and
Samelson
proved in
\cite{B-S} that isotropy representations of symmetric spaces are
variationally
complete. We will prove the following converse of their theorem in 
Section~\ref{sec:BS}.

\begin{thm} \label{thm:var}
A variationally complete
representation of a compact connected Lie group is orbit
equivalent to the isotropy representation of a symmetric
space.\end{thm}

We will in fact prove more than we have stated here and give a complete
classification of variationally complete
representations (up to image equivalence; see the beginning of 
Section~\ref{sec:O2} for this concept). As a
corollary  we get a new
proof of the
classification of polar representations due to  Dadok (\cite{Dadok}) since
these are variationally
complete (see
\cite{Conlon}) and  isotropy representations of symmetric spaces are easily
seen to be polar. It
follows that an orthogonal representation is variationally complete if and
only if
it is polar, and it is polar if and only if it is orbit equivalent to the
isotropy
representation of a symmetric space.

Another important result of Bott and Samelson in~\cite{B-S} is that the
distance functions to orbits of variationally complete
representations are perfect Morse functions. We say that submanifolds with this
property are {\it taut} and call representations
{\it taut} if all of their orbits are taut. 
In Section~\ref{sec:taut} we will prove the following
classification theorem 
which provides new examples of taut submanifolds and taut 
representations.  
\begin{thm}\label{thm:taut}
 A taut irreducible
representation~$\rho$ of a compact connected Lie group~$G$
is either orbit equivalent to the isotropy representation of a symmetric
space or it is one of the following orthogonal representations ($n\geq2$):
\[ \begin{array}{|c|c|}
\hline
 G & \rho \\
\hline
\SO2\times\Spin9 & \mbox{(standard)}\otimes_{\mathbf R}\mbox{(spin)} \\
\U2\times\SP n & \mbox{(standard)}\otimes_{\mathbf C}\mbox{(standard)} \\
\SU2\times\SP n & \mbox{(standard)}^3\otimes_{\mathbf H}\mbox{(standard)} \\
\hline
\end{array}\]
\end{thm}

As a consequence of this theorem, the classes of polar,
variationally complete and taut irreducible representations of
compact {\it simple} groups coincide. Notice that the three
exceptional representations in Theorem~\ref{thm:taut} are precisely the
irreducible representations of cohomogeneity three that are not 
polar (see~\cite{Yasukura1}). Hence, in the case
of a general compact Lie group, the taut
irreducible representations are precisely those that are 
either polar or of cohomogeneity three.
A further related result is the classification of homogeneous taut
submanifolds with flat normal bundle due to Olmos (\cite{Olmos}).
He shows that such a submanifold is a principal orbit
of the isotropy representation of a symmetric space.
It follows that the principal orbits of the three nonpolar representations
in Theorem~\ref{thm:taut} do not have flat normal bundle
which is also easy to see directly.

Bott and Samelson proved the tautness of variationally complete
representations as
an application of their $K$-cycles which are concrete representatives of a
basis for
the homology of the orbits. (In~\cite{B-S} the group acting is denoted by
$K$ 
which explains the terminology `$K$-cycle').
The three nonpolar representations in Theorem~\ref{thm:taut} admit
generalizations of the $K$-cycles of Bott and Samelson. We will not
verify this claim in the present paper since it requires
methods quite different from those employed here and would
make the paper too long; we refer instead to the forthcoming
paper~\cite{G-Th}. Instead we will prove their tautness with shorter,
although maybe
less illuminating arguments based on a theorem of Floyd (\cite{Floyd}). This
approach is adapted
from the paper \cite{Duistermaat} of Duistermaat.

Kuiper observed in \cite{Kuiper} that the second osculating space of a taut
submanifold in a Euclidean space $V$ coincides with $V$ if the submanifold
is not
contained in a proper affine subspace. In fact, he proved this more generally
for
tight submanifolds, but this is unimportant for us since an orbit is tight
if and
only if it is taut. Since the classes of representations we are dealing with
are all taut,
it follows from this observation of Kuiper that they belong to class
${\mathcal
O}^2$ if they are irreducible. The class ${\mathcal O}^2$ is much more
tractable
than the other classes of representations we are dealing with since it
involves an
infinitesimal condition. The technique of Dadok (\cite{Dadok}),
notably his invariant
$k(\lambda)$, turns out to be an extremely powerful tool to reduce the
class
${\mathcal O}^2$ in size so that the remaining cases are accessible to the
geometric methods developed in Section~\ref{sec:geom}.

We have already pointed out that Dadok's classification of polar representations
follows as
a consequence of Theorem~\ref{thm:var}.
Some other classification results can also
easily be proved with our methods as shown in 
Section~\ref{sec:other}. These include
the classification of cohomogeneity one representations
as well as the classification of cohomogeneity two
representations due to Hsiang-Lawson (\cite{H-L}). The cohomogeneity two
representations are polar and therefore included in Dadok's classification,
but
the point here is that it is very easy to see directly that they belong to
class
${\mathcal O}^2$ without referring to tautness.

One could also use our methods to classify symmetric spaces and their
isotropy
representations. This would not lead to a simpler proof than those existing
in the
literature. Insisting on the details needed to classify symmetric spaces and
their isotropy representations would make the paper much longer and more
involved than
it is. We will therefore use this classification in the
paper\footnote{Except for the classification of cohomogeneity one representations.}
and assume knowledge of the lists in~\cite{Wolf},
Tables~8.11.2 and~8.11.5.
(See also~\cite{Helgason}, but notice that the isotropy representations are
not
listed in this reference). 
Otherwise, our guiding principle in writing the paper has been to make it as
self-contained as possible. 

Previous to this paper taut representations were studied in~\cite{Galemann}
and~\cite{C-Th}.  
It is decided in~\cite{Galemann}, with some exceptions,
which
representations of $\SU n$ and ${\mathbf U} (n)$ can be taut. In~\cite{C-Th}
it is
proved among other things that a compact group admitting an almost faithful
taut
representation can have at most four simple factors. Eschenburg and Heintze
use the
theory of isoparametric submanifolds in~\cite{E-H2} to reprove Dadok's
characterization up to orbit equivalence 
of polar representations under the assumption of
cohomogeneity at
least three. In~\cite{E-H1} they use the same theory to decide which
irreducible
representations are orbit equivalent to isotropy representations of
symmetric spaces. This follows immediately from
our work in~Section~\ref{sec:BS}.
Bergmann in~\cite{Bergmann} decides which reducible representations
are polar. Her result is equivalent to
the last step in our classification of variationally complete
representations. Di Scala and Olmos give in~\cite{D-O} a direct proof 
of the fact that
a variationally complete representation is polar.

The first author wishes to thank the \emph{Alexander von Humboldt
Foundation} for its generous support and constant assistance
during the completion of this work.

\section{A review of basic definitions and results}\label{sec:review}
\setcounter{thm}{0}

Let $G$ be a compact Lie group acting on a Riemannian manifold $M$ by
isometries. A geodesic $\gamma$ in $M$  is called
{\it
$G$-transversal} if it is orthogonal to the $G$-orbit through $\gamma(t)$
for every $t$. One can show that a geodesic
$\gamma$ is $G$-transversal if there is a point $t_0$ such that
$\dot\gamma(t_0)$ is orthogonal to $G\gamma(t_0)$.
A Jacobi field along a geodesic in $M$ is called {\it $G$-transversal}
if it is the variational vector field of a
variation through
$G$-transversal geodesics. The action of $G$ on $M$ is called {\it
variationally complete} if every $G$-transversal Jacobi
field $J$ in $M$ that is tangent to the $G$-orbits at two different
parameter values is the restriction of a Killing
field on
$M$ induced by the $G$-action (see~\cite{Bott} and~\cite{B-S}).
It is proved in~\cite{B-S} on
p.~974 that instead of requiring tangency at two
different points in the definition of variational completeness it is
equivalent to require tangency at one point and
vanishing at another point.

Let $\rho$ be a variationally complete reducible representation and
$\tilde\rho$ a summand of $\rho$. Then it is easy to see
that $\tilde\rho$ is also variationally complete, see Proposition~\ref{prop:red-vc}.
 Notice though that it is not
true that the direct sum of two variationally complete representations of
a compact Lie group $G$ is variationally complete as can be seen by taking
the direct sum of of
two copies of $\SO 2$ acting on $\mathbf R^2$.

Let $N$ be a properly embedded submanifold of a Euclidean space $V$ and let
$x$ be some point in $V$. Then we define the
{\it distance function $L_x: N\to {\bf R}$ from $x$ to $N$} by setting
$L_x(p)=||p-x||^2$. It follows that $L_x$ is a
non-negative proper function since $N$ is properly embedded. Hence it is
possible to apply Morse theory to $L_x$. We say
that $L_x$ is {\it perfect with respect to a field $F$} if it is a Morse
function and the Morse inequalities for $L_x$
with respect to $F$ are equalities. Furthermore we say that $N$ is {\it
$F$-taut} or simply {\it taut} if $L_x$ is perfect
with respect to $F$ whenever $L_x$ is a Morse function, see~\cite{C-R}. The
concept of tautness can be extended to
submanifolds of complete Riemannian manifolds, see~\cite{T-Th1}, but we will
not
need that here. We will say that an orthogonal
 representation $\rho:G\to \mathbf O(V)$ of a compact Lie group $G$ is {\it
$F$-taut} or simply {\it taut} if the the orbits of $G$
are $F$-taut submanifolds of $V$.

If $\rho$ is a reducible taut representation and
$\tilde\rho$ one of its factors, then  $\tilde\rho$ is clearly taut. Notice
 that it is not true that direct sums of  taut
representations of a compact Lie group $G$ are also taut. An example of this
is the direct sum of
$n$ copies of $\SU n$ acting on $\mathbf C^n$. This representation is not
taut since one can compute that $\SU n$ is
not taut in the matrix space $M(n;\mathbf C)$.

The following theorem mentioned in the introduction was proved by Bott and
Samelson
more generally for variationally complete actions on complete
Riemannian manifolds. In fact it is an immediate corollary of  Theorem I
in~\cite{B-S}.  Notice that Bott and Samelson do not
use the concept of a taut submanifold which was introduced later.

\begin{thm}[Bott-Samelson]\label{thm:BS1}
A variationally complete
representation of a compact connected
Lie group on an Euclidean space is ${\bf Z}_2$-taut. 
\end{thm}

Bott proved in~\cite{Bott}, Propositions~7.1 and~11.6, that the action of a
compact
Lie group $G$ with a bi-invariant metric on
itself by conjugations and the adjoint representation of $G$ on its Lie
algebra $\mathfrak g$ are variationally complete.
This was generalized by Bott and Samelson in~\cite{B-S}, Theorem~II as
follows.
Let $(L,G)$ be a symmetric pair and ${\Ll}={\Lg}\oplus{\Lp}$ the
corresponding Cartan decomposition. Then
the action of $L\times L$ on $L$, the action of
$G$ on $L/G$ and the action of ${\rm Ad}_L(G)$ restricted to $\Lp$ are
variationally complete. The action of  ${\rm
Ad}_L(G)$ on $\Lp$ is equivalent to the isotropy representation of $G$
on the tangent space $T_p(L/G)$ where $p$ denotes
the coset $G$. Therefore we have the following theorem mentioned in the
introduction.

\begin{thm}[Bott-Samelson]\label{thm:BS2}
 The isotropy representation of a
symmetric space
is variationally complete. 
\end{thm}

Conlon considered in~\cite{Conlon} actions of a Lie group $G$ on a complete
Riemannian manifold $M$ with the property that there is a
connected submanifold $\Sigma$ of $M$ that meets
all orbits of $G$ in such a way that the intersections between $\Sigma$ and
the orbits of $G$ are all orthogonal. Such a submanifold
is called  a {\it section} and an action admitting a section is now usually
called {\it polar} if $\Sigma$ is properly
embedded. Notice that Conlon does not assume in~\cite{Conlon} that $\Sigma$
is
properly embedded, but it is usually required in
the recent literature on the subject.  It is easy to see that a section
$\Sigma$ is totally geodesic in
$M$. An action admitting a section that is flat in the induced metric is
called {\it hyperpolar}. There is clearly no
difference between polar and hyperpolar representations since totally
geodesic submanifolds of a Euclidean space are affine
subspaces. Moreover, the question whether $\Sigma$ should be required to be
properly embedded or not becomes redundant.  Conlon
proved the following theorem in~\cite{Conlon}.

\begin{thm}[Conlon]\label{thm:conlon}
 A hyperpolar action of a compact Lie
group on a complete
Riemannian manifold is variationally complete. \end{thm}

Polar representations were classified by Dadok in~\cite{Dadok}. We recall
that two representations $\rho_1:G_1\to O(V_1)$ and
$\rho_2:G_2\to O(V_2)$ are said to be {\it orbit equivalent} if there is an
isometry $A:V_1\to V_2$ under which the orbits of
$G_1$ and $G_2$ correspond.
As a consequence of his classification he obtained the following result.

\begin{thm}[Dadok] A polar representation of a compact connected
Lie group is orbit equivalent to the isotropy
representation of a symmetric space. \end{thm}

We can summarize this discussion as follows. Let $G$ be a compact connected
Lie group.
Denote by $\cal I$ the representations of $G$
that are isotropy representations of symmetric spaces, by ${\cal P}$ those
that are polar, by ${\cal V}$ those that  are
variationally complete, and by ${\cal T}$ those that are taut. Then
$${\cal I}\subset{\cal P}\subset {\cal V}\subset{\cal T}.$$

The starting point in the present paper will be to consider a class of
representations that is a priori larger than $\cal
T$, but easier to deal with, see~\cite{C-Th}. We recall that the {\it second
osculating space} $O_p^2(N)$ at a point $p$ of a
submanifold
$N$ in a Euclidean space $V$ is the vector space spanned by the first and
second derivatives at $p$ of the inclusion of $N$
into
$V$. It is easy to see that
$$
O_p^2(N)=T_pN\oplus\langle\{\alpha(X,Y)\,|\, X,Y\in T_pN\}\rangle.
$$
where $\alpha$ denotes the second fundamental form of $N$ and  $\langle
S\rangle$ stands for the linear hull of the set $S$.
The following is a corollary of the discussion of Kuiper in~\cite{Kuiper}.

\begin{thm}[Kuiper] \label{thm:Kuiper}
Let $N$ be a taut submanifold of a
Euclidean space. Then the affine hull of $N$
coincides with $p+O_p^2(N)$ for every point $p$ in $N$. \end{thm}

We let $\mathcal O^2$ denote the class of representations of a compact 
connected Lie group
$G$ such that the representation space coincides with $O^2_p(G p)$ for all
nonzero
$p$. The representations in $\mathcal O^2$ are clearly irreducible. If a
taut
representation is irreducible, then Theorem~\ref{thm:Kuiper} implies that
it belongs to class $\mathcal O^2$. Thus we have the inclusions
$${\cal I}_i\subset{\cal P}_i\subset {\cal V}_i\subset{\cal T}_i\subset
{\cal
O}^2,$$
where ${\cal I}_i$, ${\cal P}_i$, ${\cal V}_i$, and  ${\cal T}_i$ are the
the subclasses of ${\cal I}$, ${\cal P}$, ${\cal
V}$, and  ${\cal T}$  consisting of irreducible
representations.

\section{A collection of geometric propositions}\label{sec:geom}
\setcounter{thm}{0}

In this section we collect technical results of geometric
nature that will be needed to prove our main theorems.

\subsection{Variationally complete representations}

We first study focal points of orbits of variationally
complete representations with the aim to prove that a representation that is
orbit equivalent to a variationally complete one is also variationally
complete. We then show how reducible
variationally complete representations relate to irreducible ones. Finally
we prove an important
criterion which can be used to decide in some cases 
whether a variationally complete
representation is orbit equivalent
to the isotropy representation of a symmetric space.

Let $q$ be a focal point of an orbit $Gp$ of an orthogonal representation
of a compact connected Lie group $G$. It is clear that the orbit
$Gq$ consists of focal points of $Gp$. We therefore call $Gq$ a {\it focal
orbit of $Gp$.} It is also clear that the singular orbits of a
representation are
always focal orbits of the principal orbits. Slightly more generally,
$Gq$ is
a focal orbit of $Gp$ if the segment $\overline{qp}$ is perpendicular 
to $Gp$ in~$p$ and there is a one-parameter group in $G_q$ that does
not fix $p$. We will use the following lemma to prove
Proposition~\ref{prop:equivvc}.

\begin{lem}\label{lem:tecvc}
A representation $\rho:G\to \mathbf O(V)$ is variationally
complete if and only if every orbit $Gp$ has the following
property: if $q$ is a focal point of $Gp$ relative to~$p$ with 
multiplicity~$k>0$, then $k=\dim G_qp$.
\end{lem}

\begin{rmk}[\mbox{\rm Compare~\cite{Terng}.}]\label{rmk:tecvc}
\em 
Let $\rho:G\to
\mathbf O(V)$ be a variationally complete representation of a compact
connected Lie group $G$. Then Lemma~\ref{lem:tecvc} implies that
the focal orbits of any orbit of $G$ are singular orbits. Furthermore,
the multiplicity of the focal points of a principal orbit $Gp$ in a singular
orbit $Gq$ is equal to the difference between
the dimension of $Gp$ and the dimension of $Gq$.
In particular, it follows that any principal orbit is a \emph{totally focal
submanifold}; see~\cite{Terng} for a definition of this concept.
Terng claims in~\cite{Terng} that this fact can be used to
show that variationally complete representations are polar. However,
as Terng herself pointed out in~\cite{T-Th1}, p.~197, there is a gap in that
argument. 
\end{rmk}

{\it Proof of Lemma~\ref{lem:tecvc}.}
 We first assume that $\rho$ is variationally complete. Let $q$ be a
focal 
point of $Gp$ along the line segment $\overline{pq}$ with multiplicity $k$.
Let 
$\cal J$ be the space of Jacobi fields along $\overline{pq}$ that are the
variational 
vector fields of 
 variations $\gamma_s$ of
$\overline{pq}$ through line segments such that all $\gamma_s$ meet the
orbit $Gp$ orthogonally for
$t=0$ and the Jacobi field $J$ of the variation of $\gamma$
 vanishes for
$t=1$. The multiplicity $k$ of the focal point $q$ is equal to the dimension
of the 
space $\cal J$. The action of $G_q$ on $\overline{pq}$ induces variations
of 
$\overline{pq}$ whose variational vector fields are contained in $\cal J$. Let
us 
denote by $\hat{\cal J}$ the subspace of $\cal J$ that consists of Jacobi
fields 
induced by the action of $G_q$ on $\overline{pq}$.
It follows that
$\dim\hat{\cal J}=\dim G_qv\le k$ where $v$ is the tangent vector of
$\overline{pq}$ at $q$. We now prove
that $\dim G_qv= k$. Let  $\gamma_s$ be a variation of
$\overline{pq}$ through line segments such that all $\gamma_s$ meet the
orbit $Gp$ orthogonally for
$t=0$ and the Jacobi field $J$ of the variation along $\gamma$
 vanishes for
$t=1$. The segments
$\gamma_s$ are orthogonal to the orbit $G p$. Hence they are $G$-transversal
and $J$ is a $G$-transversal Jacobi field.
It is clear that $J$ is tangent to the orbits through $p=\gamma(0)$ and
$q=\gamma(1)$. By variational completeness, there is
a Killing field $X$ on $M$ induced by the $G$-action such that
$J=X\circ\gamma$. Hence $J\in\hat{\cal J}$ and it follows that $\hat{\cal
J}={\cal J}$.
Thus $\dim G_qv=k$.

To prove the converse, we have to show that $\dim G_qv=k$ for all focal
points $q$ implies that $\rho$ is
variationally complete.
This is equivalent to show that a Jacobi field $J$
that vanishes in 
a point $q$ and is tangent to an orbit $Gp$ in $p$ is induced by the group
action. We let the symbols $\cal J$ and $\hat{\cal J}$ have the same meaning
as in 
the first part of the proof. We have to show that $\hat{\cal J}=\cal J$.
This is clear since
$\dim {\cal J}=k$ and $\dim\hat{\cal J}=\dim G_qv=k$.
\EPf

\medskip

The following proposition that we will use to classify variationally
complete
representations is now easy.

\begin{prop}\label{prop:equivvc}
Let $\rho$ be a representation that is orbit equivalent to a variationally
complete
representation. Then $\rho$ is itself variationally complete.
\end{prop}

\Pf We can assume that we have representations $\rho$ and $\hat\rho$
of the groups $G$ and $\hat G$ with the same representation space and
the same orbits and that $\hat\rho$ is variationally complete. Clearly the
focal orbits of the two representations coincide since their orbits are
the same. Let $v$ be a normal vector in $q$ of the orbit $Gq=\hat Gq$. It
follows from the
slice theorem and the fact that orbits of $\rho$ and $\hat\rho$ are the same
that
$G_qv=\hat G_qv$. Since $\hat\rho$ is variationally complete,
Lemma~\ref{lem:tecvc} now implies that $\rho$ is also
variationally
complete.
\EPf

\medskip

The next proposition will be important in the proof of
Theorem~\ref{thm:var} in the introduction. It
is analogous to a result of Dadok on polar representations, see Theorem 4 in
\cite{Dadok}.

\begin{prop}\label{prop:red-vc}
 Let $G$ be a compact connected Lie group. Let $\rho:G\to
\mathbf O (V)$ be a nontrivial
reducible variationally complete representation and let $\rho=\rho_1\oplus\rho_2$ be
a nontrivial decomposition. Let $V_1$ and $V_2$ be the representation
spaces of $\rho_1$ and $\rho_2$ respectively. 
Then we have the following:
\begin{enumerate}
\item[(a)] The representations $\rho_1$ and $\rho_2$ are variationally complete.
\item[(b)] The restriction of $\rho_2$ to the connected component 
of the isotropy subgroup of any point in 
$V_1$ is orbit equivalent to~$\rho_2$. 
\item[(c)] The representations $\rho_1$ and $\rho_2$ are not equivalent 
as representations. 
\item[(d)] The representation $\rho$ is orbit equivalent to the representation
$\hat\rho:G\times
G\to\mathbf O(V)$ defined by $\hat\rho(g_1,g_2)=\rho_1(g_1)+\rho_2(g_2)$.
\end{enumerate}
\end{prop}

\Pf The claim in part~(a) is clear, since the orbits of $\rho_1$ 
and $\rho_2$ are also orbits of $\rho$. 

To prove part~(b), note that it is true in general, without
assuming $\rho$ is variationally complete, 
that any point $p_1\in V_1$ is a focal point 
of any orbit $Gp_2\subset V_2$ with multiplicity equal 
to~$\dim Gp_2$. Now if $\rho$ is variationally complete 
then the claim follows from Lemma~\ref{lem:tecvc} and the fact that 
the $G$-orbits are connected.

Next we prove~(c). Assume $A:V_1\to V_2$ is an isomorphism such that
$\rho_2=A\rho_1A^{-1}$. Notice that $G_p=G_{Ap}$ for all $p\in V_1$.
Part~(b) implies that $G_p$ is transitive on the orbit $G(Ap)$, 
which is impossible if we choose $p\in V_1$ such that $\dim Gp>0$.
Therefore such an isomorphism~$A$ does not exist. 

Finally we prove part~(d). Let $v=v_1+v_2$ be a point in $V$ and
$v_1$ and $v_2$ its
components in $V_1$ and $V_2$ respectively. Let $H_1$ (resp.~$H_2$)
be the isotropy group
of $v_1$ (resp.~$v_2$). Notice that $v_2$ is a focal point of $Gv_1$
with multiplicity
$\dim Gv_1$. It follows from part~(b) that
$H_2v_1=Gv_1$ and
similarly that $H_1v_2=Gv_2$. We next prove the factorization 
$G=H_1H_2$. Let $g\in G$.
Then there is an element
$h_1\in H_1$ such that $h_1v_2=gv_2$. Then $h_1^{-1}g\in H_2$. Set
$h_2=h_1^{-1}g$. Then
$g=h_1h_2$. It follows that
$$
\mathfrak g=\mathfrak h_1 +\mathfrak h_2
$$
where $\mathfrak g$, $\mathfrak h_1$ and $\mathfrak h_2$ are the Lie
algebras of $G$, $H_1$ and
$H_2$ respectively. It is clear that the orbit of $v$ under $\rho(G)$ is
contained
in its orbit under $\hat\rho(G\times G)$. The tangent space of
$\hat\rho(G\times G)v$ at $v$ is
$$
(\mathfrak g\times\mathfrak g)v=\mathfrak g v_1+\mathfrak g v_2=\mathfrak
h_2v_1+\mathfrak h_1v_2.
$$
Let $X_1\in\mathfrak h_1$ and $X_2\in\mathfrak h_2$. Then $X_1v_1=0$ and
$X_2v_2=0$ and hence
$$
X_1v_2+X_2v_1=(X_1+X_2)v\in\mathfrak gv.
$$
Hence the dimension of the tangent spaces of the orbits of $v$ under $\rho$
and $\hat\rho$
coincide. Since $G$ is connected the orbits coincide. Since $v$ was chosen
arbitrarily, all orbits
$\rho$ and $\hat\rho$ coincide.\\ \mbox{}\EPf

\begin{rmk}
\em
Part~(c) in Proposition~\ref{prop:red-vc} generalizes Proposition~2.10
in~\cite{HPTT2}.
\end{rmk}

The next lemma is a preparation for Proposition~\ref{prop:common-orbit}. 
We follow the terminology of Bredon in~\cite{Bredon}, pp.~180-181, 
and distinguish between three types of orbits, namely
principal (regular), exceptional and singular. Recall that orbits are called
\emph{exceptional} if they have maximal dimension without being principal,
and \emph{singular} if their dimension is not maximal. We also call points 
\emph{regular} when they belong to principal orbits. 

\begin{lem}\label{lem:principal-orbit}
Let $G$ be a closed subgroup of the compact connected
Lie group $K$  and let $\rho:K\to \mathbf O(V)$ be
a representation. If $\rho$ has an 
orbit of maximal dimension in common with $\rho|G$, then $\rho$ and 
$\rho|G$ are orbit equivalent. 
\end{lem}

\Pf We denote the maximal dimension of an orbit of $\rho$ by $n$.
Let us denote by ${\cal O}_G$ the union of the $G$-orbits of maximal
dimension and by ${\cal O}_K$ the corresponding set relative to~$K$.
It is well-known that both ${\cal O}_G$ and ${\cal O}_K$ contain an
open and dense subset of $V$. 
Hence the set ${\cal O}_G\cap {\cal O}_K$ also contains an open and dense
set. The dimensions of the orbits of $G$ and $K$ through
points in
${\cal O}_G\cap {\cal O}_K$ are all equal to $n$. It follows that the orbits
of $G$ and $K$ through points in ${\cal
O}_G\cap {\cal O}_K$ coincide since $K$ is connected. Let $p\in
V\setminus{\cal
O}_G\cap {\cal O}_K$. We would like to show that $G p=K
p$.  Let
$k(p)$ be some other point in $K p$. Let $(p_n)$ be a sequence in ${\cal
O}_G\cap {\cal O}_K$ that converges to $p$.
Then $(k(p_n))$ converges to $k(p)$. There is for every $n$ an element $g_n$
in $G$ such that $g_n(p_n)=k(p_n)$. Now $(g_n(p_n))$ converges to $k(p)$. 
By going to a subsequence if necessary we can assume that $g_n$ converges to an
element $g$ in~$G$, so that the continuity of $G\times V\to V$ implies that
$(g_n(p_n))$ converges to~$g(p)$. Hence $g(p)=k(p)$ 
and we have seen that all orbits of $G$ and $K$ coincide. \EPf

\begin{prop}\label{prop:common-orbit}
Let $G$ be a closed subgroup of the 
compact connected Lie group $K$ and let $\rho:K\to \mathbf O(V)$ be
the isotropy representation of an irreducible symmetric space.
If the restriction
$\rho|G$ is variationally complete and $\rho|G$ has a nontrivial orbit in common
with $\rho$, then $\rho$ and $\rho|G$ are orbit equivalent.
\end{prop}

\Pf 
Let $p\in V$ be a point such that the orbits of $G$ and $K$ through $p$
coincide and have the largest possible dimension
with this property. Assume that $Kp$ is not an orbit of maximal
dimension for the action of $K$. Let ${\La}$ be a Cartan subspace in
$V$ with respect to the symmetric space representation of $K$ such that
$p\in {\La}$. Let $T$ be a
tubular neighborhood of the common orbit $K p =G p$.
Let $q\in{\La}$ be a point in $G p=K p$ different from
$p$ and such that the segment $\overline{pq}$ from $p$ to $q$ passes through
a point $r$ in $T$ with $\dim
K_pr\ge 1$. (Here it is important that the symmetric space is
irreducible.) It
follows that $\dim Kp<\dim Kr$ by the slice theorem. Hence
$\dim Gr <\dim Kr$ since the dimension of $Gp=Kp$
is maximal among the common orbits of $G$ and $K$. It follows that $\dim
G_pr<\dim K_pr$. 

Let $X$ be an element of the Lie algebra
${\mathfrak k}_p$ of $K_p$ such that $\dot\gamma_{X,r}(0)$ is not in $T_r(G_p
r)$ where $\gamma_{X,r}(s)=\exp(sX) r$.
Now let us consider the variation $(\overline{pq})_s=\exp(sX)
(\overline{pq})$ of $\overline{pq}$. Let $J$ be the
corresponding Jacobi field along $\overline{pq}$.  This variation is
$K$-transversal and hence also $G$-transversal. The value of $J$ over $p$
and $q$ are tangent vectors of the common orbit
$G p=K p$. By the variational completeness of the action of $G$, the
Jacobi field $J$ is the restriction
of a Killing field on $V$ to $\overline{pq}$ that is induced by the action
of $G$ on $V$. Hence all values of $J$ are
tangent vectors to orbits of $G$. The value of
$J$ over $r$ is
$\dot\gamma_{X,r}(0)$ which is not a tangent vector of the orbit of $G$
through $r$. This is a contradiction. Hence $G$ and
$K$ have an orbit in common which is of maximal dimension among the 
$K$-orbits. It now follows from
Lemma~\ref{lem:principal-orbit} that all orbits of
$G$ and $K$ coincide. \EPf

\subsection{Taut representations}

We will first prove that the slice representation of a taut representation
is also taut
and end the section with an important proposition about
taut reducible representations. It also follows immediately from our methods that a
taut representation
does not have exceptional orbits.

\begin{prop}\label{prop:slice}
Let $G$ be a compact connected Lie group with a
taut representation $\rho:G\to \mathbf O(V)$. Let $p\in
V$ and let
$\nu_p:G_p\to \mathbf O(N_p(G p))$  be the slice representation of $\rho$ at~ $p$.
Then $\nu_p$ is taut.
\end{prop}

\Pf Let $N^\epsilon(G p)$ denote the bundle of normal
vectors of length less than $\epsilon$ over the
orbit $G p$. Let $\epsilon>0$ be so small that
$T=\exp(N^\epsilon(Gp))$ is a tubular neighborhood around the orbit $G p$.
Let $v\in N_p(G p)$. We want to show
that $G_p v$ is taut in the normal space $N_p(G p)$. Let $\alpha>0$ be a
number that is so small that $w=\alpha\,
v$ has length less than $\epsilon$. It is clear that $G_p v$ is taut in
$N_p(G p)$ if and only if $G_p w$ is taut in $N_p(G p)$.
It is also clear that $G_p w$ is taut in $N_p(G p)$ if and only
if $\exp_p(G_p w)=G_p\exp_p(w)$
is taut in
$\exp_p(N_p(G p))$ since $\exp_p$ is an isometry. We set $q=\exp_p(w)$. A
submanifold in an affine subspace $A$ of $V$
is taut in
$A$ if and only if it is taut in $V$. Hence $G_p v$ is taut in $N_p(G p)$ if
and only if $G_p q$ is taut in
$V$. 

Now let $L_p:G q\to {\bf R}$ be the distance function. The segment
$\overline{pq}$ is orthogonal to $G p$ and
hence also to $G q$. It follows that $q$ is a critical point of $L_p$. Since
$q$ is in the tubular neighborhood $T$ and
the length of the segment $\overline{pq}$ is smaller than $\epsilon$ it
follows that $L_p$ takes on its minimum in $q$.
Let $C$ be the subset of $G q$ on which $L_p$ takes on its minimum value.
Clearly $C=G q\cap \exp_p(N^\epsilon_p(G
p))$ since a geodesic segment between $Gp$ and $Gq$ that is orthogonal to
$Gp$ is orthogonal to $Gq$ and vice versa. By the
slice theorem
$G_p q=G q\cap
\exp_p(N^\epsilon_p(G p))$. Hence
$G_p q= C$. The set of critical points  of a distance function to a taut
submanifold is a union over taut submanifolds by a
theorem of Ozawa, see~\cite{Ozawa}. It follows that $G_pq$ is taut and hence
that the
slice representation $\nu_p$ is taut. \EPf

\begin{rmk}
\em
We are assuming in Proposition~\ref{prop:slice} that $G$ is connected.
Hence its orbits are connected and it follows from
their tautness that the set of points where a distance function takes on its
minimal value is connected. The proof of
Proposition~\ref{prop:slice} now implies that the orbits of $G_p$ are
connected even if $G_p$ is
not connected.
\end{rmk}

It is well-known that the isotropy representation of a symmetric space does
not have an exceptional orbit.
The methods used in the proof of Proposition~\ref{prop:slice}
immediately give that this is also the
case for taut representations, as we show in the next proposition.

\begin{prop}\label{prop:exceptional}
A taut representation $\rho:G\to \mathbf O(V)$ of
a compact connected Lie group $G$ does not
have exceptional orbits.
\end{prop}

\Pf Assume the orbit through $p$ in $V$ is exceptional.
Then the slice representation of $G_p$ at $p$ has
a disconnected orbit. Arguing exactly as in the proof of
Proposition~\ref{prop:slice}, we see
that the distance function from $p$ to a principal
orbit through some regular point close to $p$ assumes its minimum value on a
disconnected set, contradicting tautness.
Thus $\rho$ cannot have an exceptional orbit. \EPf

\medskip

We close this subsection with a discussion of
taut reducible representations.

\begin{prop}\label{prop:red-taut}
 Let $\rho_1$ and $\rho_2$ be representations of a
compact connected Lie group $G$ with representation
spaces
$V_1$ and $V_2$ respectively. Assume that $\rho_1\oplus\rho_2$ is
$F$-taut. Then the restriction of $\rho_2$ to the
isotropy group
$G_{v_1}$ is taut for every $v_1\in V_1$.

Furthermore, we have $p(G(v_1,v_2);F)=p(Gv_1;F)\,p(G_{v_1}v_2;F)$ where
$p(M;F)$
denotes the Poincar\'e polynomial of $M$ with respect to the field $F$. 
In particular, $G_{v_1}v_2$ is connected
and 
$$b_1(G(v_1,v_2);F)=b_1(Gv_1;F)+b_1(G_{v_1}v_2;F) $$
where $b_1(M;F)$ denotes the first Betti number of $M$ with respect to
$F$.
\end{prop}

\Pf We can work with height functions instead of distance functions in this
proof since a submanifold contained in a
round sphere in a Euclidean space is taut if and only if all height
functions are perfect Morse functions, see
\cite{PT2}. Furthermore, in this situation the set of critical
points of a distance function will also occur as
the set of critical points of a height function, and vice versa. 
Fix $(v_1,v_2)\in V_1\oplus V_2$. Let
$a\in V_1$ be such that the height function
$h_a:Gv_1\to {\bf R}$ defined by
$h_a(v)=\langle a,v\rangle$ is a Morse function. We define the height
function $h_{(a,0)}: G(v_1,v_2)\to {\bf R}$ similarly.
The point
$(u_1,u_2)$ is a critical point of $h_{(a,0)}$ with index $i$ if and only if
$u_1$ is a critical point of $h_a$ with index
$i$. Hence the critical set $C$ on the critical level
$h_{(a,0)}(u_1,u_2)=h_a(u_1)$ is
$$
C=\{(w_1,w_2)\in G(v_1,v_2)\,|\, w_1=u_1\}=\{(u_1,w_2)\,|\, w_2\in
G_{u_1}v_2\}.
$$
Ozawa proves in \cite{Ozawa} that the set of critical points of a distance
function on a taut submanifold is a union over
nondegenerate critical submanifolds that are again taut submanifolds. It
follows that
$C$ is taut. The projection of $C$ into $V_2$, which coincides with
$G_{u_1}v_2$, is then also taut. We can choose $a$ such
that
$h_a$ is a Morse function of which $v_1$ is a critical point. It follows that
$G_{v_1} v_2$ is taut for every $v_1\in V_1$ and $v_2\in V_2$ 
and hence that the restriction
of $\rho_2$ to $G_{v_1}$ is taut.

Now fix again a point $(v_1,v_2)$. If $C\subset G(v_1,v_2)$ is an arbitrary
critical submanifold of $h_{(a,0)}$, then it is
diffeomorphic to $G_{v_1}v_2$. To see this let $\pi_1:G(v_1,v_2)\to Gv_1$ be
the projection onto the first factor. Then one
easily sees that 
$$\pi_1^{-1}(gv_1)=G_{gv_1}gv_2=gG_{v_1}g^{-1}(gv_2)=gG_{v_1}v_2.$$
Hence we have that $\pi_1$ is a $G$-equivariant
fibration. The critical submanifolds of
$h_{(a,0)}$ are fibers of $\pi_1$ and hence
diffeomorphic to each other.

 Since the orbit $G(v_1,v_2)$ is taut we have by \cite{Ozawa} that the
Morse-Bott
inequalities are equalities. Let $C_1,\dots, C_k$ be the critical manifolds
of $h_{(a,0)}$. Let $i(C_j)$ be the index of
$C_j$ for $j=1,\dots, k$. Then
$$
p(G(v_1,v_2))=\sum_{j=1}^k p(C_j)t^{i(C_j)}=p(G_{v_1}v_2)\sum_{j=1}^k
t^{i(C_j)}.
$$
since the critical submanifolds are all diffeomorphic to $G_{v_1}v_2$. We
have 
$$
p(Gv_1)=\sum_{j=1}^k t^{i(C_j)}
$$
since $h_a$ is a perfect Morse functions with critical points
$p_1=\pi_1(C),\dots, p_k=\pi_1(C_k)$ and the index of $p_j$ is
equal to $i(C_j)$. It follows that
$$
p(G(v_1,v_2))=p(Gv_1)\,p(G_{v_1}v_2).
$$
Multiplying out the Poincar\'e polynomials gives $b_0(G_{v_1}v_2)=1$ and
$$
b_1(G(v_1,v_2))=b_1(Gv_1)+b_1(G_{v_1}v_2).
$$
This finishes the proof of the proposition.  \EPf

\medskip

\begin{rmk}
\em The property that $\rho_2|G_{v_1}$ is taut in
Proposition~\ref{prop:red-taut} is a weakening of the property
that $\rho_2|G_{v_1}$ is orbit equivalent to $\rho_2$ 
for a variationally complete reducible representation $\rho_1\oplus\rho_2$ 
(cf.~Proposition~\ref{prop:red-vc}, part~(b)).
\end{rmk}

We give examples of how Proposition~\ref{prop:red-taut} can be used.

\begin{exs} 
\rm (i) Let $G=\SO n$ and let $\rho_1$ be the $\SO n$-conjugation
on the space $V_1$ of real traceless symmetric $n\times n$ matrices. Then
$\rho_1$ is taut since it is the isotropy
representation of the symmetric space ${\bf SL}(n,{\bf R}) /\SO n$. Let
$\rho_2$ be any other nontrivial  representation of
$\SO n$ with representation space $V_2$. Then $\rho_1\oplus\rho_2$ cannot be
taut if $n\ge 3$.  To see this, let $v_1\in V_1$ be
a regular point. Then
$G_{v_1}$ is the discrete group consisting of all diagonal matrices with
determinant one and entries $\pm1$ on the diagonal.
The kernel of $\rho_2$ is contained in the center of $\SO n$. Since $n\ge
3$, we see that $G_{v_1}$ cannot be contained in the
center of $\rho_2$. Hence there is an element $v_2\in V_2$ that is not fixed
by $G_{v_1}$. It follows that $G_{v_1}v_2$ is
disconnected. Now Proposition~\ref{prop:red-taut}
implies that $\rho_1\oplus\rho_2$ is not taut.

(ii) Now let $G$ be a compact connected simple Lie group of rank at least two and
let $\rho_1$ denote the adjoint representation of
$G$. We assume that $G$ is simply connected. Let
$\rho_2$ be any other nontrivial representation of $G$. Then
$\rho_1\oplus\rho_2$ is not taut. To see this let
$T$ be a maximal torus in $G$. We denote the representation spaces of
$\rho_1$ and $\rho_2$ by $V_1$ and $V_2$ respectively.
There is a regular element $v_1\in V_1$ with $G_{v_1}=T$. The restriction of
$\rho_2$ to $T$ has a discrete kernel that is
contained in the center of $G$. If $v_2\in V_2$ is a $T$-regular point
then the isotropy subgroup $T_{v_2}$ coincides with the
kernel of $\rho_2|T$. Hence $G_{v_1}v_2$ is diffeomorphic to $T$ and it
follows that $b_1(G_{v_1}v_2;F)$ is equal to the rank
of
$G$. In particular $b_1(G_{v_1}v_2;F)\ge 2$.
Now notice that the isotropy group of $(v_1,v_2)$ is also $T_{v_2}$.
Hence $\pi_1(G(v_1,v_2))=T_{v_2}$ which implies $H_1(G(v_1,v_2);{\bf
Z})=T_{v_2}$ since $T_{v_2}$ is Abelian. If $G\ne \Spin {4k}$
then the center of
$G$ is a cyclic group and it follows that
$b_1(G(v_1,v_2);F)\le 1$.
If $G=\Spin{4k}$, then $k\ge2$ and we get $b_1(G_{v_1}v_2;F)=2k\ge 4$;
since the center of $\Spin{4k}$ is ${\bf Z}_2\times {\bf Z}_2$,
we have $b_1(G(v_1,v_2);F)\le 2$. In either case, 
$b_1(G_{v_1}v_2;F)>b_1(G(v_1,v_2);F)$ which implies by
Proposition~\ref{prop:red-taut} that $\rho_1\oplus\rho_2$ is not taut. 
\end{exs}

\subsection{A reduction principle}\label{subsec:red}

We now discuss a reduction principle which will be used
later to prove that certain representations are or are not taut.
Let $\rho:G\to \mathbf O(V)$ be an orthogonal linear action of a compact
Lie group $G$ which is not assumed to be connected. Denote by $H$ a fixed
principal isotropy subgroup of the $G$-action on $V$ and let $V^H$ be the
subspace of $V$ that is left pointwise fixed by the action of
$H$. Let $N$ be the normalizer of $H$ in $G$.
Then the group $N/H$ acts on $V^H$ with trivial principal
isotropy subgroup. Moreover, the following result is known
(\cite{G-S,Luna,LR,Schwartz,SS,Straume1}):

\begin{thm}[Luna-Richardson]
The inclusion $V^H\to V$ induces a stratification preserving
homeomorphism between orbit spaces 
\setcounter{equation}{\value{thm}}
\stepcounter{thm}
\begin{equation}\label{eqn:homeo}
V^H/N\to V/G.
\end{equation}
\end{thm}

The injectivity of the map~(\ref{eqn:homeo}) means that 
$Np=Gp\cap V^H$ for $p\in V^H$. In particular, the $H$-fixed point 
set of a $G$-orbit is a smooth manifold.
Observe also that for a regular point $p\in V^H$
the normal space to the principal orbit $M=Gp$ at~$p$
is contained in $V^H$, because the slice representation at $p$
is trivial.

\begin{lem}\label{lem:2}
Let $p\in V^H$ be a regular point and consider a normal vector $\xi$
to $M=Gp$ at $p$. Then the Weingarten operator of $M$ at $\xi$
restricts to the Weingarten operator of $M^H=M\cap V^H$ at $\xi$,
in symbols, $A_\xi^M|_{T_pM^H}=A_\xi^{M^H}$.
\end{lem}

\Pf Let $v\in T_pM^H$. Consider an extension $\tilde\xi(t)$ of $\xi$
to a normal vector field along a curve $\gamma(t)$ in~$M^H$ with $\dot\gamma(0)=v$.
Then $A_\xi^Mv=-\tilde\xi'(0)^T$,
where ``${}^T$'' denotes the orthogonal projection onto $T_pM$.
Now $\tilde\xi(t)\in V^H$, so the 
derivative $\tilde\xi'(0)\in V^H$. 
The normal component of $\tilde\xi'(0)$ is already in $V^H$, 
so its component in~$T_pM$ is also in $V^H$. Hence 
$A_\xi^{M^H}v=-\tilde\xi'(0)^T$ which proves the claim. \EPf

\begin{lem}\label{lem:3}
Let $p$, $q\in V^H$ and suppose that $q$ is a 
regular point for $G$. Consider $M=Gp$ and let $L_q:M\to\mathbf
R$ be the distance function. Then the critical set of $L_q$ is contained
in $M^H=M\cap V^H$, namely, $\mbox{Crit($L_q$)}=\mbox{Crit($L_q|_{M^H}$)}$.
\end{lem}

\Pf Let $x\in M$ be a critical point of $L_q:M\to\mathbf R$. Then $x=gp$ for some
$g\in G$ and $v=q-x$ is a normal vector to $M$ at $x$.
Consider first the case where $p$ is a regular point for $G$.
Then the isotropy subgroup at $x$ is principal, so $G_x=gHg^{-1}$
acts trivially on the normal space. In particular $gHg^{-1}v=v$,
which implies that $gHg^{-1}q=q$. Now $gHg^{-1}\subset G_q=H$ and therefore
$g\in N$. Thus $x=gp\in V^H$ and $x$ is a critical point for the
restriction $L_q:M^H\to\mathbf R$. Consider now the case where 
$p$ is not a regular point for $G$. Then there is $\epsilon>0$ sufficiently small
such that $y=x+\epsilon v$ is a regular point for $G$. 
Now $y$ is a critical point for $L_q:Gy\to\mathbf R$, since 
a line segment meeting one orbit orthogonally is orthogonal to 
all orbits it meets.  
By the previous case $y\in V^H$. Hence $x\in V^H$ and
again $x$ is a critical point for the restriction $L_q:M^H\to\mathbf R$.

On the other hand, let $x\in M^H$ be a critical point of
$L_q:M^H\to\mathbf R$. Consider first the case where $p$ is a regular point
for $G$. Then $x$ is also a regular point for $G$
so that the normal spaces $N_x(M^H)$ and $N_x(M)$ coincide. 
It follows that $x$ is a critical point of $L_q:M\to\mathbf R$. 
Consider now the case where $p$ 
is not a regular point for $G$. Then there is $\epsilon>0$ sufficiently small
such that $y=x+\epsilon v\in V^H$ is a regular point for $G$, where $v=q-x$. 
Since $Gy\cap V^H=Ny$, we get that $y$ is a critical point for $L_q:Ny\to\mathbf R$.
By the previous case $y$ is a critical point for $L_q:Gy\to\mathbf R$. 
Hence again $x$ is a critical point of $L_q:M\to\mathbf R$.\\ \mbox{} \EPf

\medskip

For the next lemma and proposition 
we found inspiration in~\cite{Duistermaat}.

\begin{lem}\label{lem:4}
Suppose there is a subgroup $L\subset H$ which is a 
finitely iterated $\mathbf Z_2$-extension of the identity and such that 
the fixed point sets $V^L=V^H$. Let $M=Gp$ be a principal orbit.
Then the sum of the Betti numbers (with respect to $\mathbf Z_2$
coefficients) of $M^H$ is less
than or equal the sum of the Betti numbers of $M$, in symbols
$\beta(M^H;\mathbf Z_2)\leq\beta(M;\mathbf Z_2)$.
\end{lem}

\Pf We can write $L=\{h_1,\ldots,h_n\}\subset H$
such that $h_{i+1}$ normalizes $H_i$ and $h_{i+1}^2 \in H_i$
for $i:0,\ldots,n-1$, where $H_0$ is the identity group 
and $H_i$ is the subgroup of $H$ generated by $h_1,\ldots,h_i$. 
Now each transformation $h_{i+1}$ has order $2$
in the fixed point set $M^{H_i}=M\cap V^{H_i}$, so by a theorem of 
Floyd~(\cite{Floyd}) 
\[ \beta(M^L;\mathbf Z_2)\leq\beta(M^{H_{n-1}};\mathbf Z_2)\leq
\ldots\leq\beta(M^{H_0};\mathbf Z_2). \]
But $M^L=M^H$ and $M^{H_0}=M$ and this completes the proof. \EPf

\begin{prop}\label{prop:suf-taut}
Suppose there is a subgroup $L\subset H$ which is a 
finitely iterated $\mathbf Z_2$-extension of the identity and such that 
the fixed point sets $V^L=V^H$. Suppose also
that the reduced representation 
$\bar\rho:(N/H)^0\to \mathbf  O(V^H)$ is $\mathbf Z_2$-taut
where $(N/H)^0$ denotes the connected component of the identity of
$N/H$. Then $\rho:G\to \mathbf O(V)$ is $\mathbf Z_2$-taut.
\end{prop}

\Pf Let $M=Gp$ be an arbitrary orbit. We shall prove that $M$ is $\mathbf Z_2$-taut. 
Since $V^H$ intersects all orbits, we may assume $p\in V^H$. 
Suppose $L_q:M\to\mathbf R$ is a Morse distance function.
We must show that $L_q$ is $\mathbf Z_2$-perfect, and in fact
it is enough to do that for $q$ in a dense subset of $V$, so
we may take $q$ to be a regular point.
Since $G$ acts on $V$ by isometries, we may moreover assume that $q\in V^H$. 
If $\bar\rho$ is $\mathbf Z_2$-taut, then each connected component of
$M^H$ is $\mathbf Z_2$-taut,
hence $M^H$ itself is $\mathbf Z_2$-taut. Note that 
$L_q|_{M^H}$ is a Morse function. Hence $L_q|_{M^H}$ is
$\mathbf Z_2$-perfect. Now Lemmas~\ref{lem:3} and~\ref{lem:4} and the 
Morse inequalities imply that
\[ \beta(M^H;\mathbf Z_2)\leq\beta(M;\mathbf Z_2)\leq 
    \#\mbox{Crit($L_q$)}=\#\mbox{Crit($L_q|_{M^H}$)},
 \]
where $\#A$ denotes the cardinality of the finite set $A$. 
It follows that $L_q$ is $\mathbf Z_2$-perfect. \EPf

%name of the file: intermezzo.tex

\section{A necessary condition for a representation to be of class
  $\mathcal O^2$}\label{sec:necessary}
\setcounter{thm}{0}

Let $\pi$ be a complex representation of the compact 
connected Lie group $G$ on a finite-dimensional vector space.
We say that $\pi$ is of \emph{real type} if it comes from a representation 
of $G$ on a real vector space by extension of scalars, and we say 
that $\pi$ is of \emph{quaternionic type} if it comes from a representation
of $G$ on a quaternionic vector space by restriction of scalars. 
If $\pi$ is neither of real type nor of quaternionic type, we say that
$\pi$ is of \emph{complex type}\footnote{It is useful to remark that
$\pi$ is of real type (resp., quaternionic type) if and only
if there exists an invariant real (resp., quaternionic)
structure on the representation space, i.~e.~a $\mathbf C$-conjugate
linear, $G$-invariant endomorphism $\epsilon$ of the
representation space such that $\epsilon^2=1$ (resp., 
$\epsilon^2=-1$).}. 

Now it is known that the finite-dimensional 
real irreducible representations $\rho$ of $G$ 
fall into one of the following disjoint classes:
\begin{enumerate}
\item[(a)] the complexification $\rho^c$ is irreducible and $\rho^c=\pi$ is 
a complex representation of real type;
\item[(b)] the complexification $\rho^c$ is reducible and 
$\rho^c=\pi\oplus\pi$ where $\pi$ is a complex irreducible representation of 
quaternionic type;
\item[(c)] the complexification $\rho^c$ is reducible and 
$\rho^c=\pi\oplus\pi^*$ where $\pi$ is a complex irreducible representation of 
complex type and $\pi^*$ is not equivalent to $\pi$ 
(where $\pi^*$ denotes the dual representation of $\pi$).
\end{enumerate}
The relation between $\rho$ and $\pi$ is that $\rho$ is a real form 
of $\pi$ in the first case ($\rho^c=\pi$), but $\rho$ is $\pi$ viewed as a real
representation in the other two cases ($\rho=\pi^r$). We shall call $\rho$ 
of \emph{real}, \emph{quaternionic} or \emph{complex type} according to whether the associated 
$\pi$ is of real, quaternionic or complex type.
Note also that $\pi$ is self-dual precisely in the first two cases. 

Suppose now that $G$ is semisimple, let $\mathfrak g$ denote its Lie 
algebra and $\Lg^c$ its complexification. Write $\Delta$ 
for the root system of $\Lg^c$ with respect to a chosen
Cartan subalgebra, $\Delta^+$ for the positive root system 
with respect to an ordering of the roots, and
 $\mathcal S=\{\alpha_1,\dots,\alpha_r\}$ for the 
corresponding simple root system. Let $\lambda_1,\ldots,\lambda_r$ be the fundamental
highest weights defined by the relations
$2(\lambda_i,\alpha_j)/(\alpha_j,\alpha_j)=\delta_{ij}$ where 
$(\,,\,)$ is the Cartan-Killing form. The Theorem of the Highest Weight of 
\'E.~Cartan states that the complex irreducible representations 
of $\Lg^c$ are parametrized by their highest weights, and these are exactly the linear 
combinations $\lambda=\sum_{i=1}^rm_i\lambda_i$ for $m_i\in\{0,1,2,\ldots\}$. 
The following proposition is a useful criterium of Dadok 
to decide for which 
$\lambda$ the corresponding representation $\pi_{\lambda}$ is of 
real, quaternionic or complex type. Recall that the roots $\alpha$, 
$\beta\in\Delta$ are called \emph{strongly orthogonal} if $\alpha\pm\beta$ 
is not a root. 

\begin{prop}[\cite{Dadok}]\label{prop:Dadok}
There is a maximal subset $\mathcal B=\{\beta_1,\ldots,\beta_s\}\subset\Delta^+$ 
of strongly orthogonal roots such that:
\begin{enumerate}
\item[(a)] we have $s_0=s_{\beta_1}\cdots s_{\beta_s}$ 
is the Weyl group element that maps 
the positive Weyl chamber into its negative;
\item[(b)] the representation 
$\pi_\lambda$ is of complex type if and only if $\lambda$ does not 
belong to the real span of $\mathcal B$;
\item[(c)] the representation $\pi_\lambda$ is of real type (resp.~quaternionic type)
if and only if $\lambda$ belongs to the real span of $\mathcal B$ 
and
\[ k(\lambda)=\sum_{i=1}^s\frac{(\lambda,\beta_i)-(s_0\lambda,\beta_i)}
{(\beta_i,\beta_i)} \]
is an even (resp.~odd) integer.
\end{enumerate}
\end{prop}

\begin{rmk}
\em
For a simple Lie algebra the set $\mathcal B$ can be constructed as
follows, as is explained in~\cite{Dadok} and will be assumed throughout our paper
(cf.~Appendix~A). Let $\beta_1$ be the highest root. 
The root system $\{\alpha\in\Delta:(\alpha,\beta_1)=0\}$ is either
irreducible or equals $\{\pm\zeta_1\}\cup\Delta_1$, with $\Delta_1$ 
being irreducible and $\zeta_1\in\Delta^+$. In the former case
set $\beta_2$ equal to the highest root of $\Delta_1$ (with the inherited order 
from $\Delta$), and proceed by induction.
In the latter case set $\beta_2=\zeta_1$ and $\beta_3$ equal to the highest root 
of $\Delta_1$, and proceed by induction.
\end{rmk}
 
We continue to assume that $G$ is semisimple. 
Let $H_\alpha$, $\alpha\in\Delta$, be the coroots of $\mathfrak g^c$. 
It is possible to choose root vectors $X_\alpha$ for $\mathfrak g^c$,
$\alpha\in\Delta$, such that the compact real form $\mathfrak g$ 
is spanned by
\setcounter{equation}{\value{thm}}
\stepcounter{thm}
\begin{equation}\label{eq:basis}
iH_\alpha,\quad X_\alpha-X_{-\alpha},\quad 
i(X_\alpha+X_{-\alpha}),\qquad\mbox{where $\alpha\in\Delta^+.$}
\end{equation} 
Now if $\pi$ is a complex representation of $G$,
we have that the adjoint map $\pi(X_\alpha)^*=\pi(X_{-\alpha})$
with respect to any $G$-invariant Hermitian product on the representation
space. Moreover, if $\pi$ is of real (resp.~quaternionic)
type and $\epsilon$ is an invariant real (resp.~quaternionic) structure 
on the representation space, then
$\pi(X_\alpha)\epsilon=\epsilon\pi(X_{-\alpha})$. 

By compactness of~$G$, any real representation is equivalent to
an orthogonal one. The following proposition, which is stated as
a remark in~\cite{Dadok}, p.~128, 
and the ensuing lemma,
which is a refinement of a result in~\cite{C-Th},
are precisely the ingredients we need to establish a necessary condition
for an orthogonal representation of~$G$ to be of class~$\mathcal O^2$.

\begin{prop}[\cite{Dadok}]\label{prop:D}
Let ${\mathcal U}^k(\Lg^c)$ be the $k$th level in the natural filtration 
of the universal enveloping algebra of $\Lg^c$. 
Fix $\pi_\lambda$, the irreducible 
representation of $\Lg^c$ with highest weight $\lambda$,
and fix $v_\lambda$, a highest weight vector. 
Let 
\[ n_i=\frac{(\lambda,\beta_i)-(s_0\lambda,\beta_i)}{(\beta_i,\beta_i)}, \] 
for $i=1,\ldots,s$. Then the element $s_0v_\lambda$
is in ${\mathcal U}^{k(\lambda)}(\Lg^c)v_{\lambda}$, but not
in ${\mathcal U}^{l}(\Lg^c)v_{\lambda}$ for any $l<k(\lambda)$, where
$s_0v_\lambda$ is given as
\[ s_0v_\lambda=\pi(X_{-\beta_1})^{n_1}\cdots
          \pi(X_{-\beta_s})^{n_s}v_{\lambda}. \]
\end{prop}

There is no proof of the above proposition in~\cite{Dadok}.
In~\cite{CF}, p.~270, there is an attempt to prove 
the proposition which in our view contains serious gaps.  

\medskip

\Pf It is enough to consider the case where $\mathfrak g^c$ 
is a complex simple Lie algebra. 
We first claim that the element $s_0v_\lambda$ is not zero. 
In fact, consider the complex subalgebra 
$\mathfrak k\subset \mathfrak g^c$ generated by the root spaces of $\mathfrak g^c$
corresponding to $\pm\beta_1,\ldots,\pm\beta_s$. Then $\mathfrak k$ is 
isomorphic to the direct product of $s$ copies of $\mathfrak{sl}(2,\mathbf
C)$. We restrict $\pi_\lambda$ to $\mathfrak k$ and let $U_\lambda$ be the 
unique irreducible $\mathfrak k$-module generated by $v_\lambda$. 
Since $\lambda-s_0\lambda=\sum_{i=1}^sn_i\beta_i$, it is now clear
that $U_\lambda$ is the representation space of
$\mathfrak k:\;\Ai{\!\!n_1}\otimes\cdots\otimes\Ai{\!\!n_s}$ and
that $s_0v_\lambda\in U_\lambda$. Therefore, our claim 
is reduced to the case of the $(n+1)$-dimensional 
complex irreducible representation of $\mathfrak{sl}(2,\mathbf
C)$, namely $\Ai{\!n}$, which is immediate to verify. 

It is obvious that $s_0v_\lambda\in\mathcal U^{k(\lambda)}(\mathfrak
g^c)v_\lambda$. We next prove that 
$s_0v_\lambda\notin\mathcal U^l(\mathfrak
g^c)v_\lambda$ for $l<k(\lambda)$ by contradiction. 
In fact,
enumerate the positive roots $\Delta^+=\{\alpha_1,\ldots,\alpha_t\}$ 
and suppose $s_0v_\lambda\in\mathcal U^l(\mathfrak
g^c)v_\lambda$ for $l<k(\lambda)$. It follows from the
Poincar\'e-Birkhoff-Witt theorem that we can write 
\[ s_0v_\lambda=
\sum_{p_1,\ldots,p_t\geq0} c_{p_1\cdots p_t}X_{-\alpha_1}^{p_1}\cdots
                      X_{-\alpha_t}^{p_t}v_{\lambda}, \]
for some complex constants $c_{p_1\cdots p_t}$. Now, using 
that $X_\alpha^*=X_{-\alpha}$,
\begin{eqnarray*} 
 0 &\neq & (s_0v_\lambda,s_0v_\lambda) \\
   &=& \sum_{p_1,\ldots,p_t\geq0}
     c_{p_1\cdots p_t} (X_{-\alpha_1}^{p_1}\cdots
     X_{-\alpha_t}^{p_t}v_{\lambda}, X_{-\beta_1}^{n_1}\cdots
     X_{-\beta_s}^{n_s}v_{\lambda} ) \\
   &=& \sum_{p_1,\ldots,p_t\geq0}
     c_{p_1\cdots p_t} (X_{\beta_s}^{n_s}\cdots
     X_{\beta_1}^{n_1}X_{-\alpha_1}^{p_1}\cdots
     X_{-\alpha_t}^{p_t}v_{\lambda}, v_{\lambda} ). 
\end{eqnarray*}
Recall that the highest weight space is one-dimensional; 
hence we can write
\[ X_{\beta_s}^{n_s}\cdots
     X_{\beta_1}^{n_1}X_{-\alpha_1}^{p_1}\cdots
     X_{-\alpha_t}^{p_t}v_{\lambda} = cv_{\lambda}, \]
for some nonnegative integers $p_1,\ldots,p_t$ such that
$k(\lambda)=n_1+\cdots+n_s>p_1+\cdots+p_t=l$
and for some nonzero complex constant $c$.
Note that $n_1\beta_1+\cdots+n_s\beta_s=p_1\alpha_1+\cdots+p_t\alpha_t$.
The contradiction we are aiming at now follows from the 
following claim.
\begin{claim}\label{claim:1}
If $N>M$ then 
\[ X_{\delta_N}\cdots X_{\delta_1}X_{-\gamma_1}\cdots
X_{-\gamma_M}v_\lambda=0, \]
where $\delta_1,\ldots,\delta_N\in\mathcal B$,
and $\gamma_1,\ldots,\gamma_M\in\Delta^+$.
\end{claim}
We proceed by induction on the integer $N$. 
The cases $N=1$ and $N=2$ are trivial (use that
$v_\lambda$ is a highest weight vector for $N=1$ and that
the sum of two roots in~$\mathcal B$ is not a root for $N=2$).
Assume the assertion is true for some 
$N-1\geq2$ and let us prove it for $N$. 
Since $X_{\beta_1},\ldots,X_{\beta_s}$ pairwise commute, we may assume that
$\delta_1=\beta_i$ where $i\leq j$ for every $j=1,\ldots,s$ which satisfies 
$\beta_j=\delta_k$ for some $k=1,\ldots,N$.  
\begin{claim}\label{claim:2}
For each $j=1,\ldots,M$, if $\delta_1-\gamma_j$ is a root
then it is a positive root.
\end{claim}
Indeed we have $\delta_1=\beta_i$ where $i$ is as above. 
If $i=1$, then we are done because $\beta_1$ is the highest root. 
If not, 
$0=(\beta_1,\delta_1+\cdots+\delta_N)=(\beta_1,\gamma_1+\cdots+\gamma_M)$.
Now $(\beta_1,\gamma_1)=-(\beta_1,\gamma_2+\cdots+\gamma_M)$, 
and since $\beta_1$ is the highest root we have 
$(\beta_1,\gamma_1)\geq0,\ldots,(\beta_1,\gamma_M)\geq0$. 
It follows that $(\beta_1,\gamma_1)=0$ and
$(\beta_1,\gamma_2+\cdots+\gamma_M)=0$.
An easy induction argument shows that $(\beta_1,\gamma_j)=0$ for
$j=1,\ldots,M$. Consider now the root system 
$\{\alpha\in\Delta:(\alpha,\beta_1)=0\}$. 
The first case occurs when it is irreducible.
Then $\beta_2$ is by definition its highest root 
so that $\beta_2\geq\gamma_j$ for $j=1,\ldots,M$. 
If $i=2$ we are done. If not, we see that 
$(\beta_2,\gamma_j)=0$ for $j=1,\ldots,M$ and
we proceed by induction.
The second case occurs when the above root system is of the
form $\{\pm\zeta_1\}\cup\Delta_1$, with $\Delta_1$ 
being irreducible and $\zeta_1\in\Delta^+$. 
Here $\beta_2=\zeta_1$ and $\beta_3$ is the highest root of 
$\Delta_1$. If $\gamma_j\in\Delta_1$ for $j=1,\ldots,M$ then 
$(\beta_2,\gamma_j)=0$ for $j=1,\ldots,M$ so that $i\geq3$
and we proceed by induction as above. On the other hand, 
if $\zeta_1=\gamma_{j_0}$ for some $j_0$, then 
$i=2$. In this case, for each $j=1,\ldots,M$ either $\beta_2=\gamma_j$ 
or $-\gamma_j\in\Delta_1$. So $\beta_2-\gamma_j$ is never a root. 
This completes the proof of Claim~\ref{claim:2}.
\begin{claim}\label{claim:3}
We have
\[ X_{\delta_N}\cdots X_{\delta_2}X_{-\gamma_1}\cdots
   X_{-\gamma_{j-1}}H_{\zeta}X_{-\gamma_{j+1}}\cdots
   X_{-\gamma_M}v_\lambda=0, \]
where $H_{\zeta}$ is the coroot vector corresponding to $\zeta\in\Delta^+$ 
and $j=1,\ldots,M$.
\end{claim}
In order to prove Claim~\ref{claim:3}, proceed by induction on $j$.
The initial case is $j=M$, which follows from the induction hypothesis
on~$N$, since $H_\zeta v_\lambda=(\zeta,\lambda)v_\lambda$. 
Next write
\setcounter{equation}{\value{thm}}
\stepcounter{thm}
\begin{eqnarray}\label{eqn:sum0}
\lefteqn{X_{\delta_N}\cdots X_{\delta_2}X_{-\gamma_1}\cdots
   X_{-\gamma_{j-1}}H_{\zeta}X_{-\gamma_{j+1}}\cdots
   X_{-\gamma_M}v_\lambda=} \nonumber \\ 
 & & X_{\delta_N}\cdots X_{\delta_2}X_{-\gamma_1}\cdots
   X_{-\gamma_{j-1}}X_{-\gamma_{j+1}}H_\zeta\cdots
   X_{-\gamma_M}v_\lambda \\
 &+& X_{\delta_N}\cdots X_{\delta_2}X_{-\gamma_1}\cdots
   X_{-\gamma_{j-1}}[H_\zeta,X_{-\gamma_{j+1}}]\cdots
   X_{-\gamma_M}v_\lambda. \nonumber
\end{eqnarray}
The first summand on the right hand side of~(\ref{eqn:sum0})
is zero by the induction hypothesis on~$j$, and 
the second summand is zero because 
$[H_\zeta,X_{-\gamma_{j+1}}]=-(\zeta,\gamma_{j+1})X_{-\gamma_{j+1}}$
so that we can use the induction hypothesis on~$N$.
This proves Claim~\ref{claim:3}.
\begin{claim}\label{claim:4}
We have
\[ X_{\delta_N}\cdots X_{\delta_2}X_{-\gamma_1}\cdots
   X_{-\gamma_{j-1}}X_{\zeta}X_{-\gamma_{j+1}}\cdots
   X_{-\gamma_M}v_\lambda=0, \]
where $X_{\zeta}$ is the root vector corresponding to $\zeta\in\Delta^+$ 
and $j=1,\ldots,M$.
\end{claim}
In order to prove Claim~\ref{claim:4}, proceed by induction on $j$.
The initial case is $j=M$, which is trivial because $v_\lambda$ is a 
highest weight vector. Next write
\setcounter{equation}{\value{thm}}
\stepcounter{thm}
\begin{eqnarray}\label{eqn:sum}
\lefteqn{X_{\delta_N}\cdots X_{\delta_2}X_{-\gamma_1}\cdots
   X_{-\gamma_{j-1}}X_{\zeta}X_{-\gamma_{j+1}}\cdots
   X_{-\gamma_M}v_\lambda=} \nonumber \\ 
 & & X_{\delta_N}\cdots X_{\delta_2}X_{-\gamma_1}\cdots
   X_{-\gamma_{j-1}}X_{-\gamma_{j+1}}X_\zeta\cdots
   X_{-\gamma_M}v_\lambda \\
 &+& X_{\delta_N}\cdots X_{\delta_2}X_{-\gamma_1}\cdots
   X_{-\gamma_{j-1}}[X_\zeta,X_{-\gamma_{j+1}}]\cdots
   X_{-\gamma_M}v_\lambda. \nonumber
\end{eqnarray}
The first summand on the right hand side of~(\ref{eqn:sum})
is zero by the induction hypothesis on~$j$.
The second summand is also zero for the following reasons.
If $\zeta=\gamma_{j+1}$ we use Claim~\ref{claim:3}.  
If not, then either $\zeta-\gamma_{j+1}$ is not a root and then 
$[X_\zeta,X_{-\gamma_{j+1}}]=0$, or else 
$\zeta-\gamma_{j+1}$ is a root. In the latter case, 
if it is a positive root we can use the induction hypothesis 
on $j$, and if it is a negative root then we can use 
the induction hypothesis on $N$. This completes the proof
of Claim~\ref{claim:4}. 
We finally turn to the proof of 
Claim~\ref{claim:1}. We can write
\setcounter{equation}{\value{thm}}
\stepcounter{thm}
\begin{eqnarray}\label{eqn:sum2}
\lefteqn{X_{\delta_N}\cdots X_{\delta_1}X_{-\gamma_1}\cdots
   X_{-\gamma_M}v_\lambda=} \nonumber \\
  & & \sum_{j=1}^M X_{\delta_N}\cdots X_{\delta_2}X_{-\gamma_1}\cdots
   X_{-\gamma_{j-1}}[X_{\delta_1},X_{-\gamma_j}]X_{-\gamma_{j+1}}\cdots
   X_{-\gamma_M}v_\lambda \\
  &+& X_{\delta_N}\cdots X_{\delta_2}X_{-\gamma_1}\cdots X_{-\gamma_M}X_{\delta_1}v_\lambda.  \nonumber
\end{eqnarray}
The second summand in the right hand side of~(\ref{eqn:sum2})
is zero because $v_\lambda$ is a highest weight vector.
The first summand is also zero because if $\delta_1=\gamma_j$
we can apply Claim~\ref{claim:3}. Otherwise,
either $\delta_1-\gamma_j$ is not a root and then 
$[X_{\delta_1},X_{-\gamma_j}]=0$ or else
$\delta_1-\gamma_j$ is a root and then 
it is a positive root by Claim~\ref{claim:2} so that we can apply 
Claim~\ref{claim:4}. This completes the proof of Claim~\ref{claim:1}
and the proof of the proposition. \EPf

\begin{lem}\label{lem:CT}
Let $\pi_\lambda$ be the complex irreducible representation 
of $G$ with highest weight $\lambda$ and representation space
$V_\lambda$. Let $\mu$ be a weight of $\pi_\lambda$ and 
fix a weight vector $v_\mu$. 
Let $\rho$ denote a real form of $\pi_\lambda$ in case $\pi_\lambda$ is of
real type, or the realification of $\pi_\lambda$ in case it is either of
quaternionic or of complex type. 
Suppose that $\rho$ is of class $\mathcal O^2$. 
\begin{enumerate}
\item[(a)] If $\pi_\lambda$ is of real type, then 
\[\mathcal U^2(\Lg^c)v_\mu+\mathcal U^2(\Lg^c)\epsilon(v_\mu)=V_\lambda, \]
where $\epsilon$ is the 
real structure on $V_\lambda$ defined by $\rho$. 
In particular, by taking $\mu=\lambda$, we find that
\[\mathcal U^2(\Lg^c)v_\lambda+\mathcal
U^2(\Lg^c)\epsilon(v_\lambda)=V_\lambda. \qquad(C_2)\]
\item[(b)] If $\pi_\lambda$ is of complex or quaternionic type, then 
$\mathcal U^2(\Lg^c)v_\mu=V_\lambda$.
\end{enumerate}
\end{lem}

\Pf (a)~Consider the real vector  $p=v_\mu+\epsilon v_\mu$. 
The second osculating space $\mathcal O^2_p(Gp)$ is spanned
over $\mathbf R$ by $Xp$, $XYp$, 
where $X$, $Y\in\mathfrak g$. 
Taking linear combinations with complex coefficients of these
vectors, we can write that the complexification 
\[ \mathcal O^2_p(Gp)^c\subset\mathcal U^2(\mathfrak g^c)v_\mu+
\mathcal U^2(\mathfrak g^c)\epsilon(v_\mu) \subset V_\lambda. \]
But $\rho$ of class $(\mathcal O^2)$ implies that 
$\mathcal O^2_p(Gp)^c=V_\lambda$. 

(b)~Let $p=v_\mu$. We have that $\mathcal O^2_p(Gp)\subset\mathcal U^2(\mathfrak
g^c)v_\mu\subset V_\lambda$ as real vector spaces
and the proof is similar as in~(a). \EPf

\medskip

We now state the main result of this section.
Notice that the same result is claimed in~\cite{CF}, p.~271,
with an attempt of a proof which in our opinion 
is not satisfactory. 

\begin{prop}\label{prop:k}
Let $\rho$ be a real (orthogonal) irreducible representation of a compact 
connected semisimple Lie group $G$ with complexified Lie algebra
$\Lg^c$, and let $\pi_\lambda$ be the associated complex irreducible
representation. Suppose that $\rho$ is of class $\mathcal O^2$. 
\begin{enumerate}
\item[(a)] If $\rho$ is of quaternionic type, then $k(\lambda)=1$.
\item[(b)] If $\rho$ is of complex type, then $k(\lambda)=1$, $2$. 
\item[(c)] If $\rho$ is of real type, then $k(\lambda)=2$, $4$. 
\end{enumerate}
\end{prop}

\Pf First consider $\rho$ to be of quaternionic or complex type,
i.~e.~$\rho^c=\pi_\lambda\oplus\pi_\lambda^*$. 
Then $\rho$ of class $\mathcal O^2$ forces 
$\mathcal U^2(\Lg^c)v_\lambda=V_\lambda$ (Lemma~\ref{lem:CT})
and then $s_0v_\lambda\in\mathcal U^2(\Lg^c)v_\lambda$ implies
that $k(\lambda)\leq2$ (Proposition~\ref{prop:D}). 
From this follow~(a) and~(b).

Now take $\rho$ to be of real type. 
Let $v_\lambda$ be a highest weight vector of $\rho^c=\pi_\lambda$. 
Since $\pi_\lambda$ is self-dual,
\[ \mbox{lowest weight of }\pi_\lambda=-\mbox{highest weight of
}\pi_\lambda^*=-\lambda=s_0\lambda. \] Also, there exists a $\mathbf C$-conjugate linear,
$G$-invariant involution $\epsilon$ of $V_\lambda$ and 
$\epsilon(v_\lambda)$ is a lowest weight vector. We shall assume 
$k(\lambda)\geq6$ and derive a contradiction. 
Since $k(\lambda)$ is even, it is possible to write
\[ \frac{k(\lambda)}2=n_1+\ldots+n_{i_0}+m \]
for some integers $0\leq i_0<s$, $0\leq m<n_{i_0+1}$,
where $n_i$ is as in Proposition~\ref{prop:D}. Set 
\[ \mu=\lambda-n_1\beta_1-\ldots-n_{i_0}\beta_{i_0}-m\beta_{i_0+1}. \]
Note that $\mu$ is a weight: a $\mu$-weight vector is
\[ u=X_{-\beta_1}^{n_1}\ldots
          X_{-\beta_{i_0}}^{n_{i_0}}
          X_{-\beta_{i_0+1}}^m v_\lambda. \]
Use Lemma~\ref{lem:CT} to decompose $u=u_1+u_2$ where
$u_1\in\mathcal U^2(\Lg^c)v_\lambda$ and 
$u_2\in\mathcal U^2(\Lg^c)\epsilon(v_\lambda)$. 
First, assume that $u_2\neq0$. Since $\frac{k(\lambda)}2\geq3$,
we have
\setcounter{equation}{\value{thm}}
\stepcounter{thm}
\begin{equation}\label{eq:u2}
 u_2=u-u_1\in \mathcal U^{\frac{k(\lambda)}2}(\Lg^c)v_\lambda + \mathcal U^2(\Lg^c)v_\lambda
          \subset \mathcal U^{\frac{k(\lambda)}2}(\Lg^c)v_\lambda.
\end{equation} 
Moreover, it is clear that $u_2$ may be assumed to
be a $\mu$-weight vector (because a component of $u_2$ in a different 
weight space has to cancel with the corresponding component of $u_1$ 
in the same weight space) and that
\[ u_2=\sum_{\gamma,\delta\in\Delta^+}
       c_{\gamma,\delta}X_{-\gamma}X_{-\delta}\epsilon(v_\lambda) \]
for some complex constants $c_{\gamma,\delta}$ (we do not need to consider
terms of first order in the sum because $\lambda-\mu$ cannot be a root). Therefore,
\[ 0 \neq (u_2,u_2) = (u_2, \sum_{\gamma,\delta\in\Delta^+}
       c_{\gamma,\delta}X_{-\gamma}X_{-\delta}\epsilon(v_\lambda) )
 = (\sum_{\gamma,\delta\in\Delta^+}
       \bar{c}_{\gamma,\delta}X_{\delta}X_{\gamma}u_2,\epsilon(v_\lambda)).\]
This shows that \[ \sum_{\gamma,\delta\in\Delta^+}
       \bar{c}_{\gamma,\delta}X_{\delta}X_{\gamma}u_2 \]
is a nonzero multiple of the lowest weight vector $\epsilon(v_\lambda)$.
This, combined with~(\ref{eq:u2}), gives that $\epsilon(v_\lambda)$
is in $\mathcal U^{\frac{k(\lambda)}2+2}(\Lg^c)v_\lambda$. 
But $\frac{k(\lambda)}2+2<k(\lambda)$, contradicting Proposition~\ref{prop:D}.

In case $u_2=0$ we have that $u=u_1\in\mathcal U^2(\Lg^c)v_\lambda$. Since
\[ X_{-\beta_{i_0+1}}^{n_{i_0+1}-m}X_{-\beta_{i_0+2}}^{n_{i_0+2}}
             \ldots X_{-\beta_s}^{n_s}u=s_0v_\lambda \]
is a nonzero multiple of the lowest weight vector $\epsilon(v_\lambda)$ and
$n_{i_0+1}-m+n_{i_0+2}+\ldots+n_s=\frac{k(\lambda)}2$, again  
$\epsilon(v_\lambda)\in\mathcal U^{\frac{k(\lambda)}2+2}(\Lg^c)v_\lambda$, 
contradicting the same proposition. \EPf

\section{The candidates to a position in class $\mathcal O^2$}\label{sec:O2}
\setcounter{thm}{0}

In this section we want to elaborate a list of possibilities
for real irreducible representations of class
$\mathcal O^2$ which we shall use later in Sections~\ref{sec:BS}
and~\ref{sec:taut} to classify variationally complete and taut 
representations. Our point of view, as is usually the case in other papers
on the subject, is to classify 
orthogonal representations up to the following equivalence
relation:
we call two representations
$\rho:G\to\mathbf O(V)$, $\rho':G'\to\mathbf O(V')$ 
\emph{image equivalent}
if there exists an isometry $\psi:V\to V'$ such that 
$\mathbf O(\psi)(\rho(G))=\rho'(G')$, where
$\mathbf O(\psi):\mathbf O(V)\to\mathbf O(V')$ is the 
induced conjugation map\footnote{This is finer than the notion of 
orbit equivalence which is 
explained in Section~\ref{sec:review}, but of course not as fine as the
usual notion of equivalence for representations.}.
In particular, the image equivalence class of $\rho$ always contains 
its dual $\rho^*$ (in fact any representation
obtained from $\rho$ by an automorphism of the Dynkin diagram), as well as 
the pull back representation $\tilde\rho=\rho\circ p$ of a covering group 
$p:\tilde G\to G$. 

\bigskip

Let $G$ be a compact connected Lie group. 
We have that $G$ is finitely covered by $T^n\times G_s$, where 
$T^n$ is an $n$-dimensional torus and $G_s$ is a 
compact connected semisimple Lie group (which may also be assumed to be
simply-connected, whenever convenient), and any representation of $G$
pulls back to a representation of $T^n\times G_s$. In view of
image equivalence, 
in order to study representations of $G$, we may assume that 
$G=T^n\times G_s$ and restrict to almost faithful 
representations. Next we observe that any complex irreducible
representation of $T^n$ is a character (hence of complex type).
Therefore a complex irreducible representation of $G$ 
can be of real type or of quaternionic type only if it is 
trivial on $T^n$, and in case it is of complex type, then
it has an $(n-1)$-dimensional kernel on~$T^n$. 
Now we have come to our working 
hypothesis (for a circle group $S^1$, let 
$x^n$ denote the $n$th power representation of that circle which is a 
complex representation of complex type):
\begin{quote}
Let $\rho$ be a real irreducible representation of a compact
connected Lie group $G$ on a finite-dimensional real vector space $V$. 
\begin{itemize}
\item If $\rho$ is of real type, 
we assume that $G=G_s$ is a 
compact connected semisimple Lie group, and $\rho$ is a real form of a complex
irreducible representation $\pi=\pi_\lambda$ of $G$ of real type on the 
complex 
vector space $V_\lambda$ for some highest weight $\lambda$.
\item If $\rho$ is of quaternionic type, 
we assume that $G=G_s$, where $G_s$ is as above,
and $\rho$ is the realification of a complex irreducible representation 
$\pi=\pi_\lambda$
of $G$ of quaternionic type on the complex 
vector space $V_\lambda$ for some highest weight $\lambda$.
\item If $\rho$ is of complex type, 
we assume that $G=G_s$ or $G=S^1\times G_s$, where $G_s$ is as above,
and $\rho$ is the realification of a complex irreducible representation $\pi$
of $G$ of complex type.
According to the form of $G$, we write 
$\pi=\pi_\lambda$ or $\pi=x\otimes\pi_\lambda$ as a
representation on the complex vector space 
$V_\lambda=\mathbf C\otimes_{\mathbf C}V_\lambda$ for 
some highest weight $\lambda$. 
\end{itemize}
\end{quote}

Isotropy representations of symmetric spaces and cohomogeneity one
representations of compact connected Lie groups
are of course two classes of representations 
that belong to class $\mathcal O^2$. 
In order to keep the exposition as simple as possible, 
we assume knowledge of the classification of symmetric spaces. 
So in the course of this section
we will be disregarding 
all the isotropy representations of symmetric spaces
that we encounter, according to tables~8.11.2
and~8.11.5 in~\cite{Wolf}.
Moreover, we classify in Section~\ref{sec:other}
the cohomogeneity one representations of compact connected 
Lie groups. (As a side remark, note that 
our proof of Theorem~\ref{thm:BMS} is
completely independent from the
classification of the symmetric spaces.) Therefore, in addition to 
our previous hypotheses
we assume throughout this section: 
\begin{quote}
The real irreducible representation $\rho$ is of 
class~$\mathcal O^2$, cohomogeneity greater than one 
and it is not the isotropy representation of 
a symmetric space. 
\end{quote}

The results of our investigation in this section 
(namely, Propositions~\ref{prop:quat},
\ref{prop:complex}, \ref{prop:real4}, \ref{prop:real3}, \ref{prop:real21},
\ref{prop:real22} and \ref{prop:simple})
finally imply:
the possibilities for~$\rho$ 
are those collected in the table of Proposition~\ref{prop:pairs}
and in Tables~C.1, C.2, C.3 and~C.4 in
Appendix~C.

\bigskip

We need some amount of notation (cf.~\cite{Wolf}, p.~237). 
Let $\A n=\SU{n+1}$,
$\B n=\Spin{2n+1}$, $\C n=\SP n$, $\D n=\Spin{2n}$, 
$\G$, $\F$, $\E6$, $\E7$ and $\E8$ refer
to \'E. Cartan's classification types of simple Lie groups 
and Lie algebras. For the complex semisimple Lie 
algebra $\mathfrak g^c$, we use the following notation for its 
complex irreducible representation $\pi_\lambda$:
if the integer $2(\lambda,\alpha_i)/(\alpha_i,\alpha_i)\neq0$,
then we write it next to the vertex of the Dynkin diagram of 
$\mathfrak g^c$ which corresponds to $\alpha_i$ (recall that
$\mathcal S=\{\alpha_1,\ldots,\alpha_r\}$ is the simple
root system). For example, 
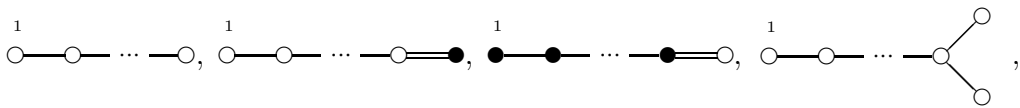
\[ \An1{}{},\;\Bn1{}{}{},\;\Cn1{}{}{},\;\Dn1{}{}{}{}, \]
denote respectively the vector representations of
$\A n$, $\B n$, $\C n$, $\D n$, and the adjoint representations are given
by:
\setlength{\extrarowheight}{0.3cm}
\[ \begin{array}{|l|c|l|c|}
\hline
\A1 & \Ai2 & \G & \Gii{}1 \\
\A n, \; n\geq2 & \An1{}1 & \F & \Fiv1{}{}{} \\
\B n, \; n\geq3 & \Bn{}1{}{} & \E6 & \Evi{}1{}{}{}{} \\
\C n, \; n\geq2 & \Cn2{}{}{} & \E7 & \Evii1{}{}{}{}{}{} \\
\D n, \; n\geq4 & \Dn{}1{}{}{} & \E8 & \Eviii{}{}{}{}{}{}{}1 \\
\hline
\end{array} \]
Note that we are using the dot convention: if there are two lengths 
of roots, then the short roots are black in the Dynkin diagram. 
And a real irreducible representation is denoted
by the diagram of its complexification.

We add that we shall refer repeatedly to  
Table~B.1 in Appendix~B, which contains
the values of the invariant $k(\lambda)$ for the fundamental
representations of the complex simple Lie algebras. 

\subsection{The case where $\rho$ is of quaternionic type}

Here there is nothing:
\begin{prop}\label{prop:quat}
The representation $\rho$ cannot be of quaternionic type.
\end{prop}

\Pf If $\rho$ is of class $\mathcal O^2$ and quaternionic type,
Proposition~\ref{prop:k} says that $k(\lambda)=1$. A glance at 
Table~B.1 now shows that $G=\SP n$ and 
$\pi_{\lambda}$ is the vector representation.
But then $\rho$ is of cohomogeneity one. \EPf

\subsection{The case where $\rho$ is of complex type}

Here the result is:
\begin{prop}\label{prop:complex}
Let $\rho$ be of complex type. Then $\rho$ is 
one of the following:
\setlength{\extrarowheight}{0.6cm}
\[ \begin{array}{|c|c|c|}
\hline
G & \rho & Conditions\\
\hline
\SU n\times\SU m & \parbox[c]{2.5in}{$
                     {}\ (\Ans1{}\otimes\Ans{}1) \\
                     {}\ \qquad \oplus(\Ans{}1\otimes\Ans1{})
                   $} & n\neq m\\
\SU n\times\SP m           & \parbox[c]{2.5in}{$
                               {}\ (\Ans1{}\oplus\Ans{}1) \\
                                {}\ \ \qquad\otimes\Cns1{}{} $}
                        & n\geq3,\; m\geq2\\
S^1\times\SU n\times\SP m  & \parbox[c]{2.5in}{$
                     (x\otimes\Ans1{}\oplus x^{-1}\otimes\Ans{}1) \\
                         {}\ \ \qquad\otimes\Cns1{}{} $} & m\geq2\\
\SU n                      & \An{}1{}\oplus\Anl{}1{} & \mbox{$n$ odd}\\
\Spin{10}                  & \parbox[c]{2.8cm}{$\Dv{}{}{}1{}$}\oplus
                             \parbox[c]{2.8cm}{$\Dv{}{}{}{}1$} & - \\
\SO2\times\Spin7            & (x\oplus x^{-1})\otimes\Biii{}{}1 & - \\
\SO2\times\Spin9            & (x\oplus x^{-1})\otimes\Biv{}{}{}1 & - \\
\SO2\times\G                & (x\oplus x^{-1})\otimes\Gii1{} & - \\
\hline
\end{array} \]
\end{prop}

\medskip

The proof of Proposition~\ref{prop:complex} will be given after the 
proof of the following two lemmas. 

\begin{lem}\label{lem:S1}
Let $\pi_\lambda$ be a complex irreducible representation 
of a compact connected semisimple Lie group $G_s$ and suppose that the
realification of the tensor product representation 
$\pi=x\otimes\pi_\lambda$ of $G=S^1\times G_s$ is of class~$\mathcal O^2$. 
We have:
\begin{enumerate}
\item[(a)] if $\pi_\lambda$ is of real type and $G_s$ is simple,
then $G_s=\Spin7$ or $G_s=\Spin9$  and $\pi_\lambda$ is the 
respective spin representation, or $G_s=\G$ and $\pi_\lambda$ is the 
$7$-dimensional representation, or 
$G_s=\SO m$ for $m\neq2$, $4$ and $\pi_\lambda$ is the vector representation;
\item[(b)] if $\pi_\lambda$ is of real type and $G_s$ is not simple,
then $G_s=\SP1\times\SP m$ and $\pi_\lambda$ is the tensor
product of the vector representations of each of the factors;
\item[(c)] if $\pi_\lambda$ is of quaternionic type or complex type,
then $k(\lambda)=1$, $2$. 
\end{enumerate}
\end{lem}

\Pf (a) and~(b)~Let $V_\lambda$ be the representation space of $\pi_\lambda$ 
and $V$ be the real subspace where a real form acts.  
Note that the representation space for $\pi$ is still $V_\lambda$, but 
$\pi$ is of complex type. Let $p=1\otimes a\in\mathbf C\otimes_{\mathbf
  R}V=V^c=V_\lambda$, $n$ the dimension of $V$ and $q$ the codimension 
of $G_s a$ in~$V$. Since the realification of $\pi$ is of class 
$\mathcal O^2$, 
we have that $\mathcal{O}_p^2(Gp)=V^{cr}\cong V\oplus iV$. But
\[ \mathcal{O}_p^2(Gp)\subset1\otimes[\mathbf Ra + \mathcal{O}_a^2(G_s
a)]+i\otimes[\mathbf Ra + T_a(G_s a)], \]
so by comparing real and imaginary parts we get that $V=\mathbf Ra +
T_a(G_s a)$. Therefore
$q=1$ and $G$ acts with cohomogeneity one on~$V$.
Now~(a) follows from Lemmas~\ref{lem:C1}
and~\ref{lem:C1-classif} below, and~(b) follows from 
Lemmas~\ref{lem:C1} and~\ref{lem:C1-classif2} below.

(c)~Consider the point $p=v_\lambda\in V_\lambda$. 
We have that $\mathcal O^2_p(Gp)\subset\mathcal U^2(\mathfrak
g^c)v_\lambda\subset V_\lambda$ as real vector spaces. 
Since the realification of~$\pi$ 
is of class~$\mathcal O^2$, we have that
$\mathcal U^2(\mathfrak g^c)v_\lambda=V_\lambda$. 
But $\mathcal U^2(\mathfrak g_s^c)v_\lambda=\mathcal U^2(\mathfrak
g^c)v_\lambda$, since the circle subgroup $S^1$ preserves the weight spaces
of $\pi_\lambda$. Therefore 
$\mathcal U^2(\mathfrak g_s^c)v_\lambda=V_\lambda$ 
and we can apply Proposition~\ref{prop:D} 
as in the first paragraph of the proof of Proposition~\ref{prop:k}. \EPf 

\medskip

Part~(a) of the following lemma also follows from~\cite{TY}.

\begin{lem}\label{lem:iv}
\begin{enumerate}
\item[(a)] Let $\pi_\lambda$ be one of the following complex irreducible 
representations of complex type of the compact connected semisimple 
Lie group~$G_s$:
\setlength{\extrarowheight}{0.35cm}
\[ \begin{array}{|c|c|c|}
\hline
G_s & \pi_\lambda & Conditions \\
\hline 
\SU n\times\SU m & \Ans1{}\otimes\Ans{}1 & n\neq m \\
\SU n  & \An{}1{} & n=2p+1 \\
\Spin{10} & \Dv{}{}{}1{} & - \\
\hline
\end{array} \]
Consider the tensor product representation
$\pi=x\otimes\pi_\lambda$ of $S^1\times G_s$. Then the realifications
of $\pi_\lambda$ and $\pi$ are orbit equivalent.
\item[(b)] The realifications of the following complex irreducible 
representations $\pi_\lambda$ 
of complex type of the compact connected semisimple 
Lie group~$G_s$ are 
\emph{not} of class $\mathcal O^2$:
\setlength{\extrarowheight}{0.35cm}
\[ \begin{array}{|c|c|c|}
\hline
G_s & \pi_\lambda  & Conditions \\
\hline 
\SU n  & \An2{}{} & n\geq3\\
\SU n  & \An{}1{} & n=2p\geq6 \\
\SU n\times\SU n & \Ans1{}\otimes\Ans{}1 & n\geq3\\
\E6       & \Evi1{}{}{}{}{} & - \\
\hline
\end{array} \]
\end{enumerate}
\end{lem}

\Pf {\sc Preliminary considerations:}
Each one of these representations $\pi_\lambda$ has
the property that the realification of $\pi=x\otimes\pi_\lambda$ is 
the isotropy 
representation of a compact irreducible Hermitian symmetric space
$X=L_0/K_0$, where $K_0$ is locally isomorphic to $S^1\times G_s$. 
We refer to~\cite{Helgason}, Chapter~VIII, $\S7$, for results 
about Hermitian symmetric spaces that we will be using in the following.
Let $\Ll_0=\Lk_0+\Lp_0$ be the decomposition of the Lie algebra
$\Ll_0$ of $L_0$ into the $\pm1$-eigenspaces of the symmetry.
Let $\Lc_0$ be the Lie algebra of $S^1$,
and let $\Lt_0$ be some Cartan subalgebra
of the Lie algebra $\Lg_s$ of~$G_s$. 
Then $\Lh_0=\Lc_0+\Lt_0$ is a Cartan subalgebra of $\Lk_0$ and of 
$\Ll_0$. Let $\Ll$ be the complexification of $\Ll_0$ and let 
$\Lc$, $\Lt$, $\Lh$, $\Lk$, $\Lp$ be the complex subspaces of $\Ll$ 
spanned by $\Lc_0$, $\Lt_0$, $\Lh_0$, $\Lk_0$, $\Lp_0$.
Let $\Delta$ denote the root system of $(\Ll,\Lh)$ and consider the
root space decomposition $\Ll=\Lh+\sum_{\alpha\in\Delta}\Ll_\alpha$. 
For each $\alpha\in\Delta$, we have that either $\Ll_\alpha\subset\Lk$
or $\Ll_\alpha\subset\Lp$, in which cases the root $\alpha$ is called 
respectively {\em compact} or {\em noncompact}. A root 
is compact if and only if it vanishes on $\Lc$. 
Let $\Delta_c$ and $\Delta_n$
denote respectively the subsets of $\Delta$ of compact and noncompact roots.
Then we have decompositions
$\Lk=\Lh+\sum_{\alpha\in\Delta_c}\Ll_\alpha$,
$\Lp=\sum_{\gamma\in\Delta_n}\Ll_\gamma$. 
Each root is real-valued on $i\Lh_0$. We introduce 
a lexicographic ordering on the dual of $i\Lh_0$ that takes $i\Lc_0$ 
before $i\Lt_0$. Let $\Delta^+$, $\Delta_c^+$ and $\Delta_n^+$ be the 
set of positive roots in $\Delta$, $\Delta_c$ and $\Delta_n$,
and define
\[ \Lp^+=\sum_{\gamma\in\Delta_n^+}\Ll_\gamma,
    \quad\Lp^-=\sum_{-\gamma\in\Delta_n^+}\Ll_\gamma.
\]
Then $\Lp=\Lp^+ + \Lp^-$ is an $\ad{\Lk}$-invariant decomposition. 
The $\ad{\Lk_0}$-invariant complex structure $J_0$ on $\Lp_0$
is given by $\ad{H_0}$, where $H_0\in\Lc_0$ is determined by 
$\gamma(H_0)=i$ for $\gamma\in\Delta_n^+$. The $\mathbf C$-linear
extension $J$ of $J_0$ to $\Lp$ has $\Lp^{\pm}$ as $\pm i$-eigenspaces. 
Let $\Gamma=\{\gamma_1,\ldots,\gamma_s\}\subset\Delta_n^+$ be a
maximal subset of strongly orthogonal roots. Then 
$\La_0=\sum_{j=1}^s\mathbf R(X_{\gamma_j}-X_{-\gamma_j})$ is a 
maximal Abelian subspace of $\Lp_0$. We view $\pi^r$ as the adjoint action 
of $K_0$ on $\Lp_0$, and then $\pi^r$ is polar and $\La_0$ is a 
section for $\pi^r$. We also view $\pi$ as 
the the adjoint action of $K_0$ on $\Lp^+$, and then the weight system of 
$\pi$ is $\Delta_n^+$. Now the root system of $(\Lg_s^c,\Lt)$ is
\[ \Delta_0=\{ \alpha|_{\Lt}:\alpha\in\Delta_c \} 
\quad(\mbox{also $\Delta_0^+=\{ \alpha|_{\Lt}:\alpha\in\Delta_c^+ \}$}), \]
and the weight system of $\pi_\lambda$ is 
\[ \Phi = \{ \gamma|_{\Lt} : \gamma \in\Delta_n^+ \}. \]
It is also useful to remark that any $\ad{\Lk_0}$-invariant
inner product $\inn{\cdot}{\cdot}$ on $\Lp_0$ is given by a 
negative multiple of the restriction of the Killing form of $\Ll_0$;
fix one and extend it to a Hermitian
product $\her{\cdot}{\cdot}$ on $\Lp$. Then $\her{\cdot}{\cdot}$ is
$\ad{\Lk_0}$-invariant and $J$ is a skew-Hermitian operator on $\Lp$
with respect to~$\her{\cdot}{\cdot}$.

\smallskip

{\sc Proof of~\textit{(a)}:} Since $\La_0$ is a section for $\pi^r$, it is 
clearly enough to check that given 
\[p=\sum_{j=1}^s a_j (X_{\gamma_j}-X_{-\gamma_j})\in\La_0, \]
where $a_j\in\mathbf R$, there exists $H\in i\Lt_0$ such that
\setcounter{equation}{\value{thm}}
\stepcounter{thm}
\begin{equation}\label{eqn:same}
 i\ad Hp=J_0p.
\end{equation} 
Write $H=\sum_{\alpha\in\Delta_0^+}c_\alpha H_\alpha$, 
where $H_\alpha$ are the coroots. Then~(\ref{eqn:same})
is equivalent to the following system of linear equations:
\[
 \sum_{\alpha\in\Delta_0^+}c_\alpha(\alpha,\gamma_j)=1,\quad j:1,\ldots,s,
\] 
which we can solve in each one of the three cases with the help of the 
tables in Appendix~A. 

Let $G_s=\SU n\times\SU m$, $n< m$,
$\pi_\lambda:\Ans1{}\otimes\Ans{}1$. Then $L_0$ is $\SU{n+m}$,
\[ \Delta_0^+=\{\theta_i-\theta_j:\mbox{$1\leq i<j\leq n$ or 
$n+1\leq i<j\leq n+m$}\}, \]
\[ \Phi=\{\theta_i-\theta_j:\mbox{$1\leq i\leq n$ and $n+1\leq j\leq
  n+m$}\} \] 
and $\Gamma=\{\theta_i-\theta_{n+m+1-i}:1\leq i\leq n\}$.
Here we can take $H$ to be a suitable multiple of 
$H_{\theta_m-\theta_{m+1}}+H_{\theta_m-\theta_{m+2}}+\ldots+
H_{\theta_m-\theta_{m+n}}$.

Let $G_s=\SU n$, $n$ odd, $\pi_\lambda:\An{}1{}$. Then $L_0$ is $\SO{2n}$,
$\Delta_0^+=\{\theta_i-\theta_j:1\leq i<j\leq n\}$, 
$\Phi=\{\theta_i+\theta_j:1\leq i<j\leq n\}$ and
$\Gamma=\{\theta_{2i-1}+\theta_{2i}:1\leq i\leq\frac{n-1}2\}$.
Here we can take
$H$ to be a suitable multiple of 
$H_{\theta_1-\theta_n}+H_{\theta_3-\theta_n}+\ldots+H_{\theta_{n-2}-\theta_n}$.

Let $G_s=\Spin{10}$, $\pi_\lambda:$\;\parbox[c]{2.8cm}{$\Dv{}{}{}1{}$}. Then $L_0$ is $E_6$,
$\Delta_0^+=\{\pm\theta_i+\theta_j:1\leq i<j\leq5\}$,
$\Phi=\{\frac12(\theta_8-\theta_7-\theta_6+
\sum_{i=1}^5\epsilon_i\theta_i):\Pi_{i=1}^5\epsilon_i=+1\}$ and
$\Gamma=\{\frac12(\theta_8-\theta_7-\theta_6+\theta_5+\theta_4+\theta_3+
\theta_2+\theta_1),\frac12(\theta_8-\theta_7-\theta_6+\theta_5-\theta_4-\theta_3-
\theta_2-\theta_1) \}$. Here we can take
$H$ to be a suitable multiple of 
$H_{\theta_4+\theta_5}-H_{\theta_4-\theta_5}$.
This completes the proof of~(a).

\smallskip

{\sc Proof of~\textit{(b)}:} We will verify directly that the second osculating 
space at the point 
\[ p = \sum_{j=1}^sX_{\gamma_j}-X_{-\gamma_j}\in\La_0\subset\Lp_0 \]
is not of maximal dimension. 
More precisely, we shall see that 
\setcounter{equation}{\value{thm}}
\stepcounter{thm}
\begin{equation}\label{eq:p}
\inn{\ad Xp}{J_0p}=\inn{\ad Y\ad Xp}{J_0p}=0,
\end{equation}
for all $X$, $Y\in\mathfrak g_s$. Notice that
$J_0p=\sum_{j=1}^s i(X_{\gamma_j}+X_{-\gamma_j})\in\Lp_0$. 

It suffices to 
prove the identities~(\ref{eq:p}) for the basis elements~(\ref{eq:basis})
of $\mathfrak g_s$. We have that $\inn{\ad Xp}{J_0p}=\Re\her{\ad Xp}{Jp}$ 
and $\inn{\ad Y\ad Xp}{J_0p}=-\Re\her{\ad Xp}{\ad YJp}$.
The only cases where the real parts of the Hermitian products are not
obviously zero are the following ($\alpha$, $\beta\in\Delta_0^+$):

\setcounter{equation}{\value{thm}}
\stepcounter{thm}
\begin{eqnarray}\label{eqn:gamma1}
\her{\ad{iH_\alpha}p}{Jp} & = & 
  \sum_{j,k=1}^s
  \her{\ad{H_\alpha}(X_{\gamma_j}-X_{-\gamma_j})}{X_{\gamma_k}+X_{-\gamma_k} } 
\nonumber \\
 & = & \sum_{j,k=1}^s \gamma_j(H_\alpha)
\her{X_{\gamma_j}+X_{-\gamma_j}}{X_{\gamma_k}+X_{-\gamma_k})} \nonumber \\
  & = & 2\sum_{j=1}^s\gamma_j(H_\alpha) = 2(\sum_{j=1}^s\gamma_j,\alpha),
\end{eqnarray}
assuming the root vectors to have length one. We will show below that the
last expression in~(\ref{eqn:gamma1}) vanishes for the 
representations we are interested. 
\setcounter{equation}{\value{thm}}
\stepcounter{thm}
\begin{eqnarray}\label{eqn:gamma2}
\her{\ad{i(X_{\alpha}+X_{-\alpha})}p}{Jp} & = & 
\her{\ad{X_{\alpha}+X_{-\alpha}}p}{\sum_{k=1}^s(X_{\gamma_k}+X_{-\gamma_k})}
\nonumber \\
 & = &
 \sum_{j,k=1}^s\her{\ad{X_{\alpha}+X_{-\alpha}}(X_{\gamma_j}-X_{-\gamma_j})}{
X_{\gamma_k}+X_{-\gamma_k} } \nonumber \\
 & = & 0,
\end{eqnarray}
because $\ad{X_\alpha}X_{\gamma_j}$ is either a $(\gamma_j+\alpha)$-root vector
or $0$; $\pm\gamma_j\pm\gamma_k\not\in\Delta$; and
root spaces corresponding to distinct roots are orthogonal. 

\begin{eqnarray*}\label{eqn:gamma3}
\her{\ad{iH_\alpha}p}{\ad{X_\beta-X_{-\beta}}Jp)} & = & 
\sum_{j,k=1}^s\gamma_j(H_\alpha)\her{(X_{\gamma_j}+X_{-\gamma_j})}{
\ad{X_\beta-X_{-\beta}}(X_{\gamma_k}+X_{-\gamma_k})} \\
& = & 0,
\end{eqnarray*}
for the same reason why~(\ref{eqn:gamma2}) is zero. 

\setcounter{equation}{\value{thm}}
\stepcounter{thm}
\begin{eqnarray}\label{eqn:gamma4}
\Re\her{\ad{X_\alpha-X_{-\alpha}}Jp}{\ad{i(X_\beta+X_{-\beta})}p} & = & 
  \Re\her{\ad{X_\alpha}Jp}{i\ad{X_\beta}p} \nonumber \\
  & &  \qquad -\Re\her{\ad{X_{-\alpha}}Jp}{i\ad{X_{-\beta}}p}   \nonumber  \\
  & &  \qquad +\Re\her{\ad{X_\alpha}Jp}{i\ad{X_{-\beta}}p} \nonumber  \\
  & &  \qquad -\Re\her{\ad{X_{-\alpha}}Jp}{i\ad{X_\beta}p}   \nonumber    \\
  &=&  \Re\her{Jp}{i\ad{[X_{-\beta},X_\alpha]}p} \nonumber \\
  & &  \qquad  +\Re\her{Jp}{i\ad{[X_{\beta},X_\alpha]}p},
\end{eqnarray}
because, for instance, 
\begin{eqnarray*}
\her{\ad{X_{\alpha}}Jp}{i\ad{X_{\beta}}p} & = & 
   \her{\ad{X_{-\beta}}\ad{X_{\alpha}}Jp}{ip} \\
  & = & \her{J\ad{X_{-\beta}}\ad{X_{\alpha}}p}{ip} \\
  & = & -\her{\ad{X_{-\beta}}\ad{X_{\alpha}}p}{Jip} \\
  & = & -\her{\ad{X_{-\beta}}\ad{X_{\alpha}}p}{iJp} \\
  & = & \her{i\ad{X_{-\beta}}\ad{X_{\alpha}}p}{Jp} \\
  & = & \overline{\her{Jp}{i\ad{X_{-\beta}}\ad{X_{\alpha}}p}}.
\end{eqnarray*}

Now~(\ref{eqn:gamma4}) shows that everything amounts to 
checking that~(\ref{eqn:gamma1}) vanishes. Next we do that in each one
of the four cases. Let
$\gamma=\sum_{j=1}^s\gamma_j$. Then the vanishing of~(\ref{eqn:gamma1})
is equivalent to
\setcounter{equation}{\value{thm}}
\stepcounter{thm}
\begin{equation}\label{eqn:gamma1'}
(\gamma,\Delta_0^+)=0. 
\end{equation}

Let $G_s=\SU n$, $\pi_\lambda:\Ans2{}$. Then $L_0$ is $\SP n$, 
$\Delta_0^+=\{\theta_i-\theta_j:1\leq i<j\leq n\}$, 
$\Phi=\{\theta_i+\theta_j:1\leq i<j\leq n\}\cup\{2\theta_i:1\leq
i\leq n\}$, $\Gamma=\{2\theta_i:1\leq i\leq n\}$ and
$\gamma=2\theta_1+\ldots+2\theta_n$. 

Let $G_s=\SU n$, $n$ even, $\pi_\lambda:\An{}1{}$. Then $L_0$ is $\SO{2n}$,
$\Delta_0^+=\{\theta_i-\theta_j:1\leq i<j\leq n\}$, 
$\Phi=\{\theta_i+\theta_j:1\leq i<j\leq n\}$,
$\Gamma=\{\theta_{2i-1}+\theta_{2i}:1\leq i\leq\frac n2\}$ and 
$\gamma=\theta_1+\theta_2+\ldots+\theta_n$. 

Let $G_s=\SU n\times\SU n$, $\pi_\lambda:\Ans1{}\otimes\Ans{}1$. Then $L_0$ 
is $\SU{2n}$, 
$\Delta_0^+=\{\theta_i-\theta_j:\mbox{$1\leq i<j\leq n$ or $n+1\leq i<j\leq 2n$}\}$,
$\Phi=\{\theta_i-\theta_j:\mbox{$1\leq i\leq n$ and $n+1\leq j\leq 2n$}\}$,
$\Gamma=\{\theta_i-\theta_{2n+1-i}:1\leq i\leq n\}$ and
$\gamma=\theta_1+\theta_2+\ldots+\theta_n-\theta_{n+1}-\theta_{n+2}-\ldots-\theta_{2n}$.

Let $G_s=\E6$, $\pi_\lambda:\parbox[c]{3.4cm}{$\Evi1{}{}{}{}{}$}$. Then $L_0=\E7$,
\[ \Delta_0^+=\{\pm\theta_i+\theta_j:1\leq
  i<j\leq5\}\cup\{\frac12(\theta_8-\theta_7-\theta_6+
                 \sum_{i=1}^5\epsilon_i\theta_i):\Pi_{i=1}^5\epsilon_i=+1\},
                 \]
\[ \Phi=\{\theta_8-\theta_7\}\cup\{\pm\theta_i+\theta_6:1\leq i\leq5\}
               \cup\{\frac12(\theta_8-\theta_7+\theta_6+
                 \sum_{i=1}^5\epsilon_i\theta_i):\Pi_{i=1}^5\epsilon_i=-1\},
                 \]
$\Gamma=\{\theta_8-\theta_7,\theta_6+\theta_5,\theta_6-\theta_5\}$ and
$\gamma=\theta_8-\theta_7+2\theta_6$. 
     
It is immediate to see that~(\ref{eqn:gamma1'}) holds in each case,
and this completes the proof of~(b).\\ \mbox{} \EPf

\medskip

{\it Proof of Proposition~\ref{prop:complex}.} 
We know from Proposition~\ref{prop:k} and Lemma~\ref{lem:S1}
that $k(\lambda)=1$, $2$.
In particular, $G_s$ can have at most two simple factors.

Suppose $k(\lambda)=1$.
This implies $G_s$ is simple. In the case 
$G=G_s$, we have $\pi$ is of complex type, so $G=\SU n$ and $\pi$ is the
vector representation. In the case $G=S^1\times G_s$, we have
$\pi_\lambda$ is of quaternionic or of complex type, so $G_s=\SP n$
or $G_s=\SU n$ and $\pi_\lambda$ is the respective vector 
representation. This only gives representations with 
cohomogeneity one.

Next suppose $k(\lambda)=2$. The first case 
is $G_s=G_1\times G_2$, where $G_i$ is a 
simple group. Here $\pi_\lambda$ decomposes
as an outer tensor product $\pi_{\mu_1}\otimes\pi_{\mu_2}$, 
where $\pi_{\mu_i}$ is the representation of $G_i$ of highest weight
$\mu_i$ with $k(\mu_i)=1$. This forces each $G_i$ to be either 
$\SP n$ or $\SU n$ for some $n$, and $\pi_{\mu_i}$ to be the 
corresponding vector representation. 
Using Lemma~\ref{lem:S1}, item (b), and Lemma~\ref{lem:iv}, item (b),
we get the first three representations in the table.

The second case is $G_s$ is a simple group. For $G=G_s$, an inspection of
Table~B.1 reveals four possibilities for~$\pi_\lambda$, namely:
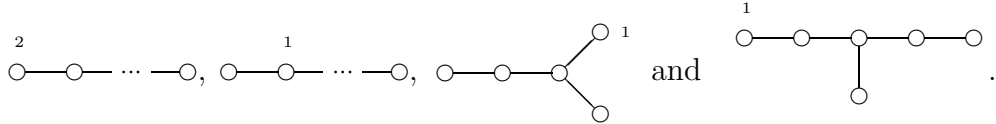
\[ \An2{}{},\;\An{}1{},\;\Dv{}{}{}1{}\;\mbox{and}\quad\Evi1{}{}{}{}{}.\] 
Now Lemma~\ref{lem:iv}, item (b), restricts these to the next two
representations in the table. 

For $G=S^1\times G_s$, if $\pi_\lambda$ is of complex type,
we only get representations associated to symmetric spaces, and if  
$\pi_\lambda$ is of real type, Lemma~\ref{lem:S1}, item (a), gives the last three
representations in our table. \EPf

\subsection{The case where $\rho$ is of real type}

The real case is much more involved than the previous two.
We know from Proposition~\ref{prop:k} that $k(\lambda)=2$, $4$.
In particular, $G$ can have at most four simple factors.

As a first step towards the classification in the case $\rho$ is
of real type, we introduce the following variations
of condition $(C_2)$ of Lemma~\ref{lem:CT} for a
representation $\pi_\lambda$ with representation space $V_\lambda$ 
and highest weight vector $v_\lambda$:
\begin{eqnarray*}
\mathcal U^2(\Lg^c)v_\lambda+\mathcal U^1(\Lg^c)\epsilon(v_\lambda)&=&V_\lambda\qquad(C_{1\frac12}) \\
\mathcal U^1(\Lg^c)v_\lambda+\mathcal U^2(\Lg^c)\epsilon(v_\lambda)&=&V_\lambda\qquad(C_{1\frac12}') \\
\mathcal U^1(\Lg^c)v_\lambda+\mathcal U^1(\Lg^c)\epsilon(v_\lambda)&=&V_\lambda\qquad(C_1) \\
\mathcal U^1(\Lg^c)v_\lambda+\mathbf C\epsilon(v_\lambda)&=&V_\lambda\qquad(C_{\frac12}) \\
\mathbf Cv_\lambda+\mathcal U^1(\Lg^c)\epsilon(v_\lambda)&=&V_\lambda\qquad(C_{\frac12}')
\end{eqnarray*}
It follows from the identity $\epsilon X_\alpha=X_{-\alpha}\epsilon$ that
conditions $(C_{1\frac12})$ and $(C_{1\frac12}')$ are equivalent, and that 
conditions $(C_{\frac12})$ and $(C_{\frac12}')$ are equivalent. 

\begin{lem}\label{lem:C1}
The following conditions are equivalent for a
complex irreducible representation $\pi_\lambda$ of \emph{real type} of a compact
connected semisimple Lie group $G$, with invariant real form $V$:
\begin{enumerate}
\item[(a)] condition $(C_1)$;
\item[(b)] $G$ is transitive on the unit sphere of $V$;
\item[(c)] $G$ acts with cohomogeneity~$1$ on $V$;
\item[(d)] $\{Xp:X\in\mathfrak g\}=(\mathbf Rp)^\perp$ for all $p\in V$;
\item[(e)] $\{Xp:X\in\mathfrak g\}=(\mathbf Rp)^\perp$ for some $p\in V$.
\end{enumerate}
\end{lem}

\Pf (b), (c), (d), and~(e)~are clearly equivalent.
We prove that (a) is equivalent to~(e). 
In fact, let $p=v_\lambda+\epsilon(v_\lambda)$. Clearly 
$\{Xp:X\in\mathfrak g\}\subset(\mathbf Rp)^{\perp}$. 
Using the basis~(\ref{eq:basis}), we find that $\{Xp:X\in\mathfrak g\}^c$ is spanned
over $\mathbf C$ by 
\[ v_\lambda-\epsilon(v_\lambda),\quad X_{-\alpha}v_\lambda,\quad
X_{\alpha}\epsilon(v_\lambda),\qquad\alpha\in\Delta^+. \]
But condition~$(C_1)$ holds precisely when this is also
a system of generators for $(\mathbf Rp)^{\perp c}$ over~$\mathbf C$.\\ \mbox{} \EPf

\begin{lem}\label{lem:C1notreal}
The following conditions are equivalent for a 
complex representation $\pi_\lambda$ of \emph{quaternionic or complex
type} of a compact connected semisimple Lie group $G$ on $V_\lambda$,
with highest weight $v_\lambda$:
\begin{enumerate}
\item[(a)] $\mathcal U^1(\mathfrak g^c) v_\lambda=V_\lambda$;
\item[(b)] $G$ is transitive on the unit sphere of $V_\lambda^r$;
\item[(c)] $G$ acts with cohomogeneity~$1$ on $V_\lambda^r$;
\item[(d)] $\{Xp:X\in\mathfrak g\}=(\mathbf Rp)^\perp$ for all $p\in V_\lambda^r$;
\item[(e)] $\{Xp:X\in\mathfrak g\}=(\mathbf Rp)^\perp$ for some $p\in V_\lambda^r$.
\end{enumerate}
\end{lem}

\Pf  The proof is analogous to the proof of Lemma~\ref{lem:C1}. \EPf

\begin{lem}\label{lem:C1/2-classif}
The self-dual complex irreducible representations 
of compact connected \emph{simple} Lie groups satisfying 
condition~$(C_\frac12)$ are precisely\,\footnote{The half-spin
representations of $\Spin8$ are image equivalent to the vector 
representation.}: the vector representation of
$\SO m$ for $m\neq2$, $4$, the vector representation of $\SP m$, 
the $7$-dimensional representation of $\G$, and the spin 
representation of $\Spin7$.  
\end{lem}

We will prove Lemma~\ref{lem:C1/2-classif} together with the
next lemma. 

\begin{lem}\label{lem:C1-classif}
The self-dual complex irreducible representations 
of compact connected \emph{simple} Lie groups satisfying 
condition~$(C_1)$ are, \emph{besides those listed in
Lemma~\ref{lem:C1/2-classif}}, precisely the following:
\setlength{\extrarowheight}{0.35cm}
\[ \begin{array}{|c|c|c|}
\hline
G & \pi_\lambda & k(\lambda) \\
\hline
\SU2  & \Ai3    & 3  \\
\SU6 & \Av{}{}1{}{} & 3 \\
\Spin9  & \Biv{}{}{}1 & 2 \\
\Spin{11}  & \Bv{}{}{}{}1 & 3 \\
\Spin{12}  & \Dvi{}{}{}{}1{} & 3 \\
\SP3   & \Ciii{}{}1 & 3 \\
\E7   &  \Evii{}{}{}{}{}{}1 & 3 \\
\hline
\end{array} \]
\end{lem}

{\it Proof of Lemmas~\ref{lem:C1/2-classif} and~\ref{lem:C1-classif}.} 
It is obvious that condition $(C_\frac12)$ implies condition $(C_1)$.
So we consider $(C_1)$. It implies the following:
\begin{enumerate}
\item[(a)] $k(\lambda)\leq3$. 
\item[(b)] The dimension of $\pi_\lambda$ is bounded by
$\dim\mathfrak{g}^c-\mbox{rank}\,\mathfrak{g}^c + 2$.
\item[(c)] If $0$ is a weight, then its multiplicity is one. 
In particular, if $\pi_\lambda$ is the adjoint representation, 
then the rank is at most one.
\end{enumerate}
In fact, (a)~follows from the fact that $k(\lambda)\geq4$ contradicts
Proposition~\ref{prop:D}.
For (b)~notice that the dimension of the left hand side
of~$(C_1)$ is bounded by $2(\mbox{number of positive roots})+2$.
Finally, if $0$ is a weight, then the $0$-weight space 
is contained in $\mathcal U^1(\mathfrak g^c)v_\lambda$ (because
$\epsilon X_\alpha=X_{-\alpha}\epsilon$),
and from this we see that $\lambda$ is a root
and the $0$-weight space is in fact spanned by
$X_{-\lambda}v_\lambda$, which implies (c).

We will use the above observations to complete the proof of the lemmas. 

Let $\pi_\lambda$ be of real type. Then $k(\lambda)=2$. 
The adjoint representations of 
$\A n$, $\B n$, $\C n$, $\D n$ with $n\geq2$,
and of $\G$, $\F$, $\E6$, $\E7$, $\E8$ are immediately excluded. 
We run through all the remaining cases of representations
of real type with $k(\lambda)=2$. 

Let $G$ be of $\A n$-type. Then $\pi_\lambda$ is 
$\Ai2$ or $\Aiii{}1{}$ which are cited in the lemmas.

Let $G$ be of $\B n$-type, $n\geq2$. Then $\pi_\lambda$ is $\Bns 1{}{}$
or $\Biii{}{}1$ or $\Biv{}{}{}1$. All
of them are cited in the lemmas.

Let $G$ be of $\C n$-type, $n\geq3$. Then $\pi_\lambda$ is 
$\Cn{}1{}{}$, which cannot occur 
because the multiplicity of the $0$-weight is $n-1$. 

Let $G$ be of $\D n$-type, $n\geq4$. Then $\pi_\lambda$ is \Dns1{}{}{},
which is cited in the lemmas.  

Let $G$ be of exceptional type. Then $\pi_\lambda$ is $\Gii1{}$ 
or $\Fiv{}{}{}1$. Only the former is cited in the lemmas;
the latter has $0$ as a weight of multiplicity $2$. 

Now let $\pi_\lambda$ be of quaternionic type. Then $k(\lambda)=1$
or $k(\lambda)=3$. The case $k(\lambda)=1$ gives the vector representation
of $\SP m$, which is cited in the lemmas. In the case $k(\lambda)=3$, the
only possibilities not cited in the lemmas
are the following four representations, which we next eliminate:
\begin{itemize}
\item $G=\B6$, $\pi_\lambda:\;\Bvi{}{}{}{}{}1$: 
here $\lambda=\frac12(\theta_1+\theta_2+\theta_3+\theta_4+\theta_5+\theta_6)$,
$\mu=\frac12(\theta_1+\theta_2+\theta_3-\theta_4-\theta_5-\theta_6)$
is a weight and $\lambda\pm\mu$ are not roots which violates
condition~$(C_1)$.
\end{itemize}
In the remaining cases we make use of the bound on the dimension of 
the representation (see e.~g.~Table~5 in the Reference Chapter
of~\cite{OV}):
\begin{itemize}
\item $G=\C n$, $\pi_\lambda:\;\Cn3{}{}{}$:
$\dim\pi_\lambda=\frac23n(n+1)(2n+1)>2n^2+2=\dim\C n-n+2$, for $n\geq2$;
\item $G=\C n$, $\pi_\lambda:\;\Cn 11{}{}$:
$\dim\pi_\lambda=\frac83n(n-1)(n+1)>2n^2+2=\dim\C n-n+2$, for $n\geq2$;
\item $G=\C n$, $\pi_\lambda:\;\Cni{}{}1{}{}$:
$\dim\pi_\lambda=\frac23n(n-2)(2n+1)>2n^2+2=\dim\C n-n+2$, for $n\geq4$;
\end{itemize}
On the other hand, the representations listed in the statements of the 
Lemmas~\ref{lem:C1/2-classif} and~\ref{lem:C1-classif} 
do satisfy $(C_1)$, as can be easily checked using the fact that
each of them has all weights with multiplicities one. 
This concludes the classification of self-dual 
representations that satisfy $(C_1)$. As for 
condition $(C_\frac12)$, we observe that $(C_\frac12)$
implies that $k(\lambda)\leq2$. This eliminates all representations
in the table of Lemma~\ref{lem:C1-classif}
except $\Biv{}{}{}1$; that one can be eliminated by a direct
check. Finally it is easy to verify that the representations
cited in Lemma~\ref{lem:C1/2-classif} do satisfy $(C_\frac12)$. \EPf

\begin{lem}\label{lem:C1-classif2}
The only complex irreducible representation of real type of a compact
connected semisimple \emph{nonsimple} Lie group $G$ satisfying 
condition~$(C_1)$ is the tensor product of the vector 
representations of $\SP 1$ and $\SP m$. 
\end{lem}

\Pf We already know that $k(\lambda)=2$ 
(cf. proof of Lemmas~\ref{lem:C1/2-classif} and~\ref{lem:C1-classif}),
so $G=G_1\times G_2$, where $G_i$ is a simple group, and 
$\pi_\lambda=\pi_{\mu_1}\otimes\pi_{\mu_2}$ with $k(\mu_i)=1$. 
Since each $\pi_{\,u_i}$ is self-dual, it must be of quaternionic type.
Therefore, each $G_i=\SP {n_i}$ and the corresponding $\pi_{\mu_i}$ is just the vector
representation. Now one of the $n_i$ must be $1$, for otherwise 
condition~$(C_1)$ for $\pi_\lambda$ is violated by consideration of a weight vector 
of type $v_1\otimes v_2$, where $v_i$ is an \emph{intermediate}
(that is, neither highest nor lowest) weight vector for $\pi_{\mu_i}$. \EPf

\begin{lem}\label{lem:C11/2-classif}
The complex irreducible representations $\pi_\lambda$ 
of quaternionic type of compact connected simple Lie groups
$G$ satisfying $k(\lambda)=3$ and 
condition~$(C_{1\frac12})$, \emph{but not condition~$(C_1)$}, are
the following:
\setlength{\extrarowheight}{0.3cm}
\[ \begin{array}{|c|c|}
\hline
G & \pi_\lambda  \\
\hline 
\Spin{13}  & \Bvi{}{}{}{}{}1 \\
\SP2   & \Cii11 \\
\hline
\end{array} \]
\end{lem}

We postpone the proof of Lemma~\ref{lem:C11/2-classif} to the end of 
Subsection~\ref{sec:tech}, since the methods used to prove it better belong
there. 

\subsubsection{The case where $G$ has four simple factors}

In this case we have:
\begin{prop}\label{prop:real4}
Let $\rho$ be of real type and $G=G_1\times G_2\times G_3\times G_4$,
where $G_i$ is a simple group. Then $G=\SU2\times\SU2\times\SU2\times\SP n$,
$n\geq2$, and $\pi$ is the tensor product of the vector representations 
of each of the factors. 
\end{prop}

\Pf Set $\pi_\lambda=\pi_{\mu_1}\otimes\pi_{\mu_2}\otimes\pi_{\mu_3}\otimes
\pi_{\mu_4}$, where $\pi_{\mu_i}$ is the representation of $G_i$ of highest 
weight ${\mu_i}$. We have $k(\mu_i)=1$ for all $i$, and since each 
$\pi_{\mu_i}$ is self-dual, it must be of quaternionic type. Therefore,
each $G_i=\SP {n_i}$ and the corresponding $\pi_{\mu_i}$ is just the vector
representation. Notice that $n_i=1$ for all $i$ gives the isotropy
representation of a symmetric space. 
Suppose now that for two indices, say $3$ and $4$, we have
$n_3>1$ and $n_4>1$. Take weight vectors $v_i$ of $\pi_{\mu_i}$ such that:
$v_1$ is the highest weight vector of $\pi_{\mu_1}$,
$v_2$ is the lowest weight vector of $\pi_{\mu_2}$, and 
$v_3$ and $v_4$ are intermediate
weight vectors of $\pi_{\mu_3}$ and $\pi_{\mu_4}$, respectively.
Then $v_1\otimes v_2\otimes v_3\otimes v_4$ is a weight vector 
of $\pi_\lambda$ which cannot be contained in 
$\mathcal U^2(\Lg^c)v_\lambda+\mathcal U^2(\Lg^c)\epsilon(v_\lambda)$, violating 
the condition~$(C_2)$ from Lemma~\ref{lem:CT}. \EPf

\subsubsection{The case where $G$ has three simple factors}

The classification is:
\begin{prop}\label{prop:real3}
Let $\rho$ be of real type and $G=G_1\times G_2\times G_3$, where $G_i$ is a simple
group. Then $\rho$ is one of the following:
\setlength{\extrarowheight}{0.35cm}
\[ \begin{array}{|c|c|c|}
\hline
G & \rho & \mbox{Conditions} \\
\hline
\SU2\times\SP n\times\SO {2m}  & \Ai1\otimes\Cns1{}{}\otimes\Dns1{}{}{}
& m\geq3,\quad n\geq2 \\ 
\SU2\times\SP n\times\SO {2m+1}  & \Ai1\otimes\Cns1{}{}\otimes\Bns1{}{}
& n\geq2 \\ 
\SU2\times\SP n\times\G        & \Ai1\otimes\Cns1{}{}\otimes\Gii1{} & - \\
\SU2\times\SP n\times\Spin7    & \Ai1\otimes\Cns1{}{}\otimes\Biii{}{}1 & - \\
\SU2\times\SU2\times\Spin9     & \Ai1\otimes\Ai1\otimes\Biv{}{}{}1 & - \\
\hline
\end{array} \]
\end{prop}

\Pf Let $\pi_\lambda=\pi_{\mu_1}\otimes\pi_{\mu_2}\otimes\pi_{\mu_3}$,
where 
$\pi_{\mu_i}$ is the representation of $G_i$ of highest 
weight ${\mu_i}$, say with $k(\mu_1)=k(\mu_2)=1$, $k(\mu_3)=2$. 
Then $\pi_{\mu_1}$, $\pi_{\mu_2}$ are each of quaternionic type and $\pi_{\mu_3}$
is of real type. So each of $\pi_{\mu_1}$, $\pi_{\mu_2}$ is the vector 
representation of $\SP {m}$ for some $m$. Moreover, one of them must have
$m=1$, say $\pi_{\mu_1}$,
for otherwise condition $(C_2)$ for $\pi_\lambda$ is violated by
consideration of a weight vector of type $v_1\otimes v_2\otimes v_3$, 
where $v_i$ is an intermediate weight vector for $\pi_{\mu_i}$.
Now we have $G_1=\SU2$, $G_2=\SP n$, and $\pi_{\mu_1}$, $\pi_{\mu_2}$
the respective vector representations. Next take a weight vector
$v_1\otimes v_2\otimes v_3$ for $\pi_\lambda$ such that $v_1$ is 
highest, $v_2$ is lowest and $v_3$ is intermediate. Condition $(C_2)$
for $\pi_\lambda$ forces $\pi_{\mu_3}$ to satisfy condition 
$(C_1)$. If $n\geq2$, we could also take $v_2$ to be intermediate, and then
condition $(C_2)$ for $\pi_\lambda$ implies something stronger, namely
that $\pi_{\mu_3}$ must satisfy condition  $(C_\frac12)$. Our table now follows 
from Lemmas~\ref{lem:C1/2-classif} and~\ref{lem:C1-classif}. \EPf

\subsubsection{The case where $G$ has two simple factors and $\rho$
has two factors of real type}

Here the result is:
\begin{prop}\label{prop:real21}
Let $\rho$ be of real type and $G=G_1\times G_2$, where $G_i$ is a simple
group. Write $\rho^c=\pi_\lambda$, decompose 
$\pi_\lambda=\pi_{\mu_1}\otimes\pi_{\mu_2}$ and suppose that 
each $\pi_{\mu_i}$ is of real type. Then each $\pi_{\mu_i}$ is one of
the following representations: the vector representation of
$\SO m$ for $m\neq2$, $4$, 
or the $7$-dimensional representation of $\G$, or the spin 
representation of $\Spin7$, or the spin representation of $\Spin9$; 
but the case where $G_1=G_2=\Spin9$ is excluded.
\end{prop}

\Pf Here $k(\mu_1)=k(\mu_2)=2$. Let $v_1\otimes v_2$ be a weight vector
for $\pi_\lambda$ such that each $v_i$ is an intermediate 
weight vector for $\pi_{\mu_i}$. Condition~$(C_2)$ for $\pi_\lambda$
applied to $v_1\otimes v_2$ forces
each $\pi_{\mu_i}$ to satisfy $(C_1)$. In fact, one of the $\pi_{\mu_i}$
must even satisfy~$(C_\frac12)$. The claim follows from 
Lemmas~\ref{lem:C1/2-classif} and~\ref{lem:C1-classif}. \EPf

\subsubsection{The case where $G$ has two simple factors and $\rho$
has two factors of quaternionic type}

The classification is:
\begin{prop}\label{prop:real22}
Let $\rho$ be of real type and $G=G_1\times G_2$, where $G_i$ is a simple
group. Write $\rho^c=\pi_\lambda$, decompose 
$\pi_\lambda=\pi_{\mu_1}\otimes\pi_{\mu_2}$ and suppose that 
each $\pi_{\mu_i}$ is of quaternionic type. Then $\rho$ is one of the 
following:
\setlength{\extrarowheight}{0.35cm}
\[ \begin{array}{|c|c|c|}
\hline
G & \rho & Conditions \\
\hline
\SP n\times\SU2      & \Cns1{}{}\otimes\Ai3                & n\geq2 \\
\SP n\times\SU6      & \Cns1{}{}\otimes\Av{}{}1{}{}        & n\geq2 \\
\SP n\times\Spin{11} & \Cns1{}{}\otimes\Bv{}{}{}{}1        & -     \\
\SP n\times\Spin{12} & \Cns1{}{}\otimes\Dvi{}{}{}{}1{}    & n\geq2 \\
\SP 1\times\Spin{13} & \Ai1\otimes\Bvi{}{}{}{}{}1           & -     \\
\SP 1\times\SP2      & \Ai1\otimes\Cii11                    & -     \\
\SP n\times\SP3      & \Cns1{}{}\otimes\Ciii{}{}1          & n\geq2 \\
\SP n\times\E7       & \Cns1{}{}\otimes\Evii{}{}{}{}{}{}1  & n\geq2 \\
\hline
\end{array} \]
\end{prop}

\Pf If $k(\lambda)=2$, then $k(\mu_1)=k(\mu_2)=1$ and each $\pi_{\mu_i}$
must be the vector representation of $\SP{n_i}$ for some $n_i$, 
so we get a representation 
associated to a symmetric space. In the case $k(\lambda)=4$ we can write 
$k(\mu_1)=1$ and $k(\mu_2)=3$. Now 
$\pi_{\mu_1}$ must be the vector representation of $\SP n$ and 
there are two cases to consider:
\begin{itemize}
\item $n=1$: here $\pi_{\mu_2}$ must satisfy condition
  $(C_{1\frac12})$;
\item $n\geq2$: here $\pi_{\mu_2}$ must satisfy the stronger condition $(C_1)$.
\end{itemize}
The claim now follows from 
Lemmas~\ref{lem:C1-classif} and~\ref{lem:C11/2-classif}. \EPf

\subsubsection{The case where $G$ is a simple group}\label{sec:tech}

The representations of real type of simple groups that are candidates
to a position in class $\mathcal O^2$ are: 

\begin{prop}\label{prop:simple}
Let $\rho$ be of real type and suppose $G$ is a simple group. 
Then $\rho$ is one of the following:
\setlength{\extrarowheight}{0.35cm}
\[ \begin{array}{|c|c|}
\hline
G & \rho \\
\hline
\Spin7      & \Biii1{}1 \\
\Spin9      & \Biv1{}{}1 \\
\Spin{15}   & \Bvii{}{}{}{}{}{}1 \\
\Spin{17}   & \Bviii{}{}{}{}{}{}{}1 \\
\hline
\end{array} \]
\end{prop}

We will accomplish the proof of the above proposition by running 
through all the cases of representations of real type of a simple 
group with $k(\lambda)\leq4$ given in Table~B.1. In fact, 
according to that table, all representations 
with $k(\lambda)=2$ are either associated to symmetric spaces or  
of cohomogeneity one.
So we only need to worry about the case
$k(\lambda)=4$. A careful inspection of Table~B.1 shows:

\begin{lem}\label{lem:k=4}
Let $\pi_\lambda$ be a representation of real type of a compact
connected simple Lie group $G$ with $k(\lambda)=4$. 
Then $(G,\pi_\lambda)$ is one of the following:
\setlength{\extrarowheight}{0.35cm}
\[\begin{array}{|r|c|c|c|}
\hline
 & G & \pi_\lambda & Conditions \\
\cline{2-4}
1 & \A1 & \Ai4 & - \\
2 & \A3 & \Aiii{}2{} & - \\
3 & \A3 & \Aiii111 & - \\
4 & \A7 & \Avii{}{}{}1{}{}{} & - \\
5 & \A n & \Ani2{}{}2 & n\geq2 \\
6 & \A n & \Ani{}11{} & n\geq4 \\
\hline
\end{array}\]

\[\begin{array}{|r|c|c|c|}
\hline
 & G & \pi_\lambda & Conditions \\
\cline{2-4}
7 & \B3 & \Biii{}{}2 & - \\
8 & \B3 & \Biii1{}1 & - \\
9 & \B3 & \Biii{}11 & - \\
10 & \B4 & \Biv{}{}{}2 & - \\
11 & \B4 & \Biv1{}{}1 & - \\
12 & \B4 & \Biv{}1{}1 & - \\
13 & \B7 & \Bvii{}{}{}{}{}{}1 & - \\
14 & \B8 & \Bviii{}{}{}{}{}{}{}1 & - \\
15 & \B n & \Bn 2{}{}{} & n\geq2 \\
16 & \B n & \Bn {}2{}{} & n\geq3 \\
17 & \B n & \Bn 11{}{} & n\geq3 \\
18 & \B n & \Bni{}{}1{}{} & n\geq4 \\
19 & \B n & \Bnii{}{}{}1{}{} & n\geq5 \\
\hline
\end{array}\]

\[\begin{array}{|r|c|c|c|}
\hline
 & G & \pi_\lambda & Conditions \\
\cline{2-4}
20 & \C n & \Cn4{}{}{} & n\geq2 \\
21 & \C n & \Cn21{}{} & n\geq2 \\
22 & \C n & \Cni1{}1{}{} & n\geq3 \\
23 & \C n & \Cni{}2{}{}{} & n\geq3 \\
24 & \C n & \Cnii{}{}{}1{}{} & n\geq4 \\
\hline
\end{array}\]

\[\begin{array}{|r|c|c|c|}
\hline
 & G & \pi_\lambda & Conditions \\
\cline{2-4}
25 & \D4 & \Div{}{}11 & - \\
26 & \D5 & \Dv{}{}{}11 & - \\
27 & \D8 & \Dviii{}{}{}{}{}{}1{}{}\quad & - \\
28 & \D n & \Dn2{}{}{}{} & n\geq4 \\
29 & \D n & \Dn{}2{}{}{} & n\geq4 \\
30 & \D n & \Dn11{}{}{} & n\geq4 \\
31 & \D n & \Dni{}{}1{}{}{} & n\geq5 \\
32 & \D n & \Dnii{}{}{}1{}{}{} & n\geq6 \\
\hline
\end{array}\]

\[\begin{array}{|r|c|c|}
\hline
 & G & \pi_\lambda  \\
\cline{2-3}
33 & \G & \Gii2{} \\
34 & \G & \Gii{}2 \\
35 & \G & \Gii11 \\
36 & \F & \Fiv2{}{}{} \\
37 & \F & \Fiv1{}{}1 \\
38 & \F & \Fiv{}{}{}2 \\
39 & \F & \Fiv{}{}1{} \\
40 & \E6 & \Evi1{}{}{}{}1 \\
41 & \E6 & \Evi{}2{}{}{}{} \\
42 & \E7 & \Evii2{}{}{}{}{}{} \\
43 & \E7 & \Evii{}{}{}{}{}1{} \\
44 & \E8 & \Eviii1{}{}{}{}{}{}{} \\
45 & \E8 & \Eviii{}{}{}{}{}{}{}2 \\
\hline
\end{array} \]
\end{lem}

To begin with, we note the representations listed in the tables
of Lemma~\ref{lem:k=4} associated to a symmetric space, namely, 
numbers $1$, $2$, $4$, $15$, $24$ ($n=4)$, $27$ and $28$.

\begin{lem}\label{lem:roots1}
Let $\pi_\lambda$ be a representation of real type with $k(\lambda)=4$.
Assume that $\lambda=\frac12(3\beta_i+\beta_j)$ where $\beta_i$, 
$\beta_j\in\mathcal B$ are distinct long roots and that
there exists a Weyl group involution $s$ such that
$s\beta_i=\beta_j$. Then a real form of 
$\pi_\lambda$ is not of class $\mathcal O^2$.
\end{lem}

%\begin{rmk}
%\em
%All the roots in $\mathcal B$ are long unless $G$ is $\G$ or
%$\SO{4m-1}$ for $m:2$, $3\ldots$
%\end{rmk}

\Pf We have that
$\mu=ss_{\beta_i}(\lambda)=\frac12(\beta_i-3\beta_j)$ is a 
weight and $\lambda-\mu=\beta_i+2\beta_j$, $\lambda+\mu=2\beta_i-\beta_j$.
Since $\beta_i$ and $\beta_j$ are orthogonal, we get that
$||\lambda\pm\mu||^2=5$ (normalizing the length of a long root
to be~$1$). On the other hand, for $\alpha'$, $\alpha''\in\Delta$
we have that
$||\alpha'+\alpha''||^2\leq(||\alpha'||+||\alpha''||)^2\leq4$,
so $\lambda\pm\mu$ is neither a root nor a sum of two roots. 
Thus a $\mu$-weight vector is not in  
$\mathcal U^2(\Lg^c)v_\lambda+\mathcal U^2(\Lg^c)\epsilon(v_\lambda)$,
violating condition~$(C_2)$.  \EPf

\medskip

One can use the last lemma together with the tables in Appendix~A 
to get rid of the representations listed in 
Lemma~\ref{lem:k=4} under the numbers $3$, $12$, $17$, $21$, $30$
and $37$. 

We next eliminate the representations $\pi_\lambda$ with
$\lambda=2(\mbox{highest root})=2\beta_1=\frac12(4\beta_1)$, using the 
following observations.

\begin{lem}\label{lem:roots2}
Let $\pi_\lambda$ be a representation of real type with $\lambda=2\beta_1$
and assume there is a nonzero weight $\mu$ such that $(\lambda,\mu)=0$.
Then a real form of $\pi_\lambda$ is not of class $\mathcal O^2$.
\end{lem}

\Pf We check that condition~$(C_2)$ is violated. In fact,
$||\lambda\pm\mu||^2=||\lambda||+||\mu||^2>4||\beta_1||^2=4$ (normalizing
the length of a long root to be~$1$), but $\alpha'$, $\alpha''\in\Delta$
implies that
$||\alpha'+\alpha''||^2\leq(||\alpha'||+||\alpha''||)^2\leq4$,
and we are done as in the proof of Lemma~\ref{lem:roots1}. \EPf

\begin{cor}\label{cor:2b}
Let $\pi_\lambda$ be a representation of real type with $\lambda=2\beta_1$.
If there is a long root $\beta_j\in\mathcal B$ with $\beta_j\neq\beta_1$,
then a real form of $\pi_\lambda$ is not of class $\mathcal O^2$.
\end{cor}

\Pf For in this case $-2\beta_1$, $-\beta_1$, $0$, $\beta_1$, $2\beta_1$ 
is the maximal $\beta_1$-string through $2\beta_1$, so $\beta_1$ is a
weight, $\beta_j$ is a weight, too ($\beta_1$ and $\beta_j$ are long roots, 
hence in the same Weyl orbit, see e.~g.~\cite{Bourbaki}, Ch.~VI, \S~1,
no.~1.3, Proposition~11) and $(\beta_1,\beta_j)=0$. \EPf 

\medskip

We use Corollary~\ref{cor:2b} to discard the representations
listed in the tables of Lemma~\ref{lem:k=4} under the numbers
$5$ ($n\geq3$), $16$, $20$, $29$, $36$, $41$, $42$ and $45$. And we use  
Lemma~\ref{lem:roots2} with $\mu=-\alpha_1+\alpha_2$
(resp.~$\mu=\beta_2=\alpha_1$)
to eliminate representation number $5$ for $n=2$ (resp.\ number~$34$). 
There still remain $20$ 
representations to be eliminated in the tables, namely numbers
$6$, $7$, $9$, $10$, $18$, $19$, $22$, $23$, $24$ ($n\geq5$), $25$, $26$, $31$, $32$,
$33$, $35$, $38$, $39$, $40$, $43$ and $44$,
which in the following we go on to
analyze case by case, but before that we want to reformulate 
an argument that has already been used and that
will frequently be used in the sequel. 

\begin{lem}\label{lem:N}
Let $\pi_\lambda$ be a representation of real type and suppose
that its highest weight $\lambda$ can be written as a
root or as a sum of two roots in $N$ different ways
(where the order of the summands is \emph{not} important).
If $0$ is a weight and its multiplicity is bigger than~$N$, then
a real form of 
$\pi_\lambda$ is not of class $\mathcal O^2$. 
\end{lem}

\Pf Write $v_\lambda$ for a highest weight vector of 
$\pi_\lambda$. If condition~$(C_2)$ is satisfied, then the $0$-weight space
must be contained in $\mathcal U^2(\Lg^c)v_\lambda+\mathcal
U^2(\Lg^c)\epsilon(v_\lambda)$. The identity 
$\epsilon X_\alpha=X_{-\alpha}\epsilon$ shows that the
$0$-weight space must in fact be contained in 
$\mathcal U^2(\Lg^c)v_\lambda$. But the complex dimension of the 
intersection of the $0$-weight space with 
$\mathcal U^2(\Lg^c)v_\lambda$ is at most $N$. \EPf

\medskip

In the remainder of this section several 
decompositions of tensor products of representations will be used,
all of which can be found in Table~5 in the Reference 
Chapter of~\cite{OV}. 

\begin{lem}\label{lem:A}
Representation number~$6$ is not of class $\mathcal O^2$. 
\end{lem}

\Pf We label the simple roots of $\A n$ as 
$\An{\!\!\alpha_1}{\!\!\alpha_2}{\!\!\alpha_n}$, where
$\alpha_i=\theta_i-\theta_{i+1}$. Then the root system
is $\Delta=\{\pm(\theta_i-\theta_j):1\leq i<j\leq n+1\}$ and the highest weight
of $\Ani{}11{}$ is $\lambda=\lambda_2+\lambda_{n-1}=\theta_1+\theta_2-
\theta_n-\theta_{n+1}$, so it can decomposed as a sum of two roots 
in exactly two ways:
\begin{eqnarray*}
 \lambda&=&(\theta_1-\theta_{n+1})+(\theta_2-\theta_{n})\\
        &=&(\theta_1-\theta_{n})+(\theta_2-\theta_{n+1}).
\end{eqnarray*}
On the other hand, we find that:
\begin{eqnarray*}
\lefteqn{\Ani{}1{}{}\otimes\Ani{}{}1{}=}\\
  & &\Ani{}11{}\oplus\Ani1{}{}1\oplus\mbox{(trivial)},
\end{eqnarray*}
from where we deduce that the multiplicity of the $0$-weight in 
$\pi_\lambda$ is $\frac{n(n+1)}2-n-1=\frac{n^2-n-2}2>2$, so we can 
apply Lemma~\ref{lem:N}. \EPf

\begin{lem}\label{lem:B1}
Representations numbers~$18$, $10$ and $19$ are not of class $\mathcal O^2$. 
\end{lem}

\Pf We have that 
\begin{eqnarray*}
\Bni{}{}1{}{}&=&\Lambda^3(\Bni1{}{}{}{}),\quad n\geq4, \\
\Biv{}{}{}2&=&\Lambda^4(\Biv1{}{}{}), \\
\Bnii{}{}{}1{}{}&=&\Lambda^4(\Bnii1{}{}{}{}{}),\quad n\geq5.
\end{eqnarray*}
This shows that the multiplicity of the $0$-weight is,
respectively, $n$, $6$ and $\binom{n}2$ (because 
the vector representation 
$\Bns1{}{}$ has weights $\theta_1$, $\theta_2,\ldots,\theta_n$,
$0$, $-\theta_n$, $-\theta_{n-1},\ldots,-\theta_1$).
But the integer $N=3$ in all cases. \EPf

\begin{lem}\label{lem:B2}
Representation number~$9$ is not of class $\mathcal O^2$. 
\end{lem}

\Pf The highest weight is $\lambda=\frac12(3\theta_1+3\theta_2+\theta_3)$.
Now $\mu=\frac12(3\theta_1-3\theta_2+\theta_3)$ is a weight and 
$\lambda+\mu=3\theta_1+\theta_3$, $\lambda-\mu=3\theta_2$ are neither
a root nor a sum of two roots, so we are done as in the proof of 
Lemma~\ref{lem:roots1}. \EPf

\begin{lem}\label{lem:C}
Representations numbers~$22$ and~$23$ are not of class $\mathcal O^2$.
Representation number~$24$ is not of class $\mathcal O^2$ if $n\geq5$.
\end{lem}

\Pf We first observe that:
\begin{itemize}
\item the vector representation $C_n:\Cn1{}{}{}$ has weights $\pm\theta_i$ 
with multiplicity $1$;
\item the representation $\Cn{}1{}{}$ is a subrepresentation of 
$\Lambda^2(\Cns1{}{})$ and has weights $\pm\theta_i\pm\theta_j$ with
multiplicity $1$, and $0$ with multiplicity $n-1$;
\item the representation $\Cni{}{}1{}{}$ is a subrepresentation of 
$\Lambda^3(\Cns1{}{})$ and has weights $\pm\theta_i\pm\theta_j\pm\theta_k$ with
multiplicity $1$, and $\pm\theta_i$ with multiplicity $n-2$;
\item the representation $\Cnii{}{}{}1{}{}$ is a subrepresentation of 
\[\Lambda^4(\Cns1{}{})\] and has weights 
$\pm\theta_i\pm\theta_j\pm\theta_k\pm\theta_l$ with
multiplicity $1$, and $\pm\theta_i\pm\theta_j$ with multiplicity $n-1$, and 
$0$ with multiplicity $\binom{n}2-n$.
\end{itemize}

If $n\geq5$ then the multiplicity of the $0$-weight in 
$\Cnii{}{}{}1{}{}$ is $\frac12n(n-3)\geq5$.
Since $\lambda=\theta_1+\theta_2+\theta_3+\theta_4$ can be written as 
a root or as a sum of two roots in $N=3$ ways only, namely,
\begin{eqnarray*}
 \lambda&=&(\theta_1+\theta_2)+(\theta_3+\theta_4)\\
        &=&(\theta_1+\theta_3)+(\theta_2+\theta_4)\\
        &=&(\theta_1+\theta_4)+(\theta_2+\theta_3),
\end{eqnarray*}
Lemma~\ref{lem:N} says that $\Cnii{}{}{}1{}{}$, $n\geq5$,
cannot be of class $\mathcal O^2$.

Next, it follows from
\begin{eqnarray*}
\lefteqn{\Ciii1{}{}\otimes\Ciii{}{}1=\Ciii1{}1\oplus\Ciii{}1{},} \\
\lefteqn{\Cni1{}{}{}{}\otimes\Cni{}{}1{}{}=}\\
 & &\Cni1{}1{}{}\oplus\Cn{}1{}{}\oplus\Cnii{}{}{}1{}{},
\quad n\geq4,
\end{eqnarray*}
that the multiplicity of the $0$-weight in $\Cni1{}1{}{}$,
$n\geq3$, is given by $\frac12(n-2)(3n-1)\geq4$. But
$\lambda=2\theta_1+\theta_2+\theta_3$ can 
be written as a root or as sum of two roots in $N=2$ ways only:
\begin{eqnarray*}
 \lambda&=&(\theta_1+\theta_2)+(\theta_1+\theta_3)\\
        &=&2\theta_1+(\theta_2+\theta_3).
\end{eqnarray*}
Lemma~\ref{lem:N} now implies that $\Cni1{}1{}{}$, $n\geq3$, is not of class
$\mathcal O^2$. 

One can also see that
$S^2(\Cn{}1{}{})$, $n\geq3$, has $0$ as a weight of multiplicity $\frac32n(n-1)$. 
Now
\begin{eqnarray*}
\lefteqn{S^2(\Ciii{}1{})=\Ciii{}2{}\oplus\Ciii{}1{}
\oplus\mbox{(trivial)},} \\
\lefteqn{S^2(\Cn{}1{}{})=}\\
 & &\Cn{}2{}{}\oplus\Cn{}1{}{}\oplus\Cnii{}{}{}1{}{}
       \oplus\mbox{(trivial)},\;n\geq4,
\end{eqnarray*}
so the multiplicity of the $0$-weight in $\Cn{}2{}{}$, $n\geq3$, is
$n(n-1)\geq6$. 
But $\lambda=2\theta_1+2\theta_2$ can be 
written as a root or as sum of two roots in $N=2$ ways only:
\begin{eqnarray*}
 \lambda&=&(\theta_1+\theta_2)+(\theta_1+\theta_2)\\
        &=&2\theta_1+2\theta_2.
\end{eqnarray*}
Lemma~\ref{lem:N} implies again that $\Cn{}2{}{}$, $n\geq3$, is not of class
$\mathcal O^2$. \EPf

\begin{lem}\label{lem:D1}
Representations numbers~$26$ and $32$ are not of class $\mathcal O^2$. 
\end{lem}

\Pf These representations are just 
\[ \D n:\Lambda^4(\Dns1{}{}{}\;),\qquad n\geq5, \]
where $\Dns1{}{}{}$ is the vector representation. Therefore
the multiplicity of the $0$-weight is $\binom{n}2\geq10$. But the
highest weight $\lambda=\theta_1+\theta_2+\theta_3+\theta_4$ decomposes
as a sum of two roots in $N=3$ ways, so we can use Lemma~\ref{lem:N}. \EPf

\begin{lem}
Representation number~$35$ is not of class $\mathcal O^2$. 
\end{lem}

\Pf Here $\lambda=\lambda_1+\lambda_2=5\alpha_1+3\alpha_2$ decomposes
only as $(2\alpha_1+\alpha_2)+(3\alpha_1+2\alpha_2)$, so that $N=1$.
Now we find that:
\[\Gii1{}\otimes\Gii{}1=\Gii11\oplus\Gii2{}\oplus\Gii1{}, \] 
where $\Gii1{}$ is the vector representation and
$S^2(\Gii1{})=\Gii2{}\oplus\mbox{(trivial)}$. Therefore 
the multiplicity of the $0$-weight in $\Gii11$ is $4$ and 
we can use Lemma~\ref{lem:N}. \EPf

\begin{lem}
Representations numbers~$38$ and~$39$ are not of class $\mathcal O^2$. 
\end{lem}

\Pf First we observe that the $26$-dimensional
representation $\F:\Fiv{}{}{}1$ has weights 
$\frac12(\pm\theta_1\pm\theta_2\pm\theta_3\pm\theta_4)$ and $\pm\theta_i$ 
with multiplicity $1$, and $0$ with multiplicity $2$. Therefore
$\Lambda^2(\Fiv{}{}{}1)$ has $0$ as a weight of multiplicity $13$. 
Now
\[ \Lambda^2(\Fiv{}{}{}1)=\Fiv{}{}1{}\oplus\Fiv1{}{}{}, \]
where $\Fiv1{}{}{}$ is the adjoint representation. It follows that
$\pi_\lambda:\Fiv{}{}1{}$ has $0$ as a weight of multiplicity $9$. 
But its highest weight
$\lambda=\frac12(3\theta_1+\theta_2+\theta_3+\theta_4)$
is not a root and decomposes as a sum of two roots as 
follows:
\begin{eqnarray*}
\lambda & = & \frac12(\theta_1+\theta_2+\theta_3+\theta_4) + \theta_1\\
        & = & \frac12(\theta_1-\theta_2+\theta_3+\theta_4) + (\theta_1+\theta_2)\\
        & = & \frac12(\theta_1+\theta_2-\theta_3+\theta_4) + (\theta_1+\theta_3)\\
        & = & \frac12(\theta_1+\theta_2+\theta_3-\theta_4) +
        (\theta_1+\theta_4).
\end{eqnarray*}
Thus the integer of Lemma~\ref{lem:N} is $N=4$ and $\Fiv{}{}1{}$ is out. 

Similarly, the multiplicity of the $0$ weight in $S^2(\Fiv{}{}{}1)$ is 15, and
\[S^2(\Fiv{}{}{}1)=\Fiv{}{}{}2\oplus\Fiv{}{}{}1\oplus\mbox{(trivial)}.  \]
It follows that 
$\pi_\lambda:\Fiv{}{}{}2$ has $0$ as a weight of multiplicity $12$. 
But its highest weight
$\lambda=2\theta_1$ is not a root and decomposes as a sum of two roots as 
follows:
\begin{eqnarray*}
\lambda & = & (\theta_1)+(\theta_1)\\
        & = & (\theta_1+\theta_2)+(\theta_1-\theta_2)\\
        & = & (\theta_1+\theta_3)+(\theta_1-\theta_3)\\
        & = & (\theta_1+\theta_4)+(\theta_1-\theta_4).
\end{eqnarray*}
Thus the integer of Lemma~\ref{lem:N} is $N=4$ and $\Fiv{}{}{}2$ is out, too. \EPf

\begin{lem}
Representation number~$40$ is not of class $\mathcal O^2$. 
\end{lem}

\Pf We have that 
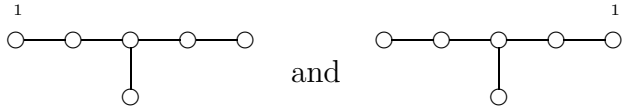
\[ \Evi1{}{}{}{}{}
\quad\mbox{and}\quad\Evi{}{}{}{}{}1 \]
are respectively the $27$-dimensional representation 
and its dual. Therefore the multiplicity of the 
$0$-weight in the tensor product of these representaions is at least $27$. 
But
\[ \Evi1{}{}{}{}{}\otimes\Evi{}{}{}{}{}1=
\Evi1{}{}{}{}1\oplus\Evi{}1{}{}{}{}
\oplus\!\mbox{(trivial)} \]
where the second summand on the right hand side is the adjoint 
representation. It follows that the multiplicity of the $0$-weight in 
$\pi_\lambda:\Evi1{}{}{}{}1$ is at least $20$. But its highest weight 
$\lambda=\theta_8-\theta_7-\theta_6+\theta_5$ is not a root and 
decomposes as a sum of two roots as follows:
\begin{eqnarray*}
\lambda & = & 
\frac12(\theta_8-\theta_7-\theta_6+\theta_5+\theta_4+\theta_3+\theta_2+\theta_1)+
\frac12(\theta_8-\theta_7-\theta_6+\theta_5-\theta_4-\theta_3-\theta_2-\theta_1) \\
 & =  &
\frac12(\theta_8-\theta_7-\theta_6+\theta_5+\theta_4+\theta_3-\theta_2-\theta_1)+
\frac12(\theta_8-\theta_7-\theta_6+\theta_5-\theta_4-\theta_3+\theta_2+\theta_1)
\\
 & = &
\frac12(\theta_8-\theta_7-\theta_6+\theta_5+\theta_4-\theta_3+\theta_2-\theta_1)+
\frac12(\theta_8-\theta_7-\theta_6+\theta_5-\theta_4+\theta_3-\theta_2+\theta_1)
\\
 & = & 
\frac12(\theta_8-\theta_7-\theta_6+\theta_5-\theta_4+\theta_3+\theta_2-\theta_1)+
\frac12(\theta_8-\theta_7-\theta_6+\theta_5+\theta_4-\theta_3-\theta_2+\theta_1).
\end{eqnarray*}
Thus the integer of Lemma~\ref{lem:N} is $N=4$ and we are done. \EPf

\begin{lem}
Representation number~$43$ is not of class $\mathcal O^2$. 
\end{lem}

\Pf We have that $\Evii{}{}{}{}{}{}1$ is the $56$-dimensional 
representation. It is self-dual and does not have $0$ as a weight. 
Therefore the multiplicity of the $0$-weight
in the exterior square of this representation is~$28$. But
\[ \Lambda^2(\;\Evii{}{}{}{}{}{}1\;) = \Evii{}{}{}{}{}1{} \oplus
\mbox{(trivial)}. \]
Therefore the multiplicity of the $0$-weight in $\pi_\lambda$ is~$27$. 
But its highest weight
$\lambda=\theta_8-\theta_7+\theta_6+\theta_5$ is not a root
and one can easily check that it 
decomposes as a sum of two roots in exactly $N=17$ different ways. \EPf 

\begin{lem}
Representation number~$44$ is not of class $\mathcal O^2$. 
\end{lem}

\Pf Here we slightly change the argument used in 
the previous lemmas. The highest weight
$\lambda=2\theta_8$ and
$(\lambda,\theta_8-\theta_7)\neq0$,
so $\nu=\lambda-(\theta_8-\theta_7)=\theta_8+\theta_7$ is a weight. 
Now $\theta_8+\theta_7$ and $\theta_6+\theta_5$ are long roots of $\E8$, 
hence in the same Weyl orbit. It follows that $\mu=\theta_6+\theta_5$
is a weight, too. Observe that each of $\lambda\pm\mu$ is not a root
and can be decomposed as a sum of two roots in only one way:
\begin{eqnarray*}
\lambda+\mu & = & (\theta_8+\theta_6)+(\theta_8+\theta_5), \\
\lambda-\mu & = & (\theta_8-\theta_6)+(\theta_8-\theta_5). 
\end{eqnarray*}
The proof will be complete if we show that the multiplicity
of $\mu$ as a weight is greater than $2$, since then we will have
that condition~$(C_2)$ is violated. 

In fact, the multiplicity of $\mu$ is the same as the multiplicity 
of $\nu$. We use Freudenthal's formula as it is stated 
in~\cite{Fulton-Harris}, namely, the multiplicity of $\nu$ as a 
weight of $\pi_\lambda$ is given by the following formula:
\[ m_\nu=\frac2{c(\nu)}\sum_{\alpha\in\Delta^+}\sum_{k\geq1}
                           (\nu+k\alpha,\alpha)m_{\nu+k\alpha}, \]
where $c(\nu)=||\lambda+\rho||^2-||\nu+\rho||^2$, 
$\rho=\frac12\sum_{\alpha\in\Delta^+}\alpha$ and $m_{\nu+k\alpha}$ is the 
multiplicity of $\nu+k\alpha$ as a weight of $\pi_\lambda$. 

We content ourselves with an estimate. Using that
$\rho=\theta_2+2\theta_3+3\theta_4+4\theta_5+5\theta_6+6\theta_7+23\theta_8$,
we compute that $c(\nu)=36||\theta_1||^2$. Since $\nu>0$, all
terms in Freudenthal's formula are nonnegative. Consider the roots
$\alpha=\frac12(\theta_8-\theta_7+\sum_{i=1}^6\epsilon_i\theta_i)$,
$\beta=\frac12(\theta_8-\theta_7-\sum_{i=1}^6\epsilon_i\theta_i)$,
where $\Pi_{i=1}^6\epsilon_i=-1$. Note that
$\alpha+\beta=\theta_8-\theta_7$. Since
$(\lambda,\alpha)\neq0$, we have that $\lambda-\alpha$ is a weight. 
Now we have two strings of weights starting at $\nu$, namely
\[ \begin{array}{l}
\nu,\quad\nu+(\theta_8-\theta_7)=\lambda\qquad\mbox{and} \\
\nu,\quad\nu+\beta=\lambda-\alpha.
\end{array} \]
Note that in fact there are $32$ different possible choices of the signs 
$\epsilon_i=\pm1$, each of which gives rise to a different string of
weights of the type $\nu$, $\nu+\beta$. 
Freudenthal's formula gives 
\begin{eqnarray*}
   m_\nu & \geq & \frac2{36||\theta_1||^2}((\lambda,\theta_8-\theta_7)m_\lambda
                 +32(\lambda-\alpha,\beta)m_{\lambda-\alpha}) \\
         &  =   & \frac1{18}(2\cdot1+32\cdot2\cdot1)=\frac{33}9>3, 
\end{eqnarray*}
so $m_\nu\geq4$ and we are done. \EPf

\begin{lem}
Representations numbers~$25$ and~$31$ are not of class $\mathcal O^2$. 
\end{lem}

\Pf Here we will use the sharper version of condition~$(C_2)$ given by 
Lemma~\ref{lem:CT}, item~(a). These representations are just
\[\pi_\lambda: \Lambda^3(\Dns1{}{}{}\;),\qquad n\geq4, \]
where $\Dns1{}{}{}$ is the vector representation. The vector representation
has weights $\pm\theta_i$, $1\leq i\leq n$; let 
$\{e_i,e_{n+i}:1\leq i\leq n\}$ be a basis of weight vectors such that
$e_i$ (resp., $e_{n+i}$) corresponds to the $\theta_i$-weight 
(resp., $-\theta_i$-weight)
weight and $\epsilon(e_i)=e_{n+i}$, where $\epsilon$ is a invariant
real structure on the representation space. 

We have that $v_\mu=e_1\wedge e_3\wedge e_{n+3}$ is a
$(\mu=\theta_1)$-weight vector for $\pi_\lambda$ and
$\epsilon(v_\mu)=e_{n+1}\wedge e_{n+3}\wedge e_3$. We want to show that
the condition~(a) of Lemma~\ref{lem:CT} is violated by proving that 
the $\theta_2$-weight 
vector $e_2\wedge e_4\wedge e_{n+4}$ is not in the complex span of
\[ v_\mu,\quad\epsilon(v_\mu),\quad X_\alpha v_\mu,\quad
X_\alpha\epsilon(v_\mu),\quad X_\beta X_\alpha(v_\mu),
\quad X_\beta X_\alpha(\epsilon(v_\mu)), \]
for $\alpha$, $\beta\in\Delta$. But this follows immediately from 
the description of the action of a root vector 
$X_\alpha$, $\alpha\in\Delta=\{\pm(\theta_i\pm\theta_j):1\leq i<j\leq n\}$,
on the weight vectors $\{e_i,e_{n+i}:1\leq i\leq n\}$ under the 
vector representation, given by the following matrices:
\begin{eqnarray*}
X_{\theta_i+\theta_j} & = & E_{i,n+j} - E_{j,n+i}, \\
X_{-\theta_i-\theta_j} & = & E_{n+i,j} - E_{n+j,i}, \\
X_{\theta_i-\theta_j} & = & E_{i,j} - E_{n+j,n+i},
\end{eqnarray*}
where $E_{ij}$ denotes the $2n\times 2n$-matrix with $1$
at the $(i,j)$-position and $0$ elsewhere. \EPf 

\begin{lem}
Representation number~$7$ is not of class $\mathcal O^2$. 
\end{lem}

\Pf We use an argument somewhat similar to the one used in the
previous lemma. The representation is 
$\pi_\lambda:\Lambda^3(\Biii1{}{})$, where $\Biii1{}{}$ is the vector 
representation. Let $e_1$, $e_2$, $e_3$, $e_4=\epsilon(e_1)$, 
$e_5=\epsilon(e_2)$, $e_6=\epsilon(e_3)$, $e_7=\epsilon(e_7)$
be weight vectors of the vector representation 
corresponding to the weights $\theta_1$, $\theta_2$, $\theta_3$, $-\theta_1$, 
$-\theta_2$, $-\theta_3$, $0$, respectively, where $\epsilon$ is the 
invariant real structure on the representation space. 

We have that $v_\mu=e_1\wedge e_4\wedge e_7$ is a ($\mu=0$)-weight vector 
for $\pi_\lambda$ and $\epsilon(v_\mu)=-v_\mu$. 
Now condition~(a) of Lemma~\ref{lem:CT} is violated
because the $\theta_2$-weight vector $e_2\wedge e_3\wedge e_6$ is not in the 
complex span of 
\[ v_\mu,\quad X_\alpha v_\mu,\quad X_\beta X_\alpha(v_\mu), \]
for $\alpha$, $\beta\in\Delta$. \EPf

\begin{lem}\label{lem:G2}
Representation number~$33$ is not of class $\mathcal O^2$. 
\end{lem}

\Pf We start with a description of the Lie algebra of $\G$ and of its
$7$-dimensional representation, as it is done in~\cite{Fulton-Harris}. 
Let $\alpha_1$, $\alpha_2$ be the simple roots. We label the other
positive roots as $\alpha_3=\alpha_1+\alpha_2$, 
$\alpha_4=2\alpha_1+\alpha_2$, $\alpha_5=3\alpha_1+\alpha_2$, 
$\alpha_6=3\alpha_1+2\alpha_2$. Choose root vectors $X_i$, $Y_i$ 
for $\alpha_i$, $-\alpha_i$, respectively, $i=1$, $2$,
such that $[H_i,X_i]=2Y_i$, $[H_i,Y_i]=-2X_i$, where
$H_i=[X_i,Y_i]$, $i=1$, $2$. Next define
\[ \begin{array}{ll}
X_3 =  [X_1,X_2], & Y_3 = - [Y_1,Y_2], \\
X_4 = \frac12[X_1,X_3], & Y_4 = - \frac12[Y_1,Y_3], \\ 
X_5 = -\frac13 [X_1,X_4], & Y_5 = \frac13[Y_1,Y_4], \\
X_6 = - [X_2,X_5], & Y_6 = [Y_2,Y_5],
\end{array} \]
and $H_i=[X_i,Y_i]$ for $i=1,\ldots,6$. Then 
$[H_i,X_i]=2Y_i$, $[H_i,Y_i]=-2X_i$ for $i=3,\ldots,6$,
and $X_i$, $Y_i$ are root vectors
for $\alpha_i$, $-\alpha_i$, respectively, $i=1,\ldots,6$.

The $7$-dimensional representation of $\G$ has weights
$\alpha_1$, $-\alpha_1$, $\alpha_3$, $-\alpha_3$, $\alpha_4$,
$-\alpha_4$, $0$ with respective weight vectors 
$v_1$, $w_1$, $v_3$, $w_3$, $v_4$, $w_4$, $u$, such that the 
action of the basis vectors of the Lie algebra of $\G$ is described in
the following table:
\setlength{\extrarowheight}{0cm}
\[ \begin{array}{|r|c|c|c|c|c|c|c|c|c|c|c|c|c|c|}
\hline
 & H_1 & H_2 & X_1 & X_2 & X_3 & X_4 & X_5 & X_6 & Y_1 & Y_2 & Y_3 & Y_4 &
 Y_5 & Y_6 \\
\hline
u & 0 & 0 & 2v_1 & 0 & 2v_3 & 2v_4 & 0 & 0 & 2w_1 & 0 & 2w_3 & 2w_4 & 0 & 0 \\
v_4 & v_4 & 0 & 0 & 0 & 0 & 0 & 0 & 0 & v_3 & 0 & -v_1 & u & w_1 & w_3 \\
v_1 & 2v_1 & -v_1 & 0 & -v_3 & -v_4 & 0 & 0 & 0 & u & 0 & 0 & w_3 & -w_4 & 0
 \\
v_3 & -v_3 & v_3 & v_4 & 0 & 0 & 0 & 0 & 0 & 0 & -v_1 & u & -w_1 & 0 & -w_4
 \\
w_4 & -w_4 & 0 & -w_3 & 0 & w_1 & u & -v_1 & -v_3 & 0 & 0 & 0 & 0 & 0 & 0
 \\
w_1 & -2w_1 & w_1 & u & 0 & 0 & -v_3 & v_4 & 0 & 0 & w_3 & w_4 & 0 & 0 &
 0\\
w_3 & w_3 & -w_3 & 0 & w_1 & u & v_1 & 0 & v_4 & -w_4 & 0 & 0 & 0 & 0 & 0 \\
\hline
\end{array} \]

Now we have the equation 
\[ S^2(\Gii1{})=\Gii2{}\oplus\mbox{(trivial)}, \]
from which we learn three things. First, $\pi_\lambda:\Gii2{}$ is a
subrepresentation of the symmetric square of the $7$-dimensional
representation, so we can read it off the table above. 
Second, the multiplicity of the weights of $\pi_\lambda$;
in particular, $\alpha_4$ is a weight of multiplicity $2$. 
And third, the trivial summand above is spanned by 
$-\frac12u^2+2v_1w_1+2v_3w_3+2v_4w_4$, so that 
an invariant real structure on the 
representation space of the $7$-dimensional 
representation is given by $\epsilon(v_i)=w_i$, $i=1$, $2$, $3$
and $\epsilon(u)=-u$. 

Let $v_\mu=2u^2-8v_1w_1+6v_3w_3+6v_4w_4$, which is a real, ($\mu=0$)-weight vector
for $\pi_\lambda$. We next show that condition~(a) of Lemma~\ref{lem:CT}
is violated because the $\alpha_4$-weight space is not contained in 
the complex span of 
\setcounter{equation}{\value{thm}}
\stepcounter{thm}
\begin{equation}\label{eqn:mu}
 v_\mu,\quad X_\alpha v_\mu,\quad X_\beta X_\alpha(v_\mu), 
\end{equation}
for $\alpha$, $\beta\in\Delta$. In fact, $\alpha_4-0$ can be written
as a root or as a sum of two roots as follows:
\begin{eqnarray*}
 \alpha_4-0 & = & 2\alpha_1 + \alpha_2 \\
            & = & \alpha_1 + (\alpha_1 + \alpha_2) \\
            & = & - \alpha_1 + (3\alpha_1 + \alpha_2 ) \\
            & = & (3\alpha_1 + 2\alpha_2) - (\alpha_1 + \alpha_2) 
\end{eqnarray*}
Now applying the corresponding root vectors to $v_\mu$ we get:
\begin{eqnarray*}
  X_4(v_\mu) & = & 14(uv_4 + v_1v_3), \\
  X_3X_1(v_\mu) & =& X_3(0) = 0, \\
  X_5Y_1(v_\mu) & = & X_5(0) = 0, \\
  Y_3X_6(v_\mu) & = & Y_3(0) = 0. 
\end{eqnarray*}
This shows that the intersection of the complex span of~(\ref{eqn:mu})
with the $\alpha_4$-weight space has dimension~$1$; but the
weight space itself has dimension~$2$. \EPf

\bigskip

{\it Proof of Lemma~\ref{lem:C11/2-classif}.} We eliminate the possibilities 
given by Table~B.1 and not listed in the table. 

Let $G=\C n$, $\pi_\lambda:\;\Cn3{}{}{}$, $n\geq2$. Then
$\lambda=3\theta_1$, $\mu=-3\theta_2$ is a weight and $\lambda\pm\mu$ are
neither roots nor a sum of two roots, so $\pi_\lambda$ does not even
satisfy~$(C_2)$.

Let $G=\C n$, $\pi_\lambda:\;\Cn11{}{}$, $n\geq3$.
We have 
\begin{eqnarray*}
\lefteqn{\Cn1{}{}{}\otimes\Cn{}1{}{}=} \\
 & & \Cn11{}{}\oplus\Cni{}{}1{}{}\oplus\Cn1{}{}{},
\end{eqnarray*}
from where we deduce that the multiplicity of $\mu=\theta_3$ as
a weight of $\Cn11{}{}$ is $2n-2\geq4$. Now $\lambda-\mu$ is not a root and can
be written as a sum of two roots in two ways,
\begin{eqnarray*}
\lambda-\mu & = & (2\theta_1)+(\theta_2-\theta_3) \\
            & = & (\theta_1+\theta_2)+(\theta_1-\theta_3),
\end{eqnarray*}
and $\lambda+\mu=2\theta_1+\theta_2+\theta_3$ is not a root. 
It follows that the dimension of the intersection 
of the $\mu$-weight space with 
$\mathcal U^2(\mathfrak g^c)v_\lambda+\mathcal U^1(\mathfrak g^c)\epsilon(v_\lambda)$
is at most~$2$, so that condition~$(C_{1\frac12})$ is not satisfied. 
 
Let $G=\C n$, $\pi_\lambda:\;\Cni{}{}1{}{}$, $n\geq4$.
Then $\lambda=\theta_1+\theta_2+\theta_3$ and 
this representation is a subrepresentation of the cubic exterior power of the 
vector representation.
Also, $\mu=-\theta_1$ is a weight of multiplicity $n-2\geq2$.  
But one can check by direct computation
 that the dimension of the intersection 
of the $\mu$-weight space with 
$\mathcal U^2(\mathfrak g^c)v_\lambda+\mathcal U^1(\mathfrak g^c)\epsilon(v_\lambda)$
is~$1$, so condition~$(C_{1\frac12})$ is not satisfied. \EPf

%name of the file: epilog.tex

\section{The proof of the converse to a Theorem of Bott and Samelson}\label{sec:BS}
\setcounter{thm}{0}

In this section we classify variationally complete representations. 
The strategy is as follows. We first exclude many irreducible
representations in class~$\mathcal O^2$ from being variationally
complete by making use of Propositions~\ref{prop:common-orbit}
and~\ref{prop:slice}. 
We then observe that the remaining irreducible representations
are all orbit equivalent to isotropy representations of symmetric
spaces, hence variationally complete by Proposition~\ref{prop:equivvc}.
Using Proposition~\ref{prop:red-vc}, item~(d),
we see that arbitrary variationally complete representations 
are precisely those that are orbit equivalent to isotropy 
representations of symmetric spaces. This result can be viewed as a 
converse to Theorem~\ref{thm:BS2} of Bott and Samelson. Finally we classify 
arbitrary variationally complete representations, mainly as a repeated
application of Proposition~\ref{prop:red-vc}, item~(b). 

Recall that the irreducible representations which according to the results
in Section~\ref{sec:O2} can belong to class~$\mathcal O^2$ without being
either the isotropy representation of a symmetric space or of
cohomogeneity one are listed in the tables of Propositions~\ref{prop:quat},
\ref{prop:complex}, \ref{prop:real4}, \ref{prop:real3}, \ref{prop:real21},
\ref{prop:real22} and~\ref{prop:simple}. For the sake of our
forthcoming arguments, these representations are rearranged in a more 
systematic way in Tables~C.1, C.2, C.3 and~C.4 in Appendix~C 
and in the table of Proposition~\ref{prop:pairs}. 

\begin{prop}\label{prop:notvc}
The representations listed in Tables~C.1, C.2 and~C.3 
are not variationally complete.
\end{prop}

\Pf For each representation $\rho$ of the compact connected
Lie group $G$ on the real vector space $V$ listed in these tables, we 
need to exhibit 
a compact symmetric space $X=L/K$ such that:
\begin{enumerate}
\item[(a)] $K$ contains $G$ as a closed subgroup; 
\item[(b)] the isotropy representation $\tilde\rho$ of $X$ restricts to $\rho$ on $G$;
\item[(c)] there exists $p\in V$ such that
the orbits of $\tilde\rho$ and $\rho$ through $p$ are the equal;
\item[(d)] there exists $q\in V$ such that
the orbits of $\tilde\rho$ and $\rho$ through $q$ are different.
\end{enumerate}
It will then follow from Proposition~\ref{prop:common-orbit} that 
$\rho$ is not variationally complete.

Let $G=\SO2\times\Spin9$. Then $X=\SO{18}/\SO2\times\SO{16}$,
$\tilde\rho$ is the tensor product of the vector representations of $\SO2$,
$\SO{16}$ on $\mathbf R^2$, $\mathbf R^{16}$, respectively, and 
we take $p=v_1\otimes v_2\in\mathbf R^2\otimes\mathbf R^{16}$. Note that the 
orbits of $\tilde\rho$ and $\rho$ through $p$ are equal because 
$\SO2$, $\Spin9$ and $\SO{16}$ are all transitive on unit spheres, but 
the orbits through a common regular point~$q$ are different because 
$\tilde\rho$ is of cohomogeneity two and $\rho$ is of cohomogeneity three.

Let $G=\U2\times\SP n$, $n\geq2$. Then
$X=\SU{2+2n}/\mathbf{S}(\U2\times\U{2n})$, 
$\tilde\rho$ can be viewed as the realification of the 
complex tensor product of the vector representations 
of $\U2$ and $\SU{2n}$ on $\mathbf C^2$, $\mathbf C^{2n}$, respectively,
and we take $p=v_1\otimes v_2\in\mathbf C^2\otimes\mathbf C^{2n}$.
The orbits of $\tilde\rho$ and $\rho$ through $p$ are equal because 
$\U2$, $\SP n$ and $\SU{2n}$ are all transitive on unit spheres, but 
the orbits through a common regular point~$q$ are different because 
$\tilde\rho$ is of cohomogeneity two and $\rho$ is of cohomogeneity three.

Let $G=\SU2\times\SP n$, $n\geq2$. Then 
$X=\SP{2+n}/\SP2\times\SP n$,
$\tilde\rho$ is the quaternionic tensor product of the 
vector representations of $\SP2$, $\SP n$ on $\mathbf H^2$, $\mathbf H^n$, 
respectively, and we take 
$p=e_1\otimes f_1 + e_2\otimes f_2\in 
\mathbf H^2\otimes_{\mathbf H}\mathbf H^n$, where 
$\{e_1,e_2\}\subset\mathbf H^2$ and 
$\{f_1,f_2\}\subset\mathbf H^n$ are orthonormal.  
The orbits of $\tilde\rho$ and $\rho$ through $p$ are equal because 
both of them are $8n-6$-dimensional, but
the orbits through a common regular point~$q$ are different because 
$\tilde\rho$ is of cohomogeneity two and $\rho$ is of cohomogeneity three.

We omit the details for the representations of Tables~C.2 and~C.3
because they are analogous to the above (For instance, 
the reason for the existence of the common orbit (item~(c)) is that as in
the first two cases above we are dealing with tensor products of
representations transitive on unit spheres.) Moreover, we will prove
something stronger later in Propositions~\ref{prop:nottaut2}
and~\ref{prop:nottaut3},
namely that these representations are not taut. \EPf

\begin{prop}\label{prop:nottaut}
The representations listed in Table C.4 are not taut.
\end{prop}

\Pf Let $\pi_\lambda$ be a representation of the compact connected
Lie group $G$ on $V_\lambda$ which is in the table. Notice that  $\pi_\lambda$ 
is of real type. Let $\rho$ be a real form acting on $V\subset V_\lambda$
and let $\epsilon$ be the invariant real structure on $V_\lambda$.  
We shall go case by case and show that the slice representation of $\rho$ 
at a specific point $p\in V$ is not taut,
and then use Proposition~\ref{prop:slice} to conclude that
$\rho$ is not taut. 
In fact, we use 
$p=v_\lambda+\epsilon(v_\lambda)\in V$ where $v_\lambda$ is
a highest weight vector for $\pi_\lambda$. 
In the following, we omit the tedious and
sometimes lengthy calculations involved and simply state the final
results.   

Let $G=\SP n\times \SU6$. Then the isotropy subalgebra at $p$
is 
\[ \mathfrak
u(1)+\mathfrak{su}(3)+\mathfrak{su}(3)+\mathfrak{sp}(n-1),\quad n\geq2 \]
The dimension of the orbit through $p$ is $\dim G(p)=4n+17$, the dimension
of the normal space at $p$ is $36n-17$, and the complexified slice representation at $p$
is given by:
\begin{eqnarray*}
\lefteqn{\mbox{(trivial)} }\\ 
 & &\oplus (x^4\otimes\Aii1{}\otimes\Aii{}1\oplus
  x^{-4}\otimes\Aii{}1\otimes\Aii1{}) \otimes \mbox{(trivial)} \\
 & & \oplus (x\otimes\Aii1{}\otimes\Aii{}1\oplus
x^{-1}\otimes\Aii{}1\otimes\Aii1{}) \otimes \Cns1{}{}. 
\end{eqnarray*}
Therefore, the slice representation contains as a summand the realification
of
\[ x\otimes\Aii1{}\otimes\Aii{}1\otimes\Cns1{}{}. \]
But this is a representation of complex type with Dadok invariant $3$.
By Proposition~\ref{prop:k}, it cannot be of class $\mathcal O^2$.
In particular, it is not taut.
                     
Let $G=\SP n\times \Spin{12}$. Then the isotropy subalgebra at $p$
is 
\[ \mathfrak
u(1)+\mathfrak{su}(6)+\mathfrak{sp}(n-1),\quad n\geq2 \]
The dimension of the orbit through $p$ is $\dim G(p)=4n+29$, the dimension
of the normal space at $p$ is $60n-29$, and the complexified slice representation at $p$
is given by:
\begin{eqnarray*}
\lefteqn{\mbox{(trivial)} }\\ 
 & &\oplus (x^4\otimes\Av{}1{}{}{}\oplus
  x^{-4}\otimes\Av{}{}{}1{}) \otimes \mbox{(trivial)} \\
 & & \oplus (x\otimes\Av{}1{}{}{}{}\oplus
x^{-1}\otimes\Av{}{}{}1{}) \otimes \Cns1{}{}. 
\end{eqnarray*}
Therefore, the slice representation contains as a summand the realification
of
\[ x\otimes\Av{}1{}{}{}\otimes\Cns1{}{}. \]
But this is a representation of complex type with Dadok invariant $3$.
By Proposition~\ref{prop:k}, it cannot be of class $\mathcal O^2$.
In particular, it is not taut.
                     
Let $G=\SP n\times \SP3$. Then the isotropy subalgebra at $p$
is 
\[ \mathfrak
u(1)+\mathfrak{su}(3)+\mathfrak{sp}(n-1),\quad n\geq2 \]
The dimension of the orbit through $p$ is $\dim G(p)=4n+11$, the dimension
of the normal space at $p$ is $24n-11$, and the complexified slice representation at $p$
is given by:
\begin{eqnarray*}
\lefteqn{\mbox{(trivial)} }\\ 
 & &\oplus (x^4\otimes\Aii2{}\oplus
  x^{-4}\otimes\Aii{}2) \otimes \mbox{(trivial)} \\
 & & \oplus (x\otimes\Aii2{}\oplus
x^{-1}\otimes\Aii{}2) \otimes \Cns1{}{}. 
\end{eqnarray*}
Therefore, the slice representation contains as a summand the realification
of
\[ x\otimes\Aii2{}\otimes\Cns1{}{}. \]
But this is a representation of complex type with Dadok invariant $3$.
By Proposition~\ref{prop:k}, it cannot be of class $\mathcal O^2$.
In particular, it is not taut.

Let $G=\SP n\times \E7$. Then the isotropy subalgebra at $p$
is 
\[ \mathfrak
u(1)+\mathfrak{e}_6+\mathfrak{sp}(n-1),\quad n\geq2 \]
The dimension of the orbit through $p$ is $\dim G(p)=4n+53$, the dimension
of the normal space at $p$ is $108n-53$, and the complexified slice representation at $p$
is given by:
\begin{eqnarray*}
\lefteqn{\mbox{(trivial)} }\\ 
 & &\oplus (x^4\otimes\Evi1{}{}{}{}{}\oplus
  x^{-4}\otimes\Evi{}{}{}{}{}1\;) \otimes \mbox{(trivial)} \\
 & & \oplus (x\otimes\Evi1{}{}{}{}{}\oplus
x^{-1}\otimes\Evi{}{}{}{}{}1\;) \otimes \Cns1{}{}. 
\end{eqnarray*}
Therefore, the slice representation contains as a summand the realification
of
\[ x\otimes\Evi1{}{}{}{}{}\otimes\Cns1{}{}. \]
But this is a representation of complex type with Dadok invariant $3$.
By Proposition~\ref{prop:k}, it cannot be of class $\mathcal O^2$.
In particular, it is not taut.

Let $G=\SP 1\times \Spin{13}$.  Then the isotropy subalgebra at $p$
is 
\[ \mathfrak
u(1)+\mathfrak{su}(6). \]
The dimension of the orbit through $p$ is $\dim G(p)=45$, the dimension
of the normal space at $p$ is $83$, and the complexified slice representation at $p$
is given by:
\begin{eqnarray*}
\lefteqn{\mbox{(trivial)} }\\ 
 & &\oplus (x^{5}\otimes\Av1{}{}{}{}\oplus
  x^{-5}\otimes\Av{}{}{}{}1)  \\
& &\oplus (x^4\otimes\Av{}1{}{}{}\oplus
  x^{-4}\otimes\Aiv{}{}{}1{})  \\
 & & \oplus (x^3\oplus x^{-3}) \otimes \Av{}{}1{}{}.
\end{eqnarray*}
Therefore, the slice representation contains as a summand the realification
of
\[ x^3\otimes\Av{}{}1{}{}. \]
But this is a representation of complex type with Dadok invariant $3$.
By Proposition~\ref{prop:k}, it cannot be of class $\mathcal O^2$.
In particular, it is not taut.                    

Let $G=\SP 1\times \SP2$.  Then the isotropy subalgebra at $p$
is 
\[ \mathfrak u(1)+\mathfrak{u}(1). \]
The dimension of the orbit through $p$ is $\dim G(p)=11$, the dimension
of the normal space at $p$ is $21$, and the slice representation at $p$
is given by the realification of:
\begin{eqnarray*}
\lefteqn{\mbox{(trivial)} }\\ 
 & & \oplus (x^{-5}\oplus x^{-2}\oplus x^1) \otimes x^2 \\
 & & \oplus (x^{5}\oplus x^{-4}\oplus 2x^{-1}\oplus 2x^2) \otimes x^4 \\
 & & \oplus \,x^3\otimes x^6.
\end{eqnarray*}
Therefore, the slice representation contains as a summand the realification
of $x^{-2}\otimes x^2\oplus x^{-4}\otimes x^4$ which is orbit equivalent to
the circle action 
\[ x^2\oplus x^4 \]
on $\mathbf R^2\oplus\mathbf R^2$.
But the orbit of the circle action 
through $((1,0),(0,1))$ is a torus-knot which is not taut. 

Let $G=\Spin{15}$. Then the isotropy subalgebra at $p$
is $\mathfrak{su}(7)$. The dimension of the orbit through $p$ is $\dim G(p)=57$, the dimension
of the normal space at $p$ is $71$, and the complexified slice representation at $p$
is given by:
\[\mbox{(trivial)}\oplus\Avi{}{}1{}{}{}\oplus\Avi{}{}{}1{}{} \]
Therefore, the slice representation contains as a summand the realification
of
\[ \Avi{}{}1{}{}{} \]
But this is a representation of complex type with Dadok invariant $3$.
By Proposition~\ref{prop:k}, it cannot be of class $\mathcal O^2$.
In particular, it is not taut.  

Let $G=\Spin{17}$. Then the isotropy subalgebra at $p$
is $\mathfrak{su}(8)$. The dimension of the orbit through $p$ is $\dim G(p)=73$, the dimension
of the normal space at $p$ is $183$, and the complexified slice representation at $p$
is given by:
\begin{eqnarray*}
\lefteqn{\mbox{(trivial)}\oplus\Avii{}{}1{}{}{}{} } \\
 & & \oplus\,\Avii{}{}{}1{}{}{}\oplus\Avii{}{}{}{}1{}{} 
\end{eqnarray*}
Therefore, the slice representation contains as a summand the realification
of
\[ \Avi{}{}1{}{}{}{} \]
But this is a representation of complex type with Dadok invariant $3$.
By Proposition~\ref{prop:k}, it cannot be of class $\mathcal O^2$.
In particular, it is not taut.  

Let $G=\Spin7$.  Then the isotropy subalgebra at $p$
is 
\[ \mathfrak u(1)+\mathfrak{su}(2) \]
The dimension of the orbit through $p$ is $\dim G(p)=17$, the dimension
of the normal space at $p$ is $31$, and the slice representation at $p$
is given by:
\begin{eqnarray*}
\lefteqn{\mbox{(trivial)} }\\ 
 & & \oplus (2(x^4\oplus x^{-4})\oplus(x^2\oplus x^{-2})) \otimes\mbox {(trivial)}\\
 & & \oplus (2(x^1\oplus x^{-1})\oplus(x^7\oplus x^{-7})) \otimes\Ai1 \\
 & & \oplus ((x^2\oplus x^{-2})\oplus (x^4\oplus x^{-4}))) \otimes\Ai2 
\end{eqnarray*}
Therefore, the slice representation contains as a summand the realification
of
\[ x^4\oplus x^2 \]
acting on $\mathbf R^2\oplus\mathbf R^2$.
But the orbit through $((1,0),(0,1))$ is a torus-knot which is  not taut. 

Let $G=\Spin9$.  Then the isotropy subalgebra at $p$
is 
\[ \mathfrak u(1)+\mathfrak{su}(3) \]
The dimension of the orbit through $p$ is $\dim G(p)=27$, the dimension
of the normal space at $p$ is $101$, and the complexified slice representation at $p$
is given by:
\begin{eqnarray*}
\lefteqn{\mbox{(trivial)} }\\ 
 & & \oplus ((x^3\oplus x^{-3})\oplus2(x^2\oplus x^{-2})\oplus(x\oplus
 x^{-1})) \otimes\mbox {(trivial)}\\
 & & \oplus ((x\oplus x^{-1})\oplus(x^2\oplus x^{-2}))\otimes\Aii11 \\
 & & \oplus ((3x\oplus2x^0\oplus x^{-3}) \otimes \Aii1{}
     \oplus (3x^{-1}\oplus2x^0\oplus x^3) \otimes \Aii{}1) \\
 & & \oplus ((x^0\oplus x^{-1}) \otimes \Aii2{} 
     \oplus (x^0\oplus x^1) \otimes \Aii{}2)
\end{eqnarray*}
Therefore, the slice representation contains as a summand the realification
of
\[ x\otimes\Aii11 \]
But this representation cannot be of class $\mathcal O^2$ (because
$\Aii11$ is of real type and we use Lemma~\ref{lem:S1}, item (a)). 
In particular, it is not taut.  

Let $G=\SP n\times \Spin{11}$. Then the isotropy subalgebra at $p$
is 
\[ \mathfrak
u(1)+\mathfrak{su}(5)+\mathfrak{sp}(n-1),\quad n\geq2 \]
The dimension of the orbit through $p$ is $\dim G(p)=4n+29$ and the dimension
of the normal space at $p$ is $60n-29$. For $n\geq2$, 
the complexified slice representation at $p$
is given by:
\begin{eqnarray*}
\lefteqn{\mbox{(trivial)} }\\ 
 & &\oplus (x^8\otimes\Aiv1{}{}{}\oplus
  x^{-8}\otimes\Aiv{}{}{}1) \otimes \mbox{(trivial)} \\
& &\oplus (x^6\otimes\Aiv{}1{}{}\oplus
  x^{-6}\otimes\Aiv{}{}1{}) \otimes \mbox{(trivial)} \\
 & & \oplus (x^3\otimes\Aiv1{}{}{}\oplus
x^{-3}\otimes\Aiv{}{}{}1) \otimes \Cns1{}{} \\
& & \oplus (x\otimes\Aiv{}1{}{}\oplus
x^{-1}\otimes\Aiv{}{}1{}) \otimes \Cns1{}{}
\end{eqnarray*}
Therefore, the slice representation contains as a summand the realification
of
\[ x\otimes\Aiv{}1{}{}\otimes\Cns1{}{}. \]
But this is a representation of complex type with Dadok invariant $3$.
By Proposition~\ref{prop:k}, it cannot be of class $\mathcal O^2$.
In particular, it is not taut.             

For $n=1$ we need a different argument. In this case the connected
component~$K$ 
of the isotropy group at $p$ is locally isomorphic to $\U5$, and 
the slice representation at $p$ contains as a summand the realification of
\[ x^8\otimes\Aiv1{}{}{}\oplus x^6\otimes\Aiv{}1{}{}. \]
Denote by $V_1$, $V_2$ the respective representation spaces
of the above summands. 
Then $V_1=\mathbf C^5$ and $V_2=\Lambda^2\mathbf C^5$. 
Choose an orthonormal basis $\{e_1,\ldots,e_5\}$ for $\mathbf C^5$. 
Let $q_2=a e_1\wedge e_2+b e_3\wedge e_4\in V_2$ be a regular point,
where $a$, $b$ are distinct positive real numbers.
We have that the isotropy $K_{q_2}$ is locally isomorphic to
$\SU2\times\SU2\times\U1$ sitting diagonally in $K$. Let $q_1=e_1+e_3\in V_1$.
Now the orbit $K_{q_2}(q_1)$ is diffeomorphic to $S^3\times S^3$, 
whereas the orbit $K(q_1,q_2)$ is diffeomorphic to $\SU5$. 
It is known that the third Betti number of a compact connected simple
Lie group is one. It follows from~Proposition~\ref{prop:red-taut} that
the slice representation at~$p$ is not taut. \EPf

\begin{prop}\label{prop:pairs}
The following is a list of pairs of irreducible representations
$\rho$, $\rho'$ of compact connected Lie groups $G$, $G'$
such that $\rho'$ is the isotropy representation of a
symmetric space, $G$ is a subgroup of $G'$ and $\rho=\rho'|G$ 
is orbit equivalent to~$\rho'$.
\setlength{\extrarowheight}{0.35cm}
\[ \begin{array}{|c|c|}
\hline
G\quad\mbox{and}\quad G' & \rho\quad\mbox{and}\quad \rho' \\
\hline\hline
\begin{array}{c}\SU n\times\SU m,\; n\neq m\\ \mbox{and}\end{array} & 
\begin{array}{c}
                       \hspace{-.7in}         (\Ans{}1\otimes\Ans1{}) \\
                  \qquad\qquad\qquad\oplus\;(\Ans1{}\otimes\Ans{}1)
\quad\mbox{and} \end{array} \\
S^1\times\SU n\times\SU m  &  \begin{array}{c} 
                \hspace{-.7in}            (\Ans{}1\otimes x\otimes\Ans1{}) \\
                          \qquad\qquad\qquad\oplus\;(\Ans1{}\otimes x^{-1}\otimes\Ans{}1) \end{array} \\
\hline
\SU n,\;\mbox{$n$ odd, and} & \An{}1{}\oplus\Anl{}1{}\quad \mbox{and} \\
\U n & (x\otimes\An{}1{})\oplus(x^{-1}\otimes\Anl{}1{}) \\
\hline
\Spin{10}\quad \mbox{and} & \Dv{}{}{}1{}\quad \mbox{and} \\
\SO2\times\Spin{10}  & (x\oplus x^{-1})\otimes\Dv{}{}{}1{} \\  
\hline
\SO2\times\Spin7\quad \mbox{and} & (x\oplus x^{-1})\otimes\Biii{}{}1\quad \mbox{and} \\
\SO2\times\SO8 & (x\oplus x^{-1})\otimes\Div1{}{}{} \\ 
\hline
\SO2\times\G \quad \mbox{and} & (x\oplus x^{-1})\otimes\Gii1{} \quad \mbox{and} \\
\SO2\times\SO7 & (x\oplus x^{-1})\otimes\Biii1{}{}  \\
\hline
\SO3\times\Spin7 \quad \mbox{and} & \Ai2\otimes\Biii{}{}1\quad\mbox{and} \\
\SO3\times\SO8 & \Ai2\otimes\Div1{}{}{} \\ 
\hline
\end{array} \]
\end{prop}

\Pf The first three cases in the table were treated 
in Lemma~\ref{lem:iv}, item~(a).
It remains to analyze the other three cases. 
It suffices to show the existence of a common principal orbit
by Lemma~\ref{lem:principal-orbit}.

Let $G=\SO2\times\Spin7$. Then $\rho$ is the tensor product of the vector
representation of $\SO 2$ on $\mathbf R^2$ and the spin representation of
$\Spin7$ on $\mathbf R^8$. Let 
$p=a_1e_1\otimes f_1+a_2e_2\otimes f_2\in
\mathbf R^2\otimes\mathbf R^8$, where $e_1$, $e_2\in\mathbf
R^2$ are orthonormal, $f_1$, $f_2\in\mathbf R^8$ are orthonormal, and
$a_1$, $a_2$ are 
distinct positive numbers. Then the orbit~$G'p$
is principal and of dimension~$14$. Since the isotropy subalgebra of $G$ 
at~$p$ is easily seen to be $\mathfrak{su}(3)$, it follows that
$Gp$ is also an orbit of dimension~$14$.
Hence, it coincides with $G'p$.
It follows from Lemma~\ref{lem:principal-orbit} that $\rho$ and $\rho'$ are 
orbit equivalent. 

Similarly, one sees that in the last two cases in the table
$\rho$ has respectively an orbit of dimension~$12$, $21$, and therefore
is orbit equivalent to $\rho'$. \EPf

\begin{thm}[Converse to a Theorem of Bott and Samelson]\label{thm:BS}
A variationally complete representation of a compact connected 
Lie group is orbit equivalent to the isotropy representation
of a symmetric space.
\end{thm}

\Pf We first assume that the given representation is irreducible.
A variationally complete, irreducible representation $\rho$ of a compact
connected Lie group~$G$ is of class $\mathcal O^2$. 
If $\rho$ is not the isotropy representation of one of the irreducible
symmetric spaces and has cohomogeneity greater than one, we have obtained
all possibilities for it in Propositions~\ref{prop:quat},
\ref{prop:complex}, \ref{prop:real4}, \ref{prop:real3}, \ref{prop:real21},
\ref{prop:real22} and \ref{prop:simple}, and 
these are listed in the table of Proposition~\ref{prop:pairs} and in
Tables~C.1, C.2, C.3 and~C.4 in Appendix~C. 
Now Propositions~\ref{prop:notvc}
and~\ref{prop:nottaut} eliminate all possibilities for $\rho$,
except those listed in the table of Proposition~\ref{prop:pairs}.
But the same proposition asserts that $\rho$ is orbit equivalent to
another representation $\rho'$ which is the isotropy representation
of a symmetric space.

Now assume that the representation is reducible. Each of its
irreducible factors is
variationally complete by part~(a) of Proposition~\ref{prop:red-vc}
and hence orbit equivalent to the isotropy
representation of an irreducible symmetric space by what we have just proved.
Part~(d) of the same proposition
now implies that the reducible representation is orbit equivalent to the
isotropy representation
of the symmetric space that is the product of the irreducible symmetric spaces
corresponding to the
irreducible summands of the representation. \EPf

\begin{cor}
A variationally complete representation of a compact connected 
Lie group is polar.
\end{cor}

The first paragraph in the proof of Theorem~\ref{thm:BS} combined with 
Theorem~\ref{thm:BMS} to be proved in Section~\ref{sec:other}
also implies:

\begin{thm}[Compare~\cite{E-H1}]\label{thm:EH}
Let $\rho:G\to\mathbf O(V)$ be a variationally complete irreducible
representation of a compact connected Lie group $G$. 
Then there are two cases:
\begin{enumerate}
\item[(a)] The image $\rho(G)$ is the maximal compact connected 
Lie group with its orbits. In this case $\rho$ is the isotropy
representation of an irreducible symmetric space.
\item[(b)] The image $\rho(G)$ is not the maximal compact connected 
Lie group with its orbits. In this case we have:
\begin{enumerate}
\item[(i)] If $\rho$ has cohomogeneity greater than one, then there exists
  a larger group $G'\supset G$ and an isotropy representation of an
irreducible symmetric space $\rho':G'\to\mathbf O(V)$ which is an
extension of $\rho$. Moreover, the possibilities for the pair 
$(\rho,\rho')$ are listed in the table of Proposition~\ref{prop:pairs}. 
\item[(ii)] If $\rho$ has cohomogeneity one, then it is listed in the
  tables of Theorem~\ref{thm:BMS}. 
\end{enumerate}
\end{enumerate}
\end{thm}

In order to complete the classification of variationally complete
representations, we finally discuss the reducible case
and state a result in Theorem~\ref{thm:bergmann}
which, in view of Theorem~\ref{thm:BS},
is equivalent to a result of Bergmann in~\cite{Bergmann}.
First note that there are some trivial ways to
produce variationally complete reducible representations
from known variationally complete representations. 
\smallskip

Let $\rho:G\to\mathbf O(V)$ be a representation of a compact 
connected Lie group $G$. 
\begin{enumerate}
\item[(a)] Assume $G=G'\times K$ where 
$G'$, $K$ are compact connected Lie groups. We say that
$K$ is \emph{nonessential} for $\rho$ if the restriction
$\rho|G'$ is orbit equivalent to the original 
representation $\rho$. Otherwise, we say that $K$ is \emph{essential} for~$\rho$.
\item[(b)] Assume $\rho$ reducible, namely $\rho=\rho_1\oplus\rho_2$ 
where $\rho_i:G\to\mathbf O(V_i)$ and $\dim V_i>0$ ($i=1$,~$2$).
We say that $\rho$ is \emph{splitting} 
if we can write $G=G_1\times G_2$ for compact
connected Lie groups $G_i$ ($i=1$, $2$) such that $\rho$ is equivalent to the 
outer direct sum representation $\rho_1|G_1\hat\oplus\rho_2|G_2:G_1\times
G_2\to\mathbf O(V_1\oplus V_2)$. 
We say that $\rho$ is \emph{almost splitting} if we can write 
$G=G_1\times G_2\times K$ for compact
connected Lie groups $K$, $G_i$ ($i=1$, $2$) such that
$K$ is nonessential at least for one of $\rho_1|G_1\times K$,
$\rho_2|G_2\times K$ and $\rho|G_1\times G_2$ is equivalent to the 
outer direct sum representation $\rho_1|G_1\hat\oplus\rho_2|G_2$.
\end{enumerate}
It is immediate from these definitions that a representation 
with a trivial summand is automatically splitting, and that
splitting implies almost splitting. We also have:

\begin{prop}
An almost splitting representation $\rho=\rho_1\oplus\rho_2$ 
is variationally complete if and only if each one of its summands
$\rho_i$ is variationally complete\footnote{Note that the isotropy 
representations of reducible symmetric spaces are splitting.}. 
\end{prop}

\Pf One direction follows from Proposition~\ref{prop:red-vc},
item~(a); we prove that the summands being variationally complete 
implies the same for the sum. We may assume that $K$ is nonessential for 
$\rho_1|G_1\times K$. Then $\rho_1$ is orbit equivalent to $\rho_1|G_1$
and $\rho_2$ is orbit equivalent to $\rho_2|G_2\times K$, so
$\rho_1|G_1$ and $\rho_2|G_2\times K$ are variationally complete. 
Now $\rho$ is orbit equivalent to $\rho_1|G_1\hat\oplus\rho_2|G_2\times K$
and hence variationally complete, because the outer direct sum of two
variationally complete representations is also variationally complete,
as is easy to see. \EPf

\medskip

We are left with the classification of
variationally complete reducible representations which are 
not almost splitting. This is given in the
following theorem.

\begin{thm}[Compare~\cite{Bergmann}]\label{thm:bergmann}
Let $\rho$ be a variationally complete reducible
representation of a compact connected Lie group $G$ which
is not almost splitting.
Then $\rho$ is one of the following orthogonal representations:
\setlength{\extrarowheight}{0.35cm}
\[ \begin{array}{|c|c|}
\hline 
 G & \rho \\
\hline
\SU4 & (\Aiii1{}{}\oplus\Aiii{}{}1)\oplus\Aiii{}1{}\\
\U4  & (x\otimes\Aiii1{}{}\oplus x^{-1}\otimes\Aiii{}{}1)\oplus\Aiii{}1{}\\
\Spin7 & \Biii1{}{}\oplus\Biii{}{}1 \\
\SO2\times\Spin7 & (x\oplus x^{-1})\otimes\Biii1{}{}\oplus\Biii{}{}1 \\
\Spin8 & \Div1{}{}{}\oplus\parbox[c]{1.8cm}{$\Div{}{}1{}$} \\
\Spin8 & \Div1{}{}{}\oplus\parbox[c]{1.8cm}{$\Div{}{}{}1$} \\
\SO2\times\Spin8 & (x\oplus x^{-1})\otimes\Div1{}{}{}\oplus\Div{}{}1{} \\
\SO2\times\Spin8 & (x\oplus x^{-1})\otimes\Div1{}{}{}\oplus\Div{}{}{}1 \\
\SO3\times\Spin8 & \Ai2\otimes\Div1{}{}{}\oplus\Div{}{}1{} \\
\SO3\times\Spin8 & \Ai2\otimes\Div1{}{}{}\oplus\Div{}{}{}1 \\ 
\hline
\end{array}\]
\end{thm}

{\it Sketch of a proof.}  Bergmann proves her theorem in~\cite{Bergmann}
using Theorem~\ref{thm:EH} as it is stated in~\cite{E-H1}.
We sketch here an argument in a similar vein, but divided into cases in a
different way. Since summands of almost faithful reducible 
representations need not be almost faithful, this is the only place in
this paper where we have to exercise some care in distinguishing between a
group $G$ and its image $\rho(G)$ under a 
representation $\rho:G\to\mathbf O(V)$. 

{\sc $G$ is a simple group.} Let $\rho=\rho_1\oplus\rho_2$ be
a decomposition where $\rho_1$ is irreducible.
Then $\rho_1$, $\rho_2$ are variationally complete. Denote by $V_1$, $V_2$
the representation spaces and by $H_1$, $H_2$ the connected 
components of the principal isotropy
subgroups. We have that $\rho_1$ is an almost faithful
variationally complete irreducible representation of a simple group,
which is orbit equivalent to the restriction $\rho_1|H_2$, where
$H_2$ is a proper subgroup of $G$; see Proposition~\ref{prop:red-vc},
item~(b). Since $\rho_1(H_2)\subset\mathbf O(V_1)$
is not the maximal compact connected group with its orbits,
$\rho_1:H_2\to\mathbf O(V_1)$ is not the isotropy representation 
of a symmetric space. It follows from
Theorem~\ref{thm:EH} that $\rho_1$ has cohomogeneity one, because 
the larger of the two groups for each pair in the table of 
Proposition~\ref{prop:pairs} is never a simple group.
Since $\rho_1$ is an arbitrary irreducible summand of $\rho$,
we conclude that every irreducible summand of $\rho$ must have
cohomogeneity one. Moreover, we know from Proposition~\ref{prop:red-vc},
item~(c), that there cannot be two equivalent summands. The only simple groups
admitting at least two inequivalent cohomogeneity one representations
are $\Spin n$ for $n\in\{3,5,6,7,8,9\}$. 
Using the fact that $\rho_2$ must be orbit equivalent to
$\rho_2|H_1$ we exclude the groups $\Spin3$, $\Spin5$, $\Spin9$ and
get the four examples corresponding to simple groups in the
table.

{\sc $G$ is not semisimple.} Let $\rho_1$ be an irreducible summand 
of $\rho$ such that $\rho_1(G)$ is locally isomorphic to 
$S^1\times G'$, where $G'$ is a nontrivial semisimple factor of $G$.
Let $\rho_2$ be another irreducible summand of $\rho$ 
with the connected component of the principal isotropy group 
being~$H_2$. If 
$\rho_1(H_2)=\rho_1(G)$ for all such $\rho_2$,
then $\rho$ is splitting. 
So we assume there exists such a $\rho_2$ 
with $\rho_1(H_2)\varsubsetneq\rho_1(G)$. 
Now we can view $\rho_1$
as an almost faithful
variationally complete irreducible representation
of $S^1\times G'$, which is orbit equivalent to the restriction 
$\rho_1|H_2'$, where $H_2'$ is a compact connected 
proper subgroup of $S^1\times G'$
and $\rho_1(H_2')=\rho_1(H_2)$.  
It follows from Theorem~\ref{thm:EH}, item~(b), 
that the possible local isomorphism classes of $H_2'$, $S^1\times G'$
and the representations $\rho_1|H_2'$, $\rho_1|S^1\times G'$ 
are given as follows:
\setlength{\extrarowheight}{0.35cm}
\[ \begin{array}{|c|c|c|}
\hline
&H_2'\quad\mbox{and}\quad S^1\times G' & \rho_1|H_2'\quad\mbox{and}\quad
\rho_1|S^1\times G' \\
\cline{2-3}
1 & \SO2\times\G \quad \mbox{and} & (x\oplus x^{-1})\otimes\Gii1{} \quad \mbox{and} \\
&\SO2\times\SO7 & (x\oplus x^{-1})\otimes\Biii1{}{}  \\
\hline
2 & \SO2\times\Spin7\quad \mbox{and} & (x\oplus x^{-1})\otimes\Biii{}{}1\quad \mbox{and} \\
&\SO2\times\SO8 & (x\oplus x^{-1})\otimes\Div1{}{}{} \\ 
\hline
3 & \Spin{10}\quad \mbox{and} & \Dv{}{}{}1{}\quad \mbox{and} \\
&\SO2\times\Spin{10}  & (x\oplus
x^{-1})\otimes\Dv{}{}{}1{} \\  
\hline
4 & \SU n,\;\mbox{$n$ odd, and} & \An{}1{}\oplus\Anl{}1{}\quad \mbox{and} \\
&S^1\times\SU n & (x\otimes\An{}1{})\oplus(x^{-1}\otimes\Anl{}1{}) \\
\hline
5 & \begin{array}{c}\SU n\times\SU m,\;n\neq m\\ \mbox{and}\end{array}  & 
\begin{array}{c}\hspace{-1in}     (\Ans{}1\otimes\Ans1{})\\
       \qquad\qquad\qquad\oplus\;(\Ans1{}\otimes\Ans{}1)
\quad\mbox{and} \end{array} \\
&S^1\times\SU n\times\SU m  & \begin{array}{c}\hspace{-1in} (\Ans{}1\otimes
  x\otimes\Ans1{})\\
               \qquad\qquad\qquad\oplus\; (\Ans1{}\otimes x^{-1}\otimes\Ans{}1) \end{array} \\
\hline
6 & \SP n \quad \mbox{and} & \Cns1{}{} \quad \mbox{and} \\
&S^1\times\SU{2n} & x\otimes\Ans1{}\oplus x^{-1}\otimes\Ans{}1 \\
\hline
7 & \SU n \quad \mbox{and} & \Ans1{}\oplus\Ans{}1 \quad \mbox{and} \\
&S^1\times\SU n & x\otimes\Ans1{}\oplus x^{-1}\otimes\Ans{}1 \\
\hline
8 & S^1\times\SP n \quad \mbox{and} & (x\oplus x^{-1})\otimes\Cns1{}{}\quad \mbox{and} \\
&S^1\times\SU{2n} & x\otimes\Ans1{}\oplus x^{-1}\otimes\Ans{}1 \\
\hline
9 & \SP n \quad \mbox{and} & \Cns1{}{} \quad \mbox{and} \\
&S^1\times\SP n & (x\oplus x^{-1})\otimes\Cns1{}{} \\
\hline
\end{array} \]

In the following we examine different possibilities for~$\rho_2$.
Recall that a variationally complete representation 
does not admit two equivalent nontrivial irreducible summands.
\begin{enumerate}
\item[(i)] $\rho_2$ is as in case~$1$. Here $\rho_2$ is trivial on
$\SO2$, so we can view $\rho_2$ as an almost faithful
representation $\Spin7\times K\to\mathbf
O(V_2)$ ($K$ is a compact connected Lie group) which is
variationally complete, irreducible 
and such that the projection of the connected component 
of its principal isotropy
subgroup in $\Spin7$ is $\G$. The only possibility
for $\rho_2$ is $\Biii{}{}1$, $K=\{1\}$. 
Any other irreducible summand $\rho_3$ of $\rho$
either would satisfy $\rho_1(H_3)=\rho_1(G)$ (where
$H_3$ is the connected component of the principal
isotropy subgroup of~$\rho_3$) or  
would have to be in the same case as~$\rho_2$.
The first possibility would make $\rho$ splitting
and the second possibility is excluded because
$\rho$ cannot have two equivalent irreducible summands.
It follows that $\rho$ 
has exactly two irreducible summands
and it is 
$(x\oplus x^{-1})\otimes\Biii1{}{}\oplus\Biii{}{}1$.
\item[(ii)] $\rho_2$ is as in case~$2$. Here $\rho_2$ is trivial on
$\SO2$, so we can view $\rho_2$ as an almost faithful representation
$\Spin8\times K\to\mathbf O(V_2)$ ($K$ is a compact connected
Lie group) which is 
variationally complete, irreducible
and such that the projection of the connected component of its  
principal isotropy subgroup in $\Spin8$ is  
$\Spin7$. Here $\Spin7$ is in the conjugacy class 
such that the vector representation of $\Spin8$ restricts 
to the spin representation of $\Spin7$. 
The only possibilities for $\rho_2$ are
the half-spin representations of $\Spin8$, $K=\{1\}$.
Similar considerations to those in~(i) show 
that $\rho$ must have exactly two irreducible summands and that
it is either one of 
\[ (x\oplus x^{-1})\otimes\Div1{}{}{}\oplus\Div{}{}1{},\quad
(x\oplus x^{-1})\otimes\Div1{}{}{}\oplus\Div{}{}{}1.\]
\item[(iii)] $\rho_2$ is as in one of the cases~$3$, $4$, $5$ or~$9$. Here
$S^1$ is not essential for $\rho_1$ and $G'\subset\ker\rho_2$.
Also, for any other irreducible summand $\rho_3$ either 
$\rho_1(H_3)=\rho_1(G)$ (where
$H_3$ is the connected component of the principal
isotropy subgroup of~$\rho_3$)
or $G'\subset\ker\rho_3$.
It follows that $\rho$ is almost splitting. 
\item[(iv)] $\rho_2$ is as in case~$6$. This case cannot occur, for otherwise
$\rho_2:\U{2n}\times K\to\mathbf O(V_2)$ ($K$~a compact connected Lie group)
is an almost faithful variationally complete irreducible representation 
such that
the projection of the connected component of 
its principal isotropy subgroup in $\U{2n}$ is 
$\SP n$. But such a $\rho_2$ does not exist. 
\item[(v)] $\rho_2$ is as in case~$7$. 
Here $S^1$ is not essential for $\rho_1$ and $G'\subset\ker\rho_2$.
Also, for any other irreducible summand $\rho_3$ we have either 
one of the following:
$\rho_1(H_3)=\rho_1(G)$ (where
$H_3$ is the connected component of the principal
isotropy subgroup of~$\rho_3$); or $G'\subset\ker\rho_3$;
or $\rho_3$ is as in case~$8$. In the first two possibilities
$\rho$ is almost splitting and in the third one we go to the next case. 
\item[(vi)] $\rho_2$ is as in case~$8$. Here $\rho_2$ is trivial on
$S^1$, so we can view $\rho_2$ as an almost faithful representation 
$\SU{2n}\times K\to\mathbf O(V_2)$ ($K$ is a compact connected Lie group)
which is variationally complete, irreducible 
and such that the projection of the connected component of its 
principal isotropy subgroup in~$\SU{2n}$ is  
$\SP n$. The only possibility
for $\rho_2$ is $\Aiii{}1{}$, $K=\{1\}$, $n=2$.
If $\rho_3$ is a third irreducible summand of $\rho$ 
then either $\rho_1(H_3)=\rho_1(G)$ (where
$H_3$ is the connected component of the principal
isotropy subgroup of~$\rho_3$) or 
$\rho_1(H_3)$ must be as in case~$7$, therefore $\rho_3$ is trivial on 
$\SU4$. So if there is any irreducible summand besides 
$\rho_1$ and $\rho_2$ then $\rho$ is almost splitting. 
Therefore there are no other irreducible summands and $\rho$ is 
\[ (x\otimes\Aiii1{}{}\oplus x^{-1}\otimes\Aiii{}{}1)\oplus\Aiii{}1{}. \]
\end{enumerate}

{\sc $G$ is semisimple but not simple.} Let $\rho_1$ be an irreducible summand 
of $\rho$ such that $\rho_1(G)$ is locally isomorphic to $G'$, a semisimple but 
not simple factor of~$G$. 
Let $\rho_2$ be another irreducible summand of $\rho$ 
with connected component of the principal isotropy group $H_2$. If 
$\rho_1(H_2)=\rho_1(G)$ for all such $\rho_2$,
then $\rho$ is splitting. 
So we assume there exists such a $\rho_2$ 
with $\rho_1(H_2)\varsubsetneq\rho_1(G)$. 
Now we can view $\rho_1$
as an almost faithful
variationally complete irreducible representation
of $G'$, which is orbit equivalent to the restriction 
$\rho_1|H_2'$, where $H_2'$ is a compact connected 
proper subgroup of $G'$
and $\rho_1(H_2')=\rho_1(H_2)$.  
It follows from Theorem~\ref{thm:EH}, item~(b),  
that the possible local isomorphism classes of $H_2'$, $G'$
and the representations $\rho_1|H_2'$, $\rho_1|G'$ 
are given as follows:
\setlength{\extrarowheight}{0.35cm}
\[ \begin{array}{|c|c|c|}
\hline
&H_2'\quad\mbox{and}\quad G' & \rho_1|H_2'\quad\mbox{and}\quad
\rho_1|G' \\
\cline{2-3}
10 & \SP n \quad\mbox{and} & \Cns1{}{} \quad\mbox{and} \\
&\SP1\times\SP n  & \Ai1\otimes\Cns1{}{} \\
\hline
11 & \SO3\times\Spin7\quad \mbox{and} & \Ai2\otimes\Biii{}{}1\quad \mbox{and} \\
&\SO3\times\SO8 & \Ai2\otimes\Div1{}{}{} \\ 
\hline
\end{array} \]
If $\rho_2$ is an in case~$10$, then $\rho$ is almost splitting 
by the same argument as in item~(iii). If $\rho_2$ is as in case~$11$,
then by an argument similar to that in item~(ii) we show that
$\rho$ must be either one of
\[ \Ai2\otimes\Div1{}{}{}\oplus\Div{}{}1{},\quad
\Ai2\otimes\Div1{}{}{}\oplus\Div{}{}{}1.\] \EPf

\section{Taut representations}\label{sec:taut}
\setcounter{thm}{0}

The main result of this section is the following classification theorem.

\begin{thm}\label{thm:taut2}
A taut irreducible 
representation~$\rho$ of a compact connected 
Lie group~$G$ is either orbit equivalent to
the isotropy representation
of a symmetric space or it is one of the following 
orthogonal representations ($n\geq2$):
\setlength{\extrarowheight}{0.35cm}
\[ \begin{array}{|c|c|}
\hline 
 G & \rho \\
\hline
\SO2\times\Spin9 & (x\oplus x^{-1})\otimes\Biv{}{}{}1  \\
\U2\times\SP n & (x\oplus x^{-1})\otimes\Ai1\otimes\Cns1{}{} \\
\SU2\times\SP n & \Ai3\otimes\Cns1{}{} \\
\hline
\end{array}\]
\end{thm}

\begin{cor}
A taut irreducible representation of a compact connected \emph{simple}
Lie group is orbit equivalent to the isotropy representation
of a symmetric space.
\end{cor}

\begin{rmk}
\em 
According to~\cite{Yasukura1}, see also~\cite{Straume2},
the representations given in the table of
Theorem~\ref{thm:taut2} are exactly the irreducible representations
of compact connected Lie groups which have cohomogeneity three
and are not orbit equivalent to the isotropy representation of a symmetric
space. 
\end{rmk}

It follows from the main results of Section~\ref{sec:BS} 
that Theorem~\ref{thm:taut2} is equivalent to the following:

\begin{thm}\label{thm:taut3}
A taut irreducible 
representation~$\rho$ of a compact connected 
Lie group~$G$ is one of the following:
\begin{enumerate}
\item[(a)] The isotropy representation of an irreducible symmetric space 
(listed in Tables~8.11.2 and~8.11.5 of~\cite{Wolf}).
\item[(b)] A cohomogeneity one representation not listed in~(a)
(cf.~Theorem~\ref{thm:BMS}). 
\item[(c)] A representation orbit equivalent to the 
isotropy representation of an irreducible symmetric space and not 
listed in~(a) or~(b)  
(cf.~Proposition~\ref{prop:pairs}).
\item[(d)] One of the three exceptional cases given in Theorem~\ref{thm:taut2}.
\end{enumerate}
\end{thm}

{\it Proof of Theorems~\ref{thm:taut2} and~\ref{thm:taut3}.} A taut 
irreducible representation $\rho$ 
of a compact connected Lie group~$G$ is of class $\mathcal O^2$. 
If $\rho$ is not orbit equivalent to the isotropy representation of 
one of the irreducible
symmetric spaces, 
we have obtained all possibilities for it, and these are listed in
Tables~C.1, C.2, C.3 and~C.4 in Appendix~C.
Now we already know from Proposition~\ref{prop:nottaut}
that the representations in Table~C.4 are not taut. 
Therefore the assertion follows from Propositions~\ref{prop:nottaut2},
\ref{prop:nottaut3} and~\ref{taut} below. \EPf

\medskip

We start by considering Table~C.3.

\begin{prop}\label{prop:nottaut2}
The representations listed in Table~C.3 are not taut.
\end{prop}

\Pf Let $G=\G\times\Spin7$. 
Then $\rho$ is the tensor product of the $7$-dimensional 
representation of $\G$ on $\mathbf R^7$ and the spin
representation of $\Spin7$ on $\mathbf R^8$. 
Let $p=v_1\otimes v_2\in\mathbf R^7\otimes\mathbf R^8$. 
The isotropy subalgebra at $p$ is $\mathfrak{su}(3)+\mathfrak g_2$.
The dimension of the orbit through $p$ is $\dim G(p)=13$, the dimension
of the normal space at $p$ is $43$, and the slice representation at $p$
minus the trivial component $\mathbf Rp$ is given by the tensor product of the
vector representation of $\mathfrak{su}(3)$ on $\mathbf R^6$ 
and the $7$-dimensional representation of $\mathfrak g_2$ 
on $\mathbf R^7$. 
Since this representation is of complex type and Dadok invariant~$3$, 
it follows from Proposition~\ref{prop:k} that it cannot be of class $\mathcal O^2$.
In particular, it is not taut. Now Proposition~\ref{prop:slice} implies
that $\rho$ is not taut. 

Let $G=\Spin7\times\Spin9$. 
Then $\rho$ is the tensor product of the spin
representation of $\Spin7$ on $\mathbf R^8$ and the spin
representation of $\Spin9$ on $\mathbf R^{16}$. 
Let $p=v_1\otimes v_2\in\mathbf R^8\otimes\mathbf R^{16}$. 
The isotropy subalgebra at $p$ is $\mathfrak g_2+\mathfrak{spin}(7)$.
The dimension of the orbit through $p$ is $\dim G(p)=22$, the dimension
of the normal space at $p$ is $106$, and the slice representation at $p$
minus the trivial component $\mathbf Rp$ is given by the real tensor product of the
$7$-dimensional representation of $\mathfrak g_2$ and the vector
representation of $\mathfrak{spin}(7)$ plus the real tensor product of the
$7$-dimensional representation of $\mathfrak g_2$ and the spin
representation of $\mathfrak{spin}(7)$.
Since we already know that the tensor product of the 
$7$-dimensional representation of $\mathfrak g_2$
and the spin representation of $\mathfrak{spin}(7)$
is not taut (see the last paragraph),
it follows from Proposition~\ref{prop:slice} that $\rho$ is not taut
either.

Let $G=\G\times\Spin9$. 
Then $\rho$ is the tensor product of the $7$-dimensional 
representation of $\G$ on $\mathbf R^7$ and the spin
representation of $\Spin9$ on $\mathbf R^{16}$. 
Let $p=v_1\otimes v_2\in\mathbf R^7\otimes\mathbf R^{16}$. 
The isotropy subalgebra at $p$ is $\mathfrak{su}(3)+\mathfrak{spin}(7)$.
The dimension of the orbit through $p$ is $\dim G(p)=21$, the dimension
of the normal space at $p$ is $91$, and the slice representation at $p$
minus the trivial component $\mathbf Rp$ is given by 
the real tensor product of the
vector representations of $\mathfrak{su}(3)$ and $\mathfrak{spin}(7)$
plus the real tensor product of the
vector representation of $\mathfrak{su}(3)$ and 
the spin representation of $\mathfrak{spin}(7)$. 
Since the second summand is a representation of complex 
type and Dadok invariant~$3$, 
it follows from Proposition~\ref{prop:k} that it cannot be of class $\mathcal O^2$.
In particular, it is not taut. Now Proposition~\ref{prop:slice} implies
that $\rho$ is not taut. 

Consider\footnote{We use the convention that $G\cdot H$ refers to a quotient
of $G\times H$ by a finite central subgroup.}
 $G=\SP1\cdot\SP n\times\G$, $n\geq2$. Let $\tau_n$ be a real form of 
$\Ai1\otimes\An1{}{}{}$. Then $\tau_n$ can be realized as the representation
of $\SP1\cdot\SP n $ on $\mathbf R^{4n}$ given by $\tau_n(q,A)x=Axq^{-1}$,
where $q\in\SP1$, $A\in\SP n$, $x\in\mathbf H^n$ and we identify 
$\mathbf H^n\cong\mathbf R^{4n}$. Now $\rho$ is the real tensor product of 
$\tau_n$ and the $7$-dimensional representation of $\G$. 
Let $p=v_1\otimes v_2\in\mathbf R^{4n}\otimes\mathbf R^7$. 
The isotropy subalgebra at~$p$ is isomorphic to
$\mathfrak{sp}(1)+\mathfrak{sp}(n-1)+\mathfrak{su}(3)$.
The dimension of the orbit through $p$ is $\dim G(p)=4n+5$, the dimension
of the normal space at $p$ is $24n-5$, and the slice representation at $p$
minus the trivial component $\mathbf Rp$ is given by
the real tensor product of $\tau_{n-1}$ and the vector representation of
$\mathfrak{su}(3)$ plus the real tensor product of the adjoint
representation of $\mathfrak{sp}(1)$ and the the vector representation of
$\mathfrak{su}(3)$. Since the second summand is a representation of complex 
type and Dadok invariant~$3$, 
it follows from Proposition~\ref{prop:k} that it cannot be of class $\mathcal O^2$.
In particular, it is not taut. Now Proposition~\ref{prop:slice} implies
that $\rho$ is not taut.  

Consider 
$G=\SP1\cdot\SP n\times\Spin7$, $n\geq2$. Let $\tau_n$ be 
as above. Now $\rho$ is the tensor product of 
$\tau_n$ and the spin representation of $\Spin7$. 
Let $p=v_1\otimes v_2\in\mathbf R^{4n}\otimes\mathbf R^8$. 
The isotropy subalgebra at~$p$ is isomorphic to
$\mathfrak{sp}(1)+\mathfrak{sp}(n-1)+\mathfrak g_2$.
The dimension of the orbit through $p$ is $\dim G(p)=4n+6$, the dimension
of the normal space at $p$ is $28n-6$, and the slice representation at $p$
minus the trivial component $\mathbf Rp$ is given by
the real tensor product of $\tau_{n-1}$ and the $7$-dimensional representation of
$\mathfrak g_2$ plus the real tensor product of the adjoint
representation of $\mathfrak{sp}(1)$ and the the $7$-dimensional representation of
$\mathfrak g_2$. Since the first summand is not a taut representation
by the above (for $n=2$ we refer to the case of $G=\SO4\times\G$ which
will be dealt with in Proposition~\ref{prop:nottaut3}), 
it follows from Proposition~\ref{prop:slice} 
that $\rho$ is not taut either. 

Let $G=\Spin7\times\Spin7$. 
Then $\rho$ is the tensor product of the spin
representations of each of the factors on $\mathbf R^8$.
Let $p=v_1\otimes v_2\in\mathbf R^8\otimes\mathbf R^8$. 
The isotropy subalgebra at $p$ is $\mathfrak g_2+\mathfrak g_2$.
The dimension of the orbit through $p$ is $\dim G(p)=14$, the dimension
of the normal space at $p$ is $50$, and the slice representation at $p$
minus the trivial component $\mathbf Rp$ is given by the real tensor product of the
$7$-dimensional representations of each of the factors. 
Since this representation is not taut (see next case),
it follows from Proposition~\ref{prop:slice} that $\rho$ is not taut
either.

Consider $G=\G\times\G$. Then $\rho$ is the tensor product of the 
$7$-dimensional representations of each of the factors on $\mathbf R^7$.
Let $\{e_1,\ldots,e_7\}$ be an orthonormal basis of $\mathbf R^7$
and take $p=e_1\otimes e_1+\ldots+e_7\otimes e_7$. Then the isotropy
subalgebra at~$p$ is 
\[ \mathfrak g_p=\{(X,X):X\in\mathfrak g_2\} \]
which is isomorphic to $\mathfrak g_2$. Now the restriction of $\rho$ to 
$\mathfrak g_p$ is a real form of 
\[ \mathfrak g_2:\Gii1{}\otimes(\Gii1{})^*=\Gii1{}\otimes\Gii1{} \]
acting on $\mathbf R^{49}$,
and $\mathbf R^{49}=T_pG(p)\oplus N_pG(p)$ is a $\mathfrak g_p$-invariant 
decomposition. We have the following decomposition into irreducible components:
\begin{eqnarray*}
 \mathfrak g_2:\otimes^2(\Gii1{}) & = & S^2(\Gii1{})\oplus\Lambda^2(\Gii1{}) \\
                    & = &
                    (\Gii2{}\oplus\mbox{(trivial)})\oplus(\Gii{}1\oplus\Gii1{}).
\end{eqnarray*}
Since $T_pG(p)$ is $14$-dimensional, it follows that the $27$-dimensional 
real representation $\Gii2{}$ must be a component of the slice
representation at~$p$. But we saw in Lemma~\ref{lem:G2}
that it is not of class $\mathcal O^2$. Now
Proposition~\ref{prop:slice} implies that $\rho$ is not taut.  \EPf

\begin{prop}\label{prop:nottaut3}
The representations listed in Table~C.2 are not taut.
\end{prop}

\Pf Let $G=\SO m\times\G$, $m\geq4$. 
Then $\rho$ is the tensor product of the vector
representation of $\SO m$ on $\mathbf R^m$ and the $7$-dimensional 
representation of $\G$ on $\mathbf R^7$. 
Let $p=v_1\otimes v_2\in\mathbf R^m\otimes\mathbf R^7$. 
The isotropy subalgebra at $p$ is $\mathfrak{so}(m-1)+\mathfrak{su}(3)$.
The dimension of the orbit through $p$ is $\dim G(p)=m+5$, the dimension
of the normal space at $p$ is $6m-5$, and the slice representation at $p$
minus the trivial component $\mathbf Rp$ is given by the tensor product of the
vector representations of $\mathfrak{so}(m-1)$ and $\mathfrak{su}(3)$
on $\mathbf R^{m-1}$ and $\mathbf R^6$, respectively. 
Since this is a representation of complex type and Dadok invariant~$3$, 
it follows from Proposition~\ref{prop:k} that it cannot be of class $\mathcal O^2$.
In particular, it is not taut. Now Proposition~\ref{prop:slice} implies
that $\rho$ is not taut. 

The proof for $G=\SO 3\times\G$ will be done in
Lemma~\ref{lem:final1} below, since this case is more involved. 

Let $G=\SO m\times\Spin7$, $m\geq4$. 
Then $\rho$ is the tensor product of the vector
representation of $\SO m$ on $\mathbf R^m$ and the spin
representation of $\Spin7$ on $\mathbf R^8$. 
Let $p=v_1\otimes v_2\in\mathbf R^m\otimes\mathbf R^8$. 
The isotropy subalgebra at $p$ is $\mathfrak{so}(m-1)+\mathfrak g_2$.
The dimension of the orbit through $p$ is $\dim G(p)=m+6$, the dimension
of the normal space at $p$ is $7m-6$, and the slice representation at $p$
minus the trivial component $\mathbf Rp$ is given by the real tensor product of the
vector representation of $\mathfrak{so}(m-1)$ and the $7$-dimensional 
representation of $\mathfrak g_2$.
Since this representation is not taut by the above,
it follows from Proposition~\ref{prop:slice} that $\rho$ is not taut
either.

Let $G=\SO m\times\Spin9$, $m\geq4$. 
Then $\rho$ is the tensor product of the vector
representation of $\SO m$ on $\mathbf R^m$ and the spin
representation of $\Spin9$ on $\mathbf R^{16}$. 
Let $p=v_1\otimes v_2\in\mathbf R^m\otimes\mathbf R^{16}$. 
The isotropy subalgebra at $p$ is $\mathfrak{so}(m-1)+\mathfrak{spin}(7)$.
The dimension of the orbit through $p$ is $\dim G(p)=m+14$, the dimension
of the normal space at $p$ is $15m-14$, and the slice representation at $p$
minus the trivial component $\mathbf Rp$ is given by the real tensor product of the
vector representations of $\mathfrak{so}(m-1)$ and $\mathfrak{spin}(7)$
plus the real tensor product of the
vector representation of $\mathfrak{so}(m-1)$ and the spin representation
of $\mathfrak{spin}(7)$. 
If $m\geq5$ then the second summand is not a taut representation by the above
and it follows from Proposition~\ref{prop:slice} that $\rho$ is not taut
either. In the case $m=4$ we give a special argument. Here we have that 
the connected component of the isotropy subgroup at $p$ is 
$K=\SO3\times\Spin7$. 
Let $\nu_1$, $\nu_2$ be respectively the irreducible summands of dimensions
$21$, $24$ of the slice representation at $p$, and denote with 
$V_1$, $V_2$ the representation spaces. Note that 
both $\nu_1$ and $\nu_2$ are taut representations of $K$,
since $\nu_1$ is the isotropy representation of a real Grassmann
manifold and $\nu_2$ is orbit equivalent to the isotropy representation 
of a real Grassmann manifold.
Let $q_1\in V_1$ be a regular point. Then the isotropy group
$K_{q_1}$ is isomorphic to $\mathbf Z_2^2\times\Spin4$. 
We can choose $q_2\in V_2$ so that
the isotropy $(K_{q_1})_{q_2}$ is contained in $\{1\}\times\Spin4$. 
Hence the orbit $K_{q_1}(q_2)$ is disconnected. Since 
the orbit $K(q_1,q_2)$ is connected, Proposition~\ref{prop:red-taut} implies that
the slice representation at $p$ is not taut.
Now Proposition~\ref{prop:slice} yields that $\rho$ is not taut either.

We postpone the proof for the case $\SO3\times\Spin9$ to Lemma~\ref{lem:final2} below.

Let $G=\SU n\times\SP m$, $n\geq3$ and $m\geq2$. Then $\rho$ is the
realification of the tensor product of the vector representations 
of $\SU n$ on $\mathbf C^n$ and of $\SP m$ on $\mathbf C^{2m}$. 
Consider complex bases $\{e_1,\ldots,e_n\}$ for $\mathbf C^n$
and $\{f_1,\ldots,f_m,f_{m+1},\ldots,f_{2m}\}$ for $\mathbf C^{2m}$
where $f_{m+j}=\epsilon f_j$, $\epsilon$ the quaternionic 
structure on $\mathbf C^{2m}$.
Let $p=e_1\otimes f_1$.
Then
the isotropy $K=G_p$ is isomorphic to $T_0\cdot\SU{n-1}\times\SP{m-1}$, 
where $T_0$ is the circle group generated by 
\begin{eqnarray*}
\lefteqn{\varphi_0(t)=(\mbox{diag}(e^{(n-1)it},e^{-it},\ldots,e^{-it}),} \\
  & &\qquad\mbox{diag}(e^{-(n-1)it},1,\ldots,1;e^{(n-1)it},1,\ldots,1))
 \in\SU n\times\SP m,
\end{eqnarray*}
and the slice representation 
$\nu_p$ minus the trivial component corresponding to
the radial direction $\mathbf Rp$
decomposes into irreducible components 
as $\nu_1\oplus\nu_2$, where $\nu_1$ restricted to 
$\SU{n-1}\times\SP{m-1}$ is 
the realification of the complex tensor product of the vector
representations of the factors
and $\nu_2$ is the realification
of the vector representation on $\SU{n-1}$ and it is trivial on $\SP{m-1}$. 
Denote the representation spaces by $V_1$, $V_2$.
Then $V_1$ is spanned by $e_\alpha\otimes f_\beta$ and $i(e_\alpha\otimes f_\beta)$,
for $2\leq\alpha\leq n$, $2\leq\beta\leq m$ or $m+2\leq\beta\leq 2m$, 
and $V_2$ is spanned by $e_\alpha\otimes f_{m+1}$ and
$i(e_\alpha\otimes f_{m+1})$, for $2\leq\alpha\leq n$.
Let $q_1=e_2\otimes f_2\in V_1$. 
Then the isotropy $K_1=K_{q_1}$ is isomorphic to 
$T_1\cdot T_2\cdot\SU{n-2}\times\SP{m-2}$, where $T_1$, $T_2$ are 
circle groups respectively generated by  
\begin{eqnarray*}
\lefteqn{\varphi_1(t)=(\mbox{diag}(1,e^{(n-2)it},e^{-it},\ldots,e^{-it}),} \\
     & &\qquad\mbox{diag}(1,e^{-(n-2)it},1,\ldots,1;1,e^{(n-2)it},1,\ldots,1))
   \in \SU n\times\SP m, 
\end{eqnarray*}
and
\begin{eqnarray*}
\lefteqn{\varphi_2(t)=
                 (\mbox{diag}(e^{(n-1)it},e^{-it},e^{-it},\ldots,e^{-it}),} \\
   & &\qquad\mbox{diag}(e^{-(n-1)it},e^{it},1,\ldots,1;e^{(n-1)it},e^{-it},
1,\ldots,1))
   \in \SU n\times\SP m.
\end{eqnarray*}
(Here we have used the hypothesis
$n\geq3$ and $m\geq2$). 
Finally let $q_2=e_2\otimes f_{m+1}\in V_2$. Then the isotropy
$K_2=(K_{q_1})_{q_2}$ is isomorphic to 
$T_3\cdot\SU{n-2}\times\SP{m-2}$, where $T_3$ is the 
circle subgroup of $T_1\cdot T_2$ generated by 
$\varphi_3(t)=\varphi_1(-t)\varphi_2(t)$. It is easy to see that 
$K_1/K_2$ is diffeomorphic to $S^1$, whereas $K/K_2$ is simply-connected.
Therefore Proposition~\ref{prop:red-taut} implies that
$\nu_p$ is not taut. Now Proposition~\ref{prop:slice} 
yields that $\rho$ is not taut either.

The case where $G=\U n\times\SP m$, $n\geq3$ and $m\geq2$, is completely
analogous to the previous case.

Let $G=\SO m\times\SP1\cdot\SP n$, $m\geq3$ and $n\geq2$. Then $\rho$
is the tensor product of
the vector representation of $\SO m$ on $\mathbf R^m$ 
and $\tau_n$, where
$\tau_n$ is a representation on $\mathbf R^{4n}$
as in Proposition~\ref{prop:nottaut2}.
We choose an orthonormal basis $\{e_1,\ldots,e_m\}$ for
$\mathbf R^m$ and an orthonormal basis 
\[ \{f_1,\ldots,f_n,f_1i,\ldots,f_ni,f_1j,\ldots,f_nj,f_1k,\ldots,f_nk\} \]
for $\mathbf R^{4n}$. 
Let $p=e_1\otimes f_1$. 
The connected component of the isotropy subgroup $K=G_p^0$ is isomorphic to
$\SO{m-1}\times\SP1'\cdot\SP{n-1}$, where
\setlength{\extrarowheight}{0cm}
\[ \SP1'=\left\{\left(q,\left(\begin{array}{cccc}
                                     q&&&\\
                                     &1&&\\
                                     &&\ddots&\\
                                     &&&1\end{array}\right)\right):
                  q\in\SP1\right\}\subset\SP1\cdot\SP n. \]
The slice representation at $p$
minus the trivial component $\mathbf Rp$ decomposes as $\nu_1\oplus\nu_2$,
where $\nu_1$ is the real tensor product of the vector representation of
$\SO{m-1}$ and $\tau_{n-1}$, and $\nu_2$ is the real tensor product of 
the vector representation of $\SO{m-1}$ and the adjoint representation of
$\SP1$ ($\nu_2$ is trivial on $\SP{n-1}$). Note that 
$\nu_2$ is the isotropy representation of a real
Grassmann manifold. Let $V_1$, $V_2$
denote the representation spaces of $\nu_1$, $\nu_2$.
Then $V_1$ is spanned by $e_\alpha\otimes f_\beta$, 
$e_\alpha\otimes f_\beta i$, $e_\alpha\otimes f_\beta j$, 
$e_\alpha\otimes f_\beta k$, for $2\leq\alpha\leq m$, $2\leq\beta\leq n$, 
and $V_2$ is spanned by $e_\alpha\otimes f_1 i$,  
$e_\alpha\otimes f_1 j$, $e_\alpha\otimes f_1 k$, for $2\leq\alpha\leq m$.
First consider the case $m\geq4$. Let $q_2\in V_2$ be a regular point,
say $q_2=ae_2\otimes f_1 i + be_3\otimes f_1 j + ce_4\otimes f_1 k$, where 
$a$, $b$, $c$ are pairwise distinct positive real numbers.  
Then the isotropy group $K_{q_2}$ is isomorphic to 
$\SO{m-4}\times Q\cdot\SP{n-1}$, where $Q$ is the quaternion subgroup
$\{\pm1,\pm i,\pm j, \pm k\}$ of $\SP1$. Consider the restriction 
of $\nu_1$ to $K_{q_2}$. We can choose $q_1\in V_1$ so that
the orbit $K_{q_2}(q_1)$ is disconnected,
say $q_1=ae_2\otimes f_2 +be_3\otimes f_2 i$, where 
$a$, $b$ are distinct positive real numbers. In fact, 
if $(A,q,B)\in\SO{m-4}\times Q\cdot\SP{n-1}$ is in the isotropy of 
$q_1$, then we must have $Bf_2q^{-1}=f_2$ and $Bf_2iq^{-1}=f_2i$,
which implies that $Bf_2q^{-1}i=Bf_2iq^{-1}$ and then $q^{-1}i=iq^{-1}$,
so that $q=\pm1$ or $q=\pm i$. In particular, no element of the 
form $(A,j,B)\in\SO{m-4}\times Q\cdot\SP{n-1}$ is in the isotropy of 
$q_1$, which implies that the orbit $K_{q_2}(q_1)$ is disconnected.
Since the orbit $K(q_1,q_2)$ is connected, 
Proposition~\ref{prop:red-taut} implies that
the slice representation at $p$ is not taut.
Now Proposition~\ref{prop:slice} yields that $\rho$ is not taut either.
Finally we analyse the case $m=3$. Then $K=\SO2\times\SP1'\cdot\SP{n-1}$.
Consider a regular point $q_2\in V_2$, say
$q_2=ae_2\otimes f_1 i + be_3\otimes f_1 j$, where
$a$, $b$ are distinct positive real numbers.
Then we have that the isotropy $K_{q_2}$
is isomorphic to $\mathbf Z_4\times\SP{n-1}$. Here $\mathbf Z_4$ is generated 
by $(-1,k)$, where $-1$ is minus the identity matrix in $\SO2$ 
and $k\in\SP1'$. Consider the restriction 
of $\nu_1$ to $K_{q_2}$. We can choose
$q_1\in V_1$ as above such that its isotropy does not contain elements of the 
form $(-1,k,B)$, for $B\in\SP{n-1}$. It follows that 
the orbit $K_{q_2}(q_1)$ is disconnected and the proof follows as above. \EPf

\begin{lem}\label{lem:final2}
The tensor product $\rho$ of the vector representation of $\SO3$ on $\mathbf R^3$
and the spin representation of $\Spin9$ on $\mathbf R^{16}$ is not taut. 
\end{lem}

\Pf Let $p=v_1\otimes v_2\in\mathbf R^3\otimes\mathbf R^{16}$. 
Then the connected component of the isotropy group $K=G_p^0$ is isomorphic to
$\SO2\times\Spin7$. The normal space at $p$ to the orbit $G(p)$ 
decomposes as $\mathbf Rp\oplus V_1\oplus V_2$, where 
$V_1=\mathbf R^2\otimes\mathbf R^7$, $V_2=\mathbf R^2\otimes\mathbf R^8$ are 
the representation spaces for the irreducible components $\nu_1$,
$\nu_2$ of the slice representation at $p$. Here $\nu_1$ is the 
tensor product of the vector representations and $\nu_2$ is the 
tensor product of the vector representation of $\SO2$ and
the spin representation of $\Spin7$. Now let $q=v_3\otimes v_4\in V_2$. 
Then the connected component of the isotropy group 
$K^0_q$ is isomorphic to $\G$ and the slice
representation at $q$ of the representation of $K$ on $V_1\oplus V_2$ is
a representation of $\G$ on $\mathbf R^{22}$ which  
decomposes as 
a trivial component $\mathbf Rq$ plus three copies of the $7$-dimensional
representation of $\G$. Next we shall show that 
the sum $\mu=\mu_1\oplus\mu_2\oplus\mu_3$ 
of three copies of $7$-dimensional representation of $H=\G$ is not
taut. It then follows from Proposition~\ref{prop:slice} that $\rho$ is not taut
either.

In fact, let $W_i$ be the representation space of $\mu_i$. 
We can select $w_i\in W_i$ such that the isotropy groups $H_{w_1}$,
$(H_{w_1})_{w_2}$, $((H_{w_1})_{w_2})_{w_3}$ are respectively
$\SU3$, $\SU2$, $\{1\}$. Since $H_{(w_1,w_2,w_3)}=((H_{w_1})_{w_2})_{w_3}$,
we get that the orbit $H(w_1,w_2,w_3)$ is diffeomorphic to $\G$. 
Now if this orbit is taut, Proposition~\ref{prop:red-taut} implies that 
it must have the product homology of the orbits through 
$w_1$, $w_2$ and $w_3$, namely, $S^6$, $S^5$ and $S^3$. But it is well
known that $\G$ has the product homology of $S^{11}$ and $S^3$
(this also follows from Satz~1 in~\cite{Hopf}). 
Thus, $\mu$ is not taut.\\ \mbox{} \EPf 

\begin{lem}\label{lem:final1}
The tensor product $\rho$ of the vector representation of $\SO3$ on $\mathbf R^3$
and the $7$-dimensional representation of $\G$ is not taut. 
\end{lem}

\Pf We shall use the reduction technique described in Subsection~\ref{subsec:red}.
Note that $G=\SO3\times\G$ is a closed subgroup of $\hat
G=\SO3\times\SO7$ and $\rho$ is the restriction of a representation
$\hat\rho$ of $\hat G$ which is the isotropy representation of a 
real Grassmann manifold.

We view the Cayley numbers $\mathbf O=\mathbf H\oplus e\mathbf H$, 
so that $\{1,e,i,j,k,ei,ej,ek\}$ is a real 
orthonormal basis for $\mathbf O$. The Lie group $\G$ is 
the automorphism group of $\mathbf O$ with respect to 
its nonassociative algebra structure.
We next compute the space of orbits
of $G$ acting by $\rho$ on $V=\mathbf R^3\otimes\mathbf R^7$.
Take $\{e,i,j,k,ei,ej,ek\}$ as a basis for $\mathbf R^7$ and 
take an orthonormal basis $\{f_1,f_2,f_3\}$ for $\mathbf R^3$.  
Recall that the subspace spanned by $\{i,j,k,ei,ej,ek\}$
has a complex structure given by left multiplication by $e$. 
Since all these representations are self-adjoint, it is equivalent to 
consider $\rho$ acting on $\mathbf R^{7*}\otimes\mathbf R^3$ which we may
identify with the space of $3\times7$ real matrices. 

Let 
\[ p_1=\left(\begin{array}{ccccccc}
            1 & 0 & 0 & 0 & 0 & 0 & 0 \\
            0 & 0 & 0 & 0 & 0 & 0 & 0 \\
            0 & 0 & 0 & 0 & 0 & 0 & 0 \end{array} \right), \]
Then the orbit $G(p_1)$ is $8$-dimensional 
and its normal space at~$p_1$ intersects all orbits. So, given an arbitrary orbit
$M$, it contains a point of the form
\[ \left( \begin{array}{ccccccc}
            a & 0 & 0 & 0 & 0 & 0 & 0 \\
            0 & * & * & * & * & * & * \\
            0 & * & * & * & * & * & * \end{array} \right), \]
for some $a\geq0$, for we can multiply on the left by 
\[ \left( \begin{array}{ccc}
           -1 & 0 & 0 \\
            0 & -1 & 0 \\
            0 & 0 & 1 \end{array} \right) \]
to change the sign of $a$, 
if necessary, and remain in the orbit $M$. Now the connected component of
the isotropy subgroup 
at $p_1$, $K_1=\SO2\times\SU3$, acts
on the normal space to $G(p_1)$ by fixing the radial
direction $\mathbf Rp_1$ and by the tensor product of the vector
representations of $\SO2$ and $\SU3$ on the orthogonal complement 
$\mathbf R^{12}$. 
The point 
\[ p_2=\left(\begin{array}{ccccccc}
            0 & 0 & 0 & 0 & 0 & 0 & 0 \\
            0 & 1 & 0 & 0 & 0 & 0 & 0 \\
            0 & 0 & 0 & 0 & 0 & 0 & 0 \end{array} \right), \]
is in~$\mathbf R^{12}$ 
and has a $6$-dimensional $K_1$-orbit. So $M$ has a representative of the form
\[ \left( \begin{array}{ccccccc}
            a & 0 & 0 & 0 & 0 & 0 & 0 \\
            0 & b & 0 & 0 & 0 & 0 & 0 \\
            0 & 0 & * & * & * & * & * \end{array} \right), \]
for some $b\geq0$, since we can also multiply on the left by 
\[ \left( \begin{array}{ccc}
            1 & 0 & 0 \\
            0 & -1 & 0 \\
            0 & 0 & -1 \end{array} \right), \]
if necessary, and remain in the same orbit.
Finally, the connected component of the isotropy subgroup of $K_1$ at $p_2$
is $\SU2$ and it acts on the normal space to $K_1(p_2)$  
with respect to $\mathbf R^{12}$ 
by fixing a  $2$-plane
and by the vector representation on the orthogonal complement. 
So $M$ has an element of the form 
\setcounter{equation}{\value{thm}}
\stepcounter{thm}
\begin{equation}\label{eqn:abcd}
  \left( \begin{array}{ccccccc}
            a & 0 & 0 & 0 & 0 & 0 & 0 \\
            0 & b & 0 & 0 & 0 & 0 & 0 \\
            0 & 0 & c & 0 & d & 0 & 0 \end{array} \right), 
\end{equation}
for some $c\geq0$. We can assume $a\geq b$ by multiplying on the left by 
\[ \left( \begin{array}{ccc}
            0 & 1 & 0 \\
            1 & 0 & 0 \\
            0 & 0 & -1 \end{array} \right), \]
and on the right by
\[ \left( \begin{array}{ccccccc}
                0 & 1 & 0 & 0 & 0 & 0 & 0 \\
                1 & 0 & 0 & 0 & 0 & 0 & 0 \\
                0 & 0 & -1 & 0 & 0 & 0 & 0 \\
                0 & 0 & 0 & 0 & 0 & -1 & 0 \\
                0 & 0 & 0 & 0 & -1 & 0 & 0 \\
                0 & 0 & 0 & -1 & 0 & 0 & 0 \\
                0 & 0 & 0 & 0 & 0 & 0 &  1  \end{array} \right), \]
if necessary. Multiplying~(\ref{eqn:abcd})
on the left by 
\[ \left( \begin{array}{ccc}
            1 & 0 & 0 \\
            0 & 0 & 1 \\
            0 & -1 & 0 \end{array} \right) \]
and on the right by 
\[ \left( \begin{array}{ccccccc}
     1 & 0 & 0 & 0 & 0 & 0 & 0 \\
     0 & 0 & \frac{-c}{\sqrt{c^2+d^2}} & 0 & \frac{-d}{\sqrt{c^2+d^2}}  & 0 & 0 \\
     0 & \frac{c}{\sqrt{c^2+d^2}} & 0 & 0 & 0 & \frac{d}{\sqrt{c^2+d^2}}  & 0 \\
     0 & 0 & 0 & 1 & 0 & 0 & 0 \\
     0 & \frac{d}{\sqrt{c^2+d^2}} & 0 & 0 & 0 & \frac{-c}{\sqrt{c^2+d^2}} & 0 \\
     0 & 0 & \frac{-d}{\sqrt{c^2+d^2}} & 0 & \frac{c}{\sqrt{c^2+d^2}} & 0 & 0 \\
     0 & 0 & 0 & 0 & 0 & 0 & 1  \end{array} \right) \]
we get 
\[ \left( \begin{array}{ccccccc}
            a & 0 & 0 & 0 & 0 & 0 & 0 \\
            0 & \sqrt{c^2+d^2} & 0 & 0 & 0 & 0 & 0 \\
            0 & 0 & \frac{bc}{\sqrt{c^2+d^2}} & 0 &
            \frac{bd}{\sqrt{c^2+d^2}} 
& 0 & 0 \end{array} \right). 
\]
We can therefore assume $b\geq\sqrt{c^2+d^2}$. It is straightforward to check that 
no further restrictions can be put on $a$, $b$, $c$, $d$.
Hence we see that the orbit space 
is described by the matrices of the form~(\ref{eqn:abcd}) subject to the
conditions $a\geq b\geq\sqrt{c^2+d^2}$ and $c\geq0$. 

Next we want to compute a principal isotropy subgroup of $G$.
Let $p$ be an element of the form~(\ref{eqn:abcd}) which is regular 
for both $G$ and $\hat G$. The principal isotropy group $\hat G_p$ 
is isomorphic to $\mathbf Z_2^2\times\SO4$, and $G_p=\hat G_p\cap G$. 
Now the $\mathbf Z_2^2$-factor is contained in $G_p$, because
\[ (\mbox{diag($-1,-1,1$)},\mbox{diag($-1,-1,1,-1,1,-1,1$)}) \]
and
\[ (\mbox{diag($1,-1,-1$)},\mbox{diag($1,-1,-1,1,-1,-1,1$)})\]
are elements of $G$ that fix~$p$. On the other hand, an element in
the $\SO4$-factor that fixes $p$ and is in $G=\SO3\times\G$ 
is in fact in $\{1\}\times\G$ and has to fix 
$e$, $i$, $cj+d(ei)\in\mathbf R^7$. Therefore it also fixes 
$j=\frac1c[(cj+d(ei))-d(ei)]$ ($c\neq0$ because $p$ is regular for $G$)
and hence it is the identity. We conclude that the principal 
isotropy group $H=G_p$ is exactly $\mathbf Z_2^2$. 

It is immediate to see that the fixed subspace $V^H$ is
\[ \left( \begin{array}{ccccccc}
            * & 0 & 0 & * & 0 & 0 & 0 \\
            0 & * & 0 & 0 & 0 & * & 0 \\
            0 & 0 & * & 0 & * & 0 & 0 \end{array} \right). \]
Let $N$ be the normalizer of $H$ in $G$ and $\bar G=N/H$ the reduced group as in 
Subsection~\ref{subsec:red}. The connected component $\bar G^0$ must be 
contained in $\{1\}\times\G$. One can use the fact that $N$ equals the normalizer 
of $V^H$ in $G$ to verify that $\bar G^0$ is the $2$-torus group
consisting of the matrices
\setcounter{equation}{\value{thm}}
\stepcounter{thm}
\begin{equation}\label{eqn:2-torus}
\left( \begin{array}{ccccccc}
                \cos t & 0 & 0 & -\sin t & 0 & 0 & 0 \\
                0 & \cos s & 0 & 0 & 0 & -\sin s & 0 \\
                0 & 0 & \cos (s-t)& 0 & -\sin(s-t) & 0 & 0 \\
                \sin t & 0 & 0 & \cos t & 0 & 0 & 0 \\
                0 & 0 & \sin(s-t) & 0 & \cos(s-t) & 0 & 0 \\
                0 & \sin s & 0 & 0 & 0 & \cos s & 0 \\
                0 & 0 & 0 & 0 & 0 & 0 & 1  \end{array}
                \right) 
\end{equation}
for $s$, $t\in\mathbf R$. 

Let $p\in V^H$ be a regular point and let $M=Gp$ be the corresponding 
principal orbit in $V$. We next prove that $M$ is not taut by
contradiction. In fact, assume $M$ is taut. 
The orbit $\bar G^0p$ is a connected component of 
$M^H$, and it is a substantial $2$-torus embedded in 
the $6$-dimensional Euclidean space $V^H$. Since its second osculating
space is at most $5$-dimensional, it cannot be taut by
Theorem~\ref{thm:Kuiper}. It is known that a surface is taut if and only if 
its lines of curvature are circles. Therefore there is a parallel
normal vector field $\xi(t)$ along a curve $\gamma(t)$ in~$\bar G^0p\subset
M^H$ such that $A^{M^H}_{\xi(t)}\dot{\gamma}(t)=\lambda(t)\dot{\gamma}(t)$
and $\lambda(t)$ is not constant. (\emph{A posteriori} it is clear that one can
find $\gamma(t)$, $\xi(t)$, $\lambda(t)$ directly by computation 
using~(\ref{eqn:2-torus}).) Since $M$ is a principal $G$-orbit,
the normal spaces of~$M$ and~$M^H$ in~$V$ coincide along~$\gamma(t)$
implying that $\xi(t)$ is also parallel with respect to~$M$. 
By Lemma~\ref{lem:2} we have 
\setcounter{equation}{\value{thm}}
\stepcounter{thm}
\begin{equation}\label{eqn:weingarten}
A^{M}_{\xi(t)}\dot{\gamma}(t)=\lambda(t)\dot{\gamma}(t).
\end{equation}
We do not claim that $\gamma(t)$ is a curvature line of~$M$
since the eigenspace corresponding to~$\lambda(t)$ 
might not be one-dimensional. Still the argument 
in~\cite{Pinkall}, Lemmas~1 and~2 (see also~\cite{Miyaoka})
carries through and shows that~(\ref{eqn:weingarten}) 
together with the tautness of~$M$ implies that $\lambda(t)$ 
is constant, which is a contradiction. Hence~$M$ is not taut.
Thus the representation $\rho$ is not taut either. \EPf

\begin{prop}\label{taut}
The representations listed in Table~C.1 are taut.
\end{prop}

\Pf We use the reduction technique discussed in
Subsection~\ref{subsec:red}.
Consider first the case of $G=\U2\times\SP n$, $n\geq2$,
acting on the realification of $\mathbf C^2\otimes\mathbf C^{2n}$
by the tensor product of the vector representations.
It is equivalent to consider $\rho$ acting on the space of 
$2\times 2n$ complex matrices, i.~e.~$\rho(A,B)X=AXB^{-1}$ 
for $A\in\U2$, $B\in\SP n\subset\SU{2n}$, $X\in M(2\times2n;\mathbf C)$.
It follows rather easily that the orbit space consists of matrices of the form
\[ \left( \begin{array}{cccccccccc}
            a & 0 & 0 & \ldots & 0 & 0 & 0 & 0 & \ldots & 0 \\
            0 & b & 0 & \ldots & 0 & c & 0 & 0 & \ldots & 0 \end{array} \right), \]
where $a\geq\sqrt{b^2+c^2}$ and $b$, $c\geq0$. Since these matrices are all 
real, it follows that the group obtained by adjoining to~$G$ 
the complex conjugation of matrices $\sigma$ still has the same orbits.
We replace $G$ by the enlarged group and continue to denote it with the
same letter. 
 
Now the principal isotropy group 
$H$ is generated by $\sigma$ and the subgroup of $G$ isomorphic to
$\U1\times\SP{n-2}$ consisting of elements of the form 
\[ \left(\left( \begin{array}{cc} e^{i\theta} & 0 \\
                                               0 & e^{-i\theta} \end{array}
                                               \right),
    \left( \begin{array}{cccccc}
           e^{i\theta} & & &0 & & \\
           & e^{-i\theta} & & &0 & \\
           & & C & & & D \\
           0& & & e^{-i\theta} & & \\
           &0 & & & e^{i\theta} & \\
           & &-\bar D & & & \bar C \end{array} \right)\right), \]
where $\theta\in\mathbf R$ and 
$\left(\begin{array}{cc}C&D\\-\bar D&\bar C\end{array}\right)\in\SP{n-2}$.
Now $V^H$ is the space of matrices of the form
\[ \left( \begin{array}{cccccccccc}
            a & 0 & 0 & \ldots & 0 & 0 & d & 0 & \ldots & 0 \\
            0 & b & 0 & \ldots & 0 & c & 0 & 0 & \ldots & 0 \end{array} \right), \]
whose entries are real, and $(N/H)^0$ is the circle group 
consisting of the coclasses defined by the matrices 
\[   \left( \begin{array}{cccccc}
           \cos\theta & 0 &   &            0 & \sin\theta &    \\
           0 & \cos\theta &   &            \sin\theta & 0 &     \\
             &            & \mathbf1 &                   &   & \mathbf0  \\
           0 & -\sin\theta &   &          \cos\theta & 0 &    \\
           -\sin\theta & 0 &   &           0 & \cos\theta&     \\
             &            & \mathbf0 &            &           & \mathbf1     \\
\end{array} \right)\in\{1\}\times\SP n \]
for $\theta\in\mathbf R$, where $\mathbf1$ (resp.~$\mathbf0$) denotes
the $(n-2)\times(n-2)$ identity (resp.~zero) block. 
We see that the reduced representation $\bar\rho:(N/H)^0\to\mathbf O(V^H)$ 
is just the direct sum 
of two copies of the vector representation of $\SO2$ on $\mathbf R^2$,
hence taut. Finally, we conclude that $\rho$ is taut
by Proposition~\ref{prop:suf-taut}, where 
we choose $L$ to be generated by the following elements:
$\sigma$, the diagonal elements in $\SP{n-2}\subset H$ with $\pm1$ 
entries, and the element 
\[ \left(\left( \begin{array}{cc} i & 0 \\
                                               0 & -i \end{array}
                                               \right),
    \left( \begin{array}{cccccc}
           i & & & & & \\
           & -i & & & & \\
           & & \mathbf 1 & & &  \\
           & & & -i & & \\
           & & & & i & \\
           & & & & & \mathbf 1 \end{array} \right)\right). \]

Next consider the case of $G=\SO2\times\Spin9$ acting on $\mathbf
R^2\otimes\mathbf R^{16}$ by the tensor product of the vector 
and spin representations. It is equivalent to consider $\rho$ acting on the space of 
$2\times 16$ real matrices, namely, $\rho(A,B)X=AXB^{-1}$ 
for $A\in\SO2$, $B\in\Spin9\subset\SO{16}$, $X\in M(2\times16;\mathbf
R)$. We identify $\mathbf R^{16}\cong\mathbf O\oplus\mathbf O$ and 
use $\{1,e,i,j,k,ei,ej,ek\}$ as a basis of $\mathbf O$ over $\mathbf R$. 
For the computations involved in this case, 
Lemmas~14.61 and~14.77 from~\cite{Harvey} are useful, which we state as follows:

\begin{lem}\label{lem:harvey1}
A matrix $g\in\OG8\cong\mathbf O(\mathbf O)$ belongs to 
$\Spin7$ if and only if $g(uv)=g(u)\chi(g)(v)$ for all $u$, $v\in\mathbf
O$, where $\chi(g)(v)=g(g^{-1}(1)v)$ for $v\in\mathbf O$ defines 
the vector representation of $\Spin7$ on $\Im\mathbf O\cong\mathbf R^7$. 
\end{lem}

\begin{lem}\label{lem:harvey2}
The subgroup $\Spin9\subset\SO{16}\cong\mathbf{SO}(\mathbf O\oplus\mathbf
O)$ is generated by the $8$-sphere 
\[ \left\{\left(\begin{array}{cc} r & R_u \\ R_{\bar u} & -r
    \end{array}\right): r\in\mathbf R, u\in\mathbf O, r^2+|u|^2=1 \right\}
\]
in $\mbox{End}(\mathbf O\oplus\mathbf O)$,
where $R_u:\mathbf O\to\mathbf O$ denotes right translation by the element 
$u\in\mathbf O$. 
\end{lem}

Mainly as a consequence of 
\setcounter{equation}{\value{thm}}
\stepcounter{thm}
\begin{equation}\label{eqn:decomp}
 \Biv{}{}{}1|_{B_3}=\mbox{(trivial)}\oplus\Biii1{}{}\oplus\Biii{}{}1 
\end{equation}
one can see that the $G$-orbit space consists of 
matrices of the form
\setcounter{equation}{\value{thm}}
\stepcounter{thm}
\begin{equation}\label{eqn:p} 
\left( \begin{array}{cccccccccccccccc}
 a & 0 & 0 & 0 & 0 & 0 & 0 & 0 & 0 & 0 & 0 & 0 & 0 & 0 & 0 & 0 \\
 0 & b & 0 & 0 & 0 & 0 & 0 & 0 & c & 0 & 0 & 0 & 0 & 0 & 0 & 0 
          \end{array} \right),
\end{equation}
and, with some more effort (using Lemma~\ref{lem:harvey2}),
that the restrictions $a\geq\sqrt{b^2+c^2}$ and
$b$, $c\geq0$ apply. Also, Lemma~\ref{lem:harvey1} shows that 
\[ g_1=\mbox{diag}(-1,1,1,1,-1-1-1,1)\in\Spin7, \]
which combined with~(\ref{eqn:decomp}) yields that 
\[ g_2=\left(\begin{array}{ccc}
              1 &  &  \\
                & \chi(g_1) &  \\
                &  &  g_1  \end{array}\right) \in\Spin9, \]
where $\chi(g_1)=\mbox{diag}(-1,-1,-1,1,1,1,-1)$.
In particular, the equality
$\left(\begin{array}{cc} 1 & 0 \\ 0 & -1 \end{array}\right)p=pg_2$
where $p$ is as in~(\ref{eqn:p}) shows that $\OG2\times\Spin9$ 
has the same orbits as $\SO2\times\Spin9$. Hereafter we
consider $G=\OG2\times\Spin9$.

The principal isotropy group 
$H$ is isomorphic to $\mathbf Z_2\times\mathbf Z_2\times\SU3$, where the first 
$\mathbf Z_2$ corresponds to the kernel of $\rho$ and it is generated by the 
element $h_1=(-1,-1)\in\OG2\times\Spin9$, the second $\mathbf Z_2$ is generated
by the element
\[ h_2=(\mbox{diag}(1,-1),g_2), \]
and the $\SU3$-factor consists of matrices of the form 
\[ \left( \begin{array}{cccccc}
           1&&&&& \\
           &1&&&& \\
           &&\begin{array}{cc} A & B \\ -B & A \end{array}&&&\\
           &&&1&&\\
           &&&&1&\\
           &&&&&\begin{array}{cc} A & B \\ -B & A \end{array}
          \end{array} \right)\in\{1\}\times\Spin9, 
\]
where $A+iB\in\SU3$. 
Now $V^H$ is the space of matrices of the form
\[ \left( \begin{array}{cccccccccccccccc}
  a & 0 & 0 & 0 & 0 & 0 & 0 & 0 & 0 & d & 0 & 0 & 0 & 0 & 0 & 0 \\
  0 & b & 0 & 0 & 0 & 0 & 0 & 0 & c & 0 & 0 & 0 & 0 & 0 & 0 & 0 
          \end{array} \right), \]
whose entries are real, and $(N/H)^0$ is the circle group 
consisting of the coclasses defined by the following matrices
of $\{1\}\times\Spin9$:
\setcounter{equation}{\value{thm}}
\stepcounter{thm}
\begin{equation}\label{eqn:matrix}
  \left( \begin{array}{cccccc}
           \cos\theta&&& &\sin\theta&\\
           &\cos\theta&& -\sin\theta&&\\
           &&\begin{array}{cc}\cos\theta\mathbf1&\\
                              &\cos\theta\mathbf1
              \end{array} &&&\begin{array}{cc}&-\sin\theta\mathbf1\\
                                                   \sin\theta\mathbf1&\\
                              \end{array}\\
           &\sin\theta&&\cos\theta&&\\
           -\sin\theta&&&&\cos\theta&\\
           &&\begin{array}{cc}&-\sin\theta\mathbf1\\
                              \sin\theta\mathbf1
              \end{array} &&&\begin{array}{cc}\cos\theta\mathbf1\\
                                              &\cos\theta\mathbf1\\
                              \end{array}
\end{array}\right)
\end{equation}
for $\theta\in\mathbf R$, where $\mathbf1$ denotes the $3\times3$
identity block (note that the matrix~(\ref{eqn:matrix}) belongs to
$\Spin9$ because it is the product of the two matrices in~$\Spin9$
which are obtained by setting $r=\cos\theta$, $u=\sin\theta e$
and $r=1$, $u=0$ in Lemma~\ref{lem:harvey2}). 
We see that the reduced representation $\bar\rho$ is just the direct sum 
of two copies of the vector representation of $\SO2$ on $\mathbf R^2$,
hence taut. We conclude that $\rho$ is taut
by Proposition~\ref{prop:suf-taut}, where 
we choose $L$ to be generated by 
$h_1$, $h_2$ and  the elements in the
$\SU3$-factor of $H$  with $A$ diagonal with $\pm1$ entries and~$B=0$. 

Finally we consider the case of $G=\SU2\times\SP n$, $n\geq2$, acting on 
a real form $V$ of $\mathbf C^4\otimes\mathbf C^{2n}$. 
It is equivalent to consider $\rho$ acting on a subspace of 
$4\times 2n$ complex matrices, namely, 
$\rho(A,B)Z=\eta(A)ZB^{-1}$ 
for $A\in\SU2$, $B\in\SP n\subset\SU{2n}$, 
$Z\in V=\left\{\left(\begin{array}{cc} X&Y\\ -\bar Y& \bar X\end{array}\right)\in 
M(4\times2n;\mathbf C) \right\}$, where $\eta:\SU2\to\SP2\subset\SU4$ is the 
$4$-dimensional complex representation of $\SU2$ given by
\[
\eta
\left(\begin{array}{cc}\alpha&\beta\\-\bar\beta&\bar\alpha\end{array}\right)=
\left(\begin{array}{cccc}
\alpha^3 & \sqrt3\alpha\beta^2 & \beta^3 & \sqrt3\alpha^2\beta \\
\sqrt3\alpha\bar\beta^2 & \alpha\bar\alpha^2-2\bar\alpha\beta\bar\beta &
\sqrt3\bar\alpha^2\beta & \beta\bar\beta^2-2\alpha\bar\alpha\bar\beta \\
-\bar\beta^3 & -\sqrt3\bar\alpha^2\bar\beta & \bar\alpha^3 &
\sqrt3\bar\alpha\bar\beta^2 \\
-\sqrt3\alpha^2\bar\beta & 2\alpha\bar\alpha\beta-\beta^2\bar\beta &
\sqrt3\bar\alpha\beta^2 & \alpha^2\bar\alpha-2\alpha\beta\bar\beta
\end{array}\right). \]

Let 
\[ p_1=\left( \begin{array}{cccccccccc}
            1 & 0 & 0 & \ldots & 0 & 0 & 0 & 0 & \ldots & 0 \\
            0 & 0 & 0 & \ldots & 0 & 0 & 0 & 0 & \ldots & 0 \\
            0 & 0 & 0 & \ldots & 0 & 1 & 0 & 0 & \ldots & 0 \\
            0 & 0 & 0 & \ldots & 0 & 0 & 0 & 0 & \ldots & 0 \\
\end{array} \right)\in V. \]
Then the orbit $G(p_1)$ is $(4n+1)$-dimensional and its normal space 
intersects all orbits at~$p_1$. So, given an arbitrary orbit $M$, it 
contains a point of the form
\setcounter{equation}{\value{thm}}
\stepcounter{thm}
\begin{equation}\label{eqn:matrix2}
 \left( \begin{array}{cccccccccc}
            a & 0 & 0 & \ldots & 0 & 0 & 0 & 0 & \ldots & 0 \\
            x_1 & x_2 & x_3 & \ldots & x_n & 0 & y_2 & y_3 & \ldots & y_n \\
            0 & 0 & 0 & \ldots & 0 & a & 0 & 0 & \ldots & 0 \\
            0 & -\bar y_2 & -\bar y_3 & \ldots & -\bar y_n & \bar x_1 &
            \bar x_2 & \bar x_3 & \ldots & \bar x_n \\
\end{array} \right) 
\end{equation}
where $x_i$, $y_i\in\mathbf C$ and $a\geq0$, for we can multiply on the 
left by the minus identity matrix to change the sign of~$a$, 
if necessary, and remain in the
orbit~$M$. The isotropy subgroup $K_1=G_{p_1}$ is isomorphic to
$S^1\times\SP{n-1}$, where the circle factor consists of the elements
\[ \left(\left(\begin{array}{cc}e^{i\theta}&\\
                                &e^{-i\theta}\end{array}\right),
   \left(\begin{array}{cccc} e^{3i\theta}&&&\\
                             &\mathbf1&&\\
                             &&e^{-3i\theta}&\\
                             &&&\mathbf1\end{array}\right)\right) \]
for $\theta\in\mathbf R$, where $\mathbf1$ denotes the $(n-1)\times(n-1)$ 
identity block. By letting $K_1$ act on the matrix~(\ref{eqn:matrix2}) 
we now see that the $G$-orbit space consists of real matrices of the
form
\[ \left( \begin{array}{cccccccccc}
            a & 0 & 0 & \ldots & 0 & 0 & 0 & 0 & \ldots & 0 \\
            c & b & 0 & \ldots & 0 & 0 & 0 & 0 & \ldots & 0 \\
            0 & 0 & 0 & \ldots & 0 & a & 0 & 0 & \ldots & 0 \\
            0 & 0 & 0 & \ldots & 0 & c & b & 0 & \ldots & 0 \\
\end{array} \right). \]
(One can also check that some restrictions on~$a$, $b$ and~$c$
apply, but this is unimportant to us.) 

The principal isotropy group 
\[ H=
\left\{ 
\left(q,\left(\begin{array}{cc}\eta(q)&\\&C\end{array}\right)\right):\;
\parbox[c]{3in}{$C\in\SP{n-2},\quad
q=\left(\begin{array}{cc}\alpha&\beta\\-\bar\beta&\bar\alpha\end{array}\right),\\
(\alpha,\beta)=\pm(1,0),\pm(i,0),\pm(0,1),\pm(0,i)$} \right\} \]
is isomorphic to $Q\times\SP{n-2}$, where 
$Q$ is quaternion group 
$\{\pm1,\pm i, \pm j, \pm k\}$. 
Now $V^H$ is the space of matrices of the form
\[ \left( \begin{array}{cccccccccc}
            a & d & 0 & \ldots & 0 & 0 & 0 & 0 & \ldots & 0 \\
            c & b & 0 & \ldots & 0 & 0 & 0 & 0 & \ldots & 0 \\
            0 & 0 & 0 & \ldots & 0 & a & d & 0 & \ldots & 0 \\
            0 & 0 & 0 & \ldots & 0 & c & b & 0 & \ldots & 0 \\
\end{array} \right), \]
where the entries are real, and $(N/H)^0$ is the circle group 
consisting of the coclasses defined by the matrices
\[ \left( \begin{array}{cccccc}
            \cos\theta & -\sin\theta &  &  &  & \\
            \sin\theta & \cos\theta &  &  &  &   \\
             & & \mathbf1 & & & \\
             & & & \cos\theta & -\sin\theta & \\
             & & & \sin\theta & \cos\theta  & \\
             & & & & & \mathbf1 
\end{array} \right)\in\{1\}\times\SP{n} \]
We see that the reduced representation $\bar\rho$ is just the direct sum 
of two copies of the vector representation of $\SO2$ on $\mathbf R^2$,
hence taut. We conclude that $\rho$ is taut
by Proposition~\ref{prop:suf-taut}, where 
we choose $L$ to be generated by 
$Q$ and the diagonal elements in 
$\SP{n-2}\subset H$ with $\pm1$ entries. \EPf

\section{Other applications}\label{sec:other}
\setcounter{thm}{0}

In this last section we collect some further 
applications of our methods.

\begin{thm}[Dadok]\label{thm:Dadok}
A polar representation of a compact connected 
Lie group is orbit equivalent to the isotropy representation
of a symmetric space.
\end{thm}

\Pf Conlon's theorem~(Theorem~\ref{thm:conlon}) 
says that a polar representation is 
variationally complete. Now use Theorem~\ref{thm:BS}. \EPf

\medskip

The next theorem, which is a special case of a result of Montgomery 
and Samelson on transitive homeomorphic actions on spheres (\cite{MS}), 
see also~\cite{Borel1,Borel2} by Borel, was used in 
Sections~\ref{sec:O2} and~\ref{sec:BS}.

\begin{thm}\label{thm:BMS}
The transitive effective linear actions of compact connected Lie 
groups $G$ on the spheres
$S^d$ are given by the following tables (in each case, the 
action is the obvious one):
\[\begin{array}{|c|c|c|c|c|c|}
\hline
G & \SO n & \G & \Spin7 & \Spin9 & \SP1\cdot\SP n \\
\hline
S^d & S^{n-1} & S^6 & S^7 & S^{15} & S^{4n-1} \\
\hline
\end{array}\]

\[\begin{array}{|c|c|c|c|c|}
\hline
G & \SU n & S^1\cdot\SU n & \SP n & S^1\cdot\SP n \\
\hline
S^d & \multicolumn{2}{c|}{S^{2n-1}} &
\multicolumn{2}{c|}{S^{4n-1}} \\
\hline
\end{array}\]
\end{thm}

\Pf We may assume that the action is given by a 
real irreducible representation $\rho$ of $G$.
It is enough to classify $\rho$ up to
image equivalence and then factor the kernel out. 
If $\rho$ is of quaternionic or complex type, then Lemma~\ref{lem:C1notreal}
implies that $k(\lambda)=1$ and we get 
the standard action of~$\SO2$ on~$S^1$ and
the actions in the
second table.
If $\rho$ is of real type, then 
Lemmas~\ref{lem:C1}, \ref{lem:C1-classif} and~\ref{lem:C1-classif2}
together give the remaining actions in the first table. \EPf

\medskip

Representations $\rho:G\to\mathbf O(V)$ with cohomogeneity equal to two were
classified up to orbit equivalence in~\cite{H-L} by Hsiang and 
Lawson\footnote{The classification of cohomogeneity two and three
representations in~\cite{H-L} contain gaps which were 
to our knowledge first corrected in the papers~\cite{Uchida,Yasukura1}.}. 
It is not difficult to prove that a representation with
cohomogeneity two is polar. Hence it is also taut and if it
is irreducible it is therefore in class $\mathcal O^2$.

We next show that it is quite easy to see directly that an irreducible
representation with cohomogeneity two is in class
$\mathcal O^2$. We formulate this in the next proposition.

\begin{prop}\label{thm:O2}
An irreducible representation $\rho:G\to\mathbf O(V)$
with cohomogeneity two is in class $\mathcal O^2$.
\end{prop}

\Pf We can clearly restrict our attention to orbits contained in the unit
sphere~$S$ of~$V$. Let $p$ in $S$ be a point whose
orbit has codimension two.
The second osculating space ${\cal O}_p^2(Gp)$ is the direct sum of
$T_p(Gp)$ and the linear span of
$\{\alpha(X,Y)\,|\,X,Y\in T_p(Gp)\}$ where $\alpha$ is the second
fundamental form of the orbit $Gp$ as a submanifold of
$V$. It is a general property of submanifolds contained in the sphere $S$
that $\langle \alpha(X,X),\xi\rangle=1$ for all unit
vectors $X\in T_p(Gp)$, where $\xi$ is the inward
pointing unit normal vector of $S$ in $p$.
We therefore need to prove
that the linear span  of
$\{\alpha(X,Y)\,|\,X,Y\in T_p(Gp)\}$ is not one dimensional.  Let $\eta$ be
a normal
vector of $Gp$ at $p$ that is tangent to the sphere and let $A_\eta^S$ be
the Weingarten map of $Gp$ at $p$ as a submanifold of $S$. The orbit $Gp$
cannot be umbilical in $S$ since then it
would be a round sphere contradicting the assumption that $\rho$ is
irreducible.
Hence
$A_\eta^S$ has at least two different eigenvalues $\lambda_1$ and
$\lambda_2$ with corresponding eigenvectors $X_1$ and $X_2$
that we assume to have length one. Then
\begin{eqnarray*}
\alpha(X_i,X_i) 
&=&\langle\alpha(X_i,X_i),\xi\rangle\xi+\langle\alpha(X_i,X_i),\eta\rangle
\eta \\
                &=&\xi+\lambda_i\eta
\end{eqnarray*}
for $i=1,2$. It follows that $\alpha(X_1,X_1)$ and $\alpha(X_2,X_2)$ are not
linearly dependent and hence that $V={\cal
O}_p^2(Gp)$.

Now assume that the orbit through~$p$ in~$S$ is singular. 
If we restrict the action of
$\rho$ to $S$, then we have a cohomogeneity one action. It is well known that
the slice representations of cohomogeneity one
actions are irreducible and transitive on spheres in the normal space. If we
now think of the orbit $Gp$ as a submanifold of
$V$, then we see that the slice representation of $G_p$ on the normal space
$N_p(Gp)$ decomposes into the trivial
representation on ${\bf R}\xi$ and an irreducible representation on
$N^S_p(Gp)$ where $\xi$ is a normal vector of the sphere $S$ as above and
$N^S_p(Gp)$ is the normal space of $Gp$ in $S$.
The linear span of $\{\alpha(X,Y)\,|\,X,Y\in T_p(Gp)\}$ is clearly invariant
under the isotropy representation of $G_p$.
It follows that it contains $N^S_p(Gp)$ since it must have a component in
$N^S_p(Gp)$. (If it had no component in $N^S_p(Gp)$, then
$Gp$ would be totally geodesic in $S$, hence a round sphere, contradicting the
assumption that
$\rho$ is irreducible.) The
linear span of $\{\alpha(X,Y)\,|\,X,Y\in T_p(Gp)\}$ 
cannot coincide with $N^S_p(Gp)$, since
$\alpha(X,X)\not\in N^S_p(Gp)$ for a nonvanishing $X$ in $T_p(Gp)$. It
follows that  $V={\cal O}_p^2(Gp)$. \EPf

\begin{thm}[Hsiang-Lawson]
The cohomogeneity two linear actions of compact connected Lie groups on the 
Euclidean spaces
are given, up to orbit equivalence,
by the isotropy representations of the rank two symmetric
spaces.
\end{thm}

\Pf A cohomogeneity two reducible representation is the direct sum of
two cohomogeneity one representations and therefore orbit equivalent to 
the isotropy representation of the product of two spheres. 

A cohomogeneity two irreducible representation of 
a compact connected Lie group
is of class $\mathcal O^2$ (Proposition~\ref{thm:O2}). Therefore the 
result follows from the fact that no representation listed in Tables~C.1,
C.2, C.3 and~C.4 is of cohomogeneity two, which is easy to check. \EPf

\appendix

\setlength{\extrarowheight}{0cm}

\section{Description of the maximal subsets of
strongly orthogonal roots of the complex simple Lie algebras}

In this appendix, for each one
of the complex simple Lie algebras, we give a description of the maximal 
subset of strongly orthogonal roots $\mathcal B=\{\beta_1,\ldots,\beta_s\}$ 
which was introduced in Proposition~\ref{prop:Dadok}. 
In each case, $\Delta^+$ is the positive root system, 
$\mathcal S=\{\alpha_1,\ldots,\alpha_n\}$ 
is the simple root system, and $\lambda_1,\ldots,\lambda_r$ 
are the fundamental highest weights ($\theta_1$, $\theta_2,\ldots$ refer 
to specific linear forms on the Cartan subalgebra; see
e.~g.~\cite{Bourbaki}).

\[\begin{array}{|l|l|l|}
\hline
\multicolumn{3}{|l|}
{\A
  n:\An{\!\!\alpha_1}{\!\!\alpha_2}{\!\!\alpha_n}} \\
\hline
\multicolumn{3}{|l|}
{\mathcal S=\{\alpha_1=\theta_1-\theta_2,\ldots,\alpha_{n}=\theta_n-\theta_{n+1},\} } \\
\hline
\multicolumn{3}{|l|}
{\Delta^+=\{\theta_i-\theta_j:1\leq i<j\leq n+1\} }
\\
\hline
\lambda_1=\theta_1                      &\beta_1=\theta_1-\theta_{n+1}&
 \beta_1=\theta_1-\theta_{n+1} \\
\lambda_2=\theta_1+\theta_2             & \beta_2=\theta_2-\theta_n  &
\beta_2=\theta_2-\theta_n  \\
\qquad\vdots                  & \qquad\vdots   & \qquad\vdots         \\
\lambda_n=\theta_1+\theta_2+\ldots+\theta_n  & 
\beta_n=\theta_{\frac n2}-\theta_{\frac n2+2} &
\beta_n=\theta_{\frac {n+1}2}-\theta_{\frac {n+1}2+1}\\
 &  \mbox{($n$ even)} & \mbox{($n$ odd)}\\
\hline
\end{array}\]

\[\begin{array}{|l|l|l|}
\hline
\multicolumn{3}{|l|}
{\B
  n:\Bn{\!\!\alpha_1}{\!\!\alpha_2}{\!\!\!\!\alpha_{n-1}}{\!\!\alpha_n}} \\
\hline
\multicolumn{3}{|l|}
{\mathcal S=\{\alpha_1=\theta_1-\theta_2,\ldots,\alpha_{n-1}=\theta_{n-1}-\theta_{n},
\alpha_{n}=\theta_{n}\} } \\
\hline
\multicolumn{3}{|l|}
{\Delta^+=\{\theta_i\pm\theta_j:1\leq i<j\leq n\}\cup\{\theta_i:1\leq n\}}
\\
\hline
\lambda_1=\theta_1 & \beta_1=\theta_1+\theta_2 & \beta_1=\theta_1+\theta_2 \\
\lambda_2=\theta_1+\theta_2 &\beta_2=\theta_1-\theta_2 & \beta_2=\theta_1-\theta_2 \\
\qquad\vdots & \qquad\vdots   & \qquad\vdots       \\
\lambda_{n-2}=\theta_1+\theta_2+\ldots+\theta_{n-2}  & 
\beta_{n-2}=\theta_{n-2}+\theta_{n-1}
& \beta_{n-1}=\theta_{n-1}+\theta_{n}\\
\lambda_{n-1}=\theta_1+\theta_2+\ldots+\theta_{n-2}+\theta_{n-1} &
\beta_{n-1}=\theta_{n-2}-\theta_{n-1} 
& \beta_{n}=\theta_{n-1}-\theta_{n}\\
\lambda_{n}=\frac12(\theta_1+\theta_2+\ldots+\theta_{n-1}+\theta_{n})   & 
\beta_n=\theta_n & \\
 & \mbox{($n$ odd)} & \mbox{($n$ even)}\\

\hline
\end{array}\]

\[\begin{array}{|l|l|}
\hline
\multicolumn{2}{|l|}
{\C
  n:\Cn{\!\!\alpha_1}{\!\!\alpha_2}{\!\!\!\!\alpha_{n-1}}{\!\!\alpha_n}} \\
\hline
\multicolumn{2}{|l|}
{\mathcal S=\{\alpha_1=\theta_1-\theta_2,\ldots,\alpha_{n-1}=\theta_{n-1}-\theta_n,
\alpha_n=2\theta_n\} } \\
\hline
\multicolumn{2}{|l|}
{\Delta^+=\{\theta_i\pm\theta_j:1\leq i<j\leq n\}\cup\{2\theta_i:1\leq n\}}
\\
\hline
\lambda_1=\theta_1                           & \beta_1=2\theta_1 \\
\lambda_2=\theta_1+\theta_2                  & \beta_2=2\theta_2 \\
\qquad\vdots                                 & \qquad\vdots            \\
\lambda_n=\theta_1+\theta_2+\ldots+\theta_n  & \beta_n=2\theta_n \\
\hline
\end{array}\]

\[\begin{array}{|l|l|l|}
\hline
\multicolumn{3}{|l|}
{\D
  n:\Dn{\!\!\alpha_1}{\!\!\alpha_2}{\!\!\!\!\!\!\!\!\alpha_{n-2}}{\alpha_{n-1}}
{\alpha_n}} \\
\hline
\multicolumn{3}{|l|}
{\mathcal S=\{\alpha_1=\theta_1-\theta_2,\ldots,\alpha_{n-2}=\theta_{n-2}-\theta_{n-1},
\alpha_{n-1}=\theta_{n-1}-\theta_{n},\alpha_{n}=\theta_{n-1}+\theta_{n}\} } \\
\hline
\multicolumn{3}{|l|}
{\Delta^+=\{\theta_i\pm\theta_j:1\leq i<j\leq n\} }
\\
\hline
\lambda_1=\theta_1 & \beta_1=\theta_1+\theta_2 & \beta_1=\theta_1+\theta_2 \\
\lambda_2=\theta_1+\theta_2 &\beta_2=\theta_1-\theta_2 & \beta_2=\theta_1-\theta_2 \\
\qquad\vdots & \qquad\vdots   & \qquad\vdots       \\
\lambda_{n-2}=\theta_1+\theta_2+\ldots+\theta_{n-2}  & 
\beta_{n-2}=\theta_{n-2}+\theta_{n-1}
& \beta_{n-1}=\theta_{n-1}+\theta_{n}\\
\lambda_{n-1}=\frac12(\theta_1+\theta_2+\ldots+\theta_{n-1}-\theta_{n}) &
\beta_{n-1}=\theta_{n-2}-\theta_{n-1} 
& \beta_{n}=\theta_{n-1}-\theta_{n}\\
\lambda_{n}=\frac12(\theta_1+\theta_2+\ldots+\theta_{n-1}+\theta_{n})   & 
\mbox{($n$ odd)} & \mbox{($n$ even)}\\

\hline
\end{array}\]

\setlength{\extrarowheight}{0cm}

\[\begin{array}{|l|l|}
\hline
\multicolumn{2}{|l|}
{\G :\Gii{\!\!\alpha_1}{\!\!\alpha_2}} \\
\hline
\multicolumn{2}{|l|}
{\Delta^+=\{\alpha_1,\alpha_2,\alpha_1+\alpha_2,2\alpha_1+\alpha_2,
3\alpha_1+\alpha_2,3\alpha_1+2\alpha_2\}}
\\
\hline
\lambda_1=2\alpha_1+\alpha_2      & \beta_1= 3\alpha_1+2\alpha_2 \\
\lambda_2=3\alpha_1+2\alpha_2     & \beta_2= \alpha_1 \\
\hline
\end{array}\]

\[\begin{array}{|l|l|}
\hline
\multicolumn{2}{|l|}
{\F :\Fiv{\!\!\alpha_1}{\!\!\alpha_2}{\!\!\alpha_3}{\!\!\alpha_4}} \\
\hline
\multicolumn{2}{|l|}
{\mathcal S=\{\alpha_1=\theta_2-\theta_3,\alpha_2=\theta_3-\theta_4,
\alpha_3=\theta_4,\alpha_4=\frac12(\theta_1-\theta_2-\theta_3-\theta_4)\}}
\\
\hline
\multicolumn{2}{|l|}
{\Delta^+=\{\theta_i:1\leq i\leq4\}\cup\{\theta_i\pm\theta_j:1\leq
  i<j\leq4\}\cup\{\frac12(\theta_1\pm\theta_2\pm\theta_3\pm\theta_4)\}}
\\
\hline
\lambda_1=\theta_1+\theta_2     & \beta_1= \theta_1+\theta_2   \\
\lambda_2=2\theta_1+\theta_2+\theta_3    & \beta_2= \theta_1-\theta_2  \\
\lambda_3=\frac12(3\theta_1+\theta_2+\theta_3+\theta_4)&\beta_3=\theta_3+\theta_4\\  
\lambda_4=\theta_1  & \beta_4=\theta_3-\theta_4 \\
\hline
\end{array}\]

\[\begin{array}{|l|l|}
\hline
\multicolumn{2}{|l|}
{\E6 :\Evi{\!\!\alpha_1}{\alpha_2}{\!\!\alpha_3}{\!\!\alpha_4}{\!\!\alpha_5}{\!\!\alpha_6}} \\
\hline
\multicolumn{2}{|l|}
{\mathcal S=\{\alpha_1=\frac12(\theta_1-\theta_2-\theta_3-\theta_4-\theta_5-\theta_6
-\theta_7+\theta_8),}\\
\multicolumn{2}{|l|}
  {\qquad\alpha_2=\theta_1+\theta_2,
   \alpha_3=\theta_2-\theta_1,
   \alpha_4=\theta_3-\theta_2,
   \alpha_5=\theta_4-\theta_3,
   \alpha_6=\theta_5-\theta_4 \} }
\\
\hline
\multicolumn{2}{|l|}
{\Delta^+=\{\pm\theta_i+\theta_j:1\leq
  i<j\leq5\}\cup\{\frac12(\theta_8-\theta_7-\theta_6+
                 \sum_{i=1}^5\epsilon_i\theta_i):\Pi_{i=1}^5\epsilon_i=+1\} }
\\
\hline
\lambda_1=\frac23(\theta_8-\theta_7-\theta_6)     
& \beta_1 =\frac12(\theta_1+\theta_2+\theta_3+\theta_4 \\
 \lambda_2=\frac12(\theta_1+\theta_2+\theta_3+\theta_4+\theta_5-\theta_6-\theta_7+\theta_8)
 &\qquad +\theta_5-\theta_6-\theta_7+\theta_8) \\ 
\lambda_3=\frac56(\theta_8-\theta_7-\theta_6)+\frac12(-\theta_1+\theta_2+\theta_3+\theta_4+\theta_5)
& \beta_2 = \frac12(-\theta_1-\theta_2-\theta_3-\theta_4 \\
\lambda_4=\theta_3+\theta_4+\theta_5-\theta_6-\theta_7+\theta_8
& \qquad +\theta_5-\theta_6-\theta_7+\theta_8)   \\
\lambda_5=\frac23(\theta_8-\theta_7-\theta_6)+\theta_4+\theta_5 
& \beta_3=\theta_4-\theta_1\\
\lambda_6=\frac13(\theta_8-\theta_7-\theta_6)+\theta_5
& \beta_4=\theta_3-\theta_2 \\
\hline
\end{array}\]

\[\begin{array}{|l|l|}
\hline
\multicolumn{2}{|l|}
{\E7 :\Evii{\!\!\alpha_1}{\alpha_2}{\!\!\alpha_3}{\!\!\alpha_4}{\!\!\alpha_5}{\!\!\alpha_6}{\!\!\alpha_7}} \\
\hline
\multicolumn{2}{|l|}
{\mathcal S=\{\alpha_1=\frac12(\theta_1-\theta_2-\theta_3-\theta_4-\theta_5-\theta_6
-\theta_7+\theta_8),}\\
\multicolumn{2}{|l|}
  {\qquad\alpha_2=\theta_1+\theta_2,
   \alpha_3=\theta_2-\theta_1,
   \alpha_4=\theta_3-\theta_2,
   \alpha_5=\theta_4-\theta_3,
   \alpha_6=\theta_5-\theta_4,
   \alpha_7=\theta_6-\theta_5\} }
\\
\hline
\multicolumn{2}{|l|}
{\Delta^+=\{\pm\theta_i+\theta_j:1\leq
  i<j\leq6\}\cup\{\theta_8-\theta_7\}\cup\{\frac12(\theta_8-\theta_7+
                 \sum_{i=1}^6\epsilon_i\theta_i):\Pi_{i=1}^6\epsilon_i=-1\} }
\\
\hline
\lambda_1=\theta_8-\theta_7    
& \beta_1=\theta_8-\theta_7\\
\lambda_2=\frac12(\theta_1+\theta_2+\theta_3+\theta_4+\theta_5+\theta_6-2\theta_7+2\theta_8)
 & \beta_2= \theta_6+\theta_5  \\
\lambda_3=\frac12(-\theta_1+\theta_2+\theta_3+\theta_4+\theta_5+\theta_6-3\theta_7+3\theta_8)
&\beta_3=\theta_6-\theta_5\\  
\lambda_4=\theta_3+\theta_4+\theta_5+\theta_6+2(\theta_8-\theta_7)
& \beta_4=\theta_4+\theta_3 \\
\lambda_5=\frac12(2\theta_4+2\theta_5+2\theta_6+3(\theta_8-\theta_7)) 
& \beta_5=\theta_4-\theta_3\\
\lambda_6=\theta_5+\theta_6-\theta_7+\theta_8 
& \beta_6=\theta_2+\theta_1\\
\lambda_7=\theta_6+\frac12(\theta_8-\theta_7)
& \beta_7=\theta_2-\theta_1\\
\hline
\end{array}\]

\[\begin{array}{|l|l|}
\hline
\multicolumn{2}{|l|}
{\E8 :\Eviii{\!\!\alpha_1}{\alpha_2}{\!\!\alpha_3}{\!\!\alpha_4}{\!\!\alpha_5}{\!\!\alpha_6}{\!\!\alpha_7}{\!\!\alpha_8}} \\
\hline
\multicolumn{2}{|l|}
{\mathcal S=\{\alpha_1=\frac12(\theta_1-\theta_2-\theta_3-\theta_4-\theta_5-\theta_6
-\theta_7+\theta_8),\alpha_2=\theta_1+\theta_2,}\\
\multicolumn{2}{|l|}
  {\qquad\alpha_3=\theta_2-\theta_1,
   \alpha_4=\theta_3-\theta_2,
   \alpha_5=\theta_4-\theta_3,
   \alpha_6=\theta_5-\theta_4,
   \alpha_7=\theta_6-\theta_5,
   \alpha_8=\theta_7-\theta_6\}}
\\
\hline
\multicolumn{2}{|l|}
{\Delta^+=\{\pm\theta_i+\theta_j:1\leq
  i<j\leq8\}\cup\{\frac12(\theta_8+
                 \sum_{i=1}^7\epsilon_i\theta_i):\Pi_{i=1}^7\epsilon_i=+1\} }
\\
\hline
\lambda_1=2\theta_8   
& \beta_1=\theta_7+\theta_8\\
\lambda_2=\frac12(\theta_1+\theta_2+\theta_3+\theta_4+\theta_5+\theta_6+\theta_7+5\theta_8)
 & \beta_2= \theta_8-\theta_7  \\
\lambda_3=\frac12(-\theta_1+\theta_2+\theta_3+\theta_4+\theta_5+\theta_6+\theta_7+7\theta_8)
&\beta_3=\theta_6+\theta_5\\  
\lambda_4=\theta_3+\theta_4+\theta_5+\theta_6+\theta_7+5\theta_8
& \beta_4=\theta_6-\theta_5 \\
\lambda_5=\theta_4+\theta_5+\theta_6+\theta_7+4\theta_8 
& \beta_5=\theta_4+\theta_3\\
\lambda_6=\theta_5+\theta_6+\theta_7+3\theta_8 
& \beta_6=\theta_4-\theta_3\\
\lambda_7=\theta_6+\theta_7+2\theta_8
& \beta_7=\theta_2+\theta_1\\
\lambda_8=\theta_7+\theta_8
& \beta_8=\theta_2-\theta_1 \\
\hline
\end{array}\]

\section{Table for Dadok's invariant $k(\lambda)$}

We next tabulate the value of the Dadok invariant $k(\lambda)$ 
for the the complex irreducible representations $\pi_\lambda$ 
of the complex simple Lie algebras, where $\lambda$ is a 
fundamental highest weight (since $k(\lambda)$ is linear on
$\lambda$, this already determines the value of 
$k(\lambda)$ for all $\lambda$). The number in parenthesis 
next to the vertex of the Dynkin diagram which specifies $\alpha_i$
is $k(\lambda_i)$. The last column tells, for each complex simple
Lie algebra, whether all representations are self-dual, or
else, the dual of $\pi_{\lambda_i}$ is given by $\pi_{\lambda_j}$,
where the vertices that specify $\alpha_i$ and $\alpha_j$ are images 
one of the other under the symmetry of the Dynkin diagram. 

\setlength{\extrarowheight}{1cm}

\[\begin{array}{|c|c|c|}
\hline
G & k(\lambda) & Conditions \\
\hline
\A{2n-1} & \Anii{\!\!(1)}{\!\!(2)}{\!\!\!\!\!\!\!(n-1)}{\!\!(n)}{\!\!\!\!\!\!\!(n-1)}{\!\!(2)}{\!\!(1)} & 
\parbox[c]{4cm}{duality given by symmetry of diagram} \\
\A{2n} & \Aniii{\!\!(1)}{\!\!(2)}{\!\!\!\!\!\!\!(n-1)}{\!\!(n)}{\!\!(n)}{\!\!\!\!\!\!\!(n-1)}{\!\!(2)}{\!\!(1)} & 
\parbox[c]{4cm}{duality given by symmetry of diagram} \\
\B{2n-1} & \Bniii{\!\!(2)}{\!\!(2)}{\!\!(4)}{\!\!(4)}{\!\!\!\!\!\!\!\!\!\!\!(2n-2)}
{\!\!\!\!\!\!\!\!(2n-2)}{\!\!(n)} &  \mbox{all self-dual} \\
\B{2n} & \Bniv{\!\!(2)}{\!\!(2)}{\!\!(4)}{\!\!(4)}{\!\!\!\!\!\!\!\!\!\!\!(2n-2)}
{\!\!\!\!\!\!\!\!(2n-2)}{\!\!\!\!(2n)}{\!\!(n)} &  \mbox{all self-dual} \\
\C n  & \Cn{\!\!(1)}{\!\!(2)}{\!\!\!\!\!\!\!(n-1)}{\!\!(n)} & \mbox{all
  self-dual} \\
\D{2n} & \Dniii{\!\!(2)}{\!\!(2)}{\!\!(4)}{\!\!(4)}
{\!\!\!\!\!\!\!\!\!\!\!\!\!\!\!\!(2n-2)}
{\!\!\!\!\!\!\!\!\!\!\!\!\!(2n-2)}{(n)}{(n)} &  \mbox{all self-dual} \\
\D{2n+1} & 
\Dniv{\!\!(2)}{\!\!(2)}{\!\!(4)}{\!\!(4)}{\!\!\!\!\!\!\!\!\!\!\!\!(2n-2)}
{\!\!\!\!\!\!\!\!\!(2n-2)}{\!\!\!\!\!\!(2n)}{(n)}{(n)} &  
\parbox[c]{4cm}{duality given by symmetry of diagram} \\
\G     & \Gii{\!\!(2)}{\!\!(2)} & \mbox{all self-dual} \\
\F     & \Fiv{\!\!(2)}{\!\!(6)}{\!\!(4)}{\!\!(2)} & \mbox{all self-dual} \\
\E6    & \Evi{\!\!(2)}{(2)}{\!\!(4)}{\!\!(6)}{\!\!(4)}{\!\!(2)} & \parbox[c]{4cm}{duality 
given by symmetry of diagram}\\
\E7    & \Evii{\!\!(2)}{(5)}{\!\!(6)}{\!\!(8)}{\!\!(7)}{\!\!(4)}{\!\!(3)} & \mbox{all self-dual} \\
\E8    &
\Eviii{\!\!(4)}{(8)}{\!\!\!(10)}{\!\!\!(14)}{\!\!\!(12)}{\!\!(8)}{\!\!(6)}{\!\!(2)}
& \mbox{all self-dual}\\
\hline
\end{array}\]
\begin{center}
\textbf{Table~B.1}: The invariant of Dadok.
\end{center}

\section{Table of candidates to class $\mathcal O^2$}

\setlength{\extrarowheight}{0.35cm}

In the tables below are compiled all the real irreducible
representations of compact connected Lie groups
that are not orbit equivalent to the isotropy representation 
of a symmetric space but that have a chance of being a representation
of class $\mathcal O^2$. These are just the representations 
mentioned in Propositions~\ref{prop:quat}, 
\ref{prop:complex}, \ref{prop:real4}, \ref{prop:real3}, \ref{prop:real21},
\ref{prop:real22} and \ref{prop:simple}, except those excluded by 
the discussion in Section~\ref{sec:BS}; note however that the presentation 
here slightly differs from the presentation in the
above mentioned propositions for the sake of convenience.

\[\begin{array}{|c|c|c|}
\hline
G & \rho & Conditions\\
\hline
\SO 2\times\Spin9  & (x\oplus x^{-1})\otimes\Biv{}{}{}1 & - \\
\U 2\times\SP n  & (x\oplus x^{-1})\otimes\Ai1\otimes\Cns1{}{} & n\geq2 \\
\SU2 \times \SP n      & \Ai3\otimes\Cns1{}{} & n\geq2 \\           
\hline
\end{array}\]
\begin{center}
\textbf{Table~C.1}: Taut representations.
\end{center}

\setlength{\extrarowheight}{0.27cm}

\[\begin{array}{|c|c|c|}
\hline
G & \rho & Conditions\\
\hline
\SU n\times\SP m           & \begin{array}{c}
                             (\Ans1{}\oplus\Ans{}1) \\
                             \qquad\otimes\;\Cns1{}{}\end{array}  & n\geq3,\;
m\geq2 \\
\U n\times\SP m  & \begin{array}{c}
                    (x\otimes\Ans1{}\oplus x^{-1}\otimes\Ans{}1) \\
                    \qquad\otimes\;\Cns1{}{} \end{array} &n\geq3,\; m\geq2 \\
 & & \\
\begin{array}{c}\SO{m}\times\SP1\cdot\SP n \\
                 (n\geq2) \end{array}
                               & \begin{array}{c}
                               \Dns1{}{}{}\otimes\Ai1\otimes\Cns1{}{}\\
                               \Bns1{}{}\otimes\Ai1\otimes\Cns1{}{}\\
                               \Ai1\otimes\Ai1\otimes\Ai1\otimes\Cns1{}{}\\
                               \Ai2\otimes\Ai1\otimes\Cns1{}{}
                              \end{array} & \begin{array}{c}
                                m=2p>4 \\
                                m=2p+1>3 \\ 
                                m=4 \\
                                m=3
                                     \end{array}  \\            
 & & \\
\SO m\times\G      & \begin{array}{c}
                        \Dns1{}{}{}\otimes\Gii1{} \\
                        \Bns1{}{}\otimes\Gii1{} \\
                        \Ai1\otimes\Ai1\otimes\Gii1{} \\
                        \Ai2\otimes\Gii1{}
                     \end{array}  &  \begin{array}{c} m=2p>4\\m=2p+1>3 \\
                                      m=4 \\ m=3 \end{array}  \\
 & & \\
\SO m\times\Spin7  & \begin{array}{c}
                       \Dns1{}{}{}\otimes\Biii{}{}1 \\
                       \Bns1{}{}\otimes\Biii{}{}1 \\
                       \Ai1\otimes\Ai1\otimes\Biii{}{}1
                     \end{array}  &  \begin{array}{c} m=2p>4 \\ m=2p+1\geq5 \\  m=4
                                     \end{array}  \\
 & & \\
\SO m\times\Spin9  & \begin{array}{c}
                       \Dns1{}{}{}\otimes\Biv{}{}{}1 \\
                       \Bns1{}{}\otimes\Biv{}{}{}1 \\
                       \Ai1\otimes\Ai1\otimes\Biv{}{}{}1 \\
                       \Ai2 \otimes\Biv{}{}{}1 
                     \end{array}  &  \begin{array}{c} m=2p>4
                       \\ m=2p+1>3 \\ m=4 \\ m=3 \end{array}  \\
\hline
\end{array}\]
\begin{center}
\textbf{Table C.2}: Nontaut representations (first group).
\end{center}

\setlength{\extrarowheight}{0.35cm}

\[\begin{array}{|c|c|c|}
\hline
 G & \rho & Conditions \\
\hline
\G\times\G         & \Gii1{}\otimes\Gii1{} & - \\
\G\times\Spin7     & \Gii1{}\otimes\Biii{}{}1 & - \\
\G\times\Spin9     & \Gii1{}\otimes\Biv{}{}{}1& - \\
\Spin7\times\Spin7 & \Biii{}{}1\otimes\Biii{}{}1 & - \\
\Spin7\times\Spin9 & \Biv{}{}{}1\otimes\Biv{}{}{}1 & - \\
\SP1\cdot\SP n\times\G     & \Ai1\otimes\Cns1{}{}\otimes\Gii1{} & n\geq2 \\
\SP1\cdot\SP n\times\Spin7 & \Ai1\otimes\Cns1{}{}\otimes\Biii{}{}1 &
n\geq2 \\
\hline
\end{array}\]
\begin{center}
\textbf{Table C.3}: Nontaut representations (second group).
\end{center}

\[\begin{array}{|c|c|c|}
\hline
 G & \rho & Conditions \\
\hline
\SP n\times\SU6      & \Cns1{}{}\otimes\Av{}{}1{}{}        & n\geq2 \\
\SP n\times\Spin{11} & \Cns1{}{}\otimes\Bv{}{}{}{}1        & -     \\
\SP n\times\Spin{12} & \Cns1{}{}\otimes\Dvi{}{}{}{}1{}     & n\geq2 \\
\SP 1\times\Spin{13} & \Ai1\otimes\Bvi{}{}{}{}{}1           & -     \\
\SP 1\times\SP2      & \Ai1\otimes\Cii11                    & -     \\
\SP n\times\SP3      & \Cns1{}{}\otimes\Ciii{}{}1          & n\geq2 \\
\SP n\times\E7       & \Cns1{}{}\otimes\Evii{}{}{}{}{}{}1  & n\geq2 \\
\Spin7      & \Biii1{}1  & - \\
\Spin9      & \Biv1{}{}1 & - \\
\Spin{15}   & \Bvii{}{}{}{}{}{}1 & - \\
\Spin{17}   & \Bviii{}{}{}{}{}{}{}1 & - \\
\hline
\end{array}\]
\begin{center}
\textbf{Table C.4}: Nontaut representations (third group).
\end{center}

\bibliographystyle{amsplain}
\bibliography{repr}

\end{document}